\documentclass[11pt]{article}
\usepackage[T1]{fontenc}
\usepackage[utf8]{inputenc}
\usepackage{lmodern}
\usepackage[a4paper,margin=27mm,headheight=14pt]{geometry}
\usepackage{amsmath,amssymb,amsthm,mathtools}
\usepackage{microtype}
\usepackage{mathrsfs}
\usepackage{enumitem}
\usepackage{needspace}
\usepackage{fancyhdr}
\usepackage[hidelinks,unicode,pagebackref]{hyperref}
\usepackage[nameinlink,noabbrev]{cleveref}
\newcommand{\manuscriptkeywords}{Optimal transport, coupled Bochner method,
moment measures, Wasserstein proximal maps, displacement convexity, functional inequalities}
\hypersetup{pdfkeywords={\manuscriptkeywords},pdftitle={Transport Stability and Wasserstein Proximal Maps: A Coupled Bochner Approach},pdfauthor={Bang-Xian Han and Zhuo-Nan Zhu},pdfsubject={Second variations of optimal transport solution maps, stability, and rigidity}}
\renewcommand*{\backrefalt}[4]{%
  \ifcase#1\relax
  \or \space(Cited on p.~#2.)%
  \else \space(Cited on pp.~#2.)%
  \fi}

\allowdisplaybreaks[2]
\numberwithin{equation}{section}
\newtheorem{theorem}{Theorem}[section]
\newtheorem{proposition}[theorem]{Proposition}
\newtheorem{lemma}[theorem]{Lemma}
\newtheorem{corollary}[theorem]{Corollary}
\theoremstyle{definition}

\theoremstyle{remark}
\newtheorem{remark}[theorem]{Remark}
\newtheorem{example}[theorem]{Example}
\crefname{theorem}{Theorem}{Theorems}
\crefname{lemma}{Lemma}{Lemmas}
\crefname{definition}{Definition}{Definitions}
\crefname{proposition}{Proposition}{Propositions}
\crefname{corollary}{Corollary}{Corollaries}
\crefname{section}{Section}{Sections}
\crefname{appendix}{appendix}{appendices}
\Crefname{appendix}{Appendix}{Appendices}

\newcommand{\R}{\mathbb R}
\newcommand{\E}{\mathbb E}
\newcommand{\N}{\mathcal N}
\newcommand{\Moment}{\mathsf M}
\newcommand{\Reconstruct}{\mathsf R}
\DeclareMathOperator{\Cov}{Cov}
\DeclareMathOperator{\range}{range}
\DeclareMathOperator{\Lip}{Lip}
\newcommand{\Ptwo}{\mathcal P_2}
\newcommand{\Pc}{\mathcal P_c}
\newcommand{\Id}{I}
\newcommand{\diff}{\mathrm{d}}
\newcommand{\D}{\mathrm{D}}
\newcommand{\dd}{\,\diff}
\newcommand{\dx}{\dd x}
\newcommand{\dy}{\dd y}
\newcommand{\dt}{\dd t}
\newcommand{\one}{\mathbf 1}
\newcommand{\HS}{\mathrm{HS}}
\newcommand{\II}{\mathrm{II}}
\newcommand{\tr}{\operatorname{tr}}
\newcommand{\diver}{\operatorname{div}}
\newcommand{\Ent}{\operatorname{Ent}}
\newcommand{\Var}{\operatorname{Var}}
\newcommand{\supp}{\operatorname{supp}}
\newcommand{\cof}{\operatorname{cof}}
\newcommand{\dist}{\operatorname{dist}}
\newcommand{\diam}{\operatorname{diam}}
\newcommand{\cJ}{\mathcal J}

\newcommand{\cE}{\mathcal E}
\newcommand{\cK}{\mathcal K}
\newcommand{\cR}{\mathcal R}

\newcommand{\cC}{\mathcal C}
\newcommand{\HH}{\mathcal H}
\newcommand{\Dcal}{\mathsf D}
\DeclareMathOperator*{\argmin}{argmin}
\newcommand{\inner}[2]{\langle #1,#2\rangle}
\newcommand{\norm}[1]{\lVert #1\rVert}

\title{Transport Stability and\\Wasserstein Proximal Maps}
\newcommand{\manuscriptsubtitle}{A Coupled Bochner Approach}

\makeatletter
\renewcommand{\@maketitle}{%
  \null
  \vspace{8pt}%
  \begingroup
    \centering
    {\fontsize{25}{31}\selectfont\bfseries\@title\par}%
    \vspace{14pt}%
    {\fontsize{16}{20}\selectfont\manuscriptsubtitle\par}%
    \vspace{24pt}%
    {\large
     \begin{tabular}[t]{c}\@author\end{tabular}\par}%
  \endgroup
  \vspace{18pt}%
}
\makeatother

\author{Bang-Xian Han\thanks{School of Mathematics, Shandong University, Jinan 250100, China. Email: \texttt{hanbx@sdu.edu.cn}}
\and Zhuo-Nan Zhu\thanks{School of Mathematical Sciences, University of Science and Technology of China, Hefei 230026, China. Email: \texttt{zhuonanzhu@mail.ustc.edu.cn}}}
\date{}
\begin{document}
\maketitle
\thispagestyle{plain}
\begin{abstract}
We develop a second-variation approach to transport stability and
Wasserstein proximal maps. For a given density variation we compare the
velocity of least Euclidean kinetic energy with the one of least
energy in the Hessian metric of the transport potential. The latter
field and its pushforward give Bochner estimates at the two ends of
the transport. In moment reconstruction and in a proximal step, the
same quadratic form gives the Hessian at a smooth minimizer for a
fixed target and controls the motion of that minimizer as the target
varies.

\smallskip
For two log-concave measures evolving under the same heat flow, we
prove that their squared Wasserstein distance is convex in time.
When the target varies in moment reconstruction, we obtain a sharp
Lipschitz estimate modulo translations for strongly log-concave
reconstructed measures. With the source fixed, we characterize
Lipschitz stability of Brenier potentials by a projected Poincar\'e
inequality.

\smallskip
For a class of strictly convex internal energies in $d\ge2$,
McCann's pressure condition characterizes nonexpansiveness of every
proximal step. For critical fast diffusion in $d\ge3$, two successive steps
preserve distance precisely for pairs of translates or pairs in the
Barenblatt family, including point masses. Penalization gives
nonexpansiveness of projection onto densities bounded by one in every
dimension, answering a question of De Philippis, M\'esz\'aros,
Santambrogio and Velichkov. We also characterize nonexpansive projections under variable density
and relative-entropy constraints and determine, for Gaussian interactions under a density
bound on a ball, the optimal Lipschitz constant over a range of step sizes.

\smallskip
For entropy with a finite $\kappa$-strongly convex reference potential,
$\kappa>0$, the proximal step propagates the Poincar\'e, Talagrand and
log-Sobolev inequalities by the same sharp law. An input bound $C$
yields the optimal universal bound $L$ on the corresponding output
constant at step size $\tau$, where
\[
 (1+\kappa\tau)L=\tau+\sqrt{CL}.
\]
On a closed connected Ricci-flat manifold, nonexpansiveness of the
entropy proximal map along a sequence of step sizes tending to zero
characterizes flatness.

\par\medskip
\noindent\textit{Keywords.} \manuscriptkeywords.
\par\smallskip
\noindent\textit{2020 Mathematics Subject Classification.}
Primary 49Q22. Secondary 49K40, 47H09, 60E15, 35K55.
\end{abstract}
\begin{NoHyper}
\renewcommand{\thefootnote}{}
\interfootnotelinepenalty=10000
\footnotetext{This research is supported in part by the National Key R\&D Program of
China (2021YFA1000900, 2021YFA1002200), the National Natural Science
Foundation of China (12201596, 12671067), the Shandong Provincial Natural
Science Foundation (ZR2025QB05), and the Taishan Scholars Program of
Shandong Province (tsqn202408059).}
\end{NoHyper}
\setcounter{tocdepth}{2} 
\clearpage
\tableofcontents
\clearpage
\part{Coupled Bochner Calculus}

\section{Introduction}\label{sec:intro}
How does the quadratic Wasserstein distance change when both measures
vary? Along two prescribed smooth curves of measures, their velocities
and the optimal map between them determine the first derivative of the
squared distance. A second differentiation requires knowing how the map
itself varies. In moment reconstruction and in a Wasserstein proximal
step, the source measure itself minimizes a functional involving transport
to the target. When the target moves, we must therefore determine
both the velocity of the minimizing measure and the variation of the map. We show that the second
variation of the functional with the target held fixed gives a way to
estimate how fast the minimizer moves when the target varies.

Let $\rho$ and $\nu$ belong to $\Ptwo(\R^d)$, the space of probability
measures with finite second moment equipped with the quadratic Wasserstein
distance $W_2$~\cite{Villani}. Suppose for now that both measures have
smooth positive densities and that the optimal map $T=\nabla\Phi$ from
$\rho$ to $\nu$ is a smooth diffeomorphism with $H=\D^2\Phi\succ0$.
If $g\rho$ and $k\nu$ are the time derivatives of the two measures,
differentiating $T_\#\rho=\nu$ gives
\[
 B\dot\Phi=k\circ T-g,\qquad B=-\diver_\rho(H^{-1}\nabla),
 \qquad \diver_\rho Z=r^{-1}\diver(rZ)\quad(\rho=r\dx).
\]
Caja-Lopez, Delgadino and Kitagawa justify this linearization for
smooth variations on uniformly convex
domains~\cite[Theorem~1.5]{CajaDelgadinoKitagawa}. For two prescribed
curves the right-hand side is known. If the source minimizes a
functional depending on the target, its density derivative $g\rho$
must be found along with $\dot\Phi$ by differentiating the
Euler--Lagrange equation.

A density derivative $g\rho$ determines a constraint on the source
velocity, namely $-\diver_\rho Z=g$. Among fields satisfying this
constraint the Wasserstein velocity $v$ minimizes $\int|Z|^2\dd\rho$.
The operator $B$ suggests a second choice. For $g\in L^2(\rho)$
of mean zero, the solution of $Bu=g$ gives the field $X=H^{-1}\nabla u$
minimizing $\int H[Z,Z]\dd\rho$ under the same constraint
(\cref{sec:compatible}). Its pushforward $Y$ satisfies
$Y\circ T=HX=\nabla u$. Since $-(\diver_\nu Y)\circ T=-\diver_\rho X=g$,
applying the Bochner identity to $X$ at the source and to $Y$ at the
target gives two expressions for the same integral $\int g^2\dd\rho$.
We used this observation in our proof of Poincar\'e
interpolation~\cite{HanZhu}.

To deduce a bound from the source identity, we need the sign of
$\tr((\D X)^2)$. This trace could be negative for a general velocity
field, and $X$ need not be a Euclidean gradient. Differentiating $HX=\nabla u$ nevertheless shows that $H\D X$ is
symmetric. Consequently $\D X$ is similar to a symmetric matrix, and
\begin{equation}\label{eq:intro-symmetry}
 H\D X=\D^2u-\D^3\Phi[X,\cdot,\cdot],\qquad
 \tr((\D X)^2)=\|H^{1/2}\D XH^{-1/2}\|_{\HS}^2 .
\end{equation}
The target identity has a corresponding square
(see~\cite[Proposition~4.8]{KolesnikovLivshytsRotem}), and the average
of the two identities gives the integrated Bochner formula for
Kolesnikov's Hessian metric, with weighted Laplacian
$-B$~\cite{KolesnikovHessian}. Suppose $\rho$ and $\nu$ are respectively
$\kappa_0$- and $\kappa_1$-strongly log-concave, with $\kappa_i>0$.
By retaining the two identities separately, we bound
$\|X\|_{L^2(\rho)}$ by $\kappa_0^{-1/2}\|g\|_{L^2(\rho)}$ and
$\|Y\|_{L^2(\nu)}$ by $\kappa_1^{-1/2}\|g\|_{L^2(\rho)}$
as in \eqref{eq:two-flux-bounds}. Together these estimates give a
spectral gap of at least $\sqrt{\kappa_0\kappa_1}$ for $B$
(\cref{mm:prop:weighted}).

The link between moment reconstruction and a proximal step is easiest
to see by expanding the cost $|x-y|^2$. If $M_2(\rho)=\int |x|^2\dd\rho$,
the maximal correlation of $\rho$ and $\nu$ is
\[
 \mathcal T(\rho,\nu)
 :=\sup_{\pi\in\Pi(\rho,\nu)}\int x\cdot y\dd\pi
 =\tfrac12\bigl(M_2(\rho)+M_2(\nu)-W_2^2(\rho,\nu)\bigr).
\]
Santambrogio~\cite{SantambrogioMoment} reconstructs a measure from
its moment measure $\nu$ by minimizing $\Ent(\rho)+\mathcal T(\rho,\nu)$,
where $\Ent(r\dx)=\int r\log r\dx$. In the Jordan--Kinderlehrer--Otto
scheme~\cite{JKO}, the measure $J_\tau^{\cE}\nu$ minimizes
$\cE(\rho)+W_2^2(\rho,\nu)/(2\tau)$, so $\mathcal T$ enters with
coefficient $-1/\tau$. Up to a constant independent of $\rho$, both
functionals thus have the form
$\mathcal F_\nu(\rho)=\cE(\rho)+cM_2(\rho)/2+a\mathcal T(\rho,\nu)$,
with $(a,c)=(1,0)$ for reconstruction and $(-1/\tau,1/\tau)$ for
a proximal step.

Consider first a variation of the source with the target fixed.
Under the hypotheses of \cref{thm:hessian-transfer}, the Wasserstein
Hessian of $\mathcal T$ in the direction $v$ is
$\int H[v-X,v-X]\dd\rho$. At a smooth critical point of the full
functional, the Euler--Lagrange equation
$\nabla q_\rho+cx+aT=0$ relates the first variation $q_\rho$ of the
energy to the same transport map. Differentiating this equation in
space gives $\D^2q_\rho+cI=-aH$. When we add the energy and
second-moment contributions to the Hessian of $a\mathcal T$, the terms
involving $v-X$ cancel and we obtain
\[
 \operatorname{Hess}_W\mathcal F_\nu(v,v)
 =\mathscr Q_\rho(X)+c\int |X|^2\dd\rho
 =\int g\,\cK_\rho g\dd\rho-a\int H[X,X]\dd\rho .
\]
Here $\cK_\rho g$ is the derivative of $q_\rho$ in the density
direction $g\rho$, and $\mathscr Q_\rho(X)=\int g\,\cK_\rho g\dd\rho+
\int\D^2q_\rho[X,X]\dd\rho$. As $g$ has mean zero, the additive
constant in $\cK_\rho g$ drops out. Testing $Bu=g$ against $u$ gives
$\int H[X,X]\dd\rho=\int gB^{-1}g\dd\rho$, so this Hessian is the
quadratic form of $\cK_\rho-aB^{-1}$ on mean-zero variations.

Now vary the target with velocity $b$ and let $g\rho$ be the
derivative along a differentiable curve of corresponding critical
points. Differentiating the scalar Euler--Lagrange equation
$q_\rho+c|x|^2/2+a\Phi=\mathrm{const}$ along this curve and
integrating against $g\rho$ removes the additive constant. The
linearized transport equation then gives
\[
 \mathscr Q_\rho(X)+c\int |X|^2\dd\rho
       =-a\int X\cdot(b\circ T)\dd\rho .
\]
The left-hand side is exactly the form we computed with the target
fixed. If it is at least $\lambda\|X\|_{L^2(\rho)}^2$ for some $\lambda>0$,
Cauchy--Schwarz and the minimality of $v$ give
$\|v\|_{L^2(\rho)}\le\|X\|_{L^2(\rho)}\le |a|\,\|b\|_{L^2(\nu)}/\lambda$.

For a proximal step, substituting $(a,c)=(-1/\tau,1/\tau)$ and
completing the square in this pairing expresses the difference of the
squared input and output speeds as two squared norms plus
$2\tau\mathscr Q_\rho(X)$ (\cref{prop:primal-contraction}). When
$\mathscr Q_\rho(X)\ge0$, the output speed is at most the input speed,
and equality forces all three terms to vanish. At critical fast
diffusion these equality conditions determine which pairs of measures
preserve their distance for two successive steps.

\subsection{Time convexity and sharp contraction along the heat flow}
Let $\mu_t$ and $\nu_t$ be the heat evolutions of two measures in
$\Ptwo(\R^d)$. To compare their distances at times $t$ and $t+h$, choose
an optimal coupling $(X,Y)$ at time $t$ and let $Z\sim\N(0,hI)$ be
independent of the pair. Convolution with the law of $Z$ advances each
heat flow by time $h$, so $X+Z$ and $Y+Z$ have laws $\mu_{t+h}$ and
$\nu_{t+h}$. Their difference is still $X-Y$. Taking the joint law of
$(X+Z,Y+Z)$ as a competitor in the infimum at time $t+h$, we obtain
\[
 W_2^2(\mu_{t+h},\nu_{t+h})
 \le \E\big|(X+Z)-(Y+Z)\big|^2
 =\E|X-Y|^2=W_2^2(\mu_t,\nu_t).
\]
This argument shows that the squared distance is nonincreasing.
When both initial measures are log-concave, we prove that over any
two successive intervals of equal length it decreases at least as
much over the first as over the second.

\begin{theorem}\label[theorem]{thm:intro-heat-convexity}
Let $\mu,\nu\in\Ptwo(\R^d)$ be log-concave, and let $\gamma_t=\N(0,tI)$
for $t>0$, with $\gamma_0=\delta_0$. Then the squared distance
$D(t)=W_2^2(\mu*\gamma_t,\nu*\gamma_t)$ is convex on $[0,\infty)$.
\end{theorem}

The comparison of losses requires following the optimal map itself.
For the regular initial densities used in the proof, let $T_t$ transport
$\mu_t$ optimally to $\nu_t$ and put $H_t=\D T_t$. Substituting the
heat velocities in the first variation and integrating by parts at both
ends gives, for almost every $t>0$,
\[
 -D'(t)=\int\tr(H_t+H_t^{-1}-2I)\dd\mu_t .
\]
This trace integral is the sum of the Dirichlet energies of the
components of $T_t-\Id$, with coefficient $H_t^{-1}$ and measure
$\mu_t$. By \eqref{heat:eq:quantitative}, the lower symmetric second derivative
of $D$ is bounded below by the sum of those same energies applied to
the entries of $H_t$, and is therefore nonnegative. This gives convexity for the regularized
measures, and the convexity inequality passes to the limit under
$W_2$ convergence. Without log-concavity the squared distance can have zero right derivative
at the origin and still fall strictly below its initial value for every
$t>0$, as an example with two-point laws on the line shows (\cref{heat:rem:two-point}).

A lower bound on the curvature of both measures gives a quantitative
estimate from the first variation alone. Convolution commutes with translations, so a pair of translates stays
at the same distance throughout the flow. We therefore center the
measures before asking for a contraction factor smaller than one.
Suppose $\mu$ and $\nu$ are centered with strong log-concavity constants
$\kappa_0,\kappa_1>0$, respectively, and write $a_i=\kappa_i^{-1}$.
After time $t$ their curvature bounds are $(a_0+t)^{-1}$ and $(a_1+t)^{-1}$.
The weighted Poincar\'e inequality of \cref{mm:prop:weighted}, applied
to the components of the displacement in the formula for $-D'(t)$,
then gives $-D'(t)\ge D(t)/\sqrt{(a_0+t)(a_1+t)}$. Integration yields
\[
 W_2(\mu_t,\nu_t)
 \le \frac{\sqrt{a_0}+\sqrt{a_1}}{\sqrt{a_0+t}+\sqrt{a_1+t}}\,
 W_2(\mu,\nu),\qquad t>0 .
\]
The factor is optimal for each $t$ and each pair of curvature bounds.
For distinct centered measures, equality requires a common Gaussian
factor with variances $a_0$ and $a_1$ and identical marginals on its
orthogonal complement (\cref{cor:heat-contraction}).

For Gaussian convolution at large times, Chen and
Niles-Weed~\cite{ChenNilesWeed} determined the asymptotics of
$W_2(\mu_t,\nu_t)$ from the moments of the initial measures under
suitable tail bounds. At small times, Fathi, Goldman and
Tsodyks~\cite[Theorem~5.10]{FathiGoldmanTsodyks} used the loss in
$W_2^2$ between a standard Gaussian and a strongly log-concave measure
to bound the squared distance of that measure to the nearest translate
of the Gaussian. 

Log-concavity also gives time convexity for certain information
quantities. For a log-concave random vector $X$ and an independent
standard Gaussian $G$, Wibisono and Jog~\cite[Theorem~1]{WibisonoJog}
showed
that the mutual information between $X$ and $X+\sqrt t\,G$ is convex
in $t>0$. A closer parallel to \cref{thm:intro-heat-convexity} arises
when both arguments of the relative entropy evolve under the heat
flow. Assuming Conjecture~4.2 of their preprint, Klartag and
Ordentlich~\cite[preprint, Theorem~1.7]{KlartagOrdentlich} proved
time convexity when the initial reference measure is log-concave.

\subsection{Inverse stability of moment measures}
In the moment measure problem the target $\nu$ is prescribed, and
one seeks a convex potential $\psi$ whose gradient transports the
density $e^{-\psi}$ to $\nu$. Thus the same potential determines both
the source and the map,
\begin{equation}\label{mm:eq:moment-definition}
 \rho=e^{-\psi}\dx,\qquad \Moment\rho=(\nabla\psi)_\#\rho=\nu .
\end{equation}
The moment measure is unchanged if the source and its potential are
translated together. Cordero-Erausquin and Klartag~\cite{CorderoKlartag}
proved that every centered probability measure of finite first moment
and full affine span has an essentially continuous convex potential in
this representation, and that translation is the only ambiguity.
Here essential continuity means lower semicontinuity with a discontinuity
set of zero $(d-1)$-dimensional Hausdorff measure. We fix the translation
by centering the source and denote the resulting measure by $\Reconstruct\nu$.

Santambrogio~\cite{SantambrogioMoment} recovers this source by
minimizing the sum of entropy and maximal correlation. We ask how far
the minimizer moves when the prescribed moment measure changes. Our bound depends on the strong convexity constants of the two
reconstructed potentials. As usual, $\kappa$-strong convexity means
that $\psi-\kappa|\cdot|^2/2$ is convex, allowing the value $+\infty$.

\begin{theorem}\label[theorem]{mm:thm:main}
Let $\nu_0,\nu_1\in\Ptwo(\R^d)$ be centered, neither supported
on a hyperplane. Write their centered reconstructions as
$\rho_i=\Reconstruct\nu_i=e^{-\psi_i}\dx$, where $\psi_i$ is
essentially continuous and $\kappa_i$-strongly convex, $\kappa_i>0$.
Then
\begin{equation}\label{mm:eq:main-estimate}
 W_2(\rho_0,\rho_1)
 \le \frac{1}{\sqrt{\kappa_0\kappa_1}}\,W_2(\nu_0,\nu_1).
\end{equation}
The constant $(\kappa_0\kappa_1)^{-1/2}$ is optimal for every
prescribed pair $\kappa_0,\kappa_1>0$, in every dimension $d\ge1$.
\end{theorem}

For a Gaussian source $\N(0,\kappa_i^{-1}I)$ the moment map is
$x\mapsto\kappa_i x$, and its image is $\N(0,\kappa_i I)$.
Two such sources attain equality for $\kappa_0\ne\kappa_1$. For a
common lower bound $\kappa$, taking precisions $\kappa$ and
$\kappa+\varepsilon$ and letting $\varepsilon\downarrow0$ gives
sharpness in the limit. In fact, \cref{mm:thm:equality} shows that
two distinct centered sources attain equality precisely when
$\kappa_0\ne\kappa_1$ and they differ only in the variance of a
Gaussian factor,
\[
 \rho_i=\N_E(0,\kappa_i^{-1}I_E)\otimes\sigma,
 \qquad i=0,1 ,
\]
where $E$ is a nonzero subspace and $\sigma$ is a centered probability
measure on $E^\perp$. For related Gaussian splitting results see
\cite{HanZhu,DePhilippisFigalli,CourtadeFathi}.

To prove the estimate, we first compare the sources by their own
optimal map $T$, writing $H=\D T$. Under the regularity assumptions of
\eqref{mm:eq:regular}, applying the transport-weighted Poincar\'e
inequality to the components of $T-\Id$ gives
\begin{equation}\label{eq:intro-moment-spectral}
 \sqrt{\kappa_0\kappa_1}\,W_2^2(\rho_0,\rho_1)
 \le\int\tr(H+H^{-1}-2I)\dd\rho_0 .
\end{equation}
The trace integral is the one that governs dissipation
along the heat flow. Here we must bound it in terms of the prescribed
moment measures. To this end we take an optimal coupling of $\nu_0$
and $\nu_1$ and pull it back by the two moment maps. Although the
resulting coupling of $\rho_0$ and $\rho_1$ need not be optimal, its
excess cost enters the four-point identity of \cref{mm:prop:finite-defect}
as a nonnegative remainder. The other terms in the deficit of the
squared estimate \eqref{mm:eq:main-estimate} are nonnegative as well
and include a positive multiple of the deficit in
\eqref{eq:intro-moment-spectral}. Thus equality in
inverse stability forces equality in the spectral estimate as well.

For the smooth centered perturbations of \cref{prop:moment-direction},
differentiating $(\nabla\psi)_\#\rho$ and projecting the resulting
velocity onto target gradients gives the Wasserstein velocity of the
moment measure.
The source speed is at most $\sqrt{C_P(\rho)/\kappa}$ times the norm
of this projected velocity,
where $C_P(\rho)$ is the source Poincar\'e constant and $\kappa$ is its
lower curvature bound. For the underlying Hessian geometry see
\cite{KlartagMoment,KolesnikovHessian}.

Bonnet and Rubinstein~\cite[Theorem~1.3]{BonnetRubinstein} estimated
the $L^1$ distance between the reconstructed densities, after a suitable
translation of one of them, in terms of the $W_1$ distance between the
prescribed moment measures.
Can one estimate the curvature of the unknown source potentials from
the prescribed moment measures? Machado and Ramos~\cite[Theorems~1.1, 5.6 and~5.12]{MachadoRamos}
used quantitative Pr\'ekopa--Leindler inequalities to study the
Brascamp--Lieb deficit and prove $W_2$-H\"older stability of moment
reconstruction modulo translations. Their estimates hold for a
nondegenerate class of compactly supported targets and for targets
with a common upper bound on the Hessians of their potentials, with
any exponent $(p-2)/(6p-4)$, $p>2$, in each class. In the latter
setting their curvature estimate~\cite[Theorem~5.10]{MachadoRamos}
provides exactly the hypothesis needed here. Indeed, if the centered
targets $\nu_i=Z_i^{-1}e^{-W_i}\dy\in\Ptwo(\R^d)$ have finite
potentials with $\D^2 W_i\preceq\Lambda_i I$ in distributions on
$\R^d$, $\Lambda_i>0$, then the reconstructed potentials are
$\Lambda_i^{-1}$-strongly convex (see also
\cite[Theorem~3.4]{KolesnikovKosov}). Substitution in
\eqref{mm:eq:main-estimate} gives
\begin{equation}\label{eq:intro-moment-inverse}
 W_2(\Reconstruct\nu_0,\Reconstruct\nu_1)
 \le\sqrt{\Lambda_0\Lambda_1}\,W_2(\nu_0,\nu_1),
\end{equation}
with optimal coefficient for every prescribed pair of upper bounds
(\cref{mm:cor:inverse}).

Another way to fix the translation is to add a confinement to
Santambrogio's functional. Let $V$ be finite and $\alpha$-strongly
convex, $\alpha>0$, and let $m_V$ be the probability measure
proportional to $e^{-V}\dx$. For every $\nu\in\Ptwo(\R^d)$,
minimizing $\Ent_{m_V}(\rho)+\mathcal T(\rho,\nu)$ over $\Ptwo(\R^d)$
gives a unique source $\rho=\Reconstruct_V\nu$, whose density is
proportional to $e^{-V-\Phi}$ for a convex potential $\Phi$ transporting
$\rho$ to $\nu$ by its gradient. For quadratic $V$, Delalande and
Farinelli~\cite[Theorem~3.1 and Remark~2.3]{DelalandeFarinelli}
proved $1/2$-H\"older continuity on $W_2$-bounded sets of targets and
also discussed strongly convex regularizations. We obtain a global
$1/\alpha$-Lipschitz estimate for every such confinement. If the two
targets have the same barycenter, the coefficient improves to
$1/(2\alpha)$ in dimension one and $1/(\sqrt{2}\,\alpha)$ in higher
dimensions. Each coefficient is sharp over the corresponding class of
confinements (\cref{mm:thm:regularized}).

\subsection{Stability of Brenier potentials}
Suppose now that the source $\rho$ is fixed. As the target varies,
we can compare the Brenier potentials themselves once we fix their
additive constants. Let $\Omega\subseteq\R^d$ be bounded, open and
convex, and suppose
\begin{equation}\label{eq:source-bounds}
 \diff\rho=r\one_\Omega\dx,\qquad
 0<m\le r\le M<\infty\ \hbox{a.e.},\qquad \int_\Omega r\dx=1.
\end{equation}
For a compactly supported target $\mu\in\Pc(\R^d)$, write $\Phi_\mu$
for the convex potential with $(\nabla\Phi_\mu)_\#\rho=\mu$ and
$\int\Phi_\mu\dd\rho=0$. We seek the least constant in the estimate
of the $L^2(\rho)$ distance between two potentials by the $W_2$
distance between their targets, namely
\begin{equation}\label{eq:S-definition}
 S(\rho)=\sup_{\substack{\mu,\nu\in\Pc(\R^d)\\\mu\ne\nu}}
 \frac{\norm{\Phi_\mu-\Phi_\nu}_{L^2(\rho)}}{W_2(\mu,\nu)}
 \in[0,\infty].
\end{equation}
The linearized equation suggests how to compute this constant. When the optimal map varies smoothly with the target,
the derivative of its normalized potential,
tested against a mean-zero function $g$, is the pairing of the target
velocity pulled back by the transport map with the source flux obtained
from a Neumann problem
(\eqref{eq:stability-adjoint}). To define this flux, let
$\cC(\Omega)$ consist of the functions smooth on a neighborhood of
$\overline\Omega$ with positive definite Hessian there. For
$\Phi\in\cC(\Omega)$ and $g\in L^2_0(\rho)$, write $u_{\Phi,g}$
for the weak solution of
\begin{equation}\label{eq:intro-Neumann}
 \int_\Omega (\D^2{\Phi})^{-1}\nabla u_{\Phi,g}\cdot\nabla h\dd\rho
       =\int_\Omega gh\dd\rho
       \qquad(h\in H^1(\Omega)).
\end{equation}
with $\int u_{\Phi,g}\dd\rho=0$. The Hessian is uniformly elliptic on
$\overline\Omega$ for each fixed $\Phi$, so the Poincar\'e inequality
gives a unique solution with this normalization. The corresponding flux is $(\D^2\Phi)^{-1}\nabla u_{\Phi,g}$.

\begin{theorem}\label{thm:stability}
Under \eqref{eq:source-bounds}, it holds that 
\begin{equation}\label{eq:intro-exact-S}
 S(\rho)=
 \sup_{\Phi\in\cC(\Omega)}\ 
 \sup_{\substack{g\in L^2_0(\rho)\\\norm g_{L^2(\rho)}=1}}
 \norm{(\D^2{\Phi})^{-1}\nabla u_{\Phi,g}}_{L^2(\rho)}.
\end{equation}
If $S(\rho)<\infty$, centered
potentials for all targets in $\Ptwo(\R^d)$ belong to $L^2(\rho)$ and
satisfy the same stability estimate. The value of $S(\rho)$ is unchanged
if both targets are required to be supported in a fixed compact set
with nonempty interior.
\end{theorem}

For smooth densities on smooth uniformly convex domains, Letrouit
and Manole~\cite[Theorem~3.2 and Lemma~3.6]{LetrouitManole} used
the adjoint pairing to deduce stability from a uniform bound on the
Neumann fluxes. For the converse, choose a smooth perturbation $f$
and use $\Phi+tf$ as the potential of a curve of targets. The
derivative of the normalized potential is then simply
$f-\int f\dd\rho$. The speed of the target at $t=0$ is bounded by
the norm of the orthogonal projection of $\nabla f$ onto the space
of fields $(\D^2\Phi)^{-1}\nabla h$. Applying the global Lipschitz
estimate to this curve therefore bounds the standard deviation of $f$
by $S(\rho)$ times the norm of its projected gradient. By duality,
this projected Poincar\'e inequality gives the required bound on the
flux with the same constant (\cref{prop:stability-equivalences}).
The argument only requires the density bounds in \eqref{eq:source-bounds}.

For a $\kappa$-strongly log-concave source, Letrouit and
Manole~\cite[Theorem~1.1]{LetrouitManole} proved the sharp universal
bound $S(\rho)\le\kappa^{-1/2}$, which Ruan~\cite[Theorem~1.1]{Ruan}
also obtained by a semi-discrete argument for finite strongly convex
source potentials. One can see why the norm is $L^2(\rho)$ by varying the energy. For a
smooth positive density $r$, the second derivative of $\int U(r)\dx$
along $(1+tg)\rho$ at $t=0$ is $\int rU''(r)g^2\dd\rho$.
Under the pressure and curvature assumptions of \cref{thm:weighted-brenier},
the dual estimate for the difference of potentials, modulo constants,
therefore carries the reciprocal weight $1/(rU''(r))$. Entropy gives
weight one. For $\rho\in\Ptwo(\R^d)$ with density proportional to
$W^{-\beta}$ on a convex domain, where $W$ is positive and strongly
convex and $\beta\ge d\ge2$, the symmetry in \eqref{eq:intro-symmetry}
gives instead a sharp estimate in $L^2(W^{-1}\rho)$
(\cref{thm:inverse-power-stability}).

An estimate for the potentials can in turn be used to control the maps.
For a uniform source on a convex body in $\R^d$, $d\ge2$, and targets
in a fixed compact set, OpenAI~\cite{OpenAIBrenier} proved that the
$L^2(\rho)$ distance between the maps is at most $C W_2(\mu,\nu)^{1/3}$. Ruan~\cite[Theorems~1.1 and~1.2]{Ruan} obtained
the same bound for finite strongly convex source potentials and targets
in a common ball. Even a Gaussian source in $d\ge2$ admits no
larger exponent. Wang and Wang~\cite[Theorem~1 and Corollary~5]{WangWang}
obtained the same exponent for isotropic log-concave sources using a
weighted estimate for the potentials. We compare their weight with
$1/(rU''(r))$ in \cref{sec:stability-curvature}. For a bounded,
compactly supported source density and targets in a fixed compact set,
M\'erigot~\cite[Theorem~3.1]{MerigotSharp} obtained the bound
$C W_1(\mu,\nu)^{1/4}$ for the maps in $L^2(\rho)$, with optimal
exponent in $d\ge2$. For earlier
H\"older estimates for maps with a fixed source, see
\cite{MerigotDelalandeChazal,DelalandeMerigot}, and for time
regularity as the target varies along a curve, see~\cite{Gigli}.

For a uniform source, whether a finite Lipschitz constant exists can
be read from the support. If $d\ge2$ and $\rho_\Omega$ is uniform
on a bounded connected open set $\Omega$ with $|\partial\Omega|=0$,
then \cref{thm:nonconvex-support} gives
\begin{equation}\label{eq:intro-essential-convexity}
 S(\rho_\Omega)<\infty\quad\Longleftrightarrow\quad
 \overline\Omega\text{ is convex}.
\end{equation}
Indeed, when $\overline\Omega$ is nonconvex, targets with at most
$2d$ atoms suffice to make the ratio in \eqref{eq:S-definition}
arbitrarily large. For nonuniform sources, convexity of the support alone cannot ensure
stability. To see why, consider a product of a bounded convex domain
and an interval. If the source density is bounded above and
below by positive constants and its marginal on the interval has a
positive $C^1$ density up to the endpoints, then $S(\rho)<\infty$
forces the conditional laws to have locally square-integrable
Wasserstein speed as the interval coordinate varies
(\cref{prop:conditional-fibres}). On a fixed square one can violate
this condition even with a positive continuous density arbitrarily
close to uniform in $L^\infty$ by choosing the conditional barycenters
to have a derivative that fails to be locally square-integrable
(\cref{prop:rough-counterexample}). The $W_1^{1/2}$ estimate of
Letrouit and M\'erigot~\cite[Theorem~1.7]{LetrouitMerigot} for
potentials with targets in a fixed compact set still applies to
these sources. Kitagawa, Letrouit and M\'erigot~\cite{KitagawaLetrouitMerigot}
proved the same estimate in $L^2(\rho)$ when the source density is
bounded above and below by positive constants on a Riemannian John
domain and the targets lie in a fixed compact set.
Under analogous density and support assumptions on bounded John domains
in finite-dimensional $\mathrm{RCD}$ spaces, we used heat kernel
regularization to obtain the same $W_1^{1/2}$ rate for normalized
Kantorovich potentials, first in $L^1(\rho)$~\cite[Theorem~1.1]{HanZhuMetricMeasure}
and then in $L^2(\rho)$~\cite[Theorem~A]{HanZhuHeatKernel}.

\subsection{The pressure criterion for proximal contraction}
An internal energy with pressure $p$ generates the diffusion
$\partial_t r=\Delta p(r)$ as its Wasserstein gradient flow
\cite{Otto,CMVContraction}. McCann's condition on $p$ makes the
energy displacement convex~\cite{McCann}, and Bolley and
Carrillo~\cite[Theorem~1 and Corollary~2]{BolleyCarrillo} showed that
the same condition characterizes $W_2$ contraction for smooth positive
solutions on a ball with Neumann boundary conditions. On $\R^d$ they
proved sufficiency for bounded initial densities with finite second
moment. Does this equivalence between the pressure condition and
$W_2$ nonexpansiveness persist at every finite step of the scheme?

In a Hilbert space every proper lower semicontinuous convex
functional has a nonexpansive proximal map~\cite{Moreau}.
Cavagnari, Savar\'e and Sodini~\cite[Definition~5.1 and
Proposition~5.4]{CavagnariSavareSodini} used this fact for
functionals of a probability law whose lift to square-integrable
random variables is convex. For such a functional the Wasserstein step is the law of the Hilbert-space
step, but convexity of the lift requires convexity along the interpolation
induced by every coupling (total convexity). For Boltzmann entropy,
coupling a standard Gaussian vector $X$ with $-X$ gives endpoints of
finite entropy but a midpoint of infinite entropy, since its law is
$\delta_0$.

Displacement convexity only tests optimal couplings, and is too weak
to imply nonexpansiveness. In \cref{ex:resolvent-scope} a step of a
displacement convex energy sends arbitrarily close inputs to unique
minimizers a unit distance apart.

Stronger convexity assumptions give related estimates. For proper
coercive lower semicontinuous energies
convex along generalized geodesics, Carlen and
Craig~\cite[Theorem~1.3]{CarlenCraig} obtained a contraction estimate
by adding to the squared distance the sum of the squared metric
slopes at its two endpoints, multiplied by $\tau^2/2$.
Di Marino, Farinelli and Naldi~\cite[Theorems~3.8 and~3.17]{DiMarino}
obtained an estimate in $W_2$ itself. Assuming $\lambda$-convexity
along two-base generalized geodesics, they glue transport plans
along an optimal coupling of the inputs and obtain the Lipschitz
bound $(1+\lambda\tau)^{-1}$ when $1+\lambda\tau>0$. This hypothesis
can also fail for entropy, since in their planar Gaussian example two
measures of finite entropy have a two-base interpolation whose
midpoint is supported on a line~\cite[Example~3.18]{DiMarino}.

For internal energies we use McCann's condition directly. It ensures
that the second variation is nonnegative along the compatible fields
selected by the optimal maps from the minimizers to their inputs.

Let $d\ge2$ and let $U:[0,\infty)\to\R$ be strictly convex, with
$U(0)=0$ and $U\in C^2((0,\infty))$. To control the negative part of
the energy at low densities, assume that $U(r)\ge-Cr^q$ for
$0<r\le1$, where $C>0$ and $d/(d+2)<q<1$. For a measure $\rho$ with singular
part $\rho^s$, the relaxed internal energy is
\begin{equation}\label{eq:internal-relaxation}
 \mathcal U(\rho)=\int U(r)\dx+U_\infty\rho^s(\R^d),\qquad
 \rho=r\dx+\rho^s,\qquad
 U_\infty=\lim_{r\to\infty}\frac{U(r)}r.
\end{equation}
The recession slope $U_\infty$ determines whether a measure with a
singular part can have finite energy. It is infinite for
superlinear $U$ and finite for the fast-diffusion powers. Write
$J_\tau^U$ for the proximal map of $\mathcal U$ wherever the
minimizer is unique, and define the pressure and its coefficient by
\begin{equation}\label{eq:pressure-coefficient}
 p(r)=rU'(r)-U(r),\qquad
 \kappa_U(r)=rp'(r)-\frac{d-1}{d}p(r).
\end{equation}
Convexity and $U(0)=0$ imply $p\ge0$.

\begin{theorem}\label{thm:pressure-global}
For $d$ and $U$ as above, the following are equivalent.
\begin{enumerate}[label=\textup{(\roman*)},leftmargin=2em]
\item $\kappa_U(r)\ge0$ for $r>0$;
\item for any $\tau>0$ and $\mu\in\Ptwo$, the relaxed proximal problem has
 a unique minimizer, and
 \begin{equation}\label{eq:global-internal-contraction}
  W_2(J_\tau^U\mu,J_\tau^U\nu)\le W_2(\mu,\nu)
       \qquad(\mu,\nu\in\Ptwo);
 \end{equation}
\item there is a sequence $\tau_j\downarrow0$ such that
 \[
  W_2(J_{\tau_j}^U\beta_R,J_{\tau_j}^U\beta_S)
       \le W_2(\beta_R,\beta_S)\qquad(R,S>0),
  \quad \beta_R=\frac{\one_{B_R}\dx}{|B_R|}.
 \]
\end{enumerate}
The minimizers in \textup{(iii)} exist and are unique even when
\textup{(i)} fails.
If $\kappa_U$ is negative at some density, each sufficiently small step
expands the distance between a suitable pair of uniform measures on concentric balls.
\end{theorem}

For sufficiency, let the input vary along a smooth curve $\eta_t$
with minimal velocity $v$, and suppose that its minimizers
$\rho_t=r_t\dx$ form a smooth curve of positive densities with minimal
velocity $w$. At a fixed time write $T$ for the optimal map from
$\rho$ to $\eta$. Solving $Bu=\dot r/r$ gives a field
$X=H^{-1}\nabla u$ representing the same density variation as $w$,
so \cref{prop:primal-contraction} gives
\[
 \|v\|_\eta^2-\|w\|_\rho^2
 =\|v\circ T-X\|_\rho^2+\|w-X\|_\rho^2+2\tau\mathscr Q_\rho(X),
\]
where $\|\cdot\|_\rho$ denotes the $L^2(\rho)$ norm. Under the pressure
condition all three terms on the right are nonnegative. Hence the
curve of minimizers is at most as long as the curve of inputs.
Approximating a geodesic between two inputs by smooth curves and
passing to the relaxed energy gives \eqref{eq:global-internal-contraction}.

For necessity we use concentric balls, as in the continuous-time
argument of Bolley and Carrillo. In radius quantiles, $W_2$ between
radial measures is the $L^2(0,1)$ distance and the inputs $\beta_R$
trace a ray $Rh$, where $h$ is the radius quantile of $\beta_1$.
Write $q_R$ for the radius quantile of $J_\tau^U\beta_R$ and
$e_\tau(R)$ for the minimum of its proximal objective. Differentiating
with respect to $R$ gives
$\tau e_\tau'(R)=\langle Rh-q_R,h\rangle_{L^2(0,1)}$.
Condition~\textup{(iii)} says that
$\|q_R-q_S\|_2\le |R-S|\,\|h\|_2$, so this derivative is nondecreasing
at each of the indicated step sizes. Along this sequence the convex functions $e_\tau(R)$ converge to
$F(R)=|B_R|U(1/|B_R|)$. The pressure condition follows because $F''(R)$
has the sign of $\kappa_U(1/|B_R|)$, and the densities $1/|B_R|$
range over $(0,\infty)$.

For a spatially dependent energy $\int U(x,r)\dx$, one can both
translate and dilate a small region of constant density. Translation
varies its center $x$, and dilation varies the volume per unit mass
$s^d=1/r$. The energy per unit mass is therefore $s^dU(x,s^{-d})$.
Under the hypotheses of \cref{thm:local-volume-criterion} in dimension
$d\ge2$, joint convexity of $s^dU(x,s^{-d})-\lambda|x|^2/2$ in $(x,s)$
is equivalent to the bound
$\Lip_{W_2}(J_\tau^U)\le(1+\lambda\tau)^{-1}$ for every $\tau>0$,
where $\lambda\ge0$ is fixed.

\subsection{Rigidity at the critical fast-diffusion exponent}
If the two speeds agree under the pressure condition, each term in
their comparison must vanish. To find the fields for which $\mathscr Q_\rho(X)=0$, take a
smooth positive density $r$ and a smooth compactly supported
compatible field $X$ on $\R^d$. The second variation is
\begin{equation}\label{eq:intro-pressure-v6}
 \mathscr Q_\rho(X)
 =\int\left\{p(r)\left\|H^{1/2}\D XH^{-1/2}
         -\frac{\diver X}{d}I\right\|_{\HS}^2
       +\kappa_U(r)(\diver X)^2\right\}\dx .
\end{equation}
The square measures the trace-free part of the symmetric matrix in
\eqref{eq:intro-symmetry}, while the last term measures the change of
volume. At the critical power $U_*(r)=-dr^{1-1/d}$ the coefficient
$\kappa_{U_*}$ vanishes identically. Thus the integrand is zero exactly when $\D X$ is a scalar multiple of $I$, allowing infinitesimal
dilations as well as translations. This accounts for the one-step
equality cases in $d\ge3$. A step preserves the distance between two
inputs precisely when its two minimizers agree up to translation and
positive dilation (\cref{thm:critical-global-equality}).

What happens if the next step also preserves distance? Write
$J_\tau^*$ for the proximal map of the relaxed critical energy, for
which $U_\infty=0$ in \eqref{eq:internal-relaxation}. The following
result identifies the profiles of pairs that preserve distance twice
without being translates. We use the Barenblatt family
\begin{equation}\label{eq:intro-barenblatt-family}
 B=c_d(1+|x|^2)^{-d}\dx,\qquad
 B_{b,R}=(b+R\Id)_\#B\quad(R\ge0),
\end{equation}
where $c_d$ makes $B$ a probability measure. Thus $B_{b,0}=\delta_b$.

\Needspace{10\baselineskip}
\begin{theorem}\label{thm:two-step-rigidity}
Let $d\ge3$, $\tau_1,\tau_2>0$, and
$\mu_0,\mu_1\in\Ptwo(\R^d)$. Set
\[
 \rho_i=J_{\tau_1}^*\mu_i,\qquad
 \sigma_i=J_{\tau_2}^*\rho_i,\qquad i=0,1.
\]
Then
\begin{equation}\label{eq:two-step-equality}
 W_2(\sigma_0,\sigma_1)=W_2(\mu_0,\mu_1)
\end{equation}
if and only if at least one of the following holds:
\begin{enumerate}[label=\textup{(\roman*)},leftmargin=2em]
\item $\mu_1=(b+\Id)_\#\mu_0$ for some $b\in\R^d$;
\item $\mu_i=B_{b_i,R_i}$ for some $b_i\in\R^d$ and $R_i\ge0$, $i=0,1$.
\end{enumerate}
Both classes preserve their distance under every common finite sequence
of positive step sizes.
\end{theorem}

To see why a second equality determines the profile, suppose the two
inputs are not translates. The one-step result gives a homothety
between $\rho_0$ and $\rho_1$ and another between $\sigma_0$ and
$\sigma_1$. We can therefore transport $\sigma_0$ optimally to
$\rho_1$ in two ways, by applying the second-step Euler map before
or after the corresponding homothety. Since the source of both maps has a density, uniqueness of optimal transport makes the maps agree almost everywhere. As the first homothety has scale different from one, this forces
the gradient of the first variation at $\sigma_0=f\dx$ to be affine
with scalar linear part. But the first variation is
$-(d-1)f^{-1/d}$, so $f$ is a Barenblatt density. Its Euler maps carry
this family back to the inputs.

Carlen and Craig~\cite[Proposition~4.7]{CarlenCraig} proved that the
discrete flow preserves the Barenblatt family above the critical
exponent. Below it, whenever the profiles have finite second moment,
every positive step strictly increases the distance between centered
Barenblatt profiles of different scales (\cref{prop:power-dilations}).

\subsection{Projections and interactions under constraints}
A density bound in models of crowd motion~\cite{MauryCrowd} leads to
a question about Wasserstein projection. Let $K_1$ consist of the
probability measures on $\R^d$ whose densities are bounded by one, and write
$P_K\mu$ for a nearest point in $K$ to $\mu$ in $W_2$. De Philippis,
M\'esz\'aros, Santambrogio and Velichkov~\cite[Remark~5.1]{DPMSV}
asked whether
\[
 W_2(P_{K_1}\mu,P_{K_1}\nu)\le W_2(\mu,\nu)
 \qquad (\mu,\nu\in\Ptwo(\R^d)).
\]
Adve and M\'esz\'aros~\cite[Theorem~3.1, Corollaries~3.2--3.3 and
Remark~3.4]{AdveMeszaros} bounded the squared distance on the left
by the cost of a particular coupling of the inputs, obtained by
transporting an optimal coupling of the projections back to $\mu$
and $\nu$. Whenever the induced coupling is optimal, its cost is
$W_2^2(\mu,\nu)$ and the desired inequality follows. This holds on
the line by monotonicity, and in any dimension when one input is a
point mass. In their planar example, however, the
induced coupling can cost strictly more than $W_2^2(\mu,\nu)$.

In higher dimensions the pressure criterion gives the bound by
approximation. Consider the power energies $(p-1)^{-1}\int r^p\dx$,
$p>1$. Their proximal maps are nonexpansive by
\cref{thm:pressure-global}, and as $p\to\infty$ their unit-step
minimizers converge narrowly to the projection onto $K_1$. Applying
this convergence to both inputs and using lower semicontinuity of
$W_2$ gives the inequality for arbitrary $\mu$ and $\nu$.

Which other constraints have nonexpansive projections? We consider
a spatially varying density bound and a sublevel set of the relative
entropy. For a reference probability measure $m$, we write
$\Ent_m(\rho)=\int\log(\diff\rho/\diff m)\dd\rho$ if $\rho\ll m$
and set $\Ent_m(\rho)=+\infty$ otherwise.

\begin{theorem}\label{thm:intro-projections}
Let $d\ge1$. For each set $K$ below, every $\mu\in\Ptwo(\R^d)$
has a unique projection $P_K\mu$, and $P_K$ is nonexpansive in $W_2$
if and only if $K$ is displacement convex.
\begin{enumerate}[label=\textup{(\roman*)},leftmargin=2em]
\item Let $f\ge0$ be locally integrable on $\R^d$, with $\int f\dx>1$,
and set $K=\{r\dx\in\Ptwo(\R^d):r\le f\text{ a.e.}\}$. Then $P_K$
is nonexpansive if and only if $f^{-1/d}$ agrees almost everywhere
with a lower semicontinuous convex function $\R^d\to[0,+\infty]$. 
\item Let $m\in\Ptwo(\R^d)$ have a Lebesgue density that is positive
almost everywhere. Fix $c>0$ and set
$K=\{\rho\in\Ptwo(\R^d):\Ent_m(\rho)\le c\}$. Then $P_K$ is
nonexpansive if and only if $m$ is log-concave.
\end{enumerate}
\end{theorem}

For $f=\one_E$ with $|E|>1$, part~\textup{(i)} says that the
projection is nonexpansive precisely when $E$ agrees up to a null set
with a convex set. Part~\textup{(ii)} gives a condition independent
of the positive level $c$. Thus nonexpansiveness at a single positive
level forces $m$ to be log-concave, in which case the optimal Lipschitz
constant is one at every positive level. At level zero the constraint
set is the singleton $\{m\}$.

For sufficiency we again approximate the constraint by energies
with nonexpansive proximal maps. In dimension $d\ge2$, a smooth
positive cap $f$ leads to a power penalty with weight
$(f^{-1/d})^{1+d(p-1)}$ (\cref{cor:weighted-power-criterion}), while
an entropy bound leads to penalties that are convex increasing
functions of $\Ent_m$. In each case the stated convexity condition
ensures that the second variation is nonnegative on compatible
fields, and the projection inherits nonexpansiveness in the limit.

For the converse, consider the unique Wasserstein midpoint of two
feasible measures. A nonexpansive projection sends it to a point at
most half the endpoint distance from either endpoint, so uniqueness
forces the midpoint itself to be feasible. For a density cap we test
this condition on measures that differ by a translation and a
dilation on a small part of their mass. The density of their midpoint
gives the convexity inequality for $f^{-1/d}$ at Lebesgue points.
For the entropy constraint, failure of log-concavity of $m$ allows us
to construct two compactly supported densities of entropy strictly
below $c$ whose midpoint has entropy strictly above $c$. The projected
midpoint must then move farther from at least one endpoint.

A density cap for which the projection is nonexpansive can still allow
expansion once an interaction energy is added. For a Gaussian kernel
we can determine the exact expansion factor. Let $\Omega$ be a ball in $\R^d$ and
equip $\Ptwo(\overline\Omega)$ with the Wasserstein distance for the
cost $|x-y|^2$. On $\Omega$, let $f=g^{-d}$, where $g$ is smooth,
positive and convex near $\overline\Omega$, and assume
$\int_\Omega f\dx>1$. For $\ell,\gamma>0$, consider the energy
\[
 \frac\gamma2\iint W_\ell(x-y)\dd\rho(x)\dd\rho(y),
 \qquad W_\ell(z)=\exp(-|z|^2/(2\ell^2)),\qquad \rho\le f\dx .
\]
Let $\Lambda_f$ be the minimum of $\D^2(W_\ell*\rho)(x)[e,e]$ and
$\alpha_{f,\ell}$ the maximum of $(W_\ell*\rho)(x)$ over admissible
probabilities $\rho$, points $x\in\overline\Omega$ and, in the first
case, unit vectors $e$. Because $-\alpha_{f,\ell}/\ell^2<\Lambda_f<0$, the quantity
$(1+\tau\gamma\Lambda_f)^{-1}$ is finite and greater than one
whenever $\tau\gamma\alpha_{f,\ell}<\ell^2$. In this range it is
the optimal Lipschitz constant of the proximal map
(\cref{thm:capped-gaussian-sharp}). Adding a confining potential with
$\D^2 V\succeq\kappa I$ gives the lower curvature bound
$\kappa+\gamma\Lambda_f$, improving the pointwise kernel bound
$\kappa-\gamma/\ell^2$. In particular, for
$-\gamma\Lambda_f<\kappa<\gamma/\ell^2$ the constrained proximal map
is strictly contractive for every $\tau>0$
(\cref{thm:capped-gaussian-bound}).

\subsection{The entropy step on a Ricci-flat manifold}
On a closed connected Ricci-flat manifold, the Wasserstein gradient flow of entropy
relative to normalized volume is the heat flow, which is nonexpansive
in $W_2$~\cite{vonRenesseSturm}. Does its variational discretization
preserve this property? On a flat manifold the entropy proximal map is
nonexpansive for every step size. If the manifold is nonflat, then for every sufficiently small $\tau>0$
two smooth positive densities can be found whose distance strictly
increases after one step. Thus nonexpansiveness along even one
sequence $\tau_j\downarrow0$ characterizes flatness
(\cref{thm:ricci-resolvent}).

To see the obstruction, compare two concentrated measures whose centers lie
nearby on a geodesic with initial velocity $v$. For an isotropic covariance the
curvature contribution to the comparison of distances is proportional
to $\tr R_v=\operatorname{Ric}(v,v)=0$, where $R_v=R(\,\cdot\,,v)v$.
At a point of nonzero sectional curvature, however, one can choose $v$
so that $R_v$ has eigenvalues of both signs. A suitable anisotropic
covariance then gives a strictly larger distance between the outputs
than between the inputs.

\subsection{Functional inequalities after an entropy step}
An entropy proximal step involves two given measures, the input
$\mu$ and the reference measure $m$. Caillet and
Coudreuse~\cite[Theorem~1.1 and Section~3]{CailletCoudreuse} used the
dissipation of entropy and Fisher information over successive steps to
derive a generalized log-Sobolev inequality for $m$. Here we fix $m$
and ask how an inequality satisfied by $\mu$ passes to the output of
one step. When $m$ and $\mu$ are centered isotropic Gaussians, the
output is again Gaussian, so its Poincar\'e, Talagrand and log-Sobolev
constants all equal its variance. This Gaussian calculation suggests a common bound,
which the following theorem proves for general inputs satisfying any
one of the three inequalities.

Write $C_{\mathrm{LS}}(\xi)$ for the best constant in
\eqref{eq:lsi-normalization}, and $C_P(\xi)$ and $C_T(\xi)$ for the
best constants in
\[
 \Var_\xi(f)\le C_P(\xi)\int |\nabla f|^2\dd\xi,
 \qquad
 W_2^2(\sigma,\xi)\le 2C_T(\xi)\Ent_\xi(\sigma),
\]
respectively, for smooth compactly supported $f$ and
$\sigma\in\Ptwo(\R^d)$ of finite entropy relative to $\xi$.

\begin{theorem}\label{thm:intro-lsi-jko}
Let $m=Z^{-1}e^{-V}\dx$, where $V:\R^d\to\R$ is finite and
$\kappa$-strongly convex for some $\kappa>0$.
For $\tau>0$ and $\mu\in\Ptwo(\R^d)$, let $\rho=J_\tau^m\mu$ be the
proximal step for $\Ent_m$. Given $0\le C<\infty$, let $L$ be the
unique positive solution of
\begin{equation}\label{eq:intro-lsi-jko}
 (1+\kappa\tau)L=\tau+\sqrt{CL}.
\end{equation}
Then
\[
 \begin{aligned}
 C_P(\mu)\le C &\ \Longrightarrow\ C_P(\rho)\le L,\\
 C_T(\mu)\le C &\ \Longrightarrow\ C_T(\rho)\le L,\\
 C_{\mathrm{LS}}(\mu)\le C &\ \Longrightarrow\ C_{\mathrm{LS}}(\rho)\le L .
 \end{aligned}
\]
All three bounds are sharp in every dimension. For the Gaussian
reference measure $m=\N(0,\kappa^{-1}I)$ and input $\mu=\N(0,CI)$,
the output is $\N(0,LI)$. For the point mass $\mu=\delta_0$, this
gives the case $C=0$ and $L=\tau/(1+\kappa\tau)$.
\end{theorem}

Given a test function $f$ of mean zero under the output law $\rho$,
we first solve $(I+\tau B)u=f$ in the smooth approximation used in the proof.
If $T$ transports $\rho$ optimally to $\mu$, then the input Poincar\'e
inequality applies to $u\circ T^{-1}$, whose gradient satisfies
$\nabla(u\circ T^{-1})\circ T=H^{-1}\nabla u$, so its energy
is controlled by the Bochner estimate. Testing the equation against $f$
then gives the bound on the variance of $f$ (\cref{thm:poincare-jko}).

For Talagrand and log-Sobolev, transporting a smooth tilt $\nu=e^h\rho$
of total mass one by $T$ gives a new input $\eta=T_\#\nu$. To use the
input inequality we must compare $\nu$ with the minimizer at $\eta$.
The excess of the proximal objective at $\nu$ over its minimum
provides this comparison. Strong convexity along a generalized
geodesic bounds the squared
distance $W_2^2(\nu,J_\tau^m\eta)$ by a multiple of this excess, while variation
at the tilted measure bounds the excess in terms of the relative
Fisher information of $\nu$ with respect to $\rho$. These two bounds
on the same excess enter the proofs of the Talagrand and log-Sobolev
estimates (\cref{lem:lsi-finite-defect,thm:talagrand-jko,thm:lsi-jko}).
Under the regularity hypotheses of \cref{thm:lsi-jko-rigidity},
a nonconstant tilt attains the output log-Sobolev bound exactly when
its log density relative to $\rho$ is affine with gradient in a
direction along which $m$ and $\mu$ have Gaussian factors of
variances $\kappa^{-1}$ and $C$, respectively.

Gaussian directions also determine equality in Talagrand interpolation
along a Wasserstein geodesic. Let $(\mu_t)_{0\le t\le1}$ join two
measures $\mu_i=e^{-V_i}\dx$ with $V_i$ finite and
$\kappa_i$-strongly convex, $\kappa_i>0$. Gozlan, Malamut and
Waldspurger~\cite[Theorem~2.1]{GMW} proved the sharp inequality
\[
 W_2^2(\sigma,\mu_t)\le
 2\left(\frac{1-t}{\sqrt{\kappa_0}}+
         \frac{t}{\sqrt{\kappa_1}}\right)^2\Ent_{\mu_t}(\sigma)
 \qquad(\sigma\in\Ptwo(\R^d)),
\]
from which they recovered the Poincar\'e interpolation of~\cite{HanZhu}
by linearization. To characterize equality we retain the remainders
in the change-of-variables formula for entropy at the two endpoints,
as in the transport proofs of
\cite{CorderoGaussian,OttoVillani,KolesnikovConvexity}. At an interior time equality forces both remainders to vanish.
Tracing this condition back through the endpoint inequalities shows
that a measure
$\sigma\ne\mu_t$ of finite relative entropy attains the bound precisely
when it is a translate of $\mu_t$ in a common Gaussian direction of
the endpoints, with variances $\kappa_0^{-1}$ and $\kappa_1^{-1}$
in that direction (\cref{cor:talagrand-rigidity}).

\subsection*{Reading guide}
Sections~\ref{sec:compatible} and~\ref{sec:tools} develop the Bochner
identities and the second-variation formula used in Part~II for heat
flow, moment reconstruction and Brenier potentials. In Part~III the pressure condition ensures that a curve of proximal
minimizers moves no faster than its curve of inputs
(Section~\ref{sec:resolvents}). Equality in that bound leads to rigidity
for critical diffusion (Section~\ref{sec:rigidity}), while passage to
the limit in penalized energies gives the projection theorems
(Section~\ref{sec:projections}). Sections~\ref{sec:ricci-resolvent}
and~\ref{sec:talagrand} treat the entropy step, first on a Ricci-flat
manifold and then in Euclidean space for functional inequalities.

\subsection*{Origins and preparation of the work}
An early question about contraction towards equilibrium arose in the
first author's 2012 master's thesis, written under the supervision of
Fran\c{c}ois Bolley~\cite[Section~3.2]{HanMasters}. Translating a
one-dimensional equilibrium eliminates the diffusion term in the
Wasserstein dissipation, leaving only the drift term, which can be
negative for a nonconvex confining potential. Bolley, Gentil and
Guillin attribute this observation to the first
author~\cite[discussion following Lemma~3.2]{BGGWJ}.
Section~\ref{sec:wj-equilibrium} returns to the case of positive
dissipation and asks whether a $WJ$ inequality, which gives a rate of
contraction towards equilibrium in continuous time, implies a
corresponding estimate after one entropy proximal step.

The immediate starting point of this work was our sharp Poincar\'e
interpolation along Euclidean Wasserstein geodesics~\cite{HanZhu}. 
Our proof drew on Kolesnikov's work on
Hessian metrics in optimal transport~\cite{KolesnikovHessian}. We then asked whether the Poincar\'e bound extends to Riemannian
manifolds and whether the Talagrand and log-Sobolev constants admit
analogous interpolation bounds in Euclidean space. At the same time,
we were studying the stability of optimal transport maps. 
Several weeks before OpenAI~\cite{OpenAIBrenier} published its counterexample to the
$1/2$-H\"older estimate for optimal transport maps conjectured by
Letrouit~\cite[Conjecture~3]{LetrouitUnstable}, we had already noticed
that the coupled Bochner identities from our interpolation proof could
also yield transport stability estimates. 

Gozlan, Malamut and Waldspurger~\cite{GMW} proved the sharp Talagrand
interpolation inequality along Wasserstein geodesics and recovered the
corresponding Poincar\'e bound by linearization. On the stability side,
Letrouit and Manole~\cite{LetrouitManole} proved a sharp
$W_2$-stability estimate for Brenier potentials by linearizing the
transport equation. In source coordinates, their equation involves
the same elliptic operator we used to construct the Hessian-weighted
velocity in our interpolation proof.

The transport proof of the sharp Sobolev inequality by
Cordero-Erausquin, Nazaret and Villani~\cite{CorderoNazaretVillani}
had already drawn our attention to the Hessian of the transport
potential through the comparison of its determinant and trace. For the present argument it is the action of that Hessian on a
velocity field that matters, since this is how we transport the field
to the second endpoint. Santambrogio's
variational formulation of moment measures~\cite{SantambrogioMoment}
suggested asking how a minimizing source changes with its prescribed
target, while our work on entropy inequalities~\cite{HanZhanZhouZhu}
and the convolution estimates of Fathi, Goldman and
Tsodyks~\cite{FathiGoldmanTsodyks} led us to follow two prescribed
measures under the heat flow. 

Our aim in this paper was to understand what the coupled Bochner approach from our work on Poincar\'e interpolation could tell us about stability and rigidity in optimal transport. The comparison between
the Wasserstein velocity and the Hessian-compatible velocity links the
second variation of a functional at a fixed target to the stability
of its minimizers as the target varies. We have presented these
applications together to bring out their common second-order structure
and developed them in some detail in the hope that the approach will
prove useful in other problems as well. 

Generative AI tools were used throughout the development and preparation
of this work as research-assistance tools. Their use included, where
appropriate, computation, testing arguments and possible
counterexamples, identifying relevant references, and supporting the drafting
and revision of the manuscript. The research program, questions, conceptual
framework, interpretation of results, mathematical
judgments, and final organization of the work were determined by the
authors. The authors take full responsibility for the mathematical
content of the manuscript. 

\section{Transport Hessians and compatible variations}\label{sec:compatible}
A change in a density determines the weighted divergence of its velocity
field, but many fields can have that divergence. When two densities are
related by an optimal map $T=\nabla\Phi$, the Hessian $\D^2\Phi$
defines a quadratic energy on source velocities. The minimizer of
this energy under the divergence constraint gives the source field,
and pushing it forward by $T$ gives the target field. The two fields
represent the same density variation in different coordinates. We show
below why the derivative terms in their Bochner identities are
nonnegative, even though the fields need not be gradients.

\subsection{The Poisson equation in source and target coordinates}
Let $\mu_0$ be a probability measure with smooth positive density, and
let $T=\nabla\Phi$ be a smooth optimal transport diffeomorphism with
$H=\D^2{\Phi}\succ0$. The infinitesimal density variation $g\mu_0$
is represented by a vector field $X$ satisfying $-\diver_{\mu_0} X=g$.
Among these fields we minimize the quadratic energy
\begin{equation}\label{eq:quadratic-flux}
 \inf_{-\diver_{\mu_0} X=g}
        \frac12\int X^{\mathsf T}HX\dd\mu_0.
\end{equation}
Suppose that $H$ is bounded and uniformly elliptic, that $\mu_0$ satisfies
a Poincar\'e inequality, and that $g\in L^2_0(\mu_0)$. We impose the
constraint in the weak sense by requiring
$\int X\cdot\nabla h\dd\mu_0=\int gh\dd\mu_0$ for every
$h\in H^1(\mu_0)$. On a bounded domain this also imposes zero normal
flux. The Poincar\'e inequality makes the bilinear form
$\int H^{-1}\nabla u\cdot\nabla h\dd\mu_0$ coercive on centered
potentials, so the Lax--Milgram theorem gives a centered $u$ for which
$H^{-1}\nabla u$ satisfies the constraint. If $X$ is any other
admissible field, testing the difference of their divergence equations
with $u$ gives
$\int (X-H^{-1}\nabla u)\cdot\nabla u\dd\mu_0=0$.
Thus expanding the square in \eqref{eq:quadratic-flux} shows that
$H^{-1}\nabla u$ is its unique minimizer.

The form $\int H^{-1}[\nabla z,\nabla z]\dd\mu_0$ also arises by
differentiating the transport semidual. Set $\mu_1=T_\#\mu_0$ and
$\mathscr K(\phi)=\int\phi\dd\mu_0+\int\phi^*\dd\mu_1$, where
$\phi^*$ denotes the Legendre transform. For $z$ smooth and compactly
supported in the source interior and $y=T(x)$, the maximizer $x_t$
in $(\Phi+t z)^*(y)$ satisfies $\nabla(\Phi+t z)(x_t)=y$.
Differentiating this relation at $t=0$ gives
$\dot x_0=-H^{-1}\nabla z(x)$. Since the first derivative of the
Legendre transform along this perturbation is $-z(x_t)$, we obtain
\[
 \D^2\mathscr K(\Phi)[z,z]=\int H^{-1}[\nabla z,\nabla z]\dd\mu_0 .
\]
In particular, the minimum in \eqref{eq:quadratic-flux} also has a
dual expression in terms of the density variation $g$. Completing the
square about the Poisson solution $u$ gives
\[
 \sup_{z\in H^1(\mu_0)}
 \left\{\int gz\dd\mu_0-\tfrac12\int H^{-1}[\nabla z,\nabla z]\dd\mu_0\right\}
 =\tfrac12\int gu\dd\mu_0 .
\]

To compare the minimizer with the field it induces on the target,
associate to each smooth scalar function $u$ the pair
\begin{equation}\label{eq:couple}
 X=H^{-1}\nabla u,\qquad Y\circ T=\nabla u.
\end{equation}
The chain rule gives $X=\nabla(u\circ T^{-1})\circ T$. Although $\D X$ need not be symmetric in Euclidean coordinates, it is
self-adjoint for the inner product defined by $H$. This observation
accounts for the sign of the trace in the next lemma.

\begin{lemma}\label[lemma]{lem:symmetry}
For the fields in \eqref{eq:couple}, $H\D X$ and $((\D Y)\circ T)H$ are
symmetric. In particular,
\begin{equation}\label{eq:trace-positive}
 \tr((\D X)^2)
 =\norm{H^{1/2}\D XH^{-1/2}}_{\HS}^2
 \ge\frac{(\diver X)^2}{d}.
\end{equation}
The corresponding target identity is
\[
 (\tr((\D Y)^2))\circ T
 =\norm{H^{-1/2}\D^2{u}H^{-1/2}}_{\HS}^2.
\]
If $T_\#\mu_0=\mu_1$, then
\begin{equation}\label{eq:common-div}
 g:=-\diver_{\mu_0}X=-(\diver_{\mu_1}Y)\circ T
\end{equation}
in the weak sense.
\end{lemma}

\begin{proof}
Differentiating $HX=\nabla u$ yields 
$H\D X=\D^2{u}-\D^3\Phi[X,\cdot,\cdot]$. The right-hand side is
symmetric because $H$ is a Hessian. The matrix $H^{1/2}\D XH^{-1/2}$
is therefore symmetric and has the same eigenvalues as $\D X$.
Applying Cauchy--Schwarz to those eigenvalues gives \eqref{eq:trace-positive}. Differentiating
$Y\circ T=\nabla u$ also yields  $((\D Y)\circ T)H=\D^2{u}$.
Finally, for a test function $\psi$, we have
\[
 \int\nabla\psi\cdot Y\dd\mu_1
 =\int\nabla\psi(T)\cdot HX\dd\mu_0
 =\int\nabla(\psi\circ T)\cdot X\dd\mu_0,
\]
which proves \eqref{eq:common-div}.
\end{proof}

The trace estimate \eqref{eq:trace-positive} can fail for a positive
definite matrix field that is not a Hessian.
For example, on $\{x>0\}\subseteq\R^2$, take
\[
 H=\begin{pmatrix}1&y/x\\y/x&x^{-1}+y^2/x^2\end{pmatrix},
 \qquad X=(-y,x).
\]
Then $H\succ0$ and $HX=\nabla y$, but $\tr((\D X)^2)=-2$.
Compatible fields need not be Euclidean gradients either. For $H=\operatorname{diag}(2,1)$
and $u(x,y)=xy$, the field $X=(y/2,x)$ is not a gradient.

When $H=\D^2\Phi$, symmetry of $H\D X$ means that the one-form with
coefficients $HX$ is closed, since the additional term in $\D(HX)$ is
symmetric. On a simply connected domain that form has a scalar
potential, whereas on a torus the existence of a periodic potential
also requires that all periods vanish.

To apply the identity at a boundary, we need the pushforward of $X$
as a vector field. Although $X$ can be written as a target gradient
composed with $T$, its pushforward includes the derivative of $T$ and is $Y=(\D T\,X)\circ T^{-1}$. Thus tangency of $X$ at a
source boundary gives tangency of $Y$ at the image boundary.

\subsection{Endpoint identities and the Hessian metric}
For a Neumann gradient, the boundary term in the integrated Bochner
identity is the second fundamental form applied to the tangential
gradient. Kolesnikov and Milman use the weighted Reilly formula and its
boundary conditions to prove Brascamp--Lieb inequalities on weighted
manifolds~\cite{KMReilly}. Our compatible field satisfies $HX=\nabla u$
and its Euclidean derivative can be nonsymmetric. We therefore retain
the term $\tr((\D X)^2)$ for a general tangent field. Let $\Omega$ be a
bounded smooth domain with outer unit
normal $n$ and second fundamental form $\II(v,v)=\inner{\D_vn}{v}$
for tangent vectors $v$. Let $\mu=r\dx=e^{-V}\dx$ on $\Omega$ and
let $X$ be a smooth field tangent to the boundary. Write
$g=-\diver_\mu X$. For a smooth scalar weight $w$,
\begin{align}\label{eq:weighted-bochner}
 \int_\Omega wg^2\dd\mu
 ={}&\int_\Omega w\tr((\D X)^2)\dd\mu
 +\int_\Omega(w\D^2{V}-\D^2{w})[X,X]\dd\mu\\
 &+2\int_\Omega g\nabla w\cdot X\dd\mu
 +\int_{\partial\Omega}wr\II(X,X)\dd S.\nonumber
\end{align}
Indeed, for $Q=-gX-(\D X)X$,
\[
 \diver_\mu Q=g^2-\tr((\D X)^2)-\D^2{V}[X,X],
 \qquad Q\cdot n=\II(X,X).
\]
Integrating the divergence of $wQ$ gives the stated boundary term.
The product rule also leaves an integral involving $\D X$, which we
evaluate by integrating the derivative of $\nabla w\cdot X$ along $X$. This gives
\[
 \int\nabla w\cdot(\D X)X\dd\mu
 =-\int \D^2{w}[X,X]\dd\mu+\int g\nabla w\cdot X\dd\mu.
\]
For compactly supported fields on $\R^d$, the boundary term is absent.

Let $\mu_i=e^{-V_i}\dx$ be probability measures with smooth positive
densities on $\R^d$, and let $T=\nabla\Phi$ be a smooth optimal transport
diffeomorphism from $\mu_0$ to $\mu_1$. Put $H=\D T$.
For $u\in C_c^\infty(\R^d)$, take the fields $X,Y$ of \eqref{eq:couple}
and set $g=-\diver_{\mu_0}X$. By \eqref{eq:common-div}, applying
\eqref{eq:weighted-bochner} with $w=1$ at both endpoints gives
\begin{equation}\label{eq:two-bochner}
 \begin{split}
 \int g^2\dd\mu_0
   &=\int\{\D^2{V_0}[X,X]+\tr((\D X)^2)\}\dd\mu_0,\\
 \int g^2\dd\mu_0
   &=\int\{\D^2{V_1}[Y,Y]+\tr((\D Y)^2)\}\dd\mu_1.
 \end{split}
\end{equation}
Kolesnikov, Livshyts and Rotem give the second identity in source
coordinates for their class of measures~\cite{KolesnikovLivshytsRotem}.
In our smooth setting, applying the target identity to the inverse
transport $T^{-1}=\nabla\Phi^*$ with test function $u\circ T^{-1}$
gives the source identity. For $u=\partial_i\Phi$ we have $X=e_i$,
$g=\partial_i V_0$ and $Y\circ T=He_i$. When the integrals are finite
and the cutoff errors vanish, the second identity then gives
Kolesnikov's directional identity for Sobolev regularity of the
transport map~\cite{KolesnikovSobolev}. 
Averaging the two identities gives the integrated formula also recorded
in~\cite{KolesnikovLivshytsRotem}.

\begin{proposition}\label[proposition]{prop:intrinsic-bochner}
In the smooth whole-space setting of \eqref{eq:two-bochner}, let
$u\in C_c^\infty(\R^d)$ and write $\nabla_H^2u$ for its covariant
Hessian in the metric $H$. For a covariant two-tensor $S$ write
$|S|_H=\|H^{-1/2}SH^{-1/2}\|_{\HS}$. With $X=H^{-1}\nabla u$ and
$Bu=-\diver_{\mu_0}X$, we have
\begin{align*}
 \int (Bu)^2\dd\mu_0
 &=\int\!\left(|\nabla_H^2u|_H^2+
       \tfrac14|\D^3\Phi[X,\cdot,\cdot]|_H^2\right)\dd\mu_0\\
 &\quad+\tfrac12\int\!\left(\D^2V_0[X,X]+\D^2V_1(T)[HX,HX]\right)\dd\mu_0 .
\end{align*}
\end{proposition}

\begin{proof}
The source and target derivative terms in \eqref{eq:two-bochner}, in
source coordinates, are $|\D^2u-\D^3\Phi[X,\cdot,\cdot]|_H^2$ and
$|\D^2u|_H^2$. The Christoffel symbols of $H$ in affine coordinates
are $\Gamma_{ij}^k=\frac12\sum_\ell (H^{-1})_{k\ell}\partial_{ij\ell}\Phi$,
so $\nabla_H^2u=\D^2u-\frac12\D^3\Phi[X,\cdot,\cdot]$.
The two tensors whose norms occur above have mean $\nabla_H^2u$ and
differ by $\D^3\Phi[X,\cdot,\cdot]$. The parallelogram identity now
gives the asserted formula.

To compare it with the usual Bochner formula, note that the volume element is $\sqrt{\det H}\dx$, so the potential of
$\mu_0$ relative to metric volume is $V_0+\tfrac12\log\det H$.
Kolesnikov's formula for the corresponding Bakry--\'Emery tensor
\cite[Corollaries~4.1--4.2]{KolesnikovHessian}, evaluated on $X$,
gives the term involving $\D^3\Phi$ and the two potential Hessians. On a bounded domain the same averaging also retains the boundary
terms from the two applications of \eqref{eq:weighted-bochner}.
\end{proof}

If $\D^2V_i\succeq\kappa_i I$ with $\kappa_i>0$, the Bakry--\'Emery
tensor is bounded below by
$(\kappa_0 I+\kappa_1 H^2)/2\succeq\sqrt{\kappa_0\kappa_1}H$,
as in \cite[Theorem~4.3]{KolesnikovHessian}. Under the ellipticity
assumptions of \cref{mm:prop:weighted}, this curvature bound gives the
spectral gap estimate \eqref{mm:eq:weighted-poincare}.
By \cref{lem:symmetry}, the trace terms in \eqref{eq:two-bochner} are
nonnegative. The two endpoint identities therefore give
\begin{equation}\label{eq:two-flux-bounds}
 \norm{X}_{L^2(\mu_0)}\le\kappa_0^{-1/2}\norm{g}_{L^2(\mu_0)},
 \qquad
 \norm{Y}_{L^2(\mu_1)}\le\kappa_1^{-1/2}\norm{g}_{L^2(\mu_0)}.
\end{equation}
Under the hypotheses of \cref{mm:prop:weighted}, these bounds extend to
the centered weak solution of $Bu=g$ for every $g\in L^2(\mu_0)$ of mean
zero~\cite[Proposition~2.4]{HanZhu}. Indeed, a function $z\in L^2(\mu_0)$
orthogonal to $B C_c^\infty$ solves $Bz=0$ distributionally. Local
elliptic regularity and testing with $\chi_R^2z$ for standard ball
cutoffs bound its energy on $B_R$ by $16R^{-2}\|H^{-1}\|_\infty\|z\|_2^2$,
so it is constant. Approximation of $g$ by $Bu_j$, $u_j\in C_c^\infty$,
then gives convergence of both endpoint fields by \eqref{eq:two-flux-bounds}.

\subsection{Finite Jacobian comparisons}\label{sec:finite-comparisons}
The sign in \eqref{eq:trace-positive} also has a finite counterpart.
To compare two transports before differentiating, let $S_i=\nabla\Psi_i$ be smooth strictly convex gradient
diffeomorphisms from a common convex domain and let $M_i=\D S_i\succ0$.
For $F=S_1\circ S_0^{-1}$ one has $\D F(S_0(y))=M_1M_0^{-1}$.
This matrix is similar to $M_0^{-1/2}M_1M_0^{-1/2}$, even when
$M_0$ and $M_1$ do not commute. Summing $\lambda-1-\log\lambda$
over its eigenvalues gives
\[
 \mathcal J(\D F):=\tr(\D F)-d-\log\det \D F\ge0 .
\]
Equality requires every eigenvalue to be one. Since $\D F$ is
diagonalizable, this is equivalent to $\D F=I$.

To see which quadratic term produces this remainder, consider
$S_s=(1-s)S_0+sS_1$, an affine interpolation of smooth convex
gradients on a convex domain, and suppose $M_s=\D S_s\succ0$.
On the image of $S_s$, define
$Z_s=\dot S_s\circ S_s^{-1}$. In the original coordinates,
\[
 \D Z_s\circ S_s=\dot M_sM_s^{-1},\qquad
 -\frac{\diff^2}{\diff s^2}\log\det M_s
 =\norm{M_s^{-1/2}\dot M_sM_s^{-1/2}}_{\HS}^2
 =\tr\bigl((\D Z_s)^2\bigr)\circ S_s.
\]
Since $M_s^{-1}(\D Z_s\circ S_s)$ is symmetric, this is the
compatibility in \cref{lem:symmetry} written for the inverse map.
Taylor's formula in $s$ thus gives a nonnegative remainder for
$-\log\det$. Section~\ref{sec:talagrand} retains this remainder
in the entropy comparison at each endpoint. On bounded domains the
images of $S_s$ may vary with $s$, and the boundary contribution
also enters the comparison.

At a fixed point where $H\D X$ is symmetric and $H\succ0$,
the same similarity argument applies to $I+\epsilon \D X$ for
sufficiently small $\epsilon$. In particular
\[
 \mathcal J(I+\epsilon \D X)
 =\tfrac12\epsilon^2\tr((\D X)^2)+o(\epsilon^2).
\]
Locally, with $S=T^{-1}$, the composition
$(\Id+sX)\circ S=S+s\nabla(u\circ S)$ remains a convex gradient
as long as its Hessian stays positive. Thus these variations can be
realized by convex gradients on the same domain. Their images may vary,
and applying them to a particular energy still requires computing its
second variation.

\section{Endpoint estimates and transport variations}\label{sec:tools}
\subsection{Poincar\'e interpolation and the transport spectral gap}
Assume first that $T$ is a smooth transport diffeomorphism on $\R^d$,
that $H$ and $H^{-1}$ are bounded and that both endpoint potentials
are smooth with $\D^2V_i\succeq\kappa_i I$, $\kappa_i>0$. The centered
Poisson solution and its two endpoint estimates are then available
by the preceding section. This gives the smooth interpolation estimate of
\cite[Proposition~2.5]{HanZhu}. Write $F_t=(1-t)\Id+tT$ and
$\mu_t=(F_t)_\#\mu_0$. For centered $h\in H^1(\mu_t)$ solve
$Bu=h\circ F_t$ and set $Z_t\circ F_t=\D F_tX=(1-t)X+tY\circ T$.
Changing variables in the weak divergence identity gives
$-\diver_{\mu_t}Z_t=h$. Testing against $h$ yields
\[
 \|h\|_2^2\le\|\nabla h\|_2\|Z_t\|_2
 \le\|\nabla h\|_2
       \big((1-t)\kappa_0^{-1/2}+t\kappa_1^{-1/2}\big)\|h\|_2 .
\]
The triangle inequality for $(1-t)X+tY\circ T$ in the last step
therefore bounds the square root of the Poincar\'e constant of $\mu_t$.
For isotropic Gaussians of variances $\kappa_i^{-1}$, the interpolant has
variance $((1-t)\kappa_0^{-1/2}+t\kappa_1^{-1/2})^2$, so the bound is exact. For general strongly
log-concave endpoints see \cite[Theorem~1.1]{HanZhu}.

The two bounds in \eqref{eq:two-flux-bounds} give a weighted
Poincar\'e inequality. When both measures are centered, applying it to
the components of $T-\Id$ gives the spectral estimate used in the
moment-measure comparison. Write $P_E$ for the orthogonal projection
onto a subspace $E$.

\begin{proposition}\label[proposition]{mm:prop:weighted}
Let $\mu_i=e^{-V_i}\dx$ be smooth probability measures on $\R^d$,
with $\D^2{V_i}\succeq\kappa_iI$, $\kappa_i>0$. Suppose that the
Brenier map $T=\nabla\Phi$ between them is a smooth diffeomorphism and
that $H=\D^2{\Phi}$ and $H^{-1}$ are bounded. Set $c=\sqrt{\kappa_0\kappa_1}$ and
$B=-\diver_{\mu_0}(H^{-1}\nabla)$, and write
\[
 \mathcal E_H(f,g)=\int\inner{H^{-1}\nabla f}{\nabla g}\dd\mu_0,
 \qquad \mathcal E_H(f)=\mathcal E_H(f,f).
\]
Then
\begin{equation}\label{mm:eq:weighted-poincare}
 \Var_{\mu_0}f\le c^{-1}\mathcal E_H(f),\qquad f\in H^1(\mu_0).
\end{equation}
If both measures are centered, then
\begin{equation}\label{mm:eq:R-lower}
 \mathfrak R:=\int\tr(H+H^{-1}-2I)\dd\mu_0\ge c\,W_2^2(\mu_0,\mu_1).
\end{equation}
Equality in \eqref{mm:eq:R-lower} holds if and only if $\mu_0=\mu_1$,
or $\kappa_0\ne\kappa_1$ and there exist a nonzero subspace $E$ and a
centered probability measure $\sigma$ on $E^\perp$ such that
\[
 \mu_i=\N_E(0,\kappa_i^{-1}I_E)\otimes\sigma,\qquad i=0,1.
\]
In the latter case $T=\sqrt{\kappa_0/\kappa_1}\,P_E+P_{E^\perp}$.
\end{proposition}

\begin{proof}
For centered $f\in H^1(\mu_0)$, solve $Bu=f$ with centered $u$ by
the Lax--Milgram theorem. The ordinary Poincar\'e inequality used here
follows from the Brascamp--Lieb variance bound~\cite{BrascampLieb}
and $\D^2V_0\succeq\kappa_0 I$.
The weak endpoint estimates \eqref{eq:two-flux-bounds} give
\[
 \mathcal E_H(u)
 \le\norm{H^{-1}\nabla u}_2\norm{\nabla u}_2
 \le c^{-1}\norm f_2^2.
\]
Testing $Bu=f$ with $f$ and applying H\"older's inequality to
$H^{-1/2}\nabla u$ and $H^{-1/2}\nabla f$, we obtain
\[
 \norm f_2^2=\mathcal E_H(u,f)
 \le\mathcal E_H(u)^{1/2}\mathcal E_H(f)^{1/2}
 \le c^{-1/2}\norm f_2\mathcal E_H(f)^{1/2}.
\]
This proves \eqref{mm:eq:weighted-poincare}. If the measures are
centered, every component of $T-\Id$ is centered and lies in
$H^1(\mu_0)$. Applying the inequality to these components and summing
proves \eqref{mm:eq:R-lower}, since
\[
 \sum_{a=1}^d\mathcal E_H((T-\Id)_a)
 =\int\tr\bigl((H-I)H^{-1}(H-I)\bigr)\dd\mu_0=\mathfrak R.
\]

Suppose $\mu_0\ne\mu_1$ and equality holds in \eqref{mm:eq:R-lower}. Each
component $f$ of $T-\Id$ then satisfies $\mathcal E_H(f)=c\norm f_2^2$.
The difference $\mathcal E_H(f+s h)-c\|f+s h\|_2^2$ is nonnegative
for centered $h$ and vanishes at $s=0$. Differentiating there gives $\mathcal E_H(f,h)=c\langle f,h\rangle$ for
every centered form-domain test $h$, and also for constant tests.
Since $cf\in L^2$, this places $f$ in the operator domain with $Bf=cf$.
Put $f_0=f$ and $f_1=f\circ T^{-1}$. By \eqref{eq:couple} and
\eqref{eq:two-flux-bounds} with $u=f$ and $g=cf$,
\[
 \int|\nabla f_i|^2\dd\mu_i\le\kappa_i\int f_i^2\dd\mu_i,\qquad i=0,1.
\]
Since each $f_i$ is centered, the reverse bound follows from the ordinary
Poincar\'e inequality. Equality gives constant gradients by
\cite[Lemma~4.1]{HanZhu}. Thus $T$ is affine. It transports one centered
measure to another, so $T(x)=Lx$ with $L$ symmetric positive definite.

Let $E=\range(L-I)$. The coefficient vectors of the components of
$(L-I)x$ span $E$, as do those of $(I-L^{-1})y$ at the target.
Polarizing the ordinary Poincar\'e equalities for these components
and taking linear combinations gives, for $i=0,1$,
\[
 \int e\cdot\nabla\zeta\dd\mu_i
 =\kappa_i\int(e\cdot x)\zeta(x)\dd\mu_i
 \qquad e\in E,\ \zeta\in C_c^\infty(\R^d).
\]
Integration by parts now gives $\partial_e V_i(x)=\kappa_i\,e\cdot x$
for every $e\in E$. Thus $V_i(x)-\kappa_i|P_E x|^2/2$ depends only on
$P_{E^\perp}x$, which gives the factorization
$\mu_i=\N_E(0,\kappa_i^{-1}I_E)\otimes(P_{E^\perp})_\#\mu_i$.
Since $L_\#\mu_0=\mu_1$ and $L$ is the identity on
$E^\perp=\ker(L-I)$, the orthogonal factors coincide. Comparing
covariances on $E$ gives $\kappa_0^{-1}(L|_E)^2=\kappa_1^{-1}I_E$.
Positivity determines $L|_E=\sqrt{\kappa_0/\kappa_1}\,I_E$.
Since $\mu_0\ne\mu_1$, $E$ is nonzero and $\kappa_0\ne\kappa_1$.

Conversely, put $\lambda=\sqrt{\kappa_0/\kappa_1}$.
For the stated product measures, $\lambda P_E+P_{E^\perp}$ pushes
$\mu_0$ to $\mu_1$ and is the gradient of a convex quadratic. Thus
\[
 c\,W_2^2(\mu_0,\mu_1)=\dim(E)\frac{(\lambda-1)^2}{\lambda}=\mathfrak R.
\]
If $\mu_0=\mu_1$, both sides vanish.
\end{proof}

\subsection{Gaussian factors under near equality}
We will also need the equality argument when the maps vary and only
their endpoint measures have limits. The following form of the argument
in \cite[proof of Theorem~1.4]{HanZhu} uses no uniform bound on the
ellipticity constants of the maps in the sequence.

\begin{lemma}\label[lemma]{lem:endpoint-saturation}
For each $j$ assume the hypotheses of \cref{mm:prop:weighted} with the
same $\kappa_i>0$ and write $c=\sqrt{\kappa_0\kappa_1}$.
Suppose $\mu_i^j\to\mu_i$ in $W_2$. Let $B_ju_j=g_j$ with $u_j,g_j$
centered and $\|g_j\|_2=1$. If
$\int X_j\cdot(Y_j\circ T_j)\dd\mu_0^j\to c^{-1}$, there is a common
unit direction along which $\mu_i$ has an independent Gaussian factor
of variance $\kappa_i^{-1}$.
\end{lemma}
\begin{proof}
We first show that both normalized fields converge to the same constant
vector along a subsequence. Take norms and inner products in $L^2(\mu_0^j)$. Since
$Y_j\circ T_j=\nabla u_j$, the hypothesis says that the inner product
of $\sqrt{\kappa_0}X_j$ and $\sqrt{\kappa_1}\nabla u_j$ tends to one.
Both norms are at most one by \eqref{eq:two-flux-bounds}. Thus both
tend to one and the norm of the difference tends to zero.

Testing $B_ju_j=g_j$ with $u_j$ and using the source Poincar\'e
inequality gives
\[
 \langle X_j,\nabla u_j\rangle=\langle g_j,u_j\rangle
 \le \|u_j\|\le\kappa_0^{-1/2}\|\nabla u_j\|\le c^{-1}.
\]
The left side tends to $c^{-1}$, so $\|u_j\|\to c^{-1}$ and
$\|\nabla u_j\|^2-\kappa_0\|u_j\|^2\to0$.
By \cite[Lemma~4.1]{HanZhu}, the $L^2$ distance from $\nabla u_j$
to its mean tends to zero. Passing to a subsequence on which the
mean converges, we obtain a unit vector $e$ such that
\[
 \|\sqrt{\kappa_0}X_j-e\|_{L^2(\mu_0^j)}\longrightarrow0,
 \qquad
 \|\sqrt{\kappa_1}Y_j-e\|_{L^2(\mu_1^j)}\longrightarrow0.
\]

To identify the limiting measures along this direction, consider the
linear function $\ell(x)=e\cdot x$, which lies in $H^1$ at both endpoints.
Testing the source divergence equation with $\ell$ and using that
$g_j$ is centered with unit norm gives
\[
 \int e\cdot X_j\dd\mu_0^j
 =\int g_j\Bigl(\ell-\int\ell\dd\mu_0^j\Bigr)\dd\mu_0^j
 \le \sqrt{\Var_{\mu_0^j}\ell}\le\kappa_0^{-1/2}.
\]
The last bound follows from the Poincar\'e inequality. At the target
the same bounds hold with $Y_j$ and $g_j\circ T_j^{-1}$. The two left
sides tend to $\kappa_0^{-1/2}$ and $\kappa_1^{-1/2}$, respectively.
Thus $\Var_{\mu_i^j}\ell\to\kappa_i^{-1}$, and convergence in $W_2$
gives $\Var_{\mu_i}\ell=\kappa_i^{-1}$ for $i=0,1$.

For $\zeta\in C_c^\infty(\R^d)$ and $s\in\R$, the Poincar\'e
inequality for $\ell+s\zeta$ passes to the limit by $W_2$ convergence.
The limiting inequality is an equality at $s=0$. Differentiating
there gives, with $b_i=\int x\dd\mu_i$,
\[
 \int\partial_e\zeta\dd\mu_i
 =\kappa_i\int \bigl(e\cdot(x-b_i)\bigr)\zeta(x)\dd\mu_i .
\]
In the coordinates $x=b_i+te+z$, $z\perp e$, this says that
$\partial_t(e^{\kappa_i t^2/2}\mu_i)=0$ in distributions. The weighted
measure is therefore Lebesgue measure in $t$ times a measure in $z$.
Normalizing gives the independent Gaussian factor of variance
$\kappa_i^{-1}$ along $e$.
\end{proof} 

\subsection{The linearized transport equation and its adjoint}\label{sec:response-equations}
Let $T_t=\nabla\Phi_t$ be a smooth transport diffeomorphism from
$\rho_t=r_t\dx$ to $\eta_t=m_t\dy$, with both densities smooth and
positive. On fixed domains we assume zero normal flux. On $\R^d$ we
assume all pairings are finite and the boundary terms vanish. At a
fixed time, suppress $t$ and write $h=\partial_t\log r_t$ and
$k=\partial_t\log m_t$ for the two density variations.
Differentiating $\int\psi(T_t)\dd\rho_t=\int\psi\dd\eta_t$ for a fixed
target test function $\psi$ and writing $u=\psi\circ T$ gives
\[
 \int u(k\circ T-h)\dd\rho
 =\int H^{-1}\nabla u\cdot\nabla\dot\Phi\dd\rho.
\]
This is the weak form of $B\dot\Phi+h=k\circ T$. To recover the same
equation directly from the change of variables, write
$\log r=\log m(T)+\log\det H$, so
$h=k\circ T+\nabla\log m(T)\cdot\nabla\dot\Phi+
\tr(H^{-1}\D^2\dot\Phi)$. The Piola identity
$\diver(\cof H)=0$ identifies the last two terms with $-B\dot\Phi$.
The three applications below differ in how their source density
depends on the potential.

For a fixed source, $h=0$. If $v_t$ is a velocity of the target and
$u=B^{-1}g$ is the centered Neumann solution for a mean-zero right-hand side $g$,
testing the target continuity equation with $u\circ T^{-1}$ gives
\begin{equation}\label{eq:fixed-source-pairing}
 \int g\dot\Phi_t\dd\rho
 =\int H^{-1}\nabla B^{-1}g\cdot v_t(T_t)\dd\rho.
\end{equation}
Taking the supremum of the flux norms over the transport potentials
gives the sharp stability constant in \cref{thm:stability}.

In the classical moment problem $\rho=e^{-\Phi}\dx$,
$H^{-1}\nabla\partial_i\Phi=e_i$, so
$B(\partial_i\Phi)=\partial_i\Phi$ in distributions.
In target coordinates these modes are the coordinate functions.
Suppose $\partial_i\Phi$ belongs to the form domain of $B$ and
$y_i\in L^2(\eta)$, where $\eta=T_\#\rho$. Since the right-hand side
$y_i\circ T$ then belongs to $L^2(\rho)$, the identity on compact
tests extends by form density to the whole-space Friedrichs realization. Indeed, translating the source translates the argument of $\Phi$, so
the derivatives $\partial_i\Phi$ are precisely the variations arising
from translations along the coordinate axes, up to sign.
Klartag used the dual Bochner formula and moment maps to obtain
weighted Poincar\'e inequalities on convex bodies~\cite{KlartagPoincare}.
For the moment-measure metric, his spectral description identifies the
coordinate functions at eigenvalue one~\cite[Section~6]{KlartagMoment}.
In the linearized reconstruction equation these modes lie in the kernel
of $B-I$, as required by translation invariance. In moment coordinates, the squared
Hessian formula for $B(B-I)$ also identifies the affine obstructions
in the even-test density argument of \cite[Eq.~(4.7) and Lemma~5.1]{OpenAILogBM}.

For regularized moment reconstruction, normalize $\Phi_t$ so that
$r_t=e^{-V-\Phi_t}$. Then $h=-\dot\Phi$ and
\[
 (B-I)\dot\Phi=k\circ T.
\]
The same Hessian $H$ occurs in $\D^2(-\log r)=\D^2V+H$ and in
$\int u Bu\dd\rho=\int H[X,X]\dd\rho$ for $X=H^{-1}\nabla u$.
Subtracting the latter term from the source Bochner identity cancels
$H$. The remaining interior terms are the curvature of $V$ and the
nonnegative derivative term in \eqref{mm:eq:regularized-bochner}.

For $\rho_t=J_\tau^{\cE}\eta_t$, the Euler equation
$T=\Id+\tau\nabla q_\rho$ gives
$\nabla\dot\Phi=\tau\nabla(\cK_\rho h)$, where $q_\rho$ is the first
variation of the energy and $\cK_\rho h=\D q_\rho[rh]$. Moving
$\cK_\rho$ to the test function by symmetry gives the adjoint
$\cR_\rho=I+\tau\cK_\rho B$ and the pairing
\[
 \int h\,\cR_\rho u\dd\rho
 =\int H^{-1}\nabla u\cdot(v\circ T)\dd\rho.
\]
For a fixed test function $f$, the derivative of $\int f\dd\rho_t$ is
$\int fh\dd\rho$. Solving $\cR_\rho u=f$ expresses this derivative
in terms of the input velocity, as needed in the proof of Poincar\'e
propagation. A bound on the Wasserstein speed of the minimizer requires
a different choice of $u$. In \cref{prop:primal-contraction} we solve
$Bu=h$, so that $X=H^{-1}\nabla u$ represents $\partial_t\rho_t$.
The Wasserstein speed of $\rho_t$ is at most $\|X\|_{L^2(\rho)}$.
The point of using $X$ is that the Euler equation pairs it with the
input velocity, so a lower bound on the second variation gives a
bound on this norm.

For the smooth variations just considered, substituting the relation
between $h$ and $\dot\Phi$ into $B\dot\Phi+h=k\circ T$ gives
\[
\begin{array}{lll}
 \text{fixed source} & h=0 & B\dot\Phi=k\circ T,\\
 \text{moment reconstruction} & h=-\dot\Phi & (B-I)\dot\Phi=k\circ T,\\
 \text{resolvent} & \dot\Phi=\tau\cK_\rho h\pmod{\R}
     & (I+\tau B\cK_\rho)h=k\circ T .
\end{array}
\]
In the first line we center the potential. In the second line the
normalization of $r$ fixes its additive constant, hence
$\int\dot\Phi\dd\rho=0$. In the third line $h$ has mean zero, and
$B$ removes the constant in the first variation. Symmetry of the
second variation in centered tests is all that the adjoint pairing
requires. In the projection problems of Part~III the corresponding
estimates pass to the constraint by penalization.

For moment measures we can also compare the endpoints directly. The four-point identity in \cref{mm:prop:finite-defect} combines
\eqref{mm:eq:R-lower} with the Bregman remainders of the two potentials.
For a smooth strictly convex potential $\Phi$, this remainder is
\[
 \Dcal_\Phi(y,x)=\Phi(y)-\Phi(x)-\nabla\Phi(x)\cdot(y-x).
\]
Along $y=x+sZ$, its second derivative at $s=0$ is
$Z^{\mathsf T}\D^2\Phi(x)Z$, the quadratic form in \eqref{eq:quadratic-flux}.

\subsection{The Hessian at a variational solution}\label{sec:variational-hessian}

For $\rho,\nu\in\Ptwo(\R^d)$, their maximal correlation is
\[
 \mathcal T(\rho,\nu)=\sup_{\pi\in\Pi(\rho,\nu)}\int x\cdot y\dd\pi
 =\tfrac12\bigl(M_2(\rho)+M_2(\nu)-W_2^2(\rho,\nu)\bigr),
 \qquad M_2(\rho)=\int|x|^2\dd\rho .
\]
Santambrogio proved displacement convexity in the first argument and
used $\Ent+\mathcal T(\cdot,\nu)$ to reconstruct a measure with
prescribed moment measure~\cite[Proposition~3.3 and Theorem~4.1]{SantambrogioMoment}.
We first differentiate $\mathcal T$ with the target fixed. Along a
differentiable branch of critical measures, the resulting quadratic
form pairs the source variation with the target velocity.

Work on $\R^d$ with smooth positive densities, or on smooth bounded
convex source and target domains with positive densities up to the
boundary. Let $F_t=\Id+tv$, where $v=\nabla f$ is smooth with bounded
derivative, and set $\rho_t=(F_t)_\#\rho$. On a bounded source take
$v$ compactly supported in its interior. For small $|t|$ this is a
displacement geodesic. Suppose the optimal transport maps $T_t=\nabla\Phi_t$
from $\rho_t$ to the fixed target are smooth diffeomorphisms,
differentiable in $t$, with $H=\D T_0$ bounded and uniformly elliptic.
We require the weighted Poisson solvability of \cref{sec:compatible},
zero conormal flux or vanishing cutoff errors, and convergence of
the differentiated pairings below, for instance by integrable domination.
Write $T=T_0$, and put
\[
 g=-\diver_\rho v,\qquad Bu=g,\qquad X=H^{-1}\nabla u,\qquad
 \|Z\|_{H,\rho}^2=\int H[Z,Z]\dd\rho .
\]
The solution $u$ is centered. For an energy $\cE$ with differentiable
first variation $q_\rho$, let $\cK_\rho g$ denote the derivative of
$q_\rho$ in the density direction $g\rho$, modulo constants. Assume
the second variation exists along the indicated smooth paths, with
symmetric density derivative on centered tests. For fields with the
same density derivative $g$, write
\[
 \mathscr Q_\rho(Z)=\langle g,\cK_\rho g\rangle_\rho+
                       \int\D^2 q_\rho[Z,Z]\dd\rho .
\]

\begin{theorem}\label{thm:hessian-transfer}
Under these assumptions, it holds that 
\begin{equation}\label{eq:correlation-hessian}
 \operatorname{Hess}_W\mathcal T(\cdot,\nu)_\rho(v,v)
       =\|v-X\|_{H,\rho}^2 .
\end{equation}
Let $a,c\in\R$ and $\mathcal F_\nu=\cE+cM_2/2+a\mathcal T(\cdot,\nu)$.
If $\nabla q_\rho+cx+aT=0$, then
\begin{equation}\label{eq:hessian-transfer}
 \operatorname{Hess}_W\mathcal F_\nu(v,v)
    =\mathscr Q_\rho(X)+c\|X\|_\rho^2
    =\langle g,\cK_\rho g\rangle_\rho-a\langle g,B^{-1}g\rangle_\rho .
\end{equation}
Here the Wasserstein Hessian denotes the second derivative along the
geodesic above, with $\nu$ fixed.

Now let $(\rho_t,\nu_t,T_t)$ be a differentiable branch of smooth optimal transport maps and critical measures with
$(T_t)_\#\rho_t=\nu_t$ and $\nabla q_{\rho_t}+cx+aT_t=0$.
Assume the same regularity for the differentiated pairings and the
same no-flux or decay conditions. At the time in question write
$g\rho=\partial_t\rho_t$, let $b$ be a velocity of $\nu_t$, and set
$X=H^{-1}\nabla B^{-1}g$. Then
\begin{equation}\label{eq:critical-branch-response}
 \mathscr Q_\rho(X)+c\|X\|_\rho^2=-a\langle X,b\circ T\rangle_\rho .
\end{equation}
\end{theorem}

\begin{proof}
Writing couplings of $\rho_t,\nu$ as pushforwards of couplings of
$\rho,\nu$ expresses $\mathcal T(\rho_t,\nu)$ as the supremum of
$\int F_t(x)\cdot y\dd\pi$. Comparison at two times, first with one
optimal plan and then with the other, gives
\[
 \frac{\diff}{\diff t}\mathcal T(\rho_t,\nu)
   =\int v(x)\cdot T_t(F_t(x))\dd\rho(x).
\]
Indeed the maximizing plans on the original source are the graphs
of $T_t\circ F_t$, and their pairings are continuous by hypothesis.
As the target stays fixed, \cref{sec:response-equations} gives
$B\dot\Phi=-g$, so $\nabla\dot\Phi=-HX$.
Differentiating $T_t(F_t(x))$ now gives $H(v-X)$, so the second
derivative is $\int H[v,v-X]\dd\rho$. Since $v$ and $X$ represent
the same density derivative, their difference has zero weighted
divergence and zero normal flux. Testing with $u$ gives
$\int H[X,v-X]\dd\rho=0$. This proves \eqref{eq:correlation-hessian}.

Differentiating the first variation of $\cE$ along the same geodesic
gives $\mathscr Q_\rho(v)$, since the material acceleration vanishes.
Differentiating the critical point equation in space gives
$\D^2q_\rho+cI=-aH$. This identifies the part of the energy Hessian
that combines with the correlation term, and the orthogonality just
proved gives
\[
 \operatorname{Hess}_W\mathcal F_\nu(v,v)
   =\langle g,\cK_\rho g\rangle_\rho
      -a\|v\|_{H,\rho}^2+a\|v-X\|_{H,\rho}^2
   =\langle g,\cK_\rho g\rangle_\rho-a\|X\|_{H,\rho}^2.
\]
Testing $Bu=g$ with $u$ and using the critical equation on $X$
gives the two expressions in \eqref{eq:hessian-transfer}.

When the target also moves, its velocity enters through the
pushforward equation. To isolate the term it contributes, differentiate
$q_\rho+c|x|^2/2+a\Phi=\mathrm{const}$ and pair with $g$.
Since $g$ is centered,
$\langle g,\cK_\rho g\rangle_\rho=-a\langle g,\dot\Phi\rangle_\rho$.
Test the differentiated pushforward with $u\circ T^{-1}$ held fixed
at the time in question, where $u=B^{-1}g$. Its pulled-back gradient
is $X$, so the weak equation for $u$ gives
\[
 \langle X,b\circ T\rangle_\rho
 =\langle u,g\rangle_\rho+\langle X,\nabla\dot\Phi\rangle_\rho
 =\langle u,g\rangle_\rho+\langle g,\dot\Phi\rangle_\rho .
\]
Using $\D^2q_\rho+cI=-aH$ and
$\int H[X,X]\dd\rho=\langle u,g\rangle_\rho$ proves
\eqref{eq:critical-branch-response}.
\end{proof}

For Santambrogio's functional $\Ent+\mathcal T(\cdot,\nu)$ we have
$(a,c)=(1,0)$ and $\cK_\rho g=g$. Thus \eqref{eq:hessian-transfer}
gives the quadratic form of $I-B^{-1}$ on centered density variations.
The eigenfunctions of $B$ at eigenvalue one generated by source
translations in \cref{sec:response-equations} lie in its kernel.
Adding a confinement $V$ changes the critical density and hence $B$,
while $\cK_\rho g=g$ still holds. For a proximal step of size $\tau$,
the coefficient of $\mathcal T$ is $-1/\tau$, so the same transport
term contributes $+\tau^{-1}\langle g,B^{-1}g\rangle_\rho$ to the Hessian.

A third example is the power entropy in the variational construction of
$q$-moment measures by Huynh Khanh and Santambrogio~\cite[Sections~2--3]{HuynhSantambrogio}.
For $q>0$, $d+q>1$, let $m=1-1/(d+q)$ and
$U(r)=r^m/(m-1)$. Up to a positive rescaling of the potential,
$-U'(r)=(d+q-1)r^{-1/(d+q)}$ is the power potential of their
reconstruction. At a smooth critical point of $\int U(r)\dx+\mathcal T$,
$T=-\nabla U'(r)$. Under the whole-space hypotheses above and the
integration-by-parts conditions of \cref{prop:volume-variation}, the
Hessian in \eqref{eq:hessian-transfer} is
\[
 \int r^m\left\{\left\|H^{1/2}\D XH^{-1/2}
            -\frac{\diver X}{d}I\right\|_{\HS}^2
       +\frac{q}{d(d+q)}(\diver X)^2\right\}\dx .
\]
Indeed $p(r)=r^m$ and $m-(1-1/d)=q/[d(d+q)]$. The coefficient of
the divergence vanishes at the critical diffusion power. The formula
concerns smooth variations of finite energy and second moment.

\subsection{Second variations of the optimal transport cost}
\label{sec:transport-second-variations}
Caja-Lopez, Delgadino and Kitagawa~\cite[Theorem~1.4]{CajaDelgadinoKitagawa}
expressed the second derivative of $W_2^2/2$ along affine perturbations
of both densities as the Dirichlet energy of the solution to the
linearized Neumann problem. With the marginals fixed, Peyr\'e, Poon
and Tron~\cite[Theorem~2.4]{PeyrePoonTron} identified the negative
cost Hessian with the squared distance of
$x\mapsto\nabla_x\dot c(x,T(x))$ from gradient fields in a weighted
$L^2$ norm. This also identifies the cost perturbations with zero
second variation.
We estimate both derivatives at the quadratic cost under the curvature
bounds $\D^2V_i\succeq\kappa_i I$, $\kappa_i>0$. The same div--curl
estimate gives the sharp coefficient $(\kappa_0\kappa_1)^{-1/2}$ in
both bounds. For marginal perturbations $HX=\nabla u$ and the skew
term vanishes. For cost perturbations the weighted divergence vanishes.

\begin{proposition} 
\label[proposition]{prop:transport-div-curl}
Let $T=\nabla\Phi$ be a smooth optimal transport diffeomorphism between
smooth probability measures $\mu_i=e^{-V_i}\dx$ on $\R^d$, with
$\D^2V_i\succeq\kappa_i I$ and $\kappa_i>0$. Assume that $H=\D T$ and
$H^{-1}$ are bounded. Write $\operatorname{Skew}A=(A-A^{\mathsf T})/2$.
For $X\in C_c^\infty(\R^d;\R^d)$ denote 
\[
 M_X=H^{1/2}\D XH^{-1/2},\qquad
 Q_H(X)=\int\bigl(|\diver_{\mu_0}X|^2+
                 \norm{\operatorname{Skew}M_X}_{\HS}^2\bigr)\dd\mu_0 .
\]
Then
\begin{equation}\label{eq:transport-div-curl}
 Q_H(X)\ge \sqrt{\kappa_0\kappa_1}\int X^{\mathsf T}HX\dd\mu_0+
 \frac12\int\big|\sqrt{\kappa_0}X-\sqrt{\kappa_1}HX\big|^2\dd\mu_0 .
\end{equation}
The estimate also holds on bounded convex $C^2$ domains if the densities
are smooth and positive up to the boundary, $T$ extends to a $C^2$
diffeomorphism of the closures, and $X$ is $C^1$ and tangent to the
source boundary.
\end{proposition}

\begin{proof}
Define the target field by $Y\circ T=HX$. Changing variables in the
weak divergence formula gives $(\diver_{\mu_1}Y)\circ T=\diver_{\mu_0}X$.
Differentiating the defining relation for $Y$ yields
\[
 H^{-1/2}(\D Y\circ T)H^{1/2}
 =H^{-1/2}\D(HX)H^{-1/2}
 =M_X+H^{-1/2}\D^3\Phi[X,\cdot,\cdot]H^{-1/2}.
\]
The last summand is symmetric, so the two conjugated derivatives have
the same skew part. Apply \eqref{eq:weighted-bochner} at each endpoint
with $w=1$, using invariance of the trace under similarity and
$\tr(M^2)=\norm{\operatorname{Sym}M}_{\HS}^2-
\norm{\operatorname{Skew}M}_{\HS}^2$. On bounded domains $\D T$ preserves
tangency, so convexity makes both boundary terms nonnegative. For $C^1$
fields, pass to the limit by $C^1$ approximation in the integrated
identity for $Q=-gX-(\D X)X$ before using tangency. This gives
\[
 Q_H(X)\ge\kappa_0\int |X|^2\dd\mu_0,\qquad
 Q_H(X)\ge\kappa_1\int |Y|^2\dd\mu_1
          =\kappa_1\int |HX|^2\dd\mu_0 .
\]
Averaging these bounds and completing the square proves
\eqref{eq:transport-div-curl}.

For $\mu_i=\N(0,\kappa_i^{-1}I)$ one has
$H=\sqrt{\kappa_0/\kappa_1}\,I$. A constant field $X=u$ has curl zero
and $Q_H(X)=\kappa_0|u|^2$. In $d\ge2$, a rotation $X(x)=Jx$ with
$J^{\mathsf T}=-J$ has weighted divergence zero and $Q_H(X)=\norm J_{\HS}^2$.
The constant field and the rotation thus give equality with,
respectively, the skew term and the divergence term equal to zero.
Gaussian integrability allows us to approximate either field in the
form norm by smooth cutoffs.
\end{proof}

For compactly supported $X$ there is an intrinsic interpretation. Let
$\omega=HX\cdot\diff x$ be the metric dual one-form. Its weighted
codifferential is $-\diver_{\mu_0}X$, and exterior differentiation
cancels the third derivatives of $\Phi$. With the usual exterior
two-form norm this gives $|\diff\omega|_H^2=2\|\operatorname{Skew}M_X\|_{\HS}^2$.
Consequently $Q_H$ is the sum of the squared codifferential norm and
half the squared exterior derivative norm. In the integrated weighted
Weitzenb\"ock formula the skew part of the covariant derivative has
squared norm $|\diff\omega|_H^2/2$. Subtracting it gives
\[
 Q_H(X)=\int\left(
   |\operatorname{Sym}\nabla^H X|_H^2+
   \operatorname{Ric}_{H,\mu_0}(X,X)\right)\dd\mu_0.
\]
One can also see this identity directly in the endpoint calculation.
In the average of the two squares, the covariant derivative is
represented by
$M_X+H^{-1/2}\D^3\Phi[X,\cdot,\cdot]H^{-1/2}/2$, and the remaining
square is Kolesnikov's cubic term. The endpoint proof also retains the
individual bounds and the two boundary contributions.

\begin{corollary} 
\label[corollary]{cor:cost-second-variation}
For $i=0,1$, let $\Omega_i\subseteq\R^d$ be a bounded uniformly convex
domain with $C^{3,\alpha}$ boundary, $0<\alpha<1$. Let
$\mu_i=e^{-V_i}\dx$ be a probability measure on $\Omega_i$ with
$V_i$ smooth up to the boundary and $\D^2V_i\succeq\kappa_i I$,
$\kappa_i>0$. Write $T$ for the Brenier map from $\mu_0$ to $\mu_1$.
For $h_i\in C^\infty(\overline{\Omega_i})$ with $\int h_i\dd\mu_i=0$,
set $\mu_{i,t}=(1+t h_i)\mu_i$ for $|t|$ small and put
$q=h_1\circ T-h_0$. Then
\begin{equation}\label{eq:cost-second-variation}
 0\le \frac12\left.\frac{\diff^2}{\diff t^2}W_2^2(\mu_{0,t},\mu_{1,t})\right|_{t=0}
 \le \frac{1}{\sqrt{\kappa_0\kappa_1}}\int_{\Omega_0} q^2\dd\mu_0 .
\end{equation}
For each $d$ and each pair $\kappa_0,\kappa_1>0$, the coefficient is
optimal over all such domains and measures.
\end{corollary}

\begin{proof}
Let $u$ solve the Neumann problem $Bu=q$, and set $H=\D T$ and
$X=H^{-1}\nabla u$. By Caja-Lopez, Delgadino and Kitagawa
\cite[Theorems~1.4--1.5 and Section~6]{CajaDelgadinoKitagawa},
$T,u\in C^2(\overline{\Omega_0})$ and the second variation in
\eqref{eq:cost-second-variation} equals $\mathcal E_H(u)\ge0$.
The conormal condition makes $X$ tangent to the boundary. Since
$HX=\nabla u$, the skew term vanishes, so \eqref{eq:transport-div-curl}
gives $\sqrt{\kappa_0\kappa_1}\,\mathcal E_H(u)\le\int q^2\dd\mu_0$.

For sharpness let $\mu_0$ be $\N(0,\kappa_0^{-1}I)$ conditioned on
$B_R$ and $\mu_1=(s\Id)_\#\mu_0$, where $s=\sqrt{\kappa_0/\kappa_1}$.
Then $H=sI$. Take $h_0=-x_1$ and $h_1=0$. Testing with $q=x_1$ and
applying H\"older's inequality gives
\[
 \frac{\mathcal E_H(u)}{\int q^2\dd\mu_0}
 \ge\frac{\int q^2\dd\mu_0}{\mathcal E_H(q)}
 =s\Var_{\mu_0}(x_1)\longrightarrow(\kappa_0\kappa_1)^{-1/2}
 \qquad(R\to\infty).
\]
\end{proof}

When both marginals are fixed, perturbing the cost produces a field
with zero weighted divergence. The other term in
\eqref{eq:transport-div-curl} will therefore control the second
variation. We use the cost $-x\cdot y$ in the next statement since
it differs from $|x-y|^2/2$ by terms whose integrals are fixed by
the marginals.

\begin{corollary} 
\label[corollary]{cor:ground-cost-variation}
Let $d\ge2$, and let $\Omega_i,\mu_i,\kappa_i$ be as in
\cref{cor:cost-second-variation}. Write $T=\nabla\Phi$ for the
Brenier map, $H=\D T$, and assume $T$ extends to a $C^{2,\alpha}$
diffeomorphism of the closures. For
$a\in C^\infty(\overline{\Omega_0}\times\overline{\Omega_1})$ define
\[
 \mathcal C(\varepsilon)=\inf_{\pi\in\Pi(\mu_0,\mu_1)}
 \int \bigl(-x\cdot y+\varepsilon a(x,y)\bigr)\dd\pi(x,y).
\]
Then $\mathcal C$ is twice differentiable at zero and
\begin{equation}\label{eq:ground-cost-variation}
 0\le-\mathcal C''(0)\le\frac{1}{\sqrt{\kappa_0\kappa_1}}
 \int_{\Omega_0}\norm{\operatorname{Skew}\!\left(
 H^{-1/2}\D^2_{xy}a(x,T(x))H^{1/2}\right)}_{\HS}^2\dd\mu_0 .
\end{equation}
The coefficient is optimal for each $d\ge2$ and each pair
$\kappa_0,\kappa_1>0$.
\end{corollary}

\begin{proof}
Put $w(x)=\nabla_xa(x,T(x))$. The potentials $\Phi,\Phi^*$ are
$C^{3,\alpha}$ up to the boundary, and uniform convexity of $\Phi$
gives a strict duality gap off the graph. The formula of
Peyr\'e, Poon and Tron~\cite[Theorem~2.4 and Remark~2.5]{PeyrePoonTron} gives
\[
 -\mathcal C''(0)=\min_{\psi\in H^1(\Omega_0)}
 \int_{\Omega_0}\inner{H^{-1}(w+\nabla\psi)}{w+\nabla\psi}\dd\mu_0 .
\]
For the $C^2$ minimizer $\psi$, the Euler equation makes
$Z=H^{-1}(w+\nabla\psi)$ divergence-free for $\mu_0$ and tangent to
the boundary. Differentiating $HZ=w+\nabla\psi$ along the graph gives
$\operatorname{Skew}M_Z=\operatorname{Skew}(H^{-1/2}\D^2_{xy}a(x,T(x))H^{1/2})$,
so \eqref{eq:transport-div-curl} yields \eqref{eq:ground-cost-variation}.

For the conditioned Gaussians in the proof of \cref{cor:cost-second-variation},
take $a(x,y)=x^{\mathsf T}Jy$ with $J^{\mathsf T}=-J\ne0$. Since
$H^{-1}w=Jx$ is tangent to the sphere and divergence-free for $\mu_0$,
$\psi=0$ attains the minimum. The squared skew norm is $\norm J_{\HS}^2$,
and $-\mathcal C''(0)=s\int|Jx|^2\dd\mu_0
\to(\kappa_0\kappa_1)^{-1/2}\norm J_{\HS}^2$ as $R\to\infty$.
\end{proof}

\clearpage
\part{Transport Stability}
\section{The squared Wasserstein distance along the heat flow}\label{sec:heat-flow}
Let $\gamma_t=\N(0,tI)$ for $t>0$ and $\gamma_0=\delta_0$. Take an
optimal coupling $(X,Y)$ of $\mu$ and $\nu$ and an independent
$G\sim\gamma_t$. Then $(X+G,Y+G)$ couples
$\mu*\gamma_t$ and $\nu*\gamma_t$ at cost $\E|X-Y|^2=W_2^2(\mu,\nu)$.
Minimizing over couplings of the convolved measures gives
$W_2(\mu*\gamma_t,\nu*\gamma_t)\le W_2(\mu,\nu)$. By the semigroup
property the same comparison holds starting at any time. For log-concave
measures we show that the drop in squared distance over an interval of
fixed length is itself nonincreasing in the starting time.

\begin{theorem}\label[theorem]{heat:thm:convexity}
Let $\mu,\nu\in\Ptwo(\R^d)$ be log-concave. Then
\[
 t\longmapsto D(t):=W_2^2(\mu*\gamma_t,\nu*\gamma_t)
\]
is convex on $[0,\infty)$.
\end{theorem}

To compare the rates of decrease at different times, we need to follow
the optimal map between the two evolving measures. Its differential
already determines the first derivative of $D$, and for strongly
log-concave data this yields a sharp contraction estimate. The proof
of time convexity will use the same differential once more, this time
through the Dirichlet energies of its entries.

\subsection{First variation and sharp contraction}
Write $\mu^c$ for the translate of $\mu$ with mean zero. The identity
\[
 W_2^2(\mu,\nu)=W_2^2(\mu^c,\nu^c)
          +\left|\int x\dd\mu-\int x\dd\nu\right|^2
\]
follows by centering both coordinates of any coupling. Since convolution
with $\gamma_t$ leaves each mean unchanged, all the dissipation comes
from the centered term. We estimate it using \cref{mm:prop:weighted}.

\begin{corollary}
\label[corollary]{cor:heat-contraction}
Let $\mu_i=e^{-V_i}\dx$ be $\kappa_i$-strongly log-concave probability
measures on $\R^d$, with $\kappa_i>0$, and put $a_i=\kappa_i^{-1}$.
Then, for every $s>0$,
\begin{equation}\label{eq:heat-contraction}
 W_2(\mu_0^c*\gamma_s,\mu_1^c*\gamma_s)
 \le \frac{\sqrt{a_0}+\sqrt{a_1}}
          {\sqrt{a_0+s}+\sqrt{a_1+s}}\,W_2(\mu_0^c,\mu_1^c).
\end{equation}
The factor is optimal for every prescribed $s,\kappa_0,\kappa_1>0$
and every $d\ge1$. Equality at some $s>0$ holds if and only if
$\mu_0^c=\mu_1^c$, or $a_0\ne a_1$ and there exist a nonzero
subspace $E$ and a centered probability measure $\sigma$ on $E^\perp$
such that $\mu_i^c=\N_E(0,a_iI_E)\otimes\sigma$ for $i=0,1$.
In either case equality holds for every $s>0$.
\end{corollary}

\begin{proof}
By translating the initial measures we may assume both are centered.
Write $\mu_{i,t}=\mu_i*\gamma_t=e^{-V_{i,t}}\dx$ and
$D(t)=W_2^2(\mu_{0,t},\mu_{1,t})$. For $t>0$, the Gaussian convolution
bound in \cite[equation~(3.14)]{HanZhu} gives
\[
 (a_i+t)^{-1}I\preceq\D^2 V_{i,t}\preceq t^{-1}I .
\]
These bounds give a smooth Brenier diffeomorphism $T_t$ between the
smoothed measures, with both $H_t=\D T_t$ and $H_t^{-1}$ bounded
\cite[Lemma~3.3]{HanZhu}. To differentiate their distance, write the
two heat equations as continuity equations with velocities
$\nabla V_{i,t}/2$. The curves are locally absolutely continuous on
$(0,\infty)$, so the first variation formula~\cite[Ch.~8]{AGS} gives,
for almost every $t>0$,
\begin{align}
 D'(t)
 &=\int(x-T_t(x))\cdot\nabla V_{0,t}(x)\dd\mu_{0,t}
   +\int(y-T_t^{-1}(y))\cdot\nabla V_{1,t}(y)\dd\mu_{1,t}\notag\\
 &=-\int\tr(H_t+H_t^{-1}-2I)\dd\mu_{0,t}.
 \label{heat:eq:first-variation}
\end{align}
To pass to the second line, integrate by parts at the source and the
target, using the bounded differentials of $T_t$ and $T_t^{-1}$ and
the Gaussian tails of the two measures. The resulting trace integral
is precisely the quantity bounded in \eqref{mm:eq:R-lower}, which gives
\[
 D'(t)\le-\frac{D(t)}{\sqrt{(a_0+t)(a_1+t)}} .
\]
The positive coefficient multiplying $-D(t)$ is the logarithmic derivative
of $A(t)^2$,
where $A(t)=\sqrt{a_0+t}+\sqrt{a_1+t}$. Hence $A(t)^2D(t)$ is
nonincreasing for $t>0$, and continuity in $W_2$ at zero gives
\eqref{eq:heat-contraction}.

Suppose $D(0)>0$ and equality holds at $s>0$. Then $A(t)^2D(t)$ is
constant on $[0,s]$, so equality holds in \eqref{mm:eq:R-lower} for
almost every $t\in(0,s)$. At any such time \cref{mm:prop:weighted} gives
$a_0\ne a_1$ and a nonzero subspace $E$ such that
\[
 \mu_{i,t}=\N_E(0,(a_i+t)I_E)\otimes\sigma_t,\qquad i=0,1.
\]
Gaussian convolution is injective, since its Fourier multiplier never
vanishes. Applied on $E^\perp$, this shows that the two initial marginals
there coincide, say with $\sigma$. On $\R^d$, both $\mu_i$ and
$\N_E(0,a_i I_E)\otimes\sigma$ convolve to $\mu_{i,t}$, giving the
asserted factorization.

For these products,
$D(t)=\dim(E)\bigl(\sqrt{a_0+t}-\sqrt{a_1+t}\bigr)^2$,
so $A(t)^2D(t)$ is constant. The case $D(0)=0$ is immediate.
Taking $E=\R^d$ proves sharpness when $a_0\ne a_1$. If $a_0=a_1=a$,
the distance ratio for $\N(0,aI)$ and $\N(0,(a-\varepsilon)I)$ tends
to $\sqrt{a/(a+s)}$ as $\varepsilon\downarrow0$.
\end{proof}

\subsection{Second variation and time convexity}
To prove \cref{heat:thm:convexity}, we first take compactly supported initial densities of the form
$e^{-\kappa|x|^2/2}g(x)$, where $\kappa>0$ and $g$ is log-concave.
The two endpoints may have different $\kappa$. Fix $t>0$ and write
\[
 \mu_t=\mu*\gamma_t=e^{-V}\dx,\qquad
 \nu_t=\nu*\gamma_t=e^{-W}\dy .
\]
Both potentials have Hessians bounded above and below by positive
constants, as in the proof of \cref{cor:heat-contraction}. On compact
subintervals of $(0,\infty)$, the densities and their space and time
derivatives have Gaussian decay by differentiation of the heat kernel.
The optimal transport map $T=\nabla\Phi$ from $\mu_t$ to $\nu_t$ and its inverse
$S=\nabla\Phi^*$ are smooth, and $H=\D T$ and $K=H^{-1}$ are bounded.

For this regularized pair, we will bound the second time difference
by the spatial derivatives of the transport Hessian. More precisely,
we claim that
\begin{equation}\label{heat:eq:quantitative}
 \liminf_{h\downarrow0}\frac{D(t+h)+D(t-h)-2D(t)}{h^2}
 \ \ge\ \sum_{i,j=1}^d \mathcal E_H(H_{ij}),\qquad
 \mathcal E_H(f)=\int\inner{K\nabla f}{\nabla f}\dd\mu_t .
\end{equation}
To interpret the right side, recall from \eqref{heat:eq:first-variation} that
$-D'(t)=\sum_i\mathcal E_H((T-\Id)_i)$ almost everywhere.
Thus $-D'(t)$ is the energy of the displacement, and the energies
of its first spatial derivatives bound the lower limit of the second
difference quotients in \eqref{heat:eq:quantitative} from below. For a single transport between smooth
densities, Kolesnikov obtained an $L^2$ bound for the Hilbert--Schmidt
norm of $\D^2\Phi$ in terms of the source Fisher information under a
uniform positive lower curvature bound on the target~\cite{KolesnikovSobolev}.

Here the heat equation acts at both endpoints. When we vary the two
Kantorovich potentials, their derivatives are linked by Legendre
duality. As \cref{heat:lem:variation} shows, optimizing over that variation amounts
to approximating a field built from the two scores by gradients on the
target. In target coordinates $y=T(x)$ we evaluate $H$ and $K$ at
$S(y)$, so that $K=\D_y^2\Phi^*$ there. In particular, the squared norm
used for this approximation is
$\|v\|_{H,\nu_t}^2=\int v^{\mathsf T}(H\circ S)v\dd\nu_t$.

\begin{lemma}\label[lemma]{heat:lem:variation}
At this fixed time denote by $\mathcal P$ the squared distance of
$q(y)=\nabla V(S(y))-\D^2\Phi^*(y)\nabla W(y)$ from gradients on the target,
\[
 \mathcal P=\inf_{\varphi\in C_c^\infty(\R^d)}
                      \norm{q-\nabla\varphi}_{H,\nu_t}^{\,2}.
\]
On the source put $M=H-K$ and
$Q=\int M:\bigl(\D^2V-\D^2W(T)\bigr)\dd\mu_t$,
where $C:E=\tr(C^{\mathsf T}E)$. Then
\begin{equation}\label{heat:eq:dual}
 2\liminf_{h\downarrow0}
 \frac{D(t+h)+D(t-h)-2D(t)}{h^2}\ge Q-\mathcal P .
\end{equation}
Moreover, $\D^3\Phi\in L^2(\mu_t)$. Taking derivatives in fixed
coordinates and then diagonalizing $H$ at each point, with eigenvalues
$\lambda_i$, we have
\begin{equation}\label{heat:eq:Q-split}
 Q=\int\left[\D^2W(T):HM^2+
 \sum_{i,j,k=1}^d\left(
 \frac{\lambda_i}{\lambda_j\lambda_k}+\frac1{\lambda_i}
 \right)|\partial_{ijk}\Phi|^2\right]\dd\mu_t .
\end{equation}
\end{lemma}

\begin{proof}
It suffices to bound $D(t+u)$ from below by a smooth function agreeing
with it at $u=0$. To do so, vary the optimal dual potential at time $t$
by setting $\Phi_u=\Phi-uf/2$ for $f\in C_c^\infty(\R^d)$. Uniform
convexity persists for small $|u|$, so Kantorovich duality gives
\[
 D(t+u)\ge L_f(u):=
 \int(|x|^2-2\Phi_u)\dd\mu_{t+u}
 +\int(|y|^2-2\Phi_u^*)\dd\nu_{t+u},
 \qquad D(t)=L_f(0).
\]
Differentiating $y=\nabla\Phi_u(x)$ at its unique solution gives
\[
 \left.\partial_u\Phi_u^*\right|_{u=0}=\tfrac12 f\circ S,\qquad
 \left.\partial_u^2\Phi_u^*\right|_{u=0}
 =\tfrac14(\nabla f^{\mathsf T}K\nabla f)\circ S .
\]
For small $|u|$ the changes in both potentials are supported in fixed
compact sets, since $\nabla\Phi_u=\nabla\Phi$ off $\supp f$. Together
with the quadratic growth of the unperturbed potentials and the Gaussian
bounds on the heat densities and their derivatives, this justifies two
differentiations under the integrals. The second moments are affine
in $u$ and hence drop out. Using the heat equation for the remaining
terms, integrating those involving $f$ once by parts and completing the
square on the target gives
\begin{align*}
 2L_f''(0)
 ={}&-\int\Phi\,\Delta^2(e^{-V})\dx-\int\Phi^*\,\Delta^2(e^{-W})\dy\\
 &+\norm q_{H,\nu_t}^2-\norm{q-\nabla(f\circ S)}_{H,\nu_t}^2 .
\end{align*}
Here $\nabla(f\circ S)=(K\nabla f)\circ S$. Two integrations by parts
give $\int\Phi\,\Delta^2(e^{-V})\dx=\int H:\D^2(e^{-V})\dx$, and
likewise at the target. To justify these integrations, take a cutoff $\chi_R$ supported on
$B_{2R}$, equal to one on $B_R$, with $0\le\chi_R\le1$ and
$|\D^j\chi_R|\le C_jR^{-j}$ for $j=1,2$. The potentials have at most
quadratic growth and their gradients at most linear growth, so the
Gaussian decay of the density derivatives makes the cutoff errors vanish.

Expanding $\D^2(e^{-V})$ and $\D^2(e^{-W})$, the first three terms
on the right sum to
\[
 \int\!\left(H:\D^2V+K:\D^2W(T)
             -2\nabla V\cdot\nabla W(T)\right)\dd\mu_t .
\]
Integration by parts at the two endpoints gives
\[
 \int H:\D^2W(T)\dd\mu_t
 =\int\nabla V\cdot\nabla W(T)\dd\mu_t
 =\int K:\D^2V\dd\mu_t .
\]
Thus $2L_f''(0)=Q-\norm{q-\nabla(f\circ S)}_{H,\nu_t}^2$. Since
$f\circ S$ ranges over $C_c^\infty(\R^d)$, taking the supremum in
the symmetric second difference bound proves \eqref{heat:eq:dual}.

To express $Q$ in terms of derivatives of $H$, differentiate the change
of variables formula $\log\det H=W(T)-V$. The operator that acts on
each entry of $H$ is $L=e^V\diver_x(e^{-V}K\nabla_x)$, since the Piola
identity $\diver_x((\det H)K)=0$ gives
$L=K:\D_x^2-\nabla W(T)\cdot\nabla_x$. Taking two derivatives of
the change of variables formula and commuting the derivatives of $\Phi$
now yields
\begin{equation}\label{heat:eq:LH}
 LH=H\D^2W(T)H-\D^2V+G,\qquad
 G_{ij}=\tr\bigl(K(\partial_{x_i}H)K(\partial_{x_j}H)\bigr).
\end{equation}
The multiplier $M=H-K$ in $Q$ is the derivative of
$H\mapsto\tfrac12\tr(H^2)-\log\det H$. We may therefore use the chain
rule in \eqref{heat:eq:LH} to integrate the term $M:LH$. With
$HMH-M=HM^2$, the resulting identity is
\begin{align*}
 &M:(\D^2V-\D^2W(T))+L\bigl(\tfrac12\tr(H^2)-\log\det H\bigr)\\
 &\qquad=\D^2W(T):HM^2+H:G
           +\sum_{i,j}K_{ij}(\partial_{x_i}H):(\partial_{x_j}H).
\end{align*}
Here the second derivative of $-\log\det H$ contributes $K:G$,
which combines with $M:G$ to give $H:G$. Each term on the right is
nonnegative. Choose the cutoffs above increasing to one. Since $\nabla W(T)$
has at most linear growth, $L\chi_R$ is bounded uniformly and supported on
$B_{2R}\setminus B_R$. Integrate against $\chi_R$ and move $L$ onto the
cutoff by symmetry. Since $\tfrac12\tr(H^2)-\log\det H$ is bounded, this
term tends to zero. Monotone convergence on the right gives
\eqref{heat:eq:Q-split}. Indeed, in a basis diagonalizing $H$,
the last two terms have coefficients $\lambda_i/(\lambda_j\lambda_k)$ and
$1/\lambda_i$, respectively, on $|\partial_{ijk}\Phi|^2$. The latter also
gives $\D^3\Phi\in L^2(\mu_t)$ since $H$ is bounded.
\end{proof}

The lower bound in \eqref{heat:eq:dual} still contains the squared
distance $\mathcal P$. We estimate how much of the positive expression
in \eqref{heat:eq:Q-split} can be lost through this term.
\begin{lemma}\label[lemma]{heat:lem:projection}
With the notation above, denote 
\[
 J=\frac12\sum_{i,j,k}
       \frac{(\lambda_i-\lambda_j)^2}{\lambda_i\lambda_j\lambda_k}
       |\partial_{ijk}\Phi|^2.
\]
Then
\begin{equation}\label{heat:eq:projection-bound}
 \mathcal P\le\int\bigl(\D^2W(T):HM^2+J\bigr)\dd\mu_t .
\end{equation}
\end{lemma}

\begin{proof}
In one dimension $\mathcal P=0$. Indeed, a field $v$ orthogonal to
gradients in $L^2(H,\nu_t)$ satisfies $(e^{-W}Hv)'=0$, and any
nonzero value of this constant would give $\int H|v|^2\dd\nu_t=\infty$.
In higher dimensions we estimate $\mathcal P$ by pairing $q$ with fields
in the orthogonal complement of gradients. To construct a dense family
of such fields, work in target coordinates and write
$\partial_j=\partial_{y_j}$ and
\[
 (\diff_y v)_{ij}=\partial_i v_j-\partial_j v_i,\qquad
 (\delta F)_i=\sum_j\delta_jF_{ij},\qquad
 \delta_j=-\partial_j+\partial_j W .
\]
For each smooth compactly supported one-form $\alpha$ take the skew field
\[
 F=H(\diff_y\alpha)H,\qquad Z=\delta F .
\]
We claim that these $KZ$ are dense in the orthogonal complement of gradients
in $L^2(H,\nu_t)$. They belong to this complement since
$\int\nabla\varphi\cdot\delta F\dd\nu_t=\int\D^2\varphi:F\dd\nu_t=0$.
If $v$ in the complement is orthogonal to every $KZ$, integration by
parts in distributions gives
\[
 \diver_{\nu_t}(Hv)=0,\qquad \delta\bigl(H(\diff_y v)H\bigr)=0 .
\]
To test these equations against a cutoff of $v$ itself, we first
establish interior smoothness. Freeze $H$ at $y$. The principal symbol
of $\delta(H\diff_y v H)$ at frequency $\xi$, acting on a vector $b$, is
$(H\xi)(\xi^{\mathsf T}Hb)-(\xi^{\mathsf T}H\xi)Hb$.
Adding $H\nabla\diver_{\nu_t}(Hv)$ cancels the first term.
The negative of the resulting operator has principal symbol
$(\xi^{\mathsf T}H\xi)H$, positive definite for $\xi\ne0$.
Interior elliptic regularity yields $v\in C^\infty$
\cite[Theorem~15.23]{DyatlovElliptic}. Test the second distributional
equation with $\chi_R^2v$. Uniform ellipticity and Cauchy--Schwarz
applied to the term containing $\diff\chi_R\wedge v$ give
\[
 \int\chi_R^2|\diff_y v|^2\dd\nu_t\le CR^{-2}\int|v|^2\dd\nu_t .
\]
Thus $\diff_y v=0$ and $v=\nabla u$ on $\R^d$. Smoothly truncating $u$
at height $N$ and multiplying by $\chi_R$ approximates $v$ by compactly
supported gradients. The cutoff term has norm $O(N/R)$, so we may send
$R\to\infty$ and then $N\to\infty$. Hence $v=0$.
The Hilbert space projection formula now gives
\begin{equation}\label{heat:eq:projection-duality}
 \mathcal P=\sup_{\alpha}
   \left(2\int q\cdot Z\dd\nu_t-\norm{Z}_{K,\nu_t}^2\right).
\end{equation}

We express both terms in \eqref{heat:eq:projection-duality} in terms
of the skew-symmetric matrix field $F$. For the linear term, differentiating $\log\det H=W(T)-V$ gives
$q=M\nabla W-(\nabla_x\log\det H)\circ S$.
By the chain rule and Hessian symmetry,
$\sum_j\partial_j H_{ij}=(\partial_{x_i}\log\det H)\circ S$ and
$\sum_j\partial_j K_{ij}=\partial_i\tr K$. Hence
$q_i=\sum_k\delta_k M_{ik}-\partial_i\tr K$.
Write $F_{ij,k}=\partial_k F_{ij}$. The gradient
of $\tr K$ pairs to zero with the compactly supported $Z$.
Using $\partial_j\delta_k=\delta_k\partial_j+\partial_{jk}W$ therefore gives
\begin{equation}\label{heat:eq:skew-numerator}
 \int q\cdot Z\dd\nu_t
 =\int\left(F:(M\D^2W)+\sum_{i,j,k}F_{ij,k}\,\partial_j M_{ik}\right)\dd\nu_t.
\end{equation}
For the quadratic term, substitute $Z=\delta F$ and integrate by parts
twice. This gives $\norm{Z}_{K,\nu_t}^2=\int(D_0+D_1)\dd\nu_t$, where
\[
 D_0=\sum_{i,j,k,l}K_{ij}(\partial_{kl}W)F_{ik}F_{jl},\qquad
 D_1=\sum_{i,j,k,l}K_{ij}F_{ik,l}F_{jl,k}.
\]
Indeed, moving $\delta_k$ from $F_{ik}$ onto $K_{ij}Z_j$ yields a term
containing $\sum_{i,k}F_{ik}\partial_k K_{ij}=0$. Next commute
$\partial_k$ past $\delta_l$ and move $\delta_l$ back onto $K_{ij}F_{ik}$.
The second coefficient derivative vanishes because
$\sum_{j,l}(\partial_lK_{ij})F_{jl,k}=0$. Both cancellations use
$K=\D_y^2\Phi^*$.

The expression for $D_1$ becomes a sum of squares because of our
choice of $F$. Indeed, $KFK=\diff_y\alpha$ is closed, and its pullback
by $T$, with coefficients $F_{ij}\circ T$, is closed as well. Writing
this condition in source coordinates and diagonalizing $H$ at $x=S(y)$ gives
\begin{equation}\label{heat:eq:closed-form}
 \lambda_iF_{jk,i}+\lambda_jF_{ki,j}+\lambda_kF_{ij,k}=0 .
\end{equation}
Multiply by $F_{ij,k}/(\lambda_i\lambda_j)$ and sum. By skew symmetry
the first two terms each sum to $-D_1$, so
\begin{equation}\label{heat:eq:skew-alternation}
 D_1=\frac12\sum_{i,j,k}\frac{\lambda_k}{\lambda_i\lambda_j}F_{ij,k}^2 .
\end{equation}
Pair $K^{1/2}F(\D^2W)^{1/2}$ and $H^{1/2}M(\D^2W)^{1/2}$ in the
Hilbert--Schmidt inner product. This gives
$|F:(M\D^2W)|^2\le D_0(\D^2W:HM^2)$.
For the other term in \eqref{heat:eq:skew-numerator}, $\partial_jK_{ik}$
is symmetric in all three indices and
$\partial_jH_{ik}=\lambda_j^{-1}\partial_{ijk}\Phi(x)$. Hence
\[
 \left|\sum_{i,j,k}F_{ij,k}\partial_jM_{ik}\right|
 =\left|\frac12\sum_{i,j,k}
    (\lambda_j^{-1}-\lambda_i^{-1})F_{ij,k}\,\partial_{ijk}\Phi(x)\right|
 \le\sqrt{D_1J(x)} .
\]
Since $\int(D_0+D_1)\dd\nu_t=\norm{Z}_{K,\nu_t}^2$, H\"older's
inequality gives
\begin{equation}\label{heat:eq:skew-estimate}
 \left|\int q\cdot Z\dd\nu_t\right|^2
 \le\norm{Z}_{K,\nu_t}^2\int\bigl(\D^2W(T):HM^2+J\bigr)\dd\mu_t .
\end{equation}
Taking the supremum in \eqref{heat:eq:projection-duality} proves
\eqref{heat:eq:projection-bound}.
\end{proof}

\begin{proof}[Proof of \cref{heat:thm:convexity}]
In the regular case, the symmetry of $\D^3\Phi$ gives
\[
 \frac12\sum_{i,j,k}
 \frac{(\lambda_i-\lambda_j)^2}{\lambda_i\lambda_j\lambda_k}|\partial_{ijk}\Phi|^2
 =\sum_{i,j,k}\left(\frac{\lambda_i}{\lambda_j\lambda_k}
                      -\frac1{\lambda_i}\right)|\partial_{ijk}\Phi|^2 .
\]
Subtract the bound of \cref{heat:lem:projection} from
\eqref{heat:eq:Q-split} and apply \eqref{heat:eq:dual}. The remaining
sum is $2\int\sum_{i,j,k}\lambda_i^{-1}|\partial_{ijk}\Phi|^2\dd\mu_t$,
which proves \eqref{heat:eq:quantitative}.

It remains to deduce convexity from the lower symmetric second
derivative, without assuming that $D$ is twice differentiable. The function $D$ is continuous by
$W_2(\mu_s,\mu_{s+h})\le\sqrt{dh}$ and the same bound for $\nu$.
If it exceeded a chord $\ell$ on $[a,b]\subseteq(0,\infty)$, then
$D(t)-\ell(t)+\varepsilon(t-a)(t-b)$ would have a positive interior
maximum for small $\varepsilon>0$. Its lower symmetric second
derivative is at least $2\varepsilon$ there, whereas every symmetric
second difference at a maximum is nonpositive.
Continuity extends convexity to $t=0$.

For general log-concave inputs, convolve first with $\gamma_\varepsilon$
so that each measure has a density, even if its original support is
lower dimensional. Multiplying each density by
$\mathbf 1_{B_R}(x)e^{-|x|^2/(2R^2)}$ and normalizing gives a pair in
the regular class. As $R\to\infty$ these densities converge in $W_2$ to
the smoothed inputs by dominated convergence of the densities and their
second moments. Letting $\varepsilon\downarrow0$ then recovers the
original measures in $W_2$. Common Gaussian convolution is nonexpansive,
so the midpoint inequality at $s,s+h,s+2h$ passes to both limits.
\end{proof}

\begin{remark}
For regular inputs, vanishing of the second time variation has a geometric
consequence. Fix $t>0$ and take $\kappa>0$ with $\D^2W\succeq\kappa I$.
In \eqref{heat:eq:skew-numerator} we may replace $\D^2W$ by
$\D^2W-\kappa I$ since $F:M=0$. This replacement only decreases $D_0$
in the Cauchy--Schwarz bound. Thus the right side of
\eqref{heat:eq:quantitative} gains the term
\[
 \frac{\kappa}{2}\int\tr\bigl(H(H-H^{-1})^2\bigr)\dd\mu_t .
\]
If the lower symmetric second derivative vanishes at this time, then
$H=I$ and $T$ is a translation. Injectivity of Gaussian convolution
gives the same translation between the initial measures. Conversely,
for translates the squared distance is constant. The criterion here
uses a second time variation at one positive time. Fathi, Goldman and
Tsodyks obtain a related rigidity estimate from the drop in squared
distance over a time interval when one endpoint is Gaussian
\cite[Theorem~5.10]{FathiGoldmanTsodyks}.
\end{remark}

\begin{remark}\label[remark]{heat:rem:two-point}
Without log-concavity, time convexity can fail even on $\R$. Take $\mu$
and $\nu$ to be the laws of $\xi$ and $2\xi$, where $\xi$ is uniform on
$\{-1,1\}$. The monotone coupling $(\xi,2\xi)$ gives $D(0)=1$.
For an independent $G\sim\N(0,1)$, the coupling
$(\xi+\sqrt t\,G,2\xi+\sqrt t\,G)$ still has cost one. For $t>0$ its
support consists of the two lines $y=x\pm1$ and contains $(0,1)$ and
$(1,0)$. An optimal coupling on the line has monotone support, so $D(t)<1$.

To rule out convexity it suffices to show that this strict decrease has
zero initial slope. Test with the $1$-Lipschitz function $x\mapsto |x|$.
By symmetry and the fact that $|2+z|-|1+z|$ equals $1$ for $z\ge-1$
and is at least $-1$ elsewhere, we obtain
\[
 1-2\mathbb P(G<-t^{-1/2})
 \le \E|2+\sqrt t\,G|-\E|1+\sqrt t\,G|
 \le \sqrt{D(t)}<1 .
\]
The Gaussian tail is $o(t)$ as $t\downarrow0$, so $D'_+(0)=0$.
Convexity would then give $D(t)\ge D(0)+tD'_+(0)=1$ for every $t>0$.
For small noise estimates with discrete inputs see~\cite{DingNilesWeed}.
\end{remark}

\section{Moment measures and inverse stability}\label{sec:moment-measures}
A centered density $\rho=e^{-\psi}\dx$ determines its moment measure
by $\nu=(\nabla\psi)_\#\rho$. We ask how a change in $\nu$ affects
$\rho$ when this relation is inverted. The connection with heat flow
comes from the optimal map $T$ between two reconstructed sources.
The integral of $\tr(\D T+(\D T)^{-1}-2I)$ against the first source
gave the instantaneous decrease of squared distance in
\eqref{heat:eq:first-variation}. Here we compare source and target couplings to bound the same integral
in terms of the distance between the moment measures.
For reconstruction in a fixed confining potential, the response to a
moving target in \cref{thm:hessian-transfer} will instead give the
stability estimate.

\subsection{Sharp inverse stability and equality cases}
\label{sec:moment-comparison}\label{sec:moment-rigidity}

Let $\nu_i$ and $\rho_i=\Reconstruct\nu_i$ be as in
\cref{mm:thm:main}. There are two natural ways to couple the sources.
One minimizes the source cost directly. The other pulls an optimal
target plan back through the moment maps. The difference between their
costs will account for one of the nonnegative terms in
\eqref{mm:eq:finite-defect}.
For a related use of excess transport cost in multi-marginal
stability, see \cite[Sections~2.3--2.4]{HanZhuMultimarginal}.

We first assume
\begin{equation}\label{mm:eq:regular}
 \psi_i\in C^\infty(\R^d),\qquad
 \kappa_iI\preceq \D^2 \psi_i\preceq MI,\qquad i=0,1,
\end{equation}
with $M<\infty$. By the transport regularity under two-sided
Hessian bounds in \cite[Lemma~3.3]{HanZhu}, the map
$T:\rho_0\to\rho_1$ is a smooth convex-gradient diffeomorphism,
and its differential and inverse differential are bounded. Set
$H=\D T$ and $c=\sqrt{\kappa_0\kappa_1}$, so that we can apply
\cref{mm:prop:weighted} to $T$. The maps $\nabla\psi_i$ also are
global diffeomorphisms by strong convexity. We may therefore pull an
optimal plan between $\nu_0$ and $\nu_1$ back to a plan
$\pi\in\Pi(\rho_0,\rho_1)$ under $(\nabla\psi_0,\nabla\psi_1)$.
For $(x,z)$ sampled from $\pi$, set $u=T(x)$ and $v=T^{-1}(z)$.
The pairs $(x,u)$ and $(v,z)$ have the same law, namely the optimal
source plan $(\Id,T)_\#\rho_0$. This gives four points at which to
compare the two costs. In particular, the cost of $(x,z)$ satisfies
\begin{equation}\label{3.2}
	 \int|x-z|^2\dd\pi\ge W_2(\rho_0,\rho_1)^2. 
\end{equation}

Two further terms arise from comparing the displacements $u-z$ and
$x-v$ and applying strong convexity to the Bregman divergences
$\Dcal_{\psi_i}$. Write
\begin{align}
 Q&=\int\left|\sqrt{\kappa_1}(u-z)-\sqrt{\kappa_0}(x-v)\right|^2\dd\pi,
 \label{mm:eq:Q-def}\\
 \mathcal B&=\int\left[
 \Dcal_{\psi_1}(u,z)-\frac{\kappa_1}{2}|u-z|^2
 +\Dcal_{\psi_0}(v,x)-\frac{\kappa_0}{2}|x-v|^2
 \right]\dd\pi, \label{mm:eq:B-def}
\end{align}
where the nonnegativity of $\mathcal{B}$ follows from the $\kappa_i$-strong convexity of $\psi_i$.

\begin{proposition}\label[proposition]{mm:prop:finite-defect}
Under \eqref{mm:eq:regular}, the deficit in the inverse stability estimate has the decomposition 
\begin{equation}\label{mm:eq:finite-defect}
 \begin{aligned}
 &W_2(\nu_0,\nu_1)^2-c^2W_2(\rho_0,\rho_1)^2\\
 &\quad=\int|\nabla \psi_0(x)-\nabla \psi_1(z)-c(u-v)|^2\dd\pi
       +2c\bigl(\mathfrak R-cW_2(\rho_0,\rho_1)^2\bigr)\\
 &\qquad+c^2\left(\int|x-z|^2\dd\pi-W_2(\rho_0,\rho_1)^2\right)+cQ+2c\mathcal B.
 \end{aligned}
\end{equation}
In particular,
\begin{equation}\label{mm:eq:control-spectral-defect}
 0\le\mathfrak R-cW_2(\rho_0,\rho_1)^2\le\frac{W_2(\nu_0,\nu_1)^2-c^2W_2(\rho_0,\rho_1)^2}{2c}.
\end{equation}
\end{proposition}

\begin{proof}
Expand the Bregman divergences under $\pi$. The potential values
cancel after integration because $u,z$ have law $\rho_1$ and $v,x$
have law $\rho_0$. The remaining terms pair $\nabla\psi_i$ with the
displacement of the optimal source map or its inverse. We can integrate
these terms by parts with cutoffs, since under \eqref{mm:eq:regular} the
potentials have at most quadratic growth, the gradients and transport
maps have at most linear growth, and the source measures have Gaussian
tails. This gives
	\begin{equation}\label{mm:eq:bregman-identity}
		\begin{aligned}
			&\int(\nabla \psi_0(x)-\nabla \psi_1(z))\cdot(u-v)\dd\pi
			-\int\left[\Dcal_{\psi_1}(u,z)+\Dcal_{\psi_0}(v,x)\right]\dd\pi\\
			=&\int\left[\nabla \psi_0(x)\cdot(u-x)+\nabla \psi_1(z)\cdot(v-z)\right]\dd\pi\\
			\overset{*}{=}&\int\tr(H-I)\dd\rho_0
			+\int\tr(\D T^{-1}-I)\dd\rho_1=\mathfrak R, 
		\end{aligned}
	\end{equation}
where $(*)$ uses the bounded differentials of $T$ and $T^{-1}$ in
the integration by parts. To recover the transport costs from the
cross term in $Q$, observe that
\begin{equation}\label{mm:eq:four-point-algebra}
 |u-v|^2-2(u-z)\cdot(x-v)
 =|u-x|^2+|z-v|^2-|x-z|^2. 
\end{equation}
Expanding $Q$ and integrating \eqref{mm:eq:four-point-algebra} gives
\begin{equation}\label{mm:eq:weighted-four-point}
	\begin{aligned}
		&\kappa_1\int|u-z|^2\dd\pi
		+\kappa_0\int|x-v|^2\dd\pi\\
		=&Q+2c\int(u-z)\cdot(x-v)\dd\pi\\
		=&Q+c\int|u-v|^2\dd\pi
		-2cW_2(\rho_0,\rho_1)^2
		+c\int|x-z|^2\dd\pi.
	\end{aligned}
\end{equation} 
Combining \eqref{mm:eq:B-def}, \eqref{mm:eq:bregman-identity} and
\eqref{mm:eq:weighted-four-point}, and using optimality of
$(\nabla \psi_0,\nabla \psi_1)_\#\pi$, yields 
\[
 \begin{aligned}
 2c\mathfrak R
 =&2c\int(\nabla \psi_0(x)-\nabla \psi_1(z))\cdot(u-v)\dd\pi-c^2\int|u-v|^2\dd\pi\\
 &-cQ+2c^2W_2(\rho_0,\rho_1)^2-c^2\int|x-z|^2\dd\pi-2c\mathcal B\\
 =&W_2(\nu_0,\nu_1)^2
 -\int\bigl|\nabla\psi_0(x)-\nabla\psi_1(z)-c(u-v)\bigr|^2\dd\pi\\
 &-cQ+2c^2W_2(\rho_0,\rho_1)^2
 -c^2\int|x-z|^2\dd\pi-2c\mathcal B, 
 \end{aligned}
\] 
which proves \eqref{mm:eq:finite-defect}. The bound
\eqref{mm:eq:control-spectral-defect} follows from
\cref{mm:prop:weighted} and \eqref{3.2}. 
\end{proof}

\begin{proof}[Proof of \cref{mm:thm:main}]
Under \eqref{mm:eq:regular}, \cref{mm:prop:finite-defect} gives
\eqref{mm:eq:main-estimate}. For general sources, set
\[
\rho_{i,t}=\rho_i*\N(0,tI),\qquad
\nu_{i,t}=\Moment\rho_{i,t},\qquad i=0,1.
\]
By \cref{mm:lem:smoothing}, the smoothed sources remain centered,
and their potentials satisfy \eqref{mm:eq:regular} with $\kappa_{i,t}={\kappa_i}/{(1+\kappa_i t)}$ and $M=t^{-1}$. By \cref{mm:prop:finite-defect}, we obtain 
\[
W_2(\rho_{0,t},\rho_{1,t})
\le \frac{1}{\sqrt{\kappa_{0,t}\kappa_{1,t}}}
W_2(\nu_{0,t},\nu_{1,t}).
\]
Passing to the limit using \eqref{mm:eq:both-converge} proves
\eqref{mm:eq:main-estimate}. For centered isotropic Gaussians with standard deviations $s_0,s_1$,
the source distance is $\sqrt d\,|s_1-s_0|$ and the moment distance is
$\sqrt d\,|s_1^{-1}-s_0^{-1}|$. Their ratio for $s_0\ne s_1$ is $s_0s_1$.
Taking $s_i\uparrow\kappa_i^{-1/2}$ proves sharpness also for equal
curvature parameters. 
\end{proof} 

To apply \cref{mm:thm:main} directly one needs lower curvature bounds
for the reconstructed potentials $\psi_i$, even though the prescribed
data are the moment measures $\nu_i$. An upper Hessian bound for the
target potential $W_i$ supplies such a hypothesis. Indeed, if
$\D^2W_i\preceq\Lambda_i I$ distributionally, the reconstruction is
$\Lambda_i^{-1}$-strongly log-concave.

\begin{corollary}\label[corollary]{mm:cor:inverse}
For $\Lambda_0,\Lambda_1>0$ and $i=0,1$, let
\[
 \nu_i=Z_i^{-1}e^{-W_i}\dy\in\Ptwo(\R^d),\qquad
 \int y\dd\nu_i(y)=0, 
\]
where each $W_i$ is finite on all of $\R^d$ and satisfies
$\D^2W_i\preceq\Lambda_i I$ in the sense of distributions. Then 
\begin{equation}\label{mm:eq:inverse-bound}
 W_2(\Reconstruct\nu_0,\Reconstruct\nu_1)
 \le\sqrt{\Lambda_0\Lambda_1}\,W_2(\nu_0,\nu_1).
\end{equation}
The coefficient is optimal for every $\Lambda_0,\Lambda_1>0$ and every $d\ge1$. 
\end{corollary} 

\begin{proof}
	The measures $\nu_i$ are centered by assumption and have full support because their
	potentials $W_i$ are finite. The moment measure theorem of
	Cordero-Erausquin and Klartag~\cite{CorderoKlartag} therefore gives
	the centered reconstructions $\rho_i=\Reconstruct\nu_i=e^{-\psi_i}\dx$. Moreover the upper Hessian bound places each target
	in the class of \cite[(5.12)]{MachadoRamos}, which allows nonconvex $W_i$.
	By \cite[Theorem~5.10]{MachadoRamos}, $\psi_i$ is
	$\Lambda_i^{-1}$-strongly convex. The moment relation gives
	\[
	\int|\nabla\psi_i|^2\dd\rho_i
	=\int|y|^2\dd\nu_i<\infty.
	\]
	Applying \cref{mm:thm:main} with $\kappa_i=\Lambda_i^{-1}$
	then proves \eqref{mm:eq:inverse-bound}. For optimality, choose 
	\[
	\nu_0=\N(0,\Lambda_0^{-1}I),\qquad
	\nu_1^\varepsilon=\N(0,(1+\varepsilon)\Lambda_1^{-1}I).
	\]
	These targets satisfy the prescribed curvature bounds. Their centered
	reconstructions are $\N(0,\Lambda_0I)$ and
	$\N(0,\Lambda_1(1+\varepsilon)^{-1}I)$, respectively. Since
	$W_2(\N(0,s^2I),\N(0,t^2I))=\sqrt d\,|s-t|$, for all sufficiently
	small $\varepsilon>0$ we obtain
	\[
	\frac{W_2(\Reconstruct\nu_0,\Reconstruct\nu_1^\varepsilon)}
	{W_2(\nu_0,\nu_1^\varepsilon)}
	=\sqrt{\frac{\Lambda_0\Lambda_1}{1+\varepsilon}}
	\longrightarrow\sqrt{\Lambda_0\Lambda_1}
	\qquad \varepsilon\rightarrow0, 
	\]
	which proves optimality in every dimension. 
\end{proof} 

The proof of \cref{mm:thm:main} used only nonnegativity of the terms in
\eqref{mm:eq:control-spectral-defect}. At equality these terms must
vanish in the limit along smooth approximations. We use this stronger
information to identify the directions in which two distinct
reconstructions can differ.

\begin{theorem}\label[theorem]{mm:thm:equality}
Under the assumptions of \cref{mm:thm:main}, suppose that $\rho_0\ne\rho_1$.
Then equality holds in \eqref{mm:eq:main-estimate} if and only if there
exist a nonzero subspace $E\subseteq\R^d$ and a centered probability
measure $\sigma$ on $E^\perp$ such that
\begin{equation}\label{mm:eq:equality-products}
 \rho_i=\N_E(0,\kappa_i^{-1}I_E)\otimes\sigma,\qquad i=0,1.
\end{equation}
In this case, $\kappa_0\ne\kappa_1$ and the pair $(\rho_0,\rho_1)$
determines $E$ uniquely.
\end{theorem}

\begin{proof}
Suppose equality holds. To compare the deficits along a smooth
approximation we keep the lower curvature bounds fixed. After convolving
$\rho_i$ with $\N(0,tI)$, dilate by $a_i=(1+\kappa_i t)^{-1/2}$.
This restores the curvature bound $\kappa_{i,t}/a_i^2=\kappa_i$ and
dilates the smoothed moment measure $\nu_{i,t}$ by $1/a_i$. Since
$a_i\to1$, both the sources and their moment measures converge to
the original pairs in $W_2$ by \cref{mm:lem:smoothing}.

For each regularized pair, take $T$ and $H=\D T$ as in
\cref{mm:prop:weighted}, with norms and inner products at its source.
By \eqref{mm:eq:control-spectral-defect}, the sum of the nonnegative
weighted Poincar\'e deficits of the components of $T-\Id$ tends to
zero. Their variances, on the other hand, sum to
$\norm{T-\Id}_2^2\to W_2^2(\rho_0,\rho_1)>0$. We can therefore select
a component of largest variance and divide it by its $L^2$ norm, which
stays bounded away from zero. The resulting centered function $g$ satisfies
\begin{equation}\label{mm:eq:normalized-displacement}
 \norm{g}_2=1,\qquad
 \mathcal E_H(g)\longrightarrow\sqrt{\kappa_0\kappa_1}.
\end{equation}

Let $X,Y$ be the fields \eqref{eq:couple} for the centered Poisson
solution with right-hand side $g$. The functions $g$ and $g\circ T^{-1}$
both have unit norm by \eqref{eq:common-div}. The coupled Bochner
bounds \eqref{eq:two-flux-bounds} give
\[
 \norm{X}_2\le\kappa_0^{-1/2},\qquad
 \norm{Y\circ T}_2\le\kappa_1^{-1/2}.
\]
Testing the source equation with $g$ gives $1=\inner{X}{\nabla g}$.
Since $Y\circ T=HX$, H\"older's inequality with weights $H$ and $H^{-1}$ yields
\begin{equation}\label{mm:eq:bochner-field-inner-product}
 \frac{1}{\mathcal E_H(g)}
 \le\inner{X}{Y\circ T}
 \le\norm{X}_2\norm{Y\circ T}_2
 \le\frac{1}{\sqrt{\kappa_0\kappa_1}} .
\end{equation}
The outer bounds have the same limit by
\eqref{mm:eq:normalized-displacement}, so the two endpoint estimates
for $X$ and $Y$ are both asymptotically equalities. By
\cref{lem:endpoint-saturation} the limiting sources have a common Gaussian
direction, with respective variances $\kappa_0^{-1}$ and $\kappa_1^{-1}$.

Let $E$ be the maximal common Gaussian subspace given by
\cite[Proposition~4.2]{HanZhu}, and write $\sigma_i$ for the
$E^\perp$-marginal of $\rho_i$. Taking the gradient of the product
potential gives
\[
 \rho_i=\N_E(0,\kappa_i^{-1}I_E)\otimes\sigma_i,\qquad
 \nu_i=\N_E(0,\kappa_i I_E)\otimes\Moment\sigma_i.
\]
Squared transport distances add over the two factors. The Gaussian
factors attain equality, so subtracting their costs gives
\[
 W_2^2(\Moment\sigma_0,\Moment\sigma_1)
   =\kappa_0\kappa_1 W_2^2(\sigma_0,\sigma_1).
\]
We can now apply the same argument on $E^\perp$ if the complementary
factors differ. To check the hypotheses, restrict the product potentials
to $E^\perp$, where they retain their curvature bounds. The product of $E$
with the discontinuity set of a complementary potential lies in the
$\mathcal H^{d-1}$-null discontinuity set of the full potential, so
slicing along translates of $E^\perp$ gives essential continuity when
the complement has positive dimension. The projected moment measures
remain centered and have finite second moments. If one were concentrated
on an affine hyperplane, the full measure would be concentrated on its
inverse image. Thus $\sigma_0\ne\sigma_1$ would give a further common
Gaussian direction in $E^\perp$, contrary to maximality. This proves
\eqref{mm:eq:equality-products}.

Conversely, under \eqref{mm:eq:equality-products} only the Gaussian
factors contribute to the two distances, so the Gaussian distance
formula gives equality. Since $\rho_0\ne\rho_1$, we have
$\kappa_0\ne\kappa_1$. The covariance difference
\[
 \Cov(\rho_0)-\Cov(\rho_1)=(\kappa_0^{-1}-\kappa_1^{-1})P_E
\]
therefore has range $E$, which proves uniqueness.
\end{proof} 

The inverse estimate also has an infinitesimal form. When a source
moves, both its density and the map given by the gradient of its
negative log density vary. Their combined variation gives a velocity
on the target, whose gradient projection determines the Wasserstein
speed of the moment measure.

\begin{proposition}\label[proposition]{prop:moment-direction}
Let $\rho=e^{-\psi}\dx$ be centered, with $\psi\in C^\infty(\R^d)$ and
$\kappa I\preceq H=\D^2\psi\preceq MI$, $0<\kappa\le M<\infty$.
Set $T=\nabla\psi$, $\nu=T_\#\rho$ and
$L_\rho=\Delta-\nabla\psi\cdot\nabla$, acting componentwise.
Let $\Pi_\psi$ be the orthogonal projection in $L^2(\rho;\R^d)$
onto the closure of $\{\nabla h\circ T:h\in C_c^\infty(\R^d)\}$.
For $f\in C_c^\infty(\R^d)$, put $v=\nabla f-\int\nabla f\dd\rho$,
$\rho_t=(\Id+tv)_\#\rho$ and $\nu_t=\Moment\rho_t$ for small $|t|$.
Then
\[
 \lim_{t\to0}\frac{W_2(\nu_t,\nu)}{|t|}
   =\|\Pi_\psi L_\rho v\|_\rho,\qquad
 \|v\|_\rho\le\sqrt{C_P(\rho)/\kappa}\,\|\Pi_\psi L_\rho v\|_\rho .
\]
Here $\|v\|_\rho=W_2(\rho_t,\rho)/|t|$ for small $t\ne0$.
\end{proposition}

\begin{proof}
For small $|t|$, $S_t=\Id+tv$ is a convex-gradient diffeomorphism,
which gives the source speed. Write $\rho_t=e^{-\psi_t}\dx$.
The identity $\psi_t\circ S_t=\psi+\log\det(\Id+t\D v)$
preserves uniform convexity for small $|t|$ by smooth convergence
on compact sets, since $S_t$ is a translation off a fixed compact set.
Differentiating in space and using $\nabla\diver v=\Delta v$ gives
\[
 \nabla\psi_t\circ S_t=T+tL_\rho v+O_{L^2(\rho)}(t^2).
\]
Thus the curve $(F_t)_\#\nu$, with $F_t=\Id+tb$ and
$b=(L_\rho v)\circ T^{-1}\in C_c^\infty(\R^d;\R^d)$, differs from
$\nu_t$ by $O(t^2)$ in $W_2$. It has the same speed at zero, which
we compute by projecting $b$ onto gradients in $L^2(\nu)$. Write
$p$ for this projection. Since $b-p$ is orthogonal to
$\nabla(h\circ F_s)$ for every test $h$, the curve $(F_s)_\#\nu$
has a velocity whose pullback is
$b-\D F_s(b-p)=p-s\D b(b-p)$. Integrating its norm gives the upper
bound $\|p\|_\nu$ for the speed at zero. For the reverse bound, Taylor
expansion of
$h\in C_c^\infty$ along an optimal coupling from $\nu$ to
$(F_t)_\#\nu$ bounds the absolute increment of its integral by
$\|\nabla h\|_\nu W_2+O(W_2^2)$, where $W_2=O(|t|)$.
The derivative of this integral at zero is
$\langle b,\nabla h\rangle_\nu$. Taking the supremum over unit
test gradients gives the lower speed $\|p\|_\nu=\|\Pi_\psi L_\rho v\|_\rho$.

To compare this target speed with the source speed, represent the
source density derivative $g=-\diver_\rho v$ by a second velocity
$X=H^{-1}\nabla u$. The centered equation
$-\diver_\rho(H^{-1}\nabla u)=g$ is solvable by ellipticity and the
Poincar\'e inequality. Its weak form makes $Z=X-v$ orthogonal to
both $v$ and $HX=\nabla u$. Moreover,
$X=\nabla(u\circ T^{-1})\circ T$, so $\Pi_\psi X=X$ by density.
The weighted Bochner identity and the orthogonality to $HX$ give
\[
 Q:=\|g\|_\rho^2-\int H[X,X]\dd\rho
   =\int\|\D v\|_{\HS}^2\dd\rho+\int H[Z,Z]\dd\rho .
\]
For Santambrogio's functional $\Ent+\mathcal T(\cdot,\nu)$, $Q$ is
the quadratic form on the right of \eqref{eq:hessian-transfer}. We
can now bound it in two ways. Testing the weak equation with $g$,
using $L_\rho v=Hv-\nabla g$ and $\Pi_\psi X=X$, gives
$Q=-\langle X,\Pi_\psi L_\rho v\rangle_\rho$. These integrations by
parts follow by cutoff from the Hessian bounds and the constancy of
$v$ off a compact set. In the other direction, componentwise
Poincar\'e gives $Q\ge C_P(\rho)^{-1}\|v\|_\rho^2+\kappa\|Z\|_\rho^2$.
Since $C_P(\rho)\le1/\kappa$ and $Z\perp v$, this also gives
$Q\ge\kappa\|X\|_\rho^2$. Cauchy--Schwarz in the pairing with the
target velocity now yields the claimed inverse estimate.
\end{proof}

\subsection{Regularized moment measures and sharp centered stability}\label{sec:regularized-moment}
In the classical inverse problem, translating a source preserves its
moment measure, so we selected the centered representative. A fixed
strongly convex confinement determines the position of the source.
With this confinement in place, equality of the target barycenters
yields a stronger stability estimate.

Fix a finite $\alpha$-strongly convex potential $V$ on $\R^d$. The
unknown transport potential $\Phi$ determines both the source
$\rho=Z^{-1}e^{-V-\Phi}\dx$ and the map $\nabla\Phi$ from $\rho$
to the prescribed target $\nu$. For smooth potentials the ordinary
moment map of this source is $\nabla V+\nabla\Phi$, whereas the
regularized problem prescribes the law of $\nabla\Phi$ alone. This
distinction is reflected in the operator $B-I$. In its Bochner estimate
the contribution $\D^2\Phi$ from the negative log density cancels,
leaving the curvature supplied by $V$.

The variational formulation allows $\Phi$ to take the value $+\infty$.
Let $m_V$ be the probability measure proportional to $e^{-V}\dx$. We minimize
\[
 \mathcal F_V(\rho;\nu)=\Ent_{m_V}(\rho)+\mathcal T(\rho,\nu)
\]
over $\rho\in\Ptwo(\R^d)$.
For quadratic confinement, Delalande and Farinelli
\cite[Theorem~3.1]{DelalandeFarinelli} obtained a $1/2$-H\"older estimate
on classes with uniformly bounded second moments.

\begin{theorem}
\label[theorem]{mm:thm:regularized}
Fix $\alpha>0$ and a finite $\alpha$-strongly convex potential $V$ on $\R^d$.
For every $\nu\in\Ptwo(\R^d)$, the functional $\mathcal F_V(\cdot;\nu)$
has a unique minimizer $\Reconstruct_V\nu$ in $\Ptwo(\R^d)$. There is a
proper lower semicontinuous convex potential $\Phi$ such that
\[
 \rho=\Reconstruct_V\nu=Z^{-1}e^{-V-\Phi}\dx,\qquad
 (\nabla\Phi)_\#\rho=\nu .
\]
Let $\nu_0,\nu_1\in\Ptwo(\R^d)$ and write $\rho_i=\Reconstruct_V\nu_i$. Then 
\begin{enumerate}[label=\textup{(\roman*)},leftmargin=2em]
\item For arbitrary targets,
\begin{equation}\label{mm:eq:regularized-stability}
 W_2(\rho_0,\rho_1)\le\frac1\alpha W_2(\nu_0,\nu_1).
\end{equation}
\item If the targets have the same barycenter, then
\begin{equation}\label{mm:eq:regularized-centered}
 W_2(\rho_0,\rho_1)\le
 \begin{cases}
  \dfrac{1}{2\alpha}\,W_2(\nu_0,\nu_1),&d=1,\\[5pt]
  \dfrac{1}{\sqrt{2}\,\alpha}\,W_2(\nu_0,\nu_1),&d\ge2.
 \end{cases}
\end{equation}
\end{enumerate}
For each $\alpha>0$ and $d\ge1$, the coefficients in both estimates
are optimal. 
\end{theorem}

\begin{proof}
We first identify the minimizer. Displacement convexity of entropy and
strong convexity of $V$ give $\alpha$-convexity of $\Ent_{m_V}$. To
handle the correlation term along a geodesic represented by
$X_t=(1-t)X_0+tX_1$, glue an optimal correlation coupling $(X_t,Y)$
to $(X_0,X_1)$ conditionally on $X_t$.
The two terms in $\E[X_t\cdot Y]=(1-t)\E[X_0\cdot Y]+t\E[X_1\cdot Y]$
are bounded by the corresponding endpoint maximal correlations. This
gives displacement convexity of $\mathcal T(\cdot,\nu)$, as in
\cite[Proposition~3.3]{SantambrogioMoment}.
The identity $2\mathcal T(\rho,\nu)=\int|x|^2\dd\rho+
\int|y|^2\dd\nu-W_2^2(\rho,\nu)$ gives $W_2$-continuity of
$\mathcal T(\cdot,\nu)$, and $\Ent_{m_V}$ is lower semicontinuous. For
$b=\int y\dd\nu(y)$, the product coupling and the entropy variational
formula give
\[
 \mathcal F_V(\rho;\nu)\ge
 \Ent_{m_V}(\rho)+\int b\cdot x\dd\rho(x)
 \ge-\log\int e^{-b\cdot x}\dd m_V(x)>-\infty .
\]
This lower bound is finite by the strong convexity of $V$, and a bounded
density of compact support gives a finite upper bound on the infimum.
For a minimizing sequence $(\rho_j)$, strong convexity at the midpoint
therefore gives
\[
 \frac{\alpha}{8} W_2^2(\rho_j,\rho_k)
 \le \frac{\mathcal F_V(\rho_j;\nu)+\mathcal F_V(\rho_k;\nu)}2
      -\inf\mathcal F_V(\cdot;\nu).
\]
Completeness and lower semicontinuity give a minimizer $\rho$, unique
by strong convexity. The argument also applies on a closed ball
$\overline\Omega$, with minimizer $\Reconstruct_V^\Omega\nu$.
In both cases \cref{mm:lem:regularized-model}\,\textup{(i)} gives the Gibbs
representation. Strong convexity along the geodesic from $\rho$
to a competitor $\sigma$ yields
\begin{equation}\label{mm:eq:regularized-growth}
 \mathcal F_V(\sigma;\nu)-\mathcal F_V(\rho;\nu)
 \ge \frac{\alpha}{2}W_2^2(\rho,\sigma).
\end{equation}

We will estimate the response to the target by differentiating a smooth
solution on a ball. To identify such a solution with the variational
minimizer, suppose that $\Phi$ is smooth and convex and that
$r=Z^{-1}e^{-V-\Phi}$ on $\overline\Omega$ solves $T_\#\rho=\nu$,
where $T=\nabla\Phi$ and $\rho=r\one_\Omega\dx$.
Given a finite-energy competitor $\sigma$ on $\overline\Omega$, let $S$
transport $\rho$ optimally to $\sigma$. Write $S_h=(1-h)\Id+hS$ and
$\rho_h=(S_h)_\#\rho$. Nonnegativity of $\Ent_\rho(\rho_h)$ gives
\[
 \Ent_{m_V}(\rho_h)-\Ent_{m_V}(\rho)
 \ge -\int[\Phi(S_h(x))-\Phi(x)]\dd\rho(x).
\]
The coupling $(S_h,T)_\#\rho$ is admissible for
$\mathcal T(\rho_h,\nu)$, and $(\Id,T)_\#\rho$ realizes
$\mathcal T(\rho,\nu)$. Adding these bounds gives
\begin{equation}\label{mm:eq:regularized-selection}
 \mathcal F_V(\rho_h;\nu)-\mathcal F_V(\rho;\nu)
 \ge -\int\Dcal_\Phi(S_h(x),x)\dd\rho(x)
 \ge -\tfrac12\|\D^2\Phi\|_\infty h^2 W_2^2(\rho,\sigma).
\end{equation}
On the other hand, $\alpha$-convexity gives
\[
 \frac{\mathcal F_V(\rho_h;\nu)-\mathcal F_V(\rho;\nu)}h
 \le \mathcal F_V(\sigma;\nu)-\mathcal F_V(\rho;\nu)
     -\frac{\alpha}{2}(1-h)W_2^2(\rho,\sigma).
\]
By \eqref{mm:eq:regularized-selection}, the lower limit on the left
is nonnegative as $h\downarrow0$. This proves
\eqref{mm:eq:regularized-growth} for the smooth solution and identifies
it with $\Reconstruct_V^\Omega\nu$.

We next obtain the estimate that permits differentiation of the moment
relation. Assume $V$ is smooth and $T$ is a diffeomorphism, and set
$H=\D^2\Phi$ and $B=-\diver_\rho(H^{-1}\nabla)$. For smooth $z$ with
$X=H^{-1}\nabla z$ tangent to $\partial\Omega$, apply
\eqref{eq:weighted-bochner} with weight one. The Hessian of the
negative log density is $\D^2V+H$, and integration by parts gives
$\int zBz\dd\rho=\int H[X,X]\dd\rho$. Subtraction cancels $H$
exactly. By \eqref{eq:trace-positive},
\begin{equation}\label{mm:eq:regularized-bochner}
\begin{aligned}
 \int (Bz-z)Bz\dd\rho
 =\int\!\bigl(\|H^{1/2}\D XH^{-1/2}\|_{\HS}^2+\D^2V[X,X]\bigr)\dd\rho
   +\int_{\partial\Omega}\II(X,X)r\dd S \ge \alpha\norm{X}_\rho^2 .
\end{aligned}
\end{equation}
The boundary term is nonnegative by convexity of the ball. In
particular, $Bz=z$ forces $X=0$ and then $z=Bz=0$. Since the
self-adjoint conormal realization of $B$ on the smooth ball has
compact resolvent, this proves that $B-I$ is invertible.
Let $\nu_t$ be a smooth positive curve on a target ball with velocity
$v_t$ tangent to the boundary. For smooth $V$,
\cref{mm:lem:regularized-model}\,\textup{(ii)} gives a smooth curve of minimizers
$\rho_t=r_t\dx$ on the source ball. Normalize $\Phi_t$ by
$r_t=e^{-V-\Phi_t}$. Fix a time, suppress $t$, and put
$g=\dot r/r=-\dot\Phi$ and $u=v\circ T$.

The source speed is the norm of the gradient velocity $p$ of $\rho$,
whereas the variation of the transport map involves $X=H^{-1}\nabla z$.
We choose $z$ as the centered conormal solution of $Bz=g$, so that
both fields represent the same density derivative. The two Neumann
problems are solvable since $\int g\dd\rho=0$, and give
\[
 -\diver_\rho p=-\diver_\rho X=g,\qquad
 p\cdot n=X\cdot n=0\quad\text{on }\partial\Omega .
\]
The difference $X-p$ pairs to zero with all gradients. Testing with
the coordinates, a velocity potential for $p$, and $z$, respectively, gives
\[
 m:=\int p\dd\rho=\int X\dd\rho,\qquad
 \langle p,X-p\rangle_\rho=0,\qquad
 \int H[X,X-p]\dd\rho=0 .
\]
The first variation of $\Ent_{m_V}$ at $\rho$ equals $-\Phi$ up to a
constant, so for a tangent field $Y$,
\[
 \mathscr Q_\rho(Y)=\norm{\diver_\rho Y}_\rho^2-\int H[Y,Y]\dd\rho .
\]
Apply \eqref{eq:critical-branch-response} with $a=1$ and $c=0$.
Since $X$ and $p$ have the same weighted divergence and
$\int H[X,X-p]\dd\rho=0$, we obtain
\begin{equation}\label{mm:eq:regularized-complete}
 -\inner uX_\rho=\mathscr Q_\rho(X)
 =\mathscr Q_\rho(p)+\int H[X-p,X-p]\dd\rho .
\end{equation}
Equation \eqref{mm:eq:regularized-bochner} already bounds
$\mathscr Q_\rho(X)$ from below by $\alpha\norm X_\rho^2$. The second
expression in \eqref{mm:eq:regularized-complete} gives additional
control of the gradient velocity. Indeed, applying
\eqref{eq:weighted-bochner} to $p$ with weight one and subtracting
$\int H[p,p]\dd\rho$ yields
\[
 \mathscr Q_\rho(p)
 =\int\!\bigl(\norm{\D p}_{\HS}^2+\D^2V[p,p]\bigr)\dd\rho
    +\int_{\partial\Omega}\II(p,p)r\dd S
 \ge\alpha\bigl(\norm p_\rho^2+\norm{p-m}_\rho^2\bigr).
\]
The last step uses the Poincar\'e inequality for the
$\alpha$-strongly log-concave density on the convex ball, applied to
each component of $p$.

These two estimates control the source speed in terms of the target
speed. Using the bound on $p$ on one side of
\eqref{mm:eq:regularized-complete} and the bound on $X$ on the other gives
\[
 \alpha\bigl(\norm p_\rho^2+\norm{p-m}_\rho^2\bigr)
 \le \mathscr Q_\rho(X)
 =-\inner uX_\rho
 \le \norm u_\rho\sqrt{\mathscr Q_\rho(X)/\alpha}.
\]
In particular, the left side is at most $\norm u_\rho^2/\alpha$.
Write $\rho_t^{\mathrm c}$ for the translate of $\rho_t$ to barycenter
zero. Since the barycenter has derivative $m_t$, the translated field
$p_t-m_t$ is a velocity for $\rho_t^{\mathrm c}$. The dynamic formula
for the two curves gives
\begin{equation}\label{mm:eq:regularized-action}
 \alpha^2\bigl[W_2^2(\rho_0,\rho_1)
             +W_2^2(\rho_0^{\mathrm c},\rho_1^{\mathrm c})\bigr]
 \le\int_0^1\norm{v_t}_{\nu_t}^2\dt .
\end{equation}

If the target barycenter is constant, then $\int u\dd\rho=0$, so the
pairing $-\inner uX_\rho$ depends only on $X-m$. To retain this gain,
average the two lower bounds for $\mathscr Q_\rho(X)$ and use
$\int X\dd\rho=m$. This gives
\[
 \alpha\norm p_\rho^2+\frac{\alpha}{2}\norm{X-m}_\rho^2
 \le \mathscr Q_\rho(X)=-\inner u{X-m}_\rho.
\]
Completing the square yields
$\norm p_\rho^2\le\norm u_\rho^2/(2\alpha^2)$.
On an interval, $(r(X-p))'=0$ and the boundary flux is zero, so
$X=p$. The lower bound for $\mathscr Q_\rho(p)$ alone now gives
\[
 \alpha\norm p_\rho^2
 \le \norm u_\rho\norm{p-m}_\rho-\alpha\norm{p-m}_\rho^2
 \le \frac{\norm u_\rho^2}{4\alpha}.
\]
Integrate the speed estimates along the approximations of
\cref{mm:lem:regularized-model}\,\textup{(iii)}, keeping the target barycenter constant
when required. The reconstructed endpoints converge in $W_2$ and the
target actions tend to $W_2^2(\nu_0,\nu_1)$. Thus
\eqref{mm:eq:regularized-action} holds on $\R^d$ with this squared
distance on the right, giving \eqref{mm:eq:regularized-stability}.
The improved speed estimates give \eqref{mm:eq:regularized-centered}.

It remains to see why the constants cannot be reduced. Take
$V(x)=\alpha|x|^2/2$. For the target $\delta_y$ the correlation is
$\int x\cdot y\dd\rho$, so completing the square gives
$\Reconstruct_V\delta_y=\N(-y/\alpha,\alpha^{-1}I)$. Two distinct
point masses thus attain equality in \eqref{mm:eq:regularized-stability}.

Keep the same $V$ and consider centered Gaussian targets. On the line, a
target of standard deviation $t\ge0$ has a centered Gaussian
reconstruction with standard deviation
\[
 s(t)=\frac{2}{t+\sqrt{t^2+4\alpha}}.
\]
Indeed, the convex potential $\Phi(x)=(s(t)^{-2}-\alpha)x^2/2$
gives the required density and sends it by its gradient to the target,
since $s(t)^{-1}-\alpha s(t)=t$. For the quadratic potentials used
here, the comparison in \eqref{mm:eq:regularized-selection} applies
on the whole space and identifies this solution as the minimizer.
Thus the distance ratio for the
targets $\delta_0$ and $\N(0,t^2)$ tends to $-s'(0)=1/(2\alpha)$
as $t\downarrow0$.

In dimension two we can vary the direction of a centered target while
keeping its variance fixed. Take the law of $(tZ,0)$, where
$Z\sim\N(0,1)$ and $t>0$, and rotate it through an angle $\theta$.
The first reconstruction is Gaussian with covariance
$\operatorname{diag}(s(t)^2,s(0)^2)$ by the product form of the preceding
potential, and the second is its rotation by invariance of $V$.
For a centered planar Gaussian with covariance
$\operatorname{diag}(a,b)$, $a,b\ge0$, the Gaussian distance
formula gives the squared distance to its rotation as
\[
 2(a+b)-2\sqrt{(a+b)^2-(a-b)^2\sin^2\theta}.
\]
Consequently its distance divided by $|\theta|$ tends to
$|a-b|/\sqrt{a+b}$ when $a+b>0$. Applying this to the two
reconstructions and to the two targets, whose covariance eigenvalues
are $t^2,0$, gives the limit as $\theta\to0$ of the distance ratio,
\[
 \frac{s(0)^2-s(t)^2}{t\sqrt{s(0)^2+s(t)^2}}
 =\frac{s(t)}{\alpha\sqrt{s(0)^2+s(t)^2}}
 \longrightarrow\frac{1}{\sqrt{2}\,\alpha}
 \qquad(t\downarrow0).
\]
In higher dimensions, tensor both targets with a point mass at the
origin on the orthogonal complement. The reconstructions acquire the
same $\N(0,\alpha^{-1}I)$ factor there, so both distances are preserved.
\end{proof} 

The unrestricted estimate admits distinct extremizing targets, as the
point masses show. The examples with a common barycenter approach the
improved constant only in a limit. We now determine when equality can
hold in either estimate.
\begin{corollary}
\label[corollary]{mm:cor:regularized-equality}
Let $\alpha$ and $V$ be as in \cref{mm:thm:regularized}.
For $\nu_0,\nu_1\in\Ptwo(\R^d)$, set $\rho_i=\Reconstruct_V\nu_i$.
\begin{enumerate}[label=\textup{(\roman*)},leftmargin=2em]
\item Equality in \eqref{mm:eq:regularized-stability} holds if and
only if there exists $h\in\R^d$ such that
\begin{equation}\label{mm:eq:regularized-translation}
 \nu_1=(\Id+h)_\#\nu_0,\qquad
 \rho_1=(\Id-h/\alpha)_\#\rho_0.
\end{equation}
For a prescribed $\nu_1=(\Id+h)_\#\nu_0$, equality holds exactly when
$V(x-h/\alpha)-V(x)+h\cdot x$ is constant on $\supp\rho_0$.
In particular, $V-\alpha|\cdot|^2/2$ is affine if and only if every pair
$(\nu,(\Id+h)_\#\nu)$ attains equality.
\item If $\nu_0$ and $\nu_1$ have the same barycenter, equality in
\eqref{mm:eq:regularized-centered} holds if and only if $\nu_0=\nu_1$.
\end{enumerate}
\end{corollary}

\begin{proof}
\emph{(i)} We start by comparing a minimizer with a translate of the
other minimizer. Let $m_i$ be the barycenter of $\nu_i$. Summing
\eqref{mm:eq:regularized-growth} at the two minimizers cancels the
entropy terms. If $\rho_1=(\Id+b)_\#\rho_0$, translation by $b$
changes the correlation with $\nu_i$ by $b\cdot m_i$, giving
\begin{equation}\label{mm:eq:regularized-translate-means}
 \alpha|b|^2\le b\cdot(m_0-m_1)\le |b|W_2(\nu_0,\nu_1).
\end{equation}
At equality in \eqref{mm:eq:regularized-stability}, the whole-space form of
\eqref{mm:eq:regularized-action} gives $\rho_0^{\mathrm c}=\rho_1^{\mathrm c}$.
Write $\rho_1=(\Id+b)_\#\rho_0$. Then $W_2(\nu_0,\nu_1)=\alpha|b|$,
and equality throughout \eqref{mm:eq:regularized-translate-means}
gives $b\cdot(m_0-m_1)=\alpha|b|^2$. For an optimal target coupling
$(Y_0,Y_1)$ it follows that
\[
 \E|Y_1-Y_0+\alpha b|^2
 =W_2^2(\nu_0,\nu_1)+2\alpha b\cdot(m_1-m_0)+\alpha^2|b|^2=0.
\]
This proves \eqref{mm:eq:regularized-translation} with $h=-\alpha b$.
Conversely, the two translations have $W_2$ lengths $|h|$ and $|h|/\alpha$.

To decide whether a prescribed translation of the target produces
the required translation of the source, we use the Euler equation.
For $\rho=g\,m_V$ of finite relative entropy and its Brenier map $T$
to the target, minimality is equivalent to
\[
 \nabla g=-Tg \qquad\text{in distributions on }\R^d.
\]
Indeed, the first variation along a smooth compactly supported field $w$
is $\int[(\nabla V+T)\cdot w-\diver w]\dd\rho$ by uniqueness of the
optimal plan. Conversely, on compact sets the density of $m_V$ is bounded above and
away from zero, so $Tg$ is locally integrable by $T\in L^2(\rho)$.
The weak equation gives $g\in W^{1,1}_{\rm loc}(\R^d)$, including
across the boundary of its positivity set. For entropy the pressure is
$g$ itself, so \cite[Theorem~10.4.9]{AGS} yields
$-T\in\partial\Ent_{m_V}(\rho)$. For any optimal transport map $S$ to a competitor, the coupling
$(S,T)_\#\rho$ bounds the correlation increment below by
$\int T\cdot(S-\Id)\dd\rho$. The two linear terms cancel.

For a prescribed target translation $h$, put $b=-h/\alpha$ and write
$\rho_0=g_0\,m_V$. Its translate by $b$ has relative density
$g_b(x+b)=e^{V(x+b)-V(x)}g_0(x)$ and map $T_b(x+b)=T_0(x)-\alpha b$
to the translated target. The product rule gives $\nabla g_b=-T_b g_b$
exactly when $\nabla V(x+b)-\nabla V(x)=\alpha b$ $\rho_0$-a.e.
Since $V$ is locally Lipschitz and the density is positive in the interior
of its convex support, this is equivalent to the existence of $k\in\R$ with
\begin{equation}\label{mm:eq:regularized-increment}
 V(x-h/\alpha)-V(x)=-h\cdot x+k \qquad(x\in\supp\rho_0).
\end{equation}
This increment also ensures that the translated source has finite relative entropy.

If $G=V-\alpha|\cdot|^2/2$ is affine, \eqref{mm:eq:regularized-increment}
holds for all $x,h$. Conversely, take $\nu_0=\delta_0$, so that
$\rho_0=m_V$ has full support. Equality for every translate then gives
$G(x+z)-G(x)=G(z)-G(0)$ for all $x,z\in\R^d$ by
\eqref{mm:eq:regularized-increment}. Thus $G-G(0)$ is continuous
and additive, hence linear.

\emph{(ii)} Suppose now that the target barycenters agree. We shall
use the remainder in the speed estimate to show that equality would
force the centered reconstructions to coincide. Write $L$ for the
coefficient in \eqref{mm:eq:regularized-centered}. Along a smooth
path with constant target barycenter as in \cref{mm:thm:regularized},
set $A=\int\D p\dd\rho$ and
$\delta=\int\norm{\D p}_{\HS}^2\dd\rho-\alpha\Var_\rho(p)$, with
variances summed over components. The Poincar\'e deficit $\delta$
controls the squared distance of $p$ from an affine field through the bounds
\[
 \frac{\alpha}{2}\Var_\rho(p-Ax)\le\delta
 \le 2\alpha\bigl(L^2\norm u_\rho^2-\norm p_\rho^2\bigr).
\]
For the upper bound, keeping $\delta$ in the Bochner estimate for
$\mathscr Q_\rho(p)$ retains $\delta/2$ after averaging with
\eqref{mm:eq:regularized-bochner}, and all of $\delta$ on the line
where $X=p$. For the lower bound let $f$ be a centered component of $p$.
On the smooth source ball solve $-\diver_\rho\nabla w=f$ with
$\partial_n w=0$ and $\int w\dd\rho=0$. Integration by parts gives
$\langle\nabla f,\nabla w\rangle_\rho=\|f\|_\rho^2$.
Apply \eqref{eq:weighted-bochner} with unit weight to $\nabla w$.
The boundary term $\int_{\partial\Omega}\II(\nabla w,\nabla w)r\dd S$
is nonnegative, and $\D^2(-\log r)\succeq\alpha I$. Applying
Poincar\'e to the components of $\nabla w$ gives
\[
 \alpha\Var_\rho(\nabla w)
 \le\int\|\D^2 w\|_{\HS}^2\dd\rho
 \le\|f\|_\rho^2-\alpha\|\nabla w\|_\rho^2 .
\]
Split $\nabla f=(\nabla f-\alpha\nabla w)+\alpha\nabla w$. Then
\[
 \begin{aligned}
 \tfrac12\Var_\rho(\nabla f)
 &\le\|\nabla f-\alpha\nabla w\|_\rho^2+
                      \alpha^2\Var_\rho(\nabla w)\\
 &\le\|\nabla f-\alpha\nabla w\|_\rho^2
      +\alpha\|f\|_\rho^2-\alpha^2\|\nabla w\|_\rho^2\\
 &=\|\nabla f\|_\rho^2-2\alpha\langle\nabla f,\nabla w\rangle_\rho
      +\alpha\|f\|_\rho^2
 =\|\nabla f\|_\rho^2-\alpha\|f\|_\rho^2 .
 \end{aligned}
\]
Summing over components and applying Poincar\'e to $p-Ax$ yields
$\alpha\Var_\rho(p-Ax)\le\int\|\D p-A\|_{\HS}^2\dd\rho\le2\delta$.

Suppose equality holds for the two given targets, and choose the smooth
curves of \cref{mm:lem:regularized-model}\,\textup{(iii)} with constant
target barycenter. To compare the bound on $\delta$ with the endpoint
distances, follow a random point of the source along
$\dot Y_s=p_s(Y_s)$, starting with the initial source law. Smoothness
and tangency of $p_s$ on each closed ball give a flow on the whole
interval with $\mathop{\rm Law}(Y_s)=\rho_s$. Write $Z_s=Y_s-\E Y_s$
and $D=Z_1-Z_0$. The speed estimate gives
\[
 W_2^2(\mathop{\rm Law}(Y_0),\mathop{\rm Law}(Y_1))
 \le \E|Y_1-Y_0|^2
 \le \int_0^1 \E|\dot Y_s|^2\dd s
 \le L^2\int_0^1\norm{v_s}_{\nu_s}^2\dd s .
\]
The two ends tend to $W_2^2(\rho_0,\rho_1)$ by the assumed equality.
The difference of the middle terms is
$\int_0^1\E|\dot Y_s-(Y_1-Y_0)|^2\dd s$. Centering decreases this
quantity, and the deficit estimate above therefore gives
\[
 \int_0^1\bigl(\delta_s+\E|\dot Z_s-D|^2\bigr)\dd s\longrightarrow0.
\]
We choose an interior time at which the velocity is close both to an
affine function of position and to the endpoint displacement $D$.
For each approximating curve, take $t\in[1/3,2/3]$ where the integrand
is at most its average over that interval. Then
$\delta_t+\E|\dot Z_t-D|^2\to0$. The interpolation
$Z=(1-t)Z_0+tZ_1$ also satisfies
$\norm{Z_t-Z}_2\le\int_0^t\norm{\dot Z_s-D}_2\dd s\to0$.
At this time $A=\int\D p_t\dd\rho_t$ satisfies
\[
 \norm A_{\HS}^2\le\int\norm{\D p_t}_{\HS}^2\dd\rho_t
 =\alpha\E|\dot Z_t|^2+\delta_t
 =\alpha\E|D|^2+o(1).
\]
In particular $A$ stays bounded. The affine estimate gives
$\norm{\dot Z_t-AZ_t}_2^2\le2\delta_t/\alpha$, hence
$\norm{D-AZ}_2\to0$.

Poincar\'e for linear functions at each approximating source endpoint
gives $\operatorname{Cov}(Z_i)\preceq\alpha^{-1}I$. Hence
\[
 \E|AZ|^2+t(1-t)\E|AD|^2
 =(1-t)\E|AZ_0|^2+t\E|AZ_1|^2\le\alpha^{-1}\norm A_{\HS}^2 .
\]
Since $\E|AZ|^2=\E|D|^2+o(1)$, this forces $\norm{AD}_2\to0$.
The covariance bounds give $\E|Z|^2\le d/\alpha$ uniformly.
As $p$ is a gradient, $A$ is symmetric and
$\E|D|^2=\E[Z\cdot AD]+o(1)\to0$.
Centering is continuous in $W_2$, since it subtracts the squared
barycenter displacement from every quadratic coupling cost. The joint
law of $(Z_0,Z_1)$ thus gives
$W_2^2(\rho_0^{\mathrm c},\rho_1^{\mathrm c})\le\liminf\E|D|^2=0$.
The sources are translates, and $m_0=m_1$ forces that translation
to vanish by \eqref{mm:eq:regularized-translate-means}. Equality in
the stability estimate now gives $\nu_0=\nu_1$, which plainly
suffices as well.
\end{proof} 

For targets with the same barycenter, \eqref{mm:eq:regularized-centered}
improves the reconstruction bound in every dimension. For Brenier
potentials from a fixed $\kappa$-strongly log-concave source $\rho$,
the improvement holds only on the line. Indeed, if $\nabla\Phi_{\nu_i}$
transports $\rho$ to $\nu_i$ and $f=\Phi_{\nu_1}-\Phi_{\nu_0}$, the
best uniform coefficient in $\Var_\rho(f)\le C W_2^2(\nu_0,\nu_1)$
is $1/(2\kappa)$ on the line and $1/\kappa$ in higher dimensions.

On the line, equality of the target barycenters gives
$\int f'\dd\rho=0$, while monotonicity of the coupling gives
$\int|f'|^2\dd\rho=W_2^2(\nu_0,\nu_1)$. The improvement therefore
follows from the Poincar\'e inequality for functions with centered
derivative \cite[Lemma~2]{CorderoFradeliziMaurey}. In a
curvature-preserving smooth approximation $\rho_j$,
replacing a smooth test $h$ by $h-x\int h'\dd\rho_j$ preserves the
derivative constraint, and Sobolev density gives the general case.
The pair $\nu_0=\delta_0$, $\nu_1=\rho=\N(0,\kappa^{-1})$ attains the bound.

For $d\ge2$ the upper bound follows from \cite[Theorem~1.1]{LetrouitManole}.
For $\rho=\N(0,\kappa^{-1}I_d)$, take the symmetric two-point target
on $\{\pm u_\theta\}$ with $u_\theta=\cos\theta\,e_1+\sin\theta\,e_2$,
as in \cite[Lemma~5.1 in the supplement]{MerigotDelalandeChazal}.
Its distance to the target at $\theta=0$
is $|u_\theta-e_1|\sim|\theta|$. Its potential is $|x\cdot u_\theta|$
up to a constant independent of $\theta$, with $L^2(\rho)$ derivative
$\operatorname{sgn}(x_1)x_2$ at zero of norm $\kappa^{-1/2}$.
This gives the matching lower bound $1/\kappa$.

\section{Stability of Brenier potentials}\label{sec:stability}
We continue to vary the target in $(\nabla\Phi)_\#\rho=\mu$, but now
the source $\rho$ is fixed and the quantity to estimate is the
normalized Brenier potential $\Phi$ in $L^2(\rho)$. This changes the
linearized equation. In a regularized reconstruction,
$r_t=e^{-V-\Phi_t}$ gives $\partial_t\log r_t=-\dot\Phi_t$.
For a fixed source $\partial_t\log r_t=0$, leaving
$B=-\diver_\rho((\D^2\Phi)^{-1}\nabla)$ in place of $B-I$ in the
equation of \cref{sec:response-equations}. Under
\eqref{eq:source-bounds}, the optimal stability constant $S(\rho)$ is the supremum over $\Phi$
of the norms of the operators $g\mapsto(\D^2\Phi)^{-1}\nabla B^{-1}g$
from \eqref{eq:fixed-source-pairing}.

This characterization distinguishes the effects of the source density
and of its support. Under a curvature assumption on the density we obtain
a weighted estimate, sharp for inverse-power densities. For the uniform
measure on a bounded connected domain with boundary of measure zero,
stability is equivalent to convexity of the closed support, as shown in \cref{sec:nonconvex}. Even on a square, however,
positive continuous densities arbitrarily close to uniform in
$L^\infty$ can have $S(\rho)=\infty$.

\subsection{A flux characterization of stability}
\label{sec:stability-tangent}\label{sec:finite-targets}\label{app:rough}
Assume \eqref{eq:source-bounds}. For $\Phi\in\cC(\Omega)$, put
$H=\D^2{\Phi}$ and let $B$ be the centered Neumann realization
of $-\diver_\rho(H^{-1}\nabla)$. Thus $u\in D(B)$ and
$B u=g$ mean that $u\in H^1(\Omega)$ is centered,
$g\in L^2_0(\rho)$, and
\begin{equation}\label{eq:Neumann-weak}
 \int H^{-1}\nabla u\cdot\nabla h\dd\rho
 =\int gh\dd\rho\qquad(h\in H^1(\Omega)).
\end{equation}
The Poincar\'e inequality \cite{PayneWeinberger} and ellipticity for
each fixed potential give a unique solution for every $g$.

For smooth positive densities on smooth uniformly convex domains,
Letrouit and Manole derived stability from a uniform bound on the
Neumann fluxes $H^{-1}\nabla B^{-1}g$
\cite[Theorem~3.2 and Lemma~3.6]{LetrouitManole}. To recover the flux
bound from stability, we perturb the potential by $tf$ and estimate
the speed of its target curve. Only the projection of $\nabla f$ onto
$H^{-1}\nabla H^1(\Omega)$ enters this speed estimate. This identifies
the optimal constant even when the source density is merely measurable
and satisfies \eqref{eq:source-bounds}.

\begin{proposition}
\label[proposition]{prop:stability-equivalences}
Let $\Omega\subseteq\R^d$ be a bounded open convex set and let
$\rho=r\one_\Omega\dx$ be a probability measure with
$0<m\le r\le M<\infty$ almost everywhere.
For $\Phi\in\cC(\Omega)$, put $H=\D^2{\Phi}$ and let $B$ be the
centered Neumann operator defined by \eqref{eq:Neumann-weak}.
Let $\Pi_\Phi$ be the orthogonal projection in $L^2(\rho;\R^d)$ onto
the closed subspace $H^{-1}\nabla H^1(\Omega)$.
For every finite $C\ge0$, the following statements are equivalent:
\begin{enumerate}[label=\textup{(\roman*)},leftmargin=2em]
\item $S(\rho)\le C$.
\item For every $\Phi\in\cC(\Omega)$ and $u\in D(B)$,
\begin{equation}\label{eq:uniform-flux}
 \norm{H^{-1}\nabla u}_{L^2(\rho)}
       \le C\norm{B u}_{L^2(\rho)}.
\end{equation}
\item For every $\Phi\in\cC(\Omega)$ and $f\in H^1(\Omega)$,
\begin{equation}\label{eq:projected-Poincare}
 \Var_\rho f\le C^2\norm{\Pi_\Phi\nabla f}_{L^2(\rho)}^2.
\end{equation}
\end{enumerate}
The best constant in each statement is $S(\rho)$, and
\begin{equation}\label{eq:exact-S}
 S(\rho)=\sup_{\Phi\in\cC(\Omega)}
 \norm{H^{-1}\nabla B^{-1}}_{L^2_0(\rho)\to L^2(\rho;\R^d)}
 \quad\text{in }[0,\infty].
\end{equation}
If $S(\rho)<\infty$, the same stability estimate holds for all
targets in $\Ptwo(\R^d)$. Restricting the compact targets to any fixed
compact set with nonempty interior leaves $S(\rho)$ unchanged.
\end{proposition}

\begin{proof}[Proof of \cref{prop:stability-equivalences,thm:stability}]
We first relate the two linear estimates for a fixed
$\Phi\in\cC(\Omega)$. Ellipticity and the Poincar\'e inequality
give closed range for $H^{-1}\nabla$ on $H^1(\Omega)$. To compute the
adjoint of $H^{-1}\nabla B^{-1}$, take $b\in L^2(\rho;\R^d)$ and
let $f_b\in H^1(\Omega)$ be the centered solution of
\begin{equation}\label{eq:source-linearized}
 \int H^{-1}\nabla f_b\cdot\nabla h\dd\rho
 =\int b\cdot H^{-1}\nabla h\dd\rho
 \qquad(h\in H^1(\Omega)).
\end{equation}
Testing with $h=B^{-1}g$, $g\in L^2_0(\rho)$, gives the weak form
of \eqref{eq:fixed-source-pairing},
\begin{equation}\label{eq:stability-adjoint}
 \inner{f_b}{g}_{L^2(\rho)}
 =\int b\cdot H^{-1}\nabla B^{-1}g\dd\rho.
\end{equation}
Thus $b\mapsto f_b$ and $H^{-1}\nabla B^{-1}$ have the same norm.
By \eqref{eq:source-linearized}, $\Pi_\Phi\nabla f_b=\Pi_\Phi b$ and
$f_{\Pi_\Phi\nabla f}=f-\int f\dd\rho$. Hence
\begin{equation}\label{eq:flux-adjoint}
 \norm{H^{-1}\nabla B^{-1}}
 =\sup_{\substack{f\in H^1(\Omega)\\\Var_\rho f>0}}
   \frac{\sqrt{\Var_\rho f}}{\norm{\Pi_\Phi\nabla f}_2},
\end{equation}
which gives \textup{(ii)}$\Leftrightarrow$\textup{(iii)}.

Now assume the global stability estimate \textup{(i)}. For
$f\in C_c^\infty(\R^d)$, the potentials $\Phi+tf$ remain in
$\cC(\Omega)$ for small $|t|$. They give an explicit target curve
$\mu_t=(T_t)_\#\rho$ with $T_t=\nabla(\Phi+tf)$ injective. To bound
its initial speed by the norm of $p=\Pi_\Phi\nabla f$, put
$Z=H^{-1}(\nabla f-p)$. By orthogonality,
$\int Z\cdot\nabla h\dd\rho=0$ for every $h\in H^1(\Omega)$.
Taking $h=\eta\circ T_t$, $\eta\in C_c^\infty(\R^d)$, gives
$\int\nabla\eta(T_t)\cdot \D T_t Z\dd\rho=0$.
Differentiating $\int\eta(T_t)\dd\rho$ now gives the continuity equation
for $\mu_t$ with velocity field $v_t\circ T_t=\nabla f-\D T_t Z=p-t\D^2 f\,Z$.
Integrating the speed bound of \cite[Theorem~8.3.1]{AGS} gives
\begin{equation}\label{eq:tangent-norm}
 W_2(\mu_s,\mu_0)
 \le |s|\norm p_2+\tfrac12s^2\norm{\D^2 f\,Z}_2.
\end{equation}
Since $\Phi_{\mu_s}-\Phi_{\mu_0}=s(f-\int f\dd\rho)$, \textup{(i)}
and \eqref{eq:tangent-norm} give \eqref{eq:projected-Poincare} after
division by $|s|$ and passage to the limit. Density extends the
estimate to $H^1(\Omega)$.

For the converse, we first integrate the flux bound between finitely
supported targets. This is the content of \cref{lem:finite-target-stability}
below. Given $\mu\in\Ptwo(\R^d)$, choose finitely supported
$\mu_n\to\mu$ in $W_2$. The lemma makes their centered potentials
$\Phi_{\mu_n}$ Cauchy in $L^2(\rho)$. Since $r\ge m$, convexity and
the $L^2(\rho)$ bound give a subsequence
converging uniformly on compact subsets of $\Omega$. Its gradients
converge almost everywhere. Passing to the limit in
$(\nabla\Phi_{\mu_n})_\#\rho=\mu_n$ and $\int\Phi_{\mu_n}\dd\rho=0$
identifies the $L^2$ limit as $\Phi_\mu$. Applying the finite-target
estimate to approximations of both endpoints gives \textup{(i)} and
its extension to $\Ptwo(\R^d)$. The constants agree by
\eqref{eq:flux-adjoint}, giving \eqref{eq:exact-S} in $[0,\infty]$.

Finally, for $F(y)=ay+b$, $a>0$, we have
$\Phi_{F_\#\mu}-\Phi_{F_\#\nu}=a(\Phi_\mu-\Phi_\nu)$, and $W_2$
scales by the same factor. Such a similarity places any two compact
supports inside any prescribed ball, proving the last assertion.
\end{proof}

For a smooth transport map $T=\nabla\Phi$, the chain rule pulls a
target gradient back to a field of the form $H^{-1}\nabla h$.
This explains the projection in the Poincar\'e inequality above.
For a smooth curve of centered potentials $\Phi_t$, the same pairing
\eqref{eq:fixed-source-pairing} and the flux bound
\eqref{eq:uniform-flux} give
$\norm{\dot\Phi_t}_{L^2(\rho)}\le C\norm{v_t}_{L^2(\mu_t)}$,
where $v_t$ is the Wasserstein velocity of the target. Integration
along a geodesic would give the desired distance bound directly.

For finitely supported targets the geodesic is explicit, but its
potentials are only piecewise affine. The next proof recovers the
speed bound by smoothing these potentials and keeping track of the
mass crossing their cell faces, with the source $\rho$ fixed throughout.

\begin{lemma}\label[lemma]{lem:finite-target-stability}
Under \eqref{eq:source-bounds}, suppose that
\eqref{eq:uniform-flux} holds with the same constant $C$ for every
$\Phi\in\cC(\Omega)$. Then, for all finitely supported probability
measures $\mu_0,\mu_1$,
\begin{equation}\label{eq:finite-target-stability}
 \norm{\Phi_{\mu_1}-\Phi_{\mu_0}}_{L^2(\rho)}
 \le C W_2(\mu_0,\mu_1).
\end{equation}
\end{lemma}

\begin{proof}
On the line, or if either target is a point mass, the two Brenier maps
induce an optimal coupling. In these cases \eqref{eq:projected-Poincare}
with $H=I$ gives the claim directly.

\emph{Step 1. The potential along the geodesic.}
Write an optimal plan as $\sum_{i=1}^N p_i\delta_{(y_i^0,y_i^1)}$ with distinct
pairs, and put $y_i(t)=y_i^0+tv_i$, $v_i=y_i^1-y_i^0$. The sites are
distinct for $0<t<1$ by $(y_i^0-y_j^0)\cdot(y_i^1-y_j^1)\ge0$.
When the source density is only measurable, its values on a given
cell face need not agree with its traces from shrinking slabs. We
therefore vary the masses first and recover the original ones by a
limit. For positive $q_i$ with $\sum_iq_i=1$, let $\Phi_t^q$ be the
centered potential for $\sum_iq_i\delta_{y_i(t)}$. Up to a constant it
is the maximum $\phi=\max_i\ell_i$ of affine functions
$\ell_i(x)=x\cdot y_i(t)+a_i(t)$. Each slope occurs on a cell
$C_i=\{\phi=\ell_i\}$ of mass $q_i$, and we fix the common additive
constant in the weights by $\sum_i a_i=0$.

For fixed $q$ the weights are Lipschitz in $t$. To see this, split the
ordered increments $a_i(s)-a_i(t)$ at any gap $\gamma$. Let $G_t$ be the
maximum over the upper group $I$ minus the maximum over its complement.
With $\delta=|s-t|\max_i\sup_\Omega|x\cdot v_i|$, we have
\[
 G_s\ge G_t+\gamma-2\delta,\qquad
 \rho(G_s>0)=\rho(G_t>0)=\sum_{i\in I}q_i .
\]
Since $G_t$ takes both signs on the connected set $\Omega$ and $r\ge m$,
a gap $\gamma>2\delta$ would increase the mass of the positive set.
The increments thus have oscillation at most $2(N-1)\delta$ and sum
zero. It follows that $\Phi_t^q$ is Lipschitz in $L^2(\rho)$. Since ties
have zero $\rho$-mass, differentiation on the cells gives for almost every $t$
\[
 \partial_t\Phi_t^q=f-\int f\dd\rho,\qquad
 f(x)=x\cdot v_i+\dot a_i(t)\quad(x\in C_i^\circ).
\]

\emph{Step 2. Conservation of the cell masses.}
We next justify differentiating the identity $\rho(C_i)=q_i$ at
almost every time. Take Borel representatives of $r$ and $1/r$, extended
by zero off $\Omega$. At fixed $t$, the offset of the plane where two
affine functions agree is a nonconstant linear function of the weights.
By Fubini and Lebesgue differentiation, for almost every $(t,a)$ the two
extensions have Lebesgue points both in $\R^d$ and along the normal
at $\HH^{d-1}$-almost every point of each plane. Coincidence of the planes
for two different pairs imposes a nontrivial linear condition on $a$.

The map $(t,a)\mapsto(t,(\rho(C_i))_i)$ is locally
Lipschitz for $0<t<1$, since a perturbation of size $h$ changes cells
on the bounded source only in slabs of width $O(h)$ and $r\le M$.
The parameter spaces have the same dimension, so the exceptional set
has null image. Thus for almost every $q$ these properties hold at
almost every $t$, where we may also differentiate the weights.

Let $F_{ij}$ be the open common face of $C_i,C_j$ in $\R^d$. Put
$D_{ij}=|y_i(t)-y_j(t)|$ and $n_{ij}=(y_i(t)-y_j(t))/D_{ij}$.
Write $[w]_{ij}=w_i-w_j$ for the jump from $C_j$ to $C_i$.
Since $\rho(C_i)=q_i$, the face variation formula gives
\begin{equation}\label{eq:mass-derivative}
 0=\sum_{j\ne i}\int_{F_{ij}\cap\Omega}
                   \frac{[f]_{ij}}{D_{ij}}\,r\dd S .
\end{equation}
Indeed, slice the mass exchanged by cells $i,j$ in time $h$ along
$n_{ij}$. Away from other ties, the equation $\ell_i=\ell_j$ gives
the signed outward displacement of $C_i$ as $h[f]_{ij}/D_{ij}+o(|h|)$
on $F_{ij}$.
Where a third branch is strictly larger, the slice contributes zero
for small $h$. The slab has width $O(|h|)$ on the bounded support of
$r$, so the normal Lebesgue point property and dominated convergence
give \eqref{eq:mass-derivative}
(see~\cite{deGournay} for smooth densities).

\emph{Step 3. The flux estimate.}
Fix masses and a time for which the preceding conclusions hold. Put
$b=v_i$ on $C_i^\circ$ and write $\diff S_{ij}$ for surface measure on
$F_{ij}$. To regularize the equation for the piecewise affine function $f$ above, we need the contributions of its jumps and of the jumps of
$\nabla\phi$. Integration by parts on the cells gives, on $\R^d$,
\[
 \D^2\phi=\sum_{i<j}D_{ij}\,n_{ij}\otimes n_{ij}\,\diff S_{ij},
 \qquad
 \D f=b\dx+\sum_{i<j}[f]_{ij}n_{ij}\,\diff S_{ij}.
\]
Take a standard even mollifier $\eta_\epsilon$ supported in $B_\epsilon$. Put
$H_\epsilon=\D^2(\phi*\eta_\epsilon+\epsilon |x|^2/2)$ and
$b_\epsilon=b*\eta_\epsilon$. The centered solution $u_\epsilon$ of
\eqref{eq:source-linearized} minimizes
\[
 E_\epsilon(w)=\int |H_\epsilon^{-1/2}(\nabla w-b_\epsilon)|^2\dd\rho .
\]
Equations \eqref{eq:stability-adjoint} and \eqref{eq:uniform-flux}
give $\norm{u_\epsilon}_2\le C\norm{b_\epsilon}_2$, with
$b_\epsilon\to b$ in $L^2(\rho)$. It remains to identify a weak
subsequential limit $u$ of $u_\epsilon$ with the centered derivative
of the potential. To compare the jumps of $u$ and $f$, define for any
function $w$ affine on each cell
\[
 E_0(w)=\sum_{i<j}\int_{F_{ij}\cap\Omega}
                  \frac{[w]_{ij}^2}{D_{ij}}\,r\dd S .
\]
Test $E_\epsilon$ with $f*\eta_\epsilon$ after subtracting its mean.
For any vector $\xi$, the inequality
$2[f]_{ij}(\xi\cdot n_{ij})-D_{ij}(\xi\cdot n_{ij})^2
 \le [f]_{ij}^2/D_{ij}$ holds on each face. Convolving the two
distributional identities above and maximizing in $\xi$ gives
\[
 E_\epsilon(u_\epsilon)
 \le \sum_{i<j}\int_{F_{ij}}\frac{[f]_{ij}^2}{D_{ij}}
                      (r*\eta_\epsilon)\dd S
 \longrightarrow E_0(f).
\]
Here $r*\eta_\epsilon\to r$ at almost every face point by Step~2.
On compact subsets of $C_i^\circ\cap\Omega$ we have
$H_\epsilon=\epsilon I$ and $b_\epsilon=v_i$ for small $\epsilon$.
The energy bound gives $\nabla u_\epsilon\to v_i$ locally in $L^2$,
so $u=f+k_i$ on this cell for some constant $k_i$.

The constants $k_i$ are determined by the jumps across the faces.
To obtain a lower bound involving those jumps, test with
$X\in C_c^\infty(\Omega;\R^d)$ and complete the square along the
weakly convergent subsequence. This gives
\[
 \liminf_{\epsilon\downarrow0} E_\epsilon(u_\epsilon)
 \ge \sum_{i<j}\int_{F_{ij}\cap\Omega}
 \left(2[u]_{ij}\,X\cdot n_{ij}
       -\frac{D_{ij}}{r}(X\cdot n_{ij})^2\right)\dd S .
\]
Indeed, weak convergence applies to
$\int\nabla u_\epsilon\cdot X\dx=-\int u_\epsilon(\diver X/r)\dd\rho$.
Integration by parts on the cells then gives the jumps of $u$. The
formula for $\D^2\phi$ and the Lebesgue points of $1/r$ give the quadratic
term. Smooth normal
data compactly supported in $F_{ij}\cap\Omega$ extend to test fields in
disjoint tubes. Since $m\le r\le M$, maximizing by density gives
$E_0(u)\le\liminf E_\epsilon(u_\epsilon)\le E_0(f)$.
By \eqref{eq:mass-derivative} the cross terms cancel in
\[
 E_0(u)-E_0(f)=\sum_{i<j}\int_{F_{ij}\cap\Omega}
                    \frac{(k_i-k_j)^2}{D_{ij}}\,r\dd S .
\]
The right side is nonnegative, while $E_0(u)\le E_0(f)$. Thus $u-f$
has zero jump on each interior face. It is constant inside each cell
and hence on the connected domain $\Omega$. Centering now gives
$u=f-\int f\dd\rho=\partial_t\Phi_t^q$, so weak lower semicontinuity
proves
\begin{equation}\label{eq:atomic-speed}
 \norm{\partial_t\Phi_t^q}_2\le C\norm{b}_2
   =C\Bigl(\sum_iq_i|v_i|^2\Bigr)^{1/2}.
\end{equation}

Integrate \eqref{eq:atomic-speed} and let $q\to p$. All targets lie
in a fixed ball, so their centered potentials are uniformly bounded
and Lipschitz on $\overline\Omega$. Under convergence of the targets,
almost everywhere convergence of convex gradients and Brenier uniqueness
identify each uniform limit, including at $t=0,1$ where sites may merge.
Since $\sum_i p_i|v_i|^2=W_2^2(\mu_0,\mu_1)$, we obtain
\eqref{eq:finite-target-stability}.
\end{proof}

\subsection{Weighted stability from the pressure condition}
\label{sec:stability-curvature}
For an isotropic log-concave source $r=Z^{-1}e^{-V}$ with $\min V=0$,
Wang and Wang~\cite[Theorem~1]{WangWang} bound
$\inf_b\|(f-b)/(4d+V)\|_{L^2(\rho)}$ by $3\sqrt e\,W_2(\mu,\nu)$
for a difference $f=\Phi_\mu-\Phi_\nu$ of Lipschitz Brenier potentials.
Their weight controls growth where $r$ is small and requires no strict
convexity of $V$. A complementary choice of weight comes from an
internal energy $U$. Indeed, along a density perturbation $(1+tg)\rho$,
its second derivative at $t=0$ is $\int rU''(r)g^2\dd\rho$.
The linearized transport equation pairs such perturbations with
changes of the potential, so the dual weight $1/(rU''(r))$ occurs
in the norm of that change. McCann's pressure condition and a
curvature bound on $-U'(r)$ give the following estimate.

\begin{theorem}\label[theorem]{thm:weighted-brenier}
Let $\Omega\subseteq\R^d$ be open and convex, and let $\rho=r\one_\Omega\dx
\in\Ptwo(\R^d)$ with $0<r<\infty$ almost everywhere on $\Omega$. Suppose
$U\in C([0,\infty))\cap C^\infty((0,\infty))$, $U(0)=0$, $U''>0$,
and $p(s)=sU'(s)-U(s)$ satisfies $p(s)\ge0$ and
$sp'(s)\ge (1-1/d)p(s)$ for $s>0$.
Assume that $1/(rU''(r))\in L^1(\rho)$ and that $-U'(r)$ agrees
almost everywhere on $\Omega$ with a finite $\kappa$-strongly convex
function for some $\kappa>0$.
Then for any two targets $\mu,\nu\in\Ptwo(\R^d)$ their Brenier
potentials satisfy
\begin{equation}\label{eq:weighted-brenier}
 \inf_{b\in\R}\int\frac{(\Phi_\mu-\Phi_\nu-b)^2}{rU''(r)}\dd\rho
 \le \frac{1}{\kappa}W_2^2(\mu,\nu).
\end{equation}
\end{theorem}
\begin{proof}
Write $A=rU''(r)$. First work on a bounded smooth uniformly convex
domain with $r$ smooth and positive up to its boundary. Given $H=\D^2\Phi>0$ smooth
on the closure and a smooth centered $g$, solve $Bu=g$ with conormal
boundary condition. The compatible flux $X=H^{-1}\nabla u$ is tangent
to $\partial\Omega$. The pressure condition makes the volume terms in
\cref{prop:volume-variation} nonnegative, and convexity of $\Omega$
gives $p(r)\II(X,X)\ge0$ on the boundary. Thus
\eqref{eq:complete-second-variation} gives
\[
 \int Ag^2\dd\rho+\int \D^2 U'(r)[X,X]\dd\rho\ge0 .
\]
Since $-U'(r)$ is $\kappa$-strongly convex,
$\kappa\|X\|_\rho^2\le\int Ag^2\dd\rho$.

To use this flux estimate, take a smooth path of positive target
densities on a fixed closed ball as in \cref{lem:smooth-paths}, with
velocity $v_t$ and potential $\Phi_t$. The smooth domains and data give
boundary regularity \cite{TrudingerWang}, and smooth dependence of the
target log density on $t$ gives time differentiability by
\cite[Corollary~3.5]{GonzalezSanzSheng}. Choose the constants in the
potentials by $\int\Phi_t A^{-1}\dd\rho=0$. Then
$g=\dot\Phi_t/A$ has zero mean, and testing the transport derivative
against this $g$ produces the desired weighted norm. Specifically,
\eqref{eq:fixed-source-pairing} gives
\[
 \int\frac{\dot\Phi_t^2}{A}\dd\rho
 =\int X\cdot(v_t\circ\nabla\Phi_t)\dd\rho
 \le\left(\frac1\kappa\int\frac{\dot\Phi_t^2}{A}\dd\rho\right)^{1/2}
                \|v_t\|_{\mu_t},
 \qquad \mu_t=(\nabla\Phi_t)_\#\rho .
\]
Integration and \cref{lem:smooth-paths} give \eqref{eq:weighted-brenier}
for compact targets. For endpoints in a fixed ball the potentials are
equi-Lipschitz, and normalization and Brenier uniqueness identify their
uniform limits.

To approximate a general source while preserving curvature, we smooth
the strongly convex representative $P$ of $-U'(r)$. Exhaust $\Omega$
by smooth uniformly convex domains $K_j\Subset\Omega$ and mollify $P$
near $\overline K_j$ with a symmetric kernel of radius
$\varepsilon_j\downarrow0$ smaller than the distance to
$\partial\Omega$. Each $P_j$ is $\kappa$-strongly convex and, being a
local average of $-U'(r)$, takes values in the open range of $-U'$.
By Jensen's inequality
$P_j\ge P$, so $\tilde r_j=(-U')^{-1}(P_j)\le r$ almost everywhere on
$K_j$, and $\tilde r_j\one_{K_j}\to r$ almost everywhere.
Thus $s_j=\int_{K_j}\tilde r_j\dx\to1$ by dominated convergence.
Apply the smooth estimate to
$\rho_j=s_j^{-1}\tilde r_j\one_{K_j}\dx$ and $U_j(z)=U(s_jz)/s_j$.
This rescaling preserves the pressure condition and gives
$-U_j'(\tilde r_j/s_j)=P_j$.

The source changes in this approximation, but we can keep both
potentials fixed by taking as targets their image laws
$\mu_j=(\nabla\Phi_\mu)_\#\rho_j$ and $\nu_j=(\nabla\Phi_\nu)_\#\rho_j$.
The gradients are bounded on each $K_j$ and belong to $L^2(\rho)$.
Since $\tilde r_j\le r$ and $s_j\to1$, dominated convergence gives
$\mu_j\to\mu$ and $\nu_j\to\nu$ in $W_2$. For
$f=\Phi_\mu-\Phi_\nu$ and minimizing constants $b_j$, the estimate on
$K_j$ reads
\[
 \frac1{s_j}\int_{K_j}\frac{(f-b_j)^2}{U''(\tilde r_j)}\dx
 \le \frac1\kappa W_2^2(\mu_j,\nu_j).
\]
The weights on the left converge almost everywhere to $1/U''(r)>0$.
Fatou's lemma rules out $|b_j|\to\infty$ along any subsequence and,
along a subsequence with $b_j\to b$, implies \eqref{eq:weighted-brenier}.
\end{proof}

Kolesnikov and Milman obtain weighted Poincar\'e and log-Sobolev
inequalities for a single reference measure by choosing an auxiliary
metric and expressing its gradient estimates in Euclidean
coordinates~\cite{KMMetrics}. In \cref{thm:weighted-brenier} the
reciprocal weight comes from the second variation of the internal
energy. The linearized transport equation then pairs the difference
of two potentials with a velocity of the target. For inverse power
densities the resulting coefficient is optimal over the class of sources
described below.

\Needspace{7\baselineskip}
\begin{theorem}
\label[theorem]{thm:inverse-power-stability}
Let $\Omega\subseteq\R^d$ be open and convex, with $d\ge2$, and let
$W:\Omega\to(0,\infty)$ be $\eta$-strongly convex, with $\eta>0$. 
Suppose that $\beta\ge d$ and that
$\diff\rho=Z^{-1}W^{-\beta}\one_\Omega\dx$ defines a probability
measure in $\Ptwo(\R^d)$. For each $\mu\in\Ptwo(\R^d)$, choose the
Brenier potential $\Phi_\mu$ from $\rho$ to $\mu$ so that
$\int_\Omega\Phi_\mu W^{-1}\dd\rho=0$. Then for all $\mu,\nu\in\Ptwo(\R^d)$,
\begin{equation}\label{eq:inverse-power-stability}
 \norm{\Phi_\mu-\Phi_\nu}_{L^2(W^{-1}\rho)}^2
 \le \frac{1}{\beta\eta}\,W_2^2(\mu,\nu).
\end{equation}
The coefficient $1/(\beta\eta)$ is optimal over this class of sources.
\end{theorem}

\begin{proof}
For the normalization of the potentials we need $W^{-1}\rho$ to have
finite mass. In fact $W$ is bounded away from zero. Otherwise fix an
interior ball $B$ on which $W\le M$. For $0<s<M$, choose
$x_s$ with $W(x_s)<s/2$ and put $\lambda=s/(2M)$. Convexity gives
$(1-\lambda)x_s+\lambda B\subset\{W<s\}$, hence
\[
 Z=\int_\Omega W^{-\beta}\dx
 \ge\beta\int_0^M |\{W<s\}|s^{-\beta-1}\dd s
 \ge\frac{\beta|B|}{(2M)^d}\int_0^M s^{d-\beta-1}\dd s=\infty .
\]
This contradiction gives $\inf_\Omega W>0$.

To apply the preceding theorem, we want $-U'(r)$ to be a multiple
of $W$ when $r=Z^{-1}W^{-\beta}$. This leads to
$m=1-1/\beta$ and $U(s)=s^m/(m-1)$. Its pressure $p(s)=s^m$
satisfies the required inequality since $\beta\ge d$, and
\[
 rU''(r)=\frac{\beta-1}{\beta}Z^{1/\beta}W,\qquad
 -U'(r)=(\beta-1)Z^{1/\beta}W .
\]
Thus \cref{thm:weighted-brenier} applies with
$\kappa=(\beta-1)Z^{1/\beta}\eta$ and gives
\[
 \inf_{b\in\R}\int_\Omega\frac{(\Phi_\mu-\Phi_\nu-b)^2}{W}\dd\rho
 \le \frac1{\beta\eta}W_2^2(\mu,\nu).
\]
For $\nu=\delta_0$ we may take $\Phi_\nu=0$. Since $W^{-1}\rho$ is
finite, the estimate gives $\Phi_\mu\in L^2(W^{-1}\rho)$ for every
$\mu$. With the normalization in the statement, $\Phi_\mu-\Phi_\nu$ has weighted mean
zero, so $b=0$ attains the infimum and \eqref{eq:inverse-power-stability}
follows.

For sharpness take $W(x)=1+\eta|x|^2/2$ on $B_R$ with normalized
source $\rho_R$ and targets $\delta_a,\delta_0$, $a\ne0$. Their
potential difference $a\cdot x$ has weighted mean zero by symmetry.
Integrating $\diver((a\cdot x)aW^{-\beta})$ gives
\[
 \frac{\beta\eta}{|a|^2}\int_{B_R}\frac{(a\cdot x)^2}{W}\dd\rho_R
 =1-\frac{|B_R|(1+\eta R^2/2)^{-\beta}}{\int_{B_R}W^{-\beta}\dx}
 \longrightarrow 1 \qquad (R\to\infty).
\]
Since $W_2(\delta_a,\delta_0)=|a|$, this proves optimality.
\end{proof} 

\begin{corollary}
\label[corollary]{cor:inverse-power-equality}
Let the source and normalized Brenier potentials be as in
\cref{thm:inverse-power-stability}. If distinct targets
$\mu,\nu\in\Ptwo(\R^d)$ attain equality in \eqref{eq:inverse-power-stability},
then for some $a\in\R^d$ and $b\in\R$ their Brenier maps satisfy
\begin{equation}\label{eq:inverse-power-equality}
 \nabla\Phi_\nu(x)-\nabla\Phi_\mu(x)=a+bx\quad \rho\text{-a.e.},
 \qquad b=0\quad\text{if }\beta>d .
\end{equation}
Thus $\nu=(\Id+a)_\#\mu$ when $\beta>d$. 
If the source domain $\Omega$ is bounded, the inequality is strict for
all distinct targets. 
\end{corollary}

\begin{proof}
For $f=\Phi_\nu-\Phi_\mu$, \cref{thm:inverse-power-stability} and the
coupling induced by the two Brenier maps give
\begin{equation}\label{eq:stability-midpoint-equality}
 \beta\eta\int f^2 W^{-1}\dd\rho
 \le W_2^2(\mu,\nu)\le\int|\nabla f|^2\dd\rho .
\end{equation}
If equality holds in the first inequality, we claim that the two
Brenier maps also induce an optimal coupling of their targets. Indeed,
apply the bound on each half of a target geodesic. Equality in the
Hilbert space triangle inequality identifies the normalized midpoint
potential with $(\Phi_\mu+\Phi_\nu)/2$. The coupling induced by the
two maps and the optimal coupling defining the geodesic therefore
have the same midpoint law. For either coupling $(Y_0,Y_1)$ we have
\[
 \E\left|\frac{Y_0+Y_1}{2}\right|^2
 =\frac{\E|Y_0|^2+\E|Y_1|^2}{2}-\frac14\E|Y_1-Y_0|^2
\]
so equality of the midpoint second moments implies equality of their
costs. Hence equality in the first bound of
\eqref{eq:stability-midpoint-equality} is equivalent to equality in
the weighted Poincar\'e inequality between its two outer terms.
For the smooth version of that inequality, see
Nguyen~\cite[Theorem~12, (5.2)]{NguyenVariance}.

On a smooth bounded uniformly convex domain with $W$ smooth and positive
up to the boundary, apply the pressure identity from the proof of
\cref{thm:weighted-brenier} with $U(s)=-\beta s^{1-1/\beta}$.
For a compatible tangent field $X$ and $g=-\diver_\rho X$, division
by $Z^{1/\beta}$ gives
\begin{equation}\label{eq:inverse-power-bochner}
\begin{aligned}
 \frac{\beta-1}{\beta}\int Wg^2\dd\rho
 ={}&\int W\left[\tr((\D X)^2)-\frac{(\diver X)^2}{\beta}\right]\dd\rho\\
 &+(\beta-1)\int \D^2W[X,X]\dd\rho
   +\frac1Z\int_{\partial\Omega}W^{1-\beta}\II(X,X)\dd S .
\end{aligned}
\end{equation}

Exhaust $\Omega$ by smooth uniformly convex domains $K_j\Subset\Omega$.
On each $K_j$, mollify $W$ and approximate $f$ in $H^1(K_j)$ by
smooth functions, subtracting the mean with
respect to the mollified weight. At fixed $j$, the weights converge
uniformly and are bounded away from zero, so the energy and weighted
variance converge to those for $f$ on $K_j$. Choose the two smoothing
parameters so that these errors are at most $1/j$ and the local $H^1$
approximation errors tend to zero. On each smooth domain, solve
$-\diver_\rho X=\beta\eta f/W$ for a gradient field with zero normal
component, using the same notation for the approximations. The datum
is centered by our choice of the additive constant in $f$. Integration
by parts and \eqref{eq:inverse-power-bochner} with $H=I$ give
\[
\begin{aligned}
 \int\!\left(|\nabla f|^2-\frac{\beta\eta f^2}{W}\right)\dd\rho
 \ge \int|\nabla f-X|^2\dd\rho+\frac1{(\beta-1)\eta}\int W
 \left(\norm{\D X}_{\HS}^2-\frac{(\diver X)^2}{\beta}\right)\dd\rho .
\end{aligned}
\]
Under the assumed equality, convergence of the energy and weighted
variance makes the left side tend to zero. The first term on the right
then yields $X\to\nabla f$ locally in $L^2$, and the second forces the
trace-free part of $\D X$ to vanish in the same limit. It also forces
$\diver X\to0$ if $\beta>d$. In distributions we thus have
$\D^2f=(\Delta f/d)I$. Since $d\ge2$, commuting mixed derivatives
shows that $\Delta f/d$ is constant. This proves
\eqref{eq:inverse-power-equality}.

Suppose now that $\Omega$ is bounded. For $h\in C_c^\infty(\R^d)$
let $\bar h$ be its weighted mean for $W^{-1}\rho$ and set
$f_t=f+t(h-\bar h)$. By \eqref{eq:inverse-power-equality}, $f$ agrees
with a quadratic polynomial on $\Omega$. For $A$ sufficiently large,
$f_t+A|x|^2$ and $A|x|^2$ are convex with bounded gradients on
$\Omega$. Applying
\eqref{eq:stability-midpoint-equality} to their gradient images gives
$\int|\nabla f_t|^2\dd\rho\ge\beta\eta\int f_t^2W^{-1}\dd\rho$.
Equality holds at $t=0$. Differentiating there and using
$\int fW^{-1}\dd\rho=0$ gives
\begin{equation}\label{eq:inverse-power-neumann}
	\int (a+bx)\cdot\nabla h\dd\rho
	=\beta\eta\int fh\,W^{-1}\dd\rho .
\end{equation}
We can now rule out each nonconstant polynomial allowed by
\eqref{eq:inverse-power-equality}. Since $\Omega$ is bounded and
$W$ is bounded away from zero, $f/W$ is bounded. We may also multiply
any smooth test by a cutoff equal to one near $\overline\Omega$.
For distinct targets, $a+bx$ is nonzero on a set of positive
$\rho$-measure. Testing first with $h(x)=|a+bx|^{2k}$ gives
$2k|b|\le\beta\eta\norm{f/W}_\infty$ for every positive integer $k$,
so $b=0$. Taking $h(x)=e^{ka\cdot x}$ then gives
$k|a|^2\le\beta\eta\norm{f/W}_\infty$, so $a=0$ as well.
The two Brenier maps therefore coincide, contradicting $\mu\ne\nu$. 
\end{proof} 

\begin{remark}\label[remark]{mm:rem:fixed-source}
Let $\diff\rho=Z^{-1}e^{-V}\dx$ with $V-\kappa|\cdot|^2/2$ convex for
some $\kappa>0$. Letrouit and Manole~\cite[Theorem~1.1]{LetrouitManole}
proved that for all $\mu,\nu\in\Ptwo(\R^d)$,
\begin{equation}\label{mm:eq:fixed-stability}
 \Var_\rho(\Phi_\nu-\Phi_\mu)\le\kappa^{-1}W_2^2(\mu,\nu).
\end{equation}
For finite $\kappa$-strongly convex source potentials, Ruan obtained the
same $L^2(\rho)$ coefficient $\kappa^{-1/2}$ for the potentials
\cite[Theorem~1.1(ii), Eq.~(1.6)]{Ruan}. For targets in a common ball,
Ruan also bounds the $L^2(\rho)$ distance between their Brenier maps by
a constant times $W_2(\mu,\nu)^{1/3}$ and proves that the exponent is
optimal in $d\ge2$ \cite[Theorems~1.1(iii) and~1.2]{Ruan}.

For $d\ge2$, normalize $\inf V=0$ and set $W_\beta=e^{V/\beta}$ on
$\operatorname{int}\{V<\infty\}$. Since $s\mapsto e^s-s$ is convex
and nondecreasing on $[0,\infty)$, $W_\beta-V/\beta$ is convex.
In \eqref{eq:inverse-power-stability} take $\eta=\kappa/\beta$.
Then $W_\beta^{-\beta}=e^{-V}$ and $W_\beta^{-1}\uparrow1$, so the
double-integral formula for variance yields \eqref{mm:eq:fixed-stability}.

In every dimension, distinct targets attain equality if and only if
$\nu=(\Id+a)_\#\mu$ for some $a\ne0$ and $\rho$ has a Gaussian factor
of variance $\kappa^{-1}$ in the direction of $a$. The midpoint argument
of \cref{cor:inverse-power-equality} gives Poincar\'e equality for
$f=\Phi_\nu-\Phi_\mu$. On the line the same conclusion follows
directly because the increasing Brenier maps induce the optimal
coupling. For the approximating Neumann fluxes $X$,
\eqref{eq:weighted-bochner} with $w=1$ gives $X\to\nabla f$ and $\D X\to0$
locally in $L^2$ on the interior of the support. Thus $\nabla f=a$
$\rho$-a.e. Polarization gives
$\int a\cdot\nabla h\dd\rho=\kappa\Cov_\rho(a\cdot x,h)$ for
$h\in C_c^\infty(\R^d)$, the Gaussian factor criterion of
\cite[Proposition~4.2]{HanZhu} applied to $(\rho,\rho)$.
Conversely, for such a translation $f=a\cdot x+\mathrm{const}$ and
$\kappa\Var_\rho(f)=|a|^2=W_2^2(\mu,\nu)$.
\end{remark}

\subsection{Support convexity and the role of the density}\label{sec:nonconvex}
Convexity of the source domain entered the preceding estimates through
the sign of the boundary term in \eqref{eq:inverse-power-bochner}.
We now ask how much can be inferred about the source from finiteness
of the unweighted stability constant $S(\rho)$. For a uniform source
this question has a geometric answer in terms of its closed support.
For a variable density, even on a square, we will also need to examine
how the conditional distributions change from one slice to the next.

Assume $d\ge2$ and let $\Omega\subseteq\R^d$ be bounded, connected and
open. For every probability measure $\rho=r\one_\Omega\dx$ with
$r>0$ almost everywhere, $S(\rho)$ is well defined.
Indeed, compact targets admit globally Lipschitz convex potentials.
Brenier uniqueness \cite{Brenier} identifies their gradients
$\rho$-almost everywhere, hence Lebesgue-almost everywhere on $\Omega$.
Connectedness then makes their difference constant, so centering
determines the potential.

\begin{theorem}\label{thm:nonconvex-support}
Suppose $|\partial\Omega|=0$, and let
$\rho_\Omega=|\Omega|^{-1}\one_\Omega\dx$. Then
\[
 S(\rho_\Omega)<\infty
 \quad\Longleftrightarrow\quad \overline\Omega\text{ is convex}.
\]
\end{theorem}

\begin{proof}
If $K=\overline\Omega$ is convex, then $|\partial\Omega|=0$ gives
$\rho_\Omega=\one_K\dx/|K|$. This source is log-concave, so
\cite[Theorem~1.4]{LetrouitManole} gives $S(\rho_\Omega)\le2R$
whenever $K\subseteq B(x_0,R)$.
For uniform sources on convex bodies, Letrouit~\cite[Conjecture~3.10]{LetrouitPeccot}
conjectured the existence of a Lipschitz bound for the potentials when
the targets lie in a fixed compact set.

Conversely, a gap in the closed support can be detected by moving the
corner of a piecewise affine transport potential across it. We carry
out the construction for every $\rho=r\one_\Omega\dx$ with $r>0$
almost everywhere on $\Omega$. If $\overline\Omega$ is nonconvex,
choose a segment joining two of its points and passing through its
open complement. We can perturb the endpoints into $\Omega$ while
retaining an interval in the complement. Taking the first coordinate
along the perturbed segment, choose $h>0$ and $s_0<s_1$ so that the tube
$C=\{(t,y):|y_i|<h,\ 1\le i\le d-1\}$ satisfies
\begin{equation}\label{eq:atomic-gap}
 \overline\Omega\cap
 \{(t,y):s_0\le t\le s_1,\ |y_i|\le h\ (1\le i\le d-1)\}
 =\varnothing.
\end{equation}
Choose balls $E_-,E_+\Subset\Omega\cap C$, respectively in $t<s_0$
and $t>s_1$, with $|y_i|<h/2$ on both. They have positive $\rho$-mass.

Letrouit and M\'erigot~\cite[Section~5]{LetrouitMerigot} move the corner
of $|t-s|$ across increasingly thin passages to disprove H\"older stability.
Here the gap \eqref{eq:atomic-gap} allows us to vary the corner between
$s_0$ and $s_1$ inside the tube. Outside the tube we must control the
masses sent to other target points. We do this by adjusting the affine
pieces of the potential using the finite-dimensional dual formulation of
semi-discrete transport~\cite[Theorem~1.1]{Kitagawa}. Start with
\[
 \Phi_0(t,y)=\max\{0,y_1-h,-y_1-h,\ldots,y_{d-1}-h,-y_{d-1}-h\}.
\]
Discard each outer branch that is the unique maximizer only on a null set. The tie
hyperplanes are null, so $\Phi_0$ and its gradient are unchanged
$\rho$-almost everywhere. The remaining branches $\ell_1,\ldots,\ell_N$
have distinct unit slopes $v_i=\nabla\ell_i\in e_1^\perp$. Put
\[
 p_0=\rho(C),\qquad
 p_i=\rho\{\ell_i>\max(0,\ell_k:k\ne i)\}>0,\qquad
 \sum_{i=0}^Np_i=1.
\]
Here $1\le N\le2(d-1)$. Indeed, a path in $\Omega$ from $E_-$
to $E_+$ must leave $C$. Openness of the domain and positivity of
the density therefore give positive mass outside the tube.

We now split the branch of slope zero into two branches with opposite
slopes. For $s\in\{s_0,s_1\}$ and $\epsilon>0$, set
\begin{equation}\label{eq:atomic-perturbed-potential}
 \Phi_{\epsilon,s,a}(t,y)
 =\max\{\epsilon|t-s|,\ell_1(t,y)+a_1,\ldots,\ell_N(t,y)+a_N\}
\end{equation}
and minimize
\begin{equation}\label{eq:atomic-dual}
 F_{\epsilon,s}(a)
 =\int\Phi_{\epsilon,s,a}\dd\rho-\sum_{i=1}^Np_i a_i,
 \qquad a\in\R^N.
\end{equation}
To see that offsets realizing the prescribed outer masses exist, put
$a_0=0$ and $M=\max_{1\le i\le N}\sup_\Omega|\ell_i|$. Then
\begin{equation}\label{eq:atomic-coercivity}
 F_{\epsilon,s}(a)
 \ge \sum_{i=0}^N p_i\bigl(\max_{0\le j\le N}a_j-a_i\bigr)-M
 \ge \bigl(\min_{0\le i\le N}p_i\bigr)\norm a_\infty-M.
\end{equation}
Indeed, $\max_i a_i-\min_i a_i\ge\norm a_\infty$ since $a_0=0$.
Thus a minimizer $a_{\epsilon,s}$ exists. The slopes
$\pm\epsilon e_1,v_1,\ldots,v_N$ are distinct, so differentiation
under the integral gives
\begin{equation}\label{eq:atomic-fixed-masses}
 \rho(\text{outer cell }i\text{ of }\Phi_{\epsilon,s,a_{\epsilon,s}})
 =p_i\qquad(1\le i\le N).
\end{equation}

At $\epsilon=0$, $a=0$ gives the prescribed masses. Any minimizing
potential transports $\rho$ to $\mu_*=p_0\delta_0+\sum_i p_i\delta_{v_i}$.
By Brenier uniqueness and connectedness it differs from $\Phi_0$ by a
constant. Both vanish on the common cell of slope zero, and each outer
cell has positive mass. Hence $a=0$ is the unique minimizer of $F_0$.
Moreover,
\begin{equation}\label{eq:atomic-uniform-dual-limit}
 \sup_a|F_{\epsilon,s}(a)-F_0(a)|
 \le\epsilon\sup_{(t,y)\in\Omega}|t-s|\longrightarrow0.
\end{equation}
Since $F_{\epsilon,s}(0)$ is bounded, \eqref{eq:atomic-coercivity}
bounds the minimizers, and every limit point minimizes $F_0$. Hence
\begin{equation}\label{eq:atomic-offset-limit}
 a_{\epsilon,s}\longrightarrow0\qquad(s=s_0,s_1).
\end{equation}

It remains to compare the changes in the targets and in their
potentials as $s$ moves across the gap. Write
$\Phi_{\epsilon,s}=\Phi_{\epsilon,s,a_{\epsilon,s}}$ and
$\mu_{\epsilon,s}=(\nabla\Phi_{\epsilon,s})_\#\rho$. These maps are
optimal by convexity. If $q_{\epsilon,s}$ is the mass sent to
$-\epsilon e_1$, the fixed outer masses give
\begin{equation}\label{eq:atomic-targets}
 \mu_{\epsilon,s}=q_{\epsilon,s}\delta_{-\epsilon e_1}
 +(p_0-q_{\epsilon,s})\delta_{\epsilon e_1}
 +\sum_{i=1}^Np_i\delta_{v_i}.
\end{equation}
Away from the tie hyperplanes and $t=s$, the indicator of the region
sent to $-\epsilon e_1$ converges to $\one_{C\cap\{t<s\}}$ by
\eqref{eq:atomic-offset-limit}. Since the portion of $C$ between
$s_0$ and $s_1$ has zero mass by \eqref{eq:atomic-gap}, dominated
convergence yields
\[
 q_{\epsilon,s}\to\rho(C\cap\{t<s\}),\qquad
 \delta_\epsilon:=|q_{\epsilon,s_1}-q_{\epsilon,s_0}|\to0.
\]
Coupling the common mass to itself leaves only $\delta_\epsilon$
to move between the two central atoms. Thus
\begin{equation}\label{eq:atomic-cost-upper}
 W_2^2(\mu_{\epsilon,s_0},\mu_{\epsilon,s_1})
 \le4\epsilon^2\delta_\epsilon=o(\epsilon^2).
\end{equation}
Both families lie in the closed unit ball for $\epsilon\le1$ and
converge in $W_2$ to $\mu_*$.

On $E_-$ and $E_+$, all outer branches are at most $-h/2$ before
their offsets are added. For small $\epsilon$,
$\Phi_{\epsilon,s}=\epsilon|t-s|$ there. Thus, with
$\Delta=s_1-s_0$ and $f_\epsilon=\Phi_{\epsilon,s_1}-\Phi_{\epsilon,s_0}$,
we have $f_\epsilon=\epsilon\Delta$ on $E_-$ and
$f_\epsilon=-\epsilon\Delta$ on $E_+$, whence
\begin{equation}\label{eq:atomic-variance-lower}
 \begin{split}
 \norm{\Phi_{\mu_{\epsilon,s_1}}-
             \Phi_{\mu_{\epsilon,s_0}}}_{L^2(\rho)}^2=\Var_\rho f_\epsilon
 =\frac12\iint|f_\epsilon(x)-f_\epsilon(y)|^2\dd\rho(x)\dd\rho(y)\ge4\epsilon^2\Delta^2\rho(E_-)\rho(E_+).
 \end{split}
\end{equation}
For small $\epsilon$, $\delta_\epsilon>0$, since otherwise the targets
coincide and uniqueness contradicts \eqref{eq:atomic-variance-lower}.
Combining the two estimates gives
\begin{equation}\label{eq:atomic-quotient}
 \frac{\norm{\Phi_{\mu_{\epsilon,s_1}}-
                  \Phi_{\mu_{\epsilon,s_0}}}_{L^2(\rho)}}
      {W_2(\mu_{\epsilon,s_1},\mu_{\epsilon,s_0})}
 \ge\frac{\Delta\sqrt{\rho(E_-)\rho(E_+)}}{\sqrt{\delta_\epsilon}}
 \longrightarrow\infty.
\end{equation}
This proves $S(\rho)=\infty$, and in particular establishes necessity
for the uniform source.

For every $a\in\R$, $|f_\epsilon-a|\ge\epsilon\Delta$ on at least one
of the two balls. If $W^{-1}\rho$ is finite and gives both balls
positive mass, then $\inf_a\norm{f_\epsilon-a}_{L^2(W^{-1}\rho)}
\ge c\epsilon$ with $c>0$ independent of $\epsilon$. Thus
\eqref{eq:atomic-cost-upper} also rules out a finite Lipschitz
constant in the weighted norm of \eqref{eq:inverse-power-stability}.

Choose a similarity $F(y)=b+\lambda y$, $\lambda>0$, mapping the unit
ball into the prescribed compact set with nonempty interior. Then
\begin{equation}\label{eq:target-similarity}
 \Phi_{F_\#\mu}(x)=\lambda\Phi_\mu(x)+b\cdot\left(x-\int z\dd\rho(z)\right),
 \qquad W_2(F_\#\mu,F_\#\nu)=\lambda W_2(\mu,\nu).
\end{equation}
This preserves the quotient, the numbers of atoms, and convergence
of both families to the same limit.
\end{proof}

The necessity argument used only positivity of the density on
$\Omega$. It therefore gives a corresponding characterization for
truncated log-concave densities. If $|\partial\Omega|=0$ and
$V:\R^d\to\R$ is convex, then $\diff\rho=Z^{-1}e^{-V}\one_\Omega\dx$
satisfies
\begin{equation}\label{eq:truncated-convex-characterization}
 S(\rho)<\infty
 \quad\Longleftrightarrow\quad \overline\Omega\text{ is convex}.
\end{equation}
If $\overline\Omega$ is convex, we may replace $\Omega$ by
$\overline\Omega$ in the density and apply
\cite[Theorem~1.4]{LetrouitManole}. In this case, if $V$ is
$\kappa$-strongly convex with $\kappa>0$, then $S(\rho)\le\kappa^{-1/2}$
by \cite[Theorem~1.1]{LetrouitManole}. If $\overline\Omega$ is nonconvex,
the preceding construction applies since $e^{-V}>0$ on $\Omega$.

For the normalized restriction $\rho$ of a source from
\cref{thm:inverse-power-stability} to a bounded connected open subset
$\Omega$ of its domain with $|\partial\Omega|=0$, the $L^2(W^{-1}\rho)$
stability estimate holds with a finite constant if and only if
$\overline\Omega$ is convex. In that case the theorem applies on
$\operatorname{int}\overline\Omega$.

The same construction also locates the failure of stability in target
space. For any source with positive density almost everywhere on a
bounded connected open domain with nonconvex closure, the individual
targets have at most $2d$ atoms. We can confine both to any prescribed
compact set with nonempty interior. They converge to a common measure
$\mu_*$ with at most $2d-1$ atoms in an affine hyperplane, so
\eqref{eq:atomic-quotient} rules out Lipschitz stability in every
$W_2$-neighborhood of $\mu_*$.

On a John domain the geometric characterization requires no separate
hypothesis on the boundary measure. Recall that this domain condition
asks for $x_*\in\Omega$ and $0<\alpha\le1$ such that each $x\in\Omega$
can be joined to $x_*$ by an arclength-parametrized curve
$\gamma:[0,L]\to\Omega$, $\gamma(0)=x$, satisfying
\begin{equation}\label{eq:john-condition}
 \dist(\gamma(t),\partial\Omega)\ge\alpha t\qquad(0\le t\le L).
\end{equation}

\begin{corollary}\label[corollary]{cor:john-domain}
Let $\rho_\Omega$ be the uniform probability measure on a bounded John
domain $\Omega\subseteq\R^d$, $d\ge2$. Then
\begin{equation}\label{eq:john-essential-convexity}
 S(\rho_\Omega)<\infty
 \quad\Longleftrightarrow\quad \overline\Omega\text{ is convex}.
\end{equation}
\end{corollary}

\begin{proof}
To apply the characterization, we show that the John condition prevents
any boundary point from being a density point of $\partial\Omega$.
For $z\in\partial\Omega$ and $0<r<\dist(x_*,\partial\Omega)$, choose
$x\in\Omega\cap B(z,r/8)$ and follow a John curve from $x$ to $x_*$.
Its length is at least $\dist(x_*,\partial\Omega)-r/8>r/4$, so
$w=\gamma(r/4)$ satisfies
\[
 |w-z|\le3r/8,\qquad \dist(w,\partial\Omega)\ge\alpha r/4.
\]
Thus $B(w,\alpha r/8)\subseteq\Omega\cap B(z,r)$, and
$|\partial\Omega\cap B(z,r)|/|B(z,r)|\le1-(\alpha/8)^d$.
The Lebesgue density theorem gives $|\partial\Omega|=0$, and
\cref{thm:nonconvex-support} applies.
\end{proof}

\begin{remark}\label[remark]{rem:punctured-ball}
The distinction between $\Omega$ and its closure already appears for
the punctured ball $\Omega=B(0,1)\setminus\{0\}$ in $d\ge2$. This
nonconvex domain satisfies the John condition. To join a point
to $e_1/2$, move radially to the sphere of radius $1/2$ and then
follow a shortest arc there. Removing the origin leaves the uniform
measure unchanged, so $S(\rho_\Omega)=S(\rho_{B(0,1)})<\infty$.
\end{remark}

On a bounded John domain, let the source density be bounded above
and below by positive constants. For targets supported in a fixed
compact set $\mathcal Y$, Letrouit and
M\'erigot~\cite[Theorem~1.7]{LetrouitMerigot} proved that
\[
 \norm{\Phi_\mu-\Phi_\nu}_{L^2(\rho)}
 \le C_{\rho,\mathcal Y}\,W_1(\mu,\nu)^{1/2}.
\]
On a fixed square, \cref{prop:rough-counterexample} gives sources
with densities arbitrarily close to uniform in $L^\infty$ for which
this estimate holds and $S(\rho)=\infty$.

To detect the effect of a varying density, we need a necessary condition
that sees individual conditional measures. Anisotropic quadratic costs
can single out the coordinate directions in a transport problem.
Carlier, Galichon and Santambrogio obtained the Knothe--Rosenblatt
rearrangement as the limit of optimal transports
for successively smaller coordinate weights~\cite{CarlierKnothe}.
Here we let the horizontal eigenvalues of
a test Hessian tend to zero in the flux bound for $S(\rho)$. The limit
satisfies a continuity equation for the conditional laws of a single
source and bounds their Wasserstein speed in the transverse variable.

\begin{proposition}\label[proposition]{prop:conditional-fibres}
Let $K\subset\R^{d-1}$ be bounded, open and convex, and $I$ a bounded
open interval. Let $\rho$ be a probability measure on $K\times I$
with density bounded above and below by positive constants. Disintegrate
$\rho(\diff x,\diff y)=\nu_y(\diff x)m(y)\dy$, with $m>0$ of class
$C^1(\overline I)$. If $S(\rho)<\infty$, the conditional laws admit a
locally absolutely continuous representative in $(\Ptwo(\overline K),W_2)$
with locally square-integrable metric speed. For every
$a\in C^1(\overline I)$ vanishing at both endpoints,
\begin{equation}\label{eq:conditional-fibre-bound}
 \int_I a^2\bigl(1+|\dot\nu_y|_{W_2}^2\bigr)m\dy
 \le S(\rho)^2\int_I\frac{|(ma)'|^2}{m}\dy .
\end{equation}
\end{proposition}
\begin{proof}
Choose $H_\epsilon=\operatorname{diag}(\epsilon I_{d-1},1)$ in
\cref{prop:stability-equivalences} to distinguish horizontal and
vertical derivatives. The datum $g=-(ma)'/m$ is centered because
$a$ vanishes at the endpoints. With this choice the vertical component of the
limiting flux is $a$. The flux estimate bounds
$X_\epsilon=H_\epsilon^{-1}\nabla u_\epsilon$ in $L^2(\rho)$ by
$S(\rho)\|g\|_\rho$.
Write $(Z,w)$ for a weak subsequential limit. Since the density is
bounded above and below, the identities
$\partial_{x_i}(X_\epsilon)_d=\epsilon\partial_y(X_\epsilon)_i$
give $\nabla_x w=0$ in distributions. Tests $h\in H^1(I)$ in the limiting
flux equation give $\int_I m(w-a)h'\dy=0$, hence $w=a$.

For $\zeta\in C_c^\infty(\R^{d-1})$ set $F(y)=\int\zeta\dd\nu_y$.
Testing with $\eta(y)\zeta(x)$, $\eta\in C_c^\infty(I)$, gives
\[
 -\int_I F(ma\eta)'\dy
 =\int_I m\eta\int\nabla\zeta\cdot Z\dd\nu_y\dy .
\]
On each compact subinterval where $a$ is nonzero, this is the
continuity equation with velocity $Z/a$ in $L^2(\nu_y(\diff x)\dy)$.
The moments of a countable dense family of smooth tests have
absolutely continuous representatives. Compactness of $\overline K$
then gives a narrowly continuous representative of the laws, and
\cite[Theorem~8.3.1]{AGS} gives
$|\dot\nu_y|_{W_2}^2\le\int|Z/a|^2\dd\nu_y$ almost everywhere there.
First choose $a$ positive throughout $I$ to obtain local absolute
continuity. For any admissible $a$, apply the speed bound on $\{a\ne0\}$.
Integration and weak lower semicontinuity of the flux norm give
\eqref{eq:conditional-fibre-bound}.
\end{proof}

\begin{proposition}\label[proposition]{prop:rough-counterexample}
Let $\sigma$ be the uniform probability measure on $Q=(-1,1)^2$.
For every $\varepsilon>0$ there is a probability measure $\rho$ on $Q$
with positive continuous density on $\overline Q$ such that
\[
 \|\frac{\diff\rho}{\diff\sigma}-1\|_{L^\infty(\sigma)}<\varepsilon,
 \qquad S(\rho)=\infty.
\]
Moreover, there exist an absolute constant $c>0$ and probability
measures $(\rho_n)_{n\ge2}$ on $Q$ with positive
$C^\infty(\overline Q)$ densities such that, for every $n\ge2$,
\[
 \|\frac{\diff\rho_n}{\diff\sigma}-1\|_{L^\infty(\sigma)}=n^{-1/2},
 \qquad S(\rho_n)\ge c\sqrt n. 
\]
\end{proposition}

\begin{proof}
For $0<\delta<1$ and an absolutely continuous $b:[-1,1]\to[-1,1]$,
set $\dd\rho=(1+\delta xb(y))\dd\sigma$. The $y$-marginal is $\dy/2$,
and the conditional law $\nu_y$ has barycenter
$\bar x(y)=\int_{-1}^1 x(1+\delta xb(y))\dx/2=\delta b(y)/3$.

Suppose $S(\rho)<\infty$. The barycenter map is $1$-Lipschitz in
$W_2$, so $|\bar x'|\le|\dot\nu_y|_{W_2}$ almost everywhere. Apply
\eqref{eq:conditional-fibre-bound} with $m=1/2$ and $a(y)=(1-y^2)/2$.
Since $-a'=y$ has squared norm $1/3$ in $L^2(m\dy)$, we obtain
\[
 S(\rho)^2\ge\frac32\int_{-1}^1 a^2(1+|\bar x'|^2)\dy
 =\frac25+\frac{\delta^2}{24}\int_{-1}^1(1-y^2)^2|b'(y)|^2\dy .
\]
To make the right side infinite while keeping the density continuous,
choose $b(y)=\operatorname{sgn}(y)\sqrt{|y|}$. Then
$|b'(y)|^2=1/(4|y|)$ for $y\ne0$, so $S(\rho)=\infty$ even for
$\delta<\varepsilon$. For $\delta\le1/2$ the density ratio lies in
$[1/2,3/2]$. Comparison
with the uniform square in the variational formulas for variance and
entropy and in the Dirichlet integral gives uniform Poincar\'e and
log-Sobolev constants, hence a uniform quadratic
transport--entropy bound~\cite{GigliLedoux}.

For the smooth examples, choose $b(y)=\sin(n\pi y)$ and
$\delta=n^{-1/2}$, $n\ge2$. This gives positive smooth densities
with $\norm{\diff\rho_n/\diff\sigma-1}_\infty=n^{-1/2}$, while the
derivative of the conditional mean grows in $L^2$. Indeed,
$\int_{-1/2}^{1/2}\cos^2(n\pi y)\dy=1/2$, so the same lower bound
gives $S(\rho_n)\ge c\sqrt n$ for an absolute $c>0$.
\end{proof}

\clearpage
\part{Wasserstein Proximal Maps}
\section{The pressure criterion for proximal maps}\label{sec:resolvents}
When does a Wasserstein proximal map decrease the distance between its
inputs? The squared distance in the objective contributes
$-\mathcal T/\tau$, and the Euler equation identifies the optimal map
from the minimizer to the input. When the input varies, testing the
differentiated pushforward relation against the potential of the
compatible field gives the velocity pairing used in
\eqref{eq:primal-contraction-deficit}. In this
identity a lower bound for the complete second variation of the energy
gives an upper bound for the output speed in terms of the input speed.

For internal energies the pressure condition of \cref{thm:weighted-brenier}
makes the complete second variation nonnegative. Uniform measures on concentric balls already suffice to prove that
the condition is necessary. At the critical power, translations
and dilations both give zero dissipation. We use this observation in
Section~\ref{sec:rigidity} to identify the pairs that preserve their
distance through two successive steps in $d\ge3$.

It matters here which variations enter the second derivative. Even a
displacement convex energy can have an expansive proximal map, as the
following example shows.

\begin{example}\label[example]{ex:resolvent-scope}
A step of a displacement convex energy can strictly increase the
distance between two inputs at which the minimizers are unique. Following the two-atom
construction of~\cite{AdveMeszaros}, put
$\mu_z=(\delta_z+\delta_{-z})/2$ in $\R^2$ and
$\rho_s=\mu_{(s,2)}$ for $|s|\le1$. Set $\mathcal F(\rho_s)=s^2/2$
and $\mathcal F=+\infty$ elsewhere. The two possible matchings give
\[
 W_2^2(\mu_z,\mu_w)=\min\{|z-w|^2,|z+w|^2\}.
\]
Between $\rho_s$ and $\rho_t$ the direct matching is uniquely
optimal, since the other costs an extra $16+4st>0$. Hence their
unique geodesic interpolates $s$ and $t$ linearly, and $\mathcal F$
is $1$-displacement convex. The energy is lower semicontinuous with
compact sublevels, so every input measure admits a minimizer.

For $\mu_\pm=\mu_{(1,\pm\varepsilon)}$ with $0<\varepsilon<1/2$, the
step-one objective at $\rho_s$ is $s^2-|s\pm2\varepsilon|$ up to a
constant. It is bounded below by $(|s|-1/2)^2-1/4-2\varepsilon$,
with equality at the minimum only at $s=\pm1/2$, respectively. Thus
\[
 W_2(\mu_+,\mu_-)=2\varepsilon,\qquad
 W_2(J_1\mu_+,J_1\mu_-)=W_2(\rho_{1/2},\rho_{-1/2})=1 .
\]
Both inputs have infinite energy, and their common limit at
$\varepsilon=0$ has two minimizers.
\end{example}

\subsection{The output velocity and the complete second variation}
\label{sec:adjoint-domains}\label{sec:volume}
Let $\eta_t$ be a smooth curve of input measures whose minimizers
$\rho_t=J_\tau\eta_t=r_t\dx$ form a smooth curve of positive densities.
At a fixed time write $q_\rho$ for the first variation of the energy
and assume that the optimal map $T$ from $\rho$ to $\eta$ is a smooth
diffeomorphism. Since $T$ is a convex gradient, its derivative is
positive definite. The Euler equation gives
\begin{equation}\label{eq:pressure-map}
 T=\Id+\tau\nabla q_\rho,\qquad
 H=\D T=I+\tau \D^2{q_\rho}\succ0.
\end{equation}
The minimal Wasserstein
velocity $w$ of the output determines its speed, whereas the field
$X=H^{-1}\nabla u$ satisfies the derivative estimate
$\tr((\D X)^2)\ge(\diver X)^2/d$ by \cref{lem:symmetry}. To use that
estimate for the motion of $\rho_t$, we choose $u$ so that $X$ and $w$
produce the same density derivative.

Write $h=\dot r/r$, which has mean zero by conservation of mass. The
first variation $q_\rho$ and its derivative $\cK_\rho h=\D q_\rho[rh]$
are defined up to spatial constants. We assume that $\cK_\rho$ is
symmetric on the centered variations in its form domain. To express
the continuity equation for $X$ and the derivative of the Euler map,
respectively, introduce
\begin{equation}\label{eq:resolvent-operators}
 B=-\diver_\rho(H^{-1}\nabla),\qquad
 \cR_\rho=I+\tau\cK_\rho B.
\end{equation}
We work on $\R^d$ or on a fixed smooth bounded domain $\Omega$, where
all functions are smooth up to the boundary and the fluxes have zero
normal component. On $\R^d$ we assume sufficient decay to integrate by
parts using cutoffs. In either setting $\int Bu\dd\rho=0$, so the
pairings with $h$ and $Bu$ are independent of the choice of constants.
We assume these pairings are finite, even when the velocity potential
itself lies outside $L^2(\rho)$. The regularity and conormal solvability
needed on a ball are proved in \cref{lem:compact-resolvent}.

For a vector field $X$ the complete form is
\begin{equation}\label{eq:complete-second-variation}
 \mathscr Q_\rho(X)
 =\inner{\diver_\rho X}{\cK_\rho(\diver_\rho X)}_\rho
       +\int \D^2{q_\rho}[X,X]\dd\rho.
\end{equation}
On $\R^d$, for compactly supported $X$, this is the affine second
variation:
\begin{equation}\label{eq:general-second-variation}
 \mathscr Q_\rho(X)=\left.\frac{\diff^2}{\diff s^2}\cE((\Id+sX)_\#\rho)\right|_{s=0}.
\end{equation}
Indeed, the first two density derivatives along $\Id+sX$ are
$-\diver(rX)$ and $\partial_{ij}(rX_iX_j)$. Integrating twice by parts
gives \eqref{eq:complete-second-variation}. When $X$ is a gradient and
$I+s\D X$ is positive definite, this affine curve is locally a
displacement geodesic. Here $\D X$ may be nonsymmetric. On a bounded domain even a tangent field can move points outside the
boundary under $\Id+sX$. Its flow stays in the domain, and
its acceleration accounts for the boundary term computed below.

\begin{example}[Interaction energy]\label[example]{ex:interaction-complete}
Let $W$ be smooth and even with bounded Hessian, and let
$\cE(\rho)=\tfrac12\iint W(x-y)\dd\rho(x)\dd\rho(y)$ on $\Ptwo$.
Its first variation is $W*\rho$ and $\cK_\rho h=W*(h\rho)$.
For $W(z)=|z|^2/2$ and a centered density variation $h$ of compact support,
$\inner{h}{\cK_\rho h}_\rho=-|\int xh\dd\rho|^2$.
To compute the complete form, differentiate the interaction along
$\Id+tX$. For smooth compactly supported $X$, \eqref{eq:general-second-variation}
gives, with $Z=X(x)-X(y)$,
\[
 \mathscr Q_\rho(X)=\tfrac12\iint\D^2 W(x-y)[Z,Z]\dd\rho(x)\dd\rho(y).
\]
Thus convexity of the kernel makes the complete form nonnegative,
even when the density variation alone contributes a negative term.
For the quadratic kernel the complete form is
$\|X\|_\rho^2-|\int X\dd\rho|^2$. Writing $b_\mu=\int x\dd\mu$,
one obtains the resolvent for every $\mu\in\Ptwo$,
\[
 J_\tau\mu=\left(b_\mu+\frac{\Id-b_\mu}{1+\tau}\right)_\#\mu.
\]
Indeed, minimizing over a transport plan with input $Y\sim\mu$ and
free output $Z$ amounts to minimizing
$\E|Z-Y|^2/(2\tau)+\E|Z-\E Z|^2/2$. Completing the square gives
$\E Z=b_\mu$ and $Z-b_\mu=(Y-b_\mu)/(1+\tau)$.
For two inputs with equal means their distance therefore contracts
by the exact factor $(1+\tau)^{-1}$, while translations preserve distance.
\end{example}

For a proximal step the coefficients in the moving-target identity
\eqref{eq:critical-branch-response} are $a=-1/\tau$ and $c=1/\tau$.
To turn this pairing into an estimate on the speed of the output
curve, we compare $X$ with the velocity field of smallest
$L^2(\rho)$ norm representing the same density derivative.

\begin{proposition}\label[proposition]{prop:primal-contraction}
Under the preceding smoothness and boundary assumptions, let $v$ be any
input velocity, $w=\nabla f$ the minimal output velocity, and
$h=\partial_t r/r$. Suppose $Bu=h$ has an admissible conormal
solution and put $X=H^{-1}\nabla u$. Then
\begin{equation}\label{eq:primal-pairing}
 \inner{X}{v\circ T}_\rho=\norm X_\rho^2+\tau\mathscr Q_\rho(X).
\end{equation}
For $a=1+\tau\lambda>0$ one has
\begin{equation}\label{eq:primal-contraction-deficit}
\begin{aligned}
 \norm v_\eta^2-a^2\norm w_\rho^2
 &=\norm{v\circ T-aX}_\rho^2+a^2\norm{X-w}_\rho^2\\
 &\quad+2a\tau\big(\mathscr Q_\rho(X)-\lambda\norm X_\rho^2\big).
\end{aligned}
\end{equation}
In particular a lower bound $\mathscr Q_\rho(X)\ge\lambda\norm X_\rho^2$
gives
\begin{equation}\label{eq:differential-contraction}
 \norm w_\rho\le (1+\tau\lambda)^{-1}\norm v_\eta .
\end{equation}
\end{proposition}

\begin{proof}
Fix the time at which the derivative is taken. Given an admissible
conormal test $z$ on the output, use $z\circ T^{-1}$ as a fixed test on
the input. Its gradient pulls back to $Y=H^{-1}\nabla z$. Thus
differentiating $T_\#\rho=\eta$ tests the input velocity against $Y$
and contributes $\dot T=\tau\nabla(\cK_\rho h)$ on the output side.
Integration by parts gives
\begin{equation}\label{eq:proximal-weak-pairing}
 \inner{Y}{v\circ T}_\rho
 =\inner{z}{h}_\rho+\tau\inner{Bz}{\cK_\rho h}_\rho .
\end{equation}
Taking $z=u$ with $Bu=h$ makes $Y=X$. Integration by parts then
gives $\inner{u}{h}_\rho=\int H[X,X]\dd\rho$ and inserting
$H=I+\tau\D^2q_\rho$ yields \eqref{eq:primal-pairing}. To recover the
ordinary output speed, recall that $X$ and $w$ represent the same
density derivative. Their difference has zero weighted divergence and
is therefore orthogonal to $w=\nabla f$. Substituting
$\norm X_\rho^2=\norm w_\rho^2+\norm{X-w}_\rho^2$ when expanding the
first square in \eqref{eq:primal-contraction-deficit} now proves the
identity by \eqref{eq:primal-pairing}.
\end{proof}

Otto and Westdickenberg~\cite{OttoWestdickenberg} prove contraction of
the porous medium flow by controlling the kinetic energy of an
infinitesimal perturbation as a curve of densities evolves.
Daneri and Savar\'e~\cite{DaneriSavare} use this Eulerian description
of the flow to establish displacement convexity of integral energies
on manifolds. For a single proximal step the Euler map supplies an
additional metric $H=I+\tau\D^2 q_\rho$. The term
$\|X-w\|_\rho^2$ in \eqref{eq:primal-contraction-deficit} is the
excess ordinary kinetic energy of the $H$-minimizer over the
Wasserstein velocity representing the same density change.

\begin{remark}[Varying the input and the energy]\label[remark]{rem:varying-energy}
Suppose the energy $\cE_t$ also varies smoothly on a fixed domain and
the hypotheses of \cref{prop:primal-contraction} hold along the branch.
Write $p_t=(\partial_t q_{t,\rho})|_{\rho=\rho_t}$ for the derivative
at fixed density. Differentiating the Euler map now gives
$\dot T=\tau\nabla(\cK_{\rho_t}h+p_t)$. The same test $Bu=h$ yields
\[
 \langle X,v\circ T-\tau\nabla p_t\rangle_\rho
       =\|X\|_\rho^2+\tau\mathscr Q_{t,\rho}(X).
\]
If $\mathscr Q_{t,\rho}(X)\ge\lambda\|X\|_\rho^2$ uniformly and
$1+\tau\lambda>0$, the output speed is at most
$(\|v\|_\eta+\tau\|\nabla p_t\|_\rho)/(1+\tau\lambda)$.
For instance, keep the input fixed and add the linear interpolation
of two smooth potentials $V_0,V_1$ to the same energy. Integration on
$[0,1]$ bounds the distance between the two outputs by
$\tau\|\nabla(V_1-V_0)\|_\infty/(1+\tau\lambda)$ under these branch
and curvature hypotheses. This differentiation requires a common
energy domain throughout the interpolation.
\end{remark}

On each regularized ball, $r$ is bounded above and below by positive
constants and $H$ is uniformly elliptic there. The Neumann form
$\int\inner{H^{-1}\nabla u}{\nabla z}\dd\rho$ is coercive on centered
potentials by the Poincar\'e inequality. Since $\int h\dd\rho=0$,
it yields the unique centered solution of $Bu=h$, smooth by elliptic
boundary regularity. On $\R^d$ the proposition assumes the indicated
solution and vanishing boundary terms. The global theorems below
apply it on balls and pass the distance bound to the limit.

A different choice of test determines how the integral of a prescribed
function $f$ changes along the output curve. Indeed, eliminating the
transport potential in \cref{thm:hessian-transfer} gives the density
Hessian $\cK_\rho+\tau^{-1}B^{-1}$. Its composition on the right with
$\tau B$ is $\cR_\rho$ modulo constants, so the equation for this test
is $\cR_\rho u=f$. The resulting estimate will also give conormal
invertibility and enter the proof of propagation of the Poincar\'e
constant.

\begin{proposition}\label[proposition]{prop:abstract-contraction}
Under the smoothness and boundary assumptions above, let $v$ be an
input velocity and let $u$ be an admissible smooth conormal solution
of $\cR_\rho u=f$. Set $Y=H^{-1}\nabla u$. Then
\begin{equation}\label{eq:resolvent-duality}
 \inner{f}{h}_\rho=\inner{Y}{v\circ T}_\rho ,
\end{equation}
and
\begin{equation}\label{eq:general-adjoint-energy}
 \inner{Y}{\nabla f}_\rho=\|Y\|_\rho^2+\tau\mathscr Q_\rho(Y).
\end{equation}
In particular, if $\mathscr Q_\rho(Y)\ge\lambda\|Y\|_\rho^2$ and
$1+\tau\lambda>0$, then $(1+\tau\lambda)\|Y\|_\rho\le\|\nabla f\|_\rho$.
\end{proposition}
\begin{proof}
Take $z=u$ in \eqref{eq:proximal-weak-pairing}. By symmetry of
$\cK_\rho$, the right-hand side is $\inner{\cR_\rho u}{h}_\rho
=\inner{f}{h}_\rho$, which proves \eqref{eq:resolvent-duality}.
Since $-\diver_\rho Y=Bu$ and $\nabla u=HY$, integration by parts gives
\[
 \inner{Y}{\nabla\cR_\rho u}_\rho
 =\int H[Y,Y]\dd\rho+\tau\inner{Bu}{\cK_\rho Bu}_\rho.
\]
Equation~\eqref{eq:pressure-map} identifies this expression with
the right-hand side of \eqref{eq:general-adjoint-energy}. The norm
bound follows by Cauchy--Schwarz.
\end{proof}

When $\cK_\rho$ is multiplication by a smooth positive function $a$
on the closed ball, the adjoint has a variational description.
For a field $Z$ with zero normal flux write $d_Z=-\diver_\rho Z$.
Instead of prescribing $d_Z$ as in \eqref{eq:quadratic-flux}, maximize
\[
 \int Z\cdot\nabla f\dd\rho-\tfrac12\int H[Z,Z]\dd\rho
                     -\tfrac\tau2\int a\,d_Z^2\dd\rho
\]
over $Z,d_Z\in L^2(\rho)$ with the indicated flux condition.
For $u+\tau aBu=f$ the field $Y=H^{-1}\nabla u$ satisfies
$HY=\nabla(f-\tau a d_Y)$. Integration by parts gives the first-order
condition, and the strictly negative quadratic part gives uniqueness.
In \eqref{eq:quadratic-flux} the divergence is prescribed. Here it is
part of the maximization, and the test function $f$ determines the
linear term.

For a local energy $\cE(r\dx)=\int_\Omega U(x,r)\dx$, write
\[
 s=r^{-1/d},\qquad \widehat U(x,s)=s^dU(x,s^{-d}),\qquad
 p(x,r)=rU_r(x,r)-U(x,r).
\]
The advantage of McCann's variable $s$ is that a map with Jacobian
determinant $j$ changes it to $sj^{1/d}$. Thus $\widehat U$, the energy
per unit mass, records both the displacement of a point and the change
in the volume it occupies~\cite{McCann}. Its derivative in $s$ involves
the pressure $p$. The formula below will give the contraction criterion
in \cref{thm:local-volume-criterion} and applies in particular to the
weighted powers used in Section~\ref{sec:projections}.

\begin{proposition}\label[proposition]{prop:volume-variation}
Let $\Omega=\R^d$ or let $\Omega$ be a bounded smooth domain. Let
$\rho=r\dx$ with $r>0$, and assume that $U$ and $r$ are smooth.
For a smooth vector field $X$, compactly supported on $\R^d$ or
tangent to $\partial\Omega$ on a bounded domain, assume that the
energy and the displayed variations are finite. Then
\begin{equation}\label{eq:volume-second}
 \begin{aligned}
 \mathscr Q_\rho(X)
 ={}&\int_\Omega \left.\frac{\diff^2}{\diff t^2}\right|_{t=0}
     \widehat U\!\left(x+tX,s+\frac{ts}{d}\diver X\right)\dd\rho\\
 &+\int_\Omega p\left(\tr((\D X)^2)-\frac{(\diver X)^2}{d}\right)\dx
   +\int_{\partial\Omega}p\,\II(X,X)\dd S ,
 \end{aligned}
\end{equation}
where the boundary integral is omitted on $\R^d$. If $\Omega$ is convex,
$X$ is compatible as in \eqref{eq:couple}, and
$\widehat U(x,s)-\lambda|x|^2/2$ is jointly convex and
nonincreasing in $s$, then $\mathscr Q_\rho(X)\ge\lambda\norm X_\rho^2$. 
\end{proposition}

\begin{proof}
Let $F_t$ be the flow of $X$. Since it preserves $\Omega$, a change
of variables expresses the energy of $(F_t)_\#\rho$ as
$\int_\Omega\widehat U(F_t,s_t)\dd\rho$, where
$s_t=s(\det\D F_t)^{1/d}$. To differentiate this expression twice,
write $b=\ddot F_0=(\D X)X$ and use the identity
$\partial_t\log\det\D F_t=(\diver X)\circ F_t$. It gives
\[
 \dot s_0=\frac{s}{d}\diver X,\qquad
 \ddot s_0=\frac{s}{d}
   \left(\diver b-\tr((\D X)^2)+\frac{(\diver X)^2}{d}\right).
\]
In the chain rule, the Hessian of $\widehat U$ gives the first
integral in \eqref{eq:volume-second}. Since
$r\nabla_x\widehat U=\nabla_x U$ and $\widehat U_s=-ds^{d-1}p$,
the two acceleration terms give
\[
 \int_\Omega p\left(\tr((\D X)^2)-\frac{(\diver X)^2}{d}\right)\dx
 +\int_\Omega(\nabla_x U\cdot b-p\diver b)\dx ,
\]
where $\nabla_x U$ is taken at fixed density.

It remains to remove the contribution of the acceleration $b$.
For a test function $f$, the second derivative of $f\circ F_t$ at zero
is $\D^2f[X,X]+\nabla f\cdot b$, so the second derivative of the
energy is $\mathscr Q_\rho(X)+\int_\Omega r\nabla q_\rho\cdot b\dx$.
After subtraction, the identity $\nabla p=r\nabla q_\rho-\nabla_x U$
turns the last integral above into
\[
 \int_\Omega\big((\nabla_x U-r\nabla q_\rho)\cdot b-p\diver b\big)\dx
 =-\int_\Omega\diver(pb)\dx
 =\int_{\partial\Omega}p\,\II(X,X)\dd S .
\]
Indeed, differentiating $X\cdot n=0$ along the boundary gives
$b\cdot n=-\II(X,X)$. On $\R^d$ the divergence integrates to zero.
This proves \eqref{eq:volume-second}.

Joint convexity of $\widehat U(x,s)-\lambda|x|^2/2$ gives a lower bound
of $\lambda\norm X_\rho^2$ for the first integral. Since
$\widehat U_s\le0$, we have $p\ge0$. For a compatible field the
second integral is nonnegative by \eqref{eq:trace-positive}, and
convexity of $\Omega$ gives the sign of the boundary term.
On $\R^d$, the affine variation
\eqref{eq:general-second-variation} also satisfies
\eqref{eq:volume-second} for $U\in C^2$, by local $C^2$ approximation
on the support of $X$.
\end{proof}

Integrating the velocity estimate \eqref{eq:differential-contraction}
along a smooth curve gives the corresponding bound on its action
\cite{BenamouBrenier,AGS}. On a ball the minimal velocity
$\nabla f_t$ is obtained by solving
$-\diver(r_t\nabla f_t)=\partial_t r_t$ with Neumann boundary condition.
By \cref{lem:smooth-paths} and continuity of the minimizer, the distance
estimate follows for all probability measures on the ball. We use
\cref{lem:resolvent-contraction-limit} to pass to $\R^d$.

For the converse estimates we will construct inputs for which the
minimizers are known. Given a finite-energy measure $\rho$, take $q$
supporting $\cE$ at $\rho$ and set
$\mu=(\Id+\tau\nabla q)_\#\rho\in\Ptwo$. If $-\tau q$ is a source
Kantorovich potential for $(\rho,\mu)$ and the dual pairings are
finite, then for every finite-energy competitor $\sigma\in\Ptwo$,
\begin{equation}\label{eq:resolvent-support}
 \cE(\sigma)-\cE(\rho)
 \ \ge\ \int q\dd(\sigma-\rho)
 \ \ge\ \frac{W_2^2(\rho,\mu)-W_2^2(\sigma,\mu)}{2\tau}.
\end{equation}
The first inequality is the support condition on $q$, and the second
is Kantorovich duality. Comparing the two ends shows that $\rho$
minimizes $\cE+W_2^2(\cdot,\mu)/(2\tau)$, uniquely if the first
inequality is strict for $\sigma\ne\rho$.

\subsection{Sufficiency and necessity of the pressure condition}
\label{sec:internal}\label{sec:relaxation}\label{app:internal}
Let $d$ and $U$ satisfy the hypotheses of \cref{thm:pressure-global}.
For the relaxed internal energy \eqref{eq:internal-relaxation}, the
volume transform in \cref{prop:volume-variation} depends only on $s$.
Differentiating
$\widehat U(s)=s^dU(s^{-d})$ gives
\[
 \widehat U'(s)=-ds^{d-1}p(r),\qquad
 \widehat U''(s)=d^2s^{d-2}\kappa_U(r),\qquad r=s^{-d}.
\]
Since $p\ge0$, the condition $\kappa_U\ge0$ is equivalent to
$\widehat U$ being convex and nonincreasing, as in McCann's criterion
\cite{McCann}.
For a compatible field $X=H^{-1}\nabla u$, let
$S_X^0=H^{1/2}\D XH^{-1/2}-(\diver X)I/d$.
On $\R^d$, for compactly supported $X$, \cref{prop:volume-variation} gives
\begin{equation}\label{eq:pressure-decomposition}
 \mathscr Q_\rho(X)
 =\int\left\{p(r)\norm{S_X^0}_{\HS}^2+
                      \kappa_U(r)(\diver X)^2\right\}\dx.
\end{equation}
By the coupled Bochner symmetry, $S_X^0$ is symmetric. The
coefficients $p$ and $\kappa_U$ thus weight the trace-free and scalar
parts of the same matrix. For entropy, $U(r)=r\log r$, their values
$r$ and $r/d$ recover the full squared norm
$r\|H^{1/2}\D XH^{-1/2}\|_{\HS}^2$.

To differentiate the Euler equation, assume for now that the integrand
and the positive minimizing density are smooth. Then
$q_\rho=U'(r)$, $H=I+\tau \D^2{q_\rho}$, and $\cR=I+\tau aB$ with
$a=rU''(r)$. For $g\in L^2(\rho)$ the operator $\cK_\rho g=ag$ is
defined when $ag\in L^2(\rho)$. Its quadratic form
$\inner{g}{\cK_\rho g}_\rho=\int a g^2\dd\rho$ is defined whenever
this integral is finite. At the critical power, for instance,
$a=(d-1)r^{-1/d}/d$ is unbounded on Barenblatt profiles.

On a ball the complete form on a compatible tangent field is
\begin{equation}\label{eq:internal-adjoint-energy}
 \mathscr Q_\rho(X)=
 \int_\Omega\{p(r)\norm{S_X^0}_{\HS}^2+\kappa_U(r)(\diver X)^2\}\dx
       +\int_{\partial\Omega}p(r)\II(X,X)\dd S .
\end{equation}
Thus the pressure and the convex boundary give
$\mathscr Q_\rho(X)\ge0$ when $\kappa_U\ge0$.

Suppose $\eta_t$ is a constant-speed geodesic with minimal velocity
$v_t$ and its positive minimizing branch satisfies the hypotheses
of \cref{prop:primal-contraction} on $\R^d$, with the terms in
\eqref{eq:primal-contraction-deficit} integrable in time.
When $\kappa_U\ge0$, endpoint equality forces equality in the action
bound for $\rho_t$ and in the integrated velocity estimate. Hence
$w_t=X_t=v_t\circ T_t$ and $\mathscr Q_{\rho_t}(X_t)=0$ almost
everywhere in time. If $\kappa_U>0$ on $(0,\infty)$, both terms in
\eqref{eq:pressure-decomposition} vanish and $\D X_t=0$. The spatially
constant velocity translates both curves by the same vector.
Conversely, uniqueness and translation invariance of the resolvent
give equality for every pair of translates.
\begin{remark}[Centered Gaussian pairs with small deficit]\label[remark]{rem:entropy-scale-escape}
Consider the entropy with quadratic confinement
$\cE(r\dx)=\int r\log r\dx+\frac\kappa2\int |x|^2r\dx$,
with $\kappa\ge0$, and put $a=1+\kappa\tau$.
For $s>0$ and $\mu_s=\N(0,s^2I)$ the output is $\mu_{R(s)}$, where
\[
 aR(s)-\frac{\tau}{R(s)}=s,\qquad
 R(s)=\frac{s+\sqrt{s^2+4a\tau}}{2a}.
\]
Indeed the Euler map of $\mu_{R(s)}$ is the positive homothety with
slope $s/R(s)$, so \eqref{eq:resolvent-support} applies.
Since $R'(s)\to 1/a$ at infinity, the deficit
$W_2^2(\mu_s,\mu_{s+1})-a^2 W_2^2(J_\tau\mu_s,J_\tau\mu_{s+1})$
tends to zero. Yet $W_2(\mu_s,\mu_{s+1})=\sqrt d$, even after
minimizing over translations of either input. Thus the deficit alone
admits no modulus controlling the distance modulo translations on
all of $\Ptwo$.
\end{remark}

\begin{proof}[Proof of \cref{thm:pressure-global}]
\emph{Sufficiency.}
Fix $\tau>0$ and assume $\kappa_U\ge0$. By \cref{lem:radial-wellposed}
the relaxed problem has a unique minimizer on $\R^d$ and on each closed
ball. To apply the differential estimate on a ball we first smooth $U$
while preserving the pressure condition. This is possible because the
rescaling $U(r)\mapsto U(ar)/a$ takes $\kappa_U(r)$ to
$\kappa_U(ar)/a$. Choose a nonnegative $\chi\in C_c^\infty((0,1))$
with integral one and set, for $n\ge1$,
\[
 U_n(r)=\int_0^1e^{s/n}U(e^{-s/n}r)\chi(s)\dd s,\qquad
 F_n(r)=U_n(r)+\frac1n r\log r,
\]
with $F_n(0)=0$. Averaging in $\log r$ makes $U_n$ smooth for $r>0$.
Since $U(r)/r$ is nondecreasing, $U_n\uparrow U$ and
$(U_n)_\infty=U_\infty$. Differentiating the rescaled integrands gives
\[
 \kappa_{F_n}(r)=
 \int_0^1e^{s/n}\kappa_U(e^{-s/n}r)\chi(s)\dd s
 +\frac{r}{nd}\ge0.
\]
The added entropy ensures that $F_n''(r)\ge1/(nr)$ and that $F_n'$
ranges from $-\infty$ to $+\infty$. This permits us to apply
\cref{lem:compact-resolvent} on each $B_R$ and differentiate the
minimizer along a smooth curve of positive input densities. By
\cref{prop:volume-variation}, the pressure condition and the convex
boundary give $\mathscr Q_\rho(X)\ge0$. Inserting this bound into
\eqref{eq:primal-contraction-deficit} with $\lambda=0$ shows that the
minimizer has no greater speed than the input. Integrating along the
curves from \cref{lem:smooth-paths} and using continuity of the
minimizer gives nonexpansivity on $\mathcal P(\overline B_R)$.

We first pass to the limit on $B_R$. Write $\mathcal F_n$ and
$\mathcal U_k$ for the relaxed energies with integrands $F_n$ and
$U_k$ on this ball. The functions $U_k$ have a common recession slope
and are bounded below by an affine minorant of $U_1$. Hence monotone
convergence gives $\mathcal U_k(\rho)\uparrow\mathcal U(\rho)$.
If $\rho_n\to\rho$ in $W_2$, the entropy bound on $B_R$ gives, for
each fixed $k$ and all $n\ge k$,
\[
 \mathcal F_n(\rho_n)\ge
 \mathcal U_k(\rho_n)-\frac{\log|B_R|}{n},\qquad
 \liminf_{n\to\infty}\mathcal F_n(\rho_n)\ge\mathcal U_k(\rho).
\]
The second inequality follows from lower semicontinuity. Letting
$k\to\infty$ gives the $\Gamma$-liminf inequality. For the recovery
sequence, use $U_n\le U$ and apply
\cref{lem:internal-recovery}\textup{(i)} with $\varepsilon_n=1/n$.

It remains to let $R\to\infty$. For the point mass $\delta_0$,
comparison with a fixed compactly supported competitor and
\cref{lem:radial-wellposed} give a uniform second moment bound for
the minimizers. The same lemma gives lower semicontinuity under this
bound, and \cref{lem:internal-recovery}\textup{(ii)} gives compactly
supported recovery sequences. Thus \cref{lem:resolvent-contraction-limit}
applies with $C=1$ and gives \eqref{eq:global-internal-contraction}
for all $\mu,\nu\in\Ptwo$. This proves
\textup{(i)}$\Rightarrow$\textup{(ii)}. The implication
\textup{(ii)}$\Rightarrow$\textup{(iii)} is immediate.

\emph{Necessity.}
Assume \textup{(iii)}. Write
$\beta_R=\one_{B_R}\dx/(\omega_dR^d)$, where $\omega_d=|B_1|$, and set
\begin{equation}\label{eq:radial-envelope}
 e_\tau(R)=\min_\rho\left\{
       \mathcal U(\rho)+\frac{W_2^2(\rho,\beta_R)}{2\tau}\right\}.
\end{equation}
The minimizer is radial by uniqueness in \cref{lem:radial-wellposed}
and rotation invariance. For radial measures $\mu,\nu$ with radius
quantiles $Q_\mu,Q_\nu$, the reverse triangle inequality bounds the
cost below by the cost between the radii. Coupling the radii monotonically
and using the same uniform direction attains the bound. Hence
\begin{equation}\label{eq:radial-hilbert}
 W_2^2(\mu,\nu)=\int_0^1|Q_\mu(v)-Q_\nu(v)|^2\dd v.
\end{equation}
In these coordinates the uniform balls trace the ray $Rv_0$, where
$v_0(v)=v^{1/d}$ and $\norm{v_0}_2^2=d/(d+2)$. Write $P_\tau(R)$ for
the quantile of $J_\tau^U\beta_R$. At each $\tau=\tau_j$, condition
\textup{(iii)} bounds the $L^2(0,1)$ Lipschitz constant of $P_\tau$
by $\norm{v_0}_2$. To see what this implies for the minimum value,
use the minimizer at $R$ as a competitor at $R+h$. Only the cost changes,
and expanding its square gives
\[
 \tau\big(e_\tau(R+h)-e_\tau(R)\big)
 \le h\inner{Rv_0-P_\tau(R)}{v_0}+\tfrac12h^2\norm{v_0}_2^2 .
\]
Using the minimizer at $R+h$ as a competitor at $R$ bounds the same
difference from below. Continuity of $P_\tau$ allows us to pass to the
limit in both difference quotients and obtain
\begin{equation}\label{eq:radial-envelope-derivative}
 e_\tau'(R)=\tau^{-1}\inner{Rv_0-P_\tau(R)}{v_0}.
\end{equation}
For $R>S$, the Lipschitz bound on $P_\tau$ and Cauchy--Schwarz now give
\[
 \tau(e_\tau'(R)-e_\tau'(S))
 =\frac{d}{d+2}(R-S)-\inner{P_\tau(R)-P_\tau(S)}{v_0}\ge0.
\]
Thus $e_{\tau_j}$ is convex, and by \cref{lem:radial-wellposed}
so is its pointwise limit
\[
 F(R)=\mathcal U(\beta_R)
     =\omega_dR^dU((\omega_dR^d)^{-1}),
\]
whose derivatives, with $r_R=(\omega_dR^d)^{-1}$, are
\begin{equation}\label{eq:radial-pressure-curvature}
 F'(R)=-d\omega_dR^{d-1}p(r_R),\qquad
 F''(R)=d^2\omega_dR^{d-2}\kappa_U(r_R).
\end{equation}
Since $r_R$ ranges over $(0,\infty)$, this proves
\textup{(iii)}$\Rightarrow$\textup{(i)}.

Finally, if $\kappa_U(r_0)<0$, choose three nearby radii at which
$F$ strictly violates convexity. The convergence
$e_\tau\to F$ gives the same violation for every sufficiently small
$\tau$. Nonexpansivity on that interval would make $e_\tau$ convex
by the preceding argument. Thus each such step expands the distance
between a suitable pair of uniform ball measures.
\end{proof}

For Boltzmann entropy, $p(r)=r$, so \cref{thm:pressure-global} gives
nonexpansiveness for every $\tau>0$. In the planar Gaussian example
of Di Marino, Farinelli and Naldi~\cite[Example~3.18]{DiMarino}, the
midpoint of a two-base interpolation is supported on a line although
both endpoints have finite entropy. The pressure criterion uses the
symmetry of $H\D X$ along variations of the actual minimizer. It
therefore applies to entropy even in the setting of their counterexample
to two-base convexity.

\subsection{A criterion for spatially dependent local energies}

When $U$ depends on position, a variation must test both the
$x$ and the $s=r^{-1/d}$ variables in \eqref{eq:volume-second}.
On a small ball of constant density we can vary them independently
by translating the center and dilating the ball to change the volume
per unit mass. Combining the two motions realizes each direction of the
joint Hessian and gives a converse to \cref{prop:volume-variation}.

\begin{theorem}\label[theorem]{thm:local-volume-criterion}
Let $d\ge2$ and $U\in C^\infty(\R^d\times[0,\infty))$ satisfy
\[
 U(x,0)=0,\qquad U_{rr}(x,r)>0\quad(r>0),\qquad
 U(x,r)\ge c r^\ell-Cr\quad(c>0,\ C\ge0,\ \ell>1).
\]
Extend $\cE_U(r\dx)=\int U(x,r)\dx$ by $+\infty$ on measures with a
singular part. Its proximal map $J_\tau^U$ is single-valued on $\Ptwo$
for every $\tau>0$. Fix $\lambda\ge0$. The following are equivalent.
\begin{enumerate}[label=\textup{(\roman*)},leftmargin=2em]
\item $\widehat U(x,s)-\lambda|x|^2/2$ is jointly convex on
$\R^d\times(0,\infty)$, where $\widehat U(x,s)=s^dU(x,s^{-d})$.
\item $\Lip_{W_2}(J_\tau^U)\le(1+\lambda\tau)^{-1}$ for every $\tau>0$.
\item There exists a sequence $\tau_j\downarrow0$ such that for every $j$
the estimate in \textup{(ii)} holds for all pairs of smooth compactly
supported probability densities.
\end{enumerate}
If the Hessian in \textup{(i)} has a negative direction at $(x_0,s_0)$,
every neighborhood of $x_0$ contains a pair of smooth minimizing
densities which, for each sufficiently small $\tau$, arise from smooth
inputs in the same neighborhood violating \textup{(ii)}.
\end{theorem}

\begin{proof}
Adding $Cr$ to $U$ changes the energy by a constant. After this shift
$U(x,r)\ge c r^\ell$, so convexity and $U(x,0)=0$ imply $U_r(x,0)\ge0$
and $p\ge0$. Comparison with a fixed bounded compactly supported density
bounds the second moments and the $L^\ell$ norms along a minimizing
sequence. Weak lower semicontinuity of the convex integral on each ball
and exhaustion give existence. Strict convexity of the energy under
mixtures gives uniqueness, both on balls and on $\R^d$.

Assume \textup{(i)}. On a ball add $\varepsilon r\log r$ to $U$.
The resulting density derivative is an increasing diffeomorphism onto
$\R$, and the volume transform acquires the convex term $-d\varepsilon\log s$.
Since $p\ge0$, the transform is nonincreasing in $s$.
By \cref{lem:compact-resolvent,prop:volume-variation,prop:primal-contraction},
the speed of the minimizer along a smooth positive input curve is at most
$(1+\lambda\tau)^{-1}$ times the input speed. Integration and
\cref{lem:smooth-paths} give this Lipschitz bound on the ball.

The entropy on $B_R$ is bounded below by $-\log|B_R|$, and every
finite-energy density there belongs to $L^\ell$ and has finite entropy.
It is therefore its own recovery sequence as $\varepsilon\downarrow0$.
For whole-space recovery, truncate a finite-energy density $r$ to
$\min(r,n)\one_{B_n}$ and restore the missing mass uniformly on a
fixed bounded set of positive volume where $r$ is bounded.
The removed second moment tends to zero by dominated convergence, and the restored mass
tends to zero, so these densities converge in $W_2$. Off the fixed set,
monotone convergence applies to the energy since $U(x,\cdot)$ is
nondecreasing. On the set, $U$ is bounded on the relevant compact range.
Finally the nonnegative energy and a common compactly supported competitor
bound the second moments of the ball minimizers at input $\delta_0$.
\Cref{lem:resolvent-contraction-limit} gives \textup{(ii)}, which
implies \textup{(iii)}.

For necessity, choose a negative direction $(v,s_0 b)$ for the Hessian
of $\widehat U-\lambda|x|^2/2$ at $(x_0,s_0)$ and put $r_*=s_0^{-d}$.
Take a smooth compactly supported gradient field $X$ agreeing with
$v+b(x-x_0)$ near $x_0$. On a small ball $B$ about $x_0$ the affine
map $F_t=\Id+tX$ has Jacobian determinant $(1+tb)^d$, so
\[
 \left.\frac{\diff^2}{\diff t^2}\right|_{t=0}
       \widehat U(F_t(x),s_0(1+tb))-\lambda|X(x)|^2<0
       \qquad(x\in\overline B).
\]
Choose a smooth density $r_0$ equal to $r_*$ on $B$, tapering to zero
in a thin shell on which $X$ is still affine. Shrink $B$ first so that
the mass is less than one, and place the balance in a disjoint ball where
$X=0$. The two balls and the cutoff for $X$ can all be taken in the
prescribed neighborhood. Write $\rho_t=(F_t)_\#(r_0\dx)=r_t\dx$,
$M=\int|X|^2r_0\dx$ and $A(t)=\cE_U(\rho_t)$. With $j_t=\det\D F_t$,
the integrand for $A(t)$ in the original coordinates is
$j_t U(F_t,r_0/j_t)$. Its second time derivative is bounded on the shell
independently of its thickness, since it involves no derivatives of
$r_0$ and $U$ is smooth up to vacuum. The second ball is stationary.
Thus $A''(0)<\lambda M$ for a sufficiently thin shell. For small
fixed $h>0$, $F_t$ is a convex gradient diffeomorphism on $[0,h]$
and $A'(h)-A'(0)<h\lambda M$.

To turn the strict inequality for $A$ into a test of the proximal
map, we need inputs whose minimizers are $\rho_0$ and $\rho_h$.
Choose a smooth compactly supported $0\le\chi\le1$ equal to one near
their supports and put
\[
 q_t(x)=U_r(x,r_t(x))-(1-\chi(x))U_r(x,0).
\]
This smooth compactly supported function supports the energy at $r_t$. It equals
$U_r(x,r_t)$ where $r_t>0$ and lies below $U_r(x,0)$ where $r_t=0$.
For $t=0,h$ and small $\tau$, $T_t=\Id+\tau\nabla q_t$ is a convex
gradient diffeomorphism and $-\tau q_t$ is a source cost potential.
Hence $J_\tau^U\mu_t^\tau=\rho_t$ for $\mu_t^\tau=(T_t)_\#\rho_t$ by
\eqref{eq:resolvent-support}. Both input densities are smooth and
supported in the prescribed neighborhood. Set $Z_t=\nabla q_t\circ F_t$.
Differentiating along $F_t$ gives $A'(t)=\inner{X}{Z_t}_{\rho_0}$.
The common-source coupling gives
\[
 W_2^2(\mu_0^\tau,\mu_h^\tau)
 \le h^2 M + 2h\tau\big(A'(h)-A'(0)\big)
                  +\tau^2\|Z_h-Z_0\|_{\rho_0}^2.
\]
The coefficient of $\tau$ in this upper bound is strictly smaller
than $2\lambda h^2 M$, by the choice of $h$. Since
$W_2^2(\rho_0,\rho_h)=h^2 M$, the bound is strictly less than
$(1+\lambda\tau)^2W_2^2(\rho_0,\rho_h)$ for every sufficiently small
$\tau$. This proves both the local assertion and the converse.
\end{proof}

\subsection{Dissipation towards a fixed equilibrium}\label{sec:wj-equilibrium}

The preceding criteria compare two arbitrary inputs. If one input
is fixed at equilibrium, can a weaker hypothesis give contraction?
To formulate this question, let $m=Z^{-1}e^{-V}\dx\in\Ptwo(\R^d)$
with $V$ smooth and consider its reversible Fokker--Planck equation
$\partial_t r_t=\Delta r_t+\diver(r_t\nabla V)$ for $\rho_t=r_t\dx$.
For a smooth positive $\rho$ with finite relative Fisher information,
let $T=\nabla\Phi$ be a smooth optimal diffeomorphism from $m$ to
$\rho$, and put $H=\D T$. With sufficient decay for integration by
parts, a smooth flow satisfies
$-\tfrac12\frac{\diff}{\diff t}W_2^2(\rho_t,m)=J(\rho_t\mid m)$.
In the notation of Bolley, Gentil and
Guillin~\cite[Corollary~2.4 and Definition~3.1]{BGGWJ}, this functional is
\[
 J(\rho\mid m)=\int\tr(H+H^{-1}-2I)\dd m
  +\int\big\langle\nabla V(T)-\nabla V,T-\Id\big\rangle\dd m .
\]
The diffusion term is precisely the sum of the Dirichlet energies of
the components of $T-\Id$ for the metric $H$, as in
\eqref{heat:eq:first-variation}. Near the identity map, the second-order
expansion of $J$ yields the integrated Bochner form used to deduce a
Poincar\'e inequality in~\cite[Proposition~3.9]{BGGWJ}.

Suppose $m$ satisfies $WJ(\lambda)$, that is,
$J(\rho\mid m)\ge\lambda W_2^2(\rho,m)$ with $\lambda>0$ for every
$\rho\in\Ptwo(\R^d)$ with positive $C^1$ density, with $J$ understood
as in their definition. In continuous time this gives exponential decay
towards $m$ at rate $\lambda$ in $W_2$. Since $J_\tau^m m=m$ for the
proximal map of $\Ent_m$, the corresponding finite-step question in
dimension $d\ge2$ is whether $WJ(\lambda)$ alone implies
\[
 W_2(J_\tau^m\mu,m)\le\frac{1}{1+\lambda\tau}W_2(\mu,m)
 \qquad\text{for every }\mu\in\Ptwo(\R^d),\ \tau>0 .
\]
When $V$ is $\lambda$-strongly convex this estimate follows from
convexity of $\Ent_m$ along generalized geodesics by taking one input
to be its minimizer~\cite[Lemma~3.19]{DiMarino}.

Under the preceding regularity and decay assumptions, the obstruction
in dimension $d\ge2$ can be seen directly from the two optimal maps out
of the proximal minimizer $\rho=J_\tau^m\mu$. The map $S=T^{-1}$ sends
$\rho$ to $m$, while the Euler map $U=\Id+\tau b$, with
$b=\nabla\log(\diff\rho/\diff m)$, sends it to $\mu$. Integration by
parts gives $J(\rho\mid m)=\int b\cdot(\Id-S)\dd\rho$, hence
\[
 \int (U-S)\cdot(\Id-S)\dd\rho
   =W_2^2(\rho,m)+\tau J(\rho\mid m).
\]
On the line, $U$ and $S$ are nondecreasing. Their common-source coupling
is optimal between $\mu$ and $m$, so Cauchy--Schwarz gives the stated
bound under these regularity assumptions. In higher dimensions
$U\circ S^{-1}$ can fail to be a gradient, even for linear maps, since
its derivative is a product of symmetric matrices which need not
commute. The coupling cost can then exceed $W_2^2(\mu,m)$, and the
same argument no longer gives the finite-step bound.

\section{Equality and rigidity for critical diffusion}
\label{sec:rigidity}\label{sec:global-critical-equality}
What can a proximal map preserve when the pressure condition
$\kappa_U\ge0$ holds? Suppose the input follows a constant-speed
$W_2$ geodesic with its minimal velocity and the corresponding
minimizers form a regular positive curve. If the input and output
endpoint distances agree, the integrated speed estimate forces
both terms in \eqref{eq:pressure-decomposition} to vanish for
almost every time. Since $p>0$, the first
gives $S_{X_t}^0=0$, or equivalently $\D X_t=b_t(x)I$. In dimension
$d\ge2$, commutation of mixed derivatives makes $b_t$ independent
of $x$, so the only possible velocities are $X_t(x)=a_t+b_tx$.
The second term permits $b_t\ne0$ only if $\kappa_U$ vanishes on the
density values attained. On any interval of such values the equation
$rp'=(1-1/d)p$ gives $U(r)=Ar-dCr^{1-1/d}$ for some $C>0$.

For the powers $U_m(r)=r^m/(m-1)$, $m\ne1$, we have $p(r)=r^m$ and
$\kappa_U(r)=(m-m_*)r^m$, where $m_*=1-1/d$. Above $m_*$ a dilation
therefore incurs a loss which can be estimated from the endpoint
distances. Indeed, let the input follow a constant-speed $W_2$ geodesic
with its minimal velocity and assume that the minimizing branch is
regular. Integrating \eqref{eq:primal-contraction-deficit} with
$\lambda=0$ bounds $2\tau(m-m_*)\int_0^1\!\int r_t^m(\diver X_t)^2\dx\dt$
by the difference of the squared endpoint distances, since the squared
output distance is at most the output action.

The same distinction occurs in the stability estimate for the density
$r=Z^{-1}W^{-\beta}$ in \cref{thm:inverse-power-stability}. Taking
$m=1-1/\beta$ gives $r^m=Z^{1/\beta}Wr$, so the terms involving $\D X$
in \eqref{eq:pressure-decomposition} and \eqref{eq:inverse-power-bochner}
agree up to this constant factor. At $\beta=d$, equivalently $m=m_*$,
only the trace-free square $\norm{S_X^0}_{\HS}^2$ remains in these
terms, and it vanishes on translations and dilations $X(x)=a+bx$.

Homogeneity explains why precisely this power allows equality under
dilations. If $\rho_s=(s\Id)_\#\rho$, then
$\mathcal U_m(\rho_s)=s^{d(1-m)}\mathcal U_m(\rho)$ for
$\mathcal U_m(r\dx)=(m-1)^{-1}\int r^m\dx$. At $m=m_*$ the right-hand
side is linear in $s$. Moreover, $q_{\rho_s}(sx)=s q_\rho(x)$, so
the gradient entering the Euler map does not change with the scale.
For distance to be preserved through a second step, the images under
the Euler maps must themselves differ by a translation and dilation. We derive that restriction in \eqref{eq:two-step-compatibility}.

The sign of $\kappa_{U_m}$ also determines how a single step changes
the distance between two Barenblatt profiles. Carlen and
Craig~\cite[Proposition~4.7]{CarlenCraig} obtained the scale recursion
above $m_*$. A supporting potential gives the recursion throughout the
range of finite second moment.

\begin{proposition}\label[proposition]{prop:power-dilations}
Let $d\ge3$, $m_2=d/(d+2)<m<1$, and $\tau>0$. Let
$h_m(x)=c_m(1+|x|^2)^{-1/(1-m)}$ have unit mass, and put
$B^{(m)}_s=(s\Id)_\#(h_m\dx)$ for $s>0$. For the relaxed energy
\eqref{eq:internal-relaxation} with $U_m(r)=r^m/(m-1)$,
\[
 J_\tau^{U_m}B^{(m)}_s=B^{(m)}_y,\qquad
 s=y-\tau k_m y^\alpha,\qquad
 k_m=\frac{2m c_m^{m-1}}{1-m},\qquad \alpha=d(1-m)-1 ,
\]
where $y>(\tau k_m)^{1/(1-\alpha)}$. For any two distinct scales
the step strictly contracts their distance when $m>m_*$, preserves
it when $m=m_*$ and strictly expands it when $m<m_*$.
In the last case its Lipschitz constant on this family is
$[2-d(1-m)]^{-1}$ for every $\tau>0$.
\end{proposition}

\begin{proof}
Write $r_y$ for the density of $B^{(m)}_y$. The condition
$1/(1-m)>(d+2)/2$ gives both $\int |x|^2h_m\dx<\infty$ and
$\int h_m^m\dx<\infty$. For $q_y=U_m'(r_y)$ we have
$\nabla q_y(x)=-k_m y^{\alpha-1}x$. Thus the Euler map is
$x\mapsto G(y)x/y$, where $G(y)=y-\tau k_m y^\alpha$.
On the indicated interval $0<\tau k_m y^{\alpha-1}<1$. As $\alpha<1$,
$G'(y)=1-\alpha\tau k_m y^{\alpha-1}>0$, and $G$ maps this
interval onto $(0,\infty)$.

For a given input scale $s$ we can thus choose the unique $y$ with
$G(y)=s$. The corresponding positive homothety is optimal from
$B^{(m)}_y$ to $B^{(m)}_s$, with source potential $-\tau q_y$ up to a
constant for the cost $|x-z|^2/2$. To verify minimality it remains to
check that $q_y$ supports the relaxed energy. The quadratic $q_y$ is
integrable against every competitor in $\Ptwo$. Strict convexity of
$U_m$ gives the support inequality for densities, and the inequality
is strict on any nonzero singular part because $q_y<0=(U_m)_\infty$.
It follows from \eqref{eq:resolvent-support} that $B^{(m)}_y$ is the
unique minimizer.

Since the homothety between any two profiles is optimal, their distance
is $W_2(B^{(m)}_s,B^{(m)}_t)=|s-t|(\int |x|^2 h_m\dx)^{1/2}$.
It suffices, therefore, to compare the scales before and after the step.
The derivative of the output scale is $1/G'(y)$, which is less than,
equal to, or greater than one according as $\alpha<0$, $\alpha=0$,
or $\alpha>0$. When $\alpha>0$ its supremum is $(1-\alpha)^{-1}$,
approached as $s\downarrow0$.
\end{proof}

The threshold $m_2$ is also a threshold for the proximal problem
itself. In $d\ge3$, if $0<m\le m_2$, the infimum is $-\infty$ for
every input and step size, even among smooth compactly supported
probability densities. To see this, start with a smooth positive
probability density proportional to
$|x|^{-(d+2)}(\log |x|)^{-b}$ for $|x|>2$, where $1<b<1/m_2$.
Its second moment is finite while the integral of its $m$th power
diverges. Smooth truncations equal to the density on $B_R$, followed
by mass normalization, have bounded second moments and hence bounded
transport costs to the input. Their energies tend to $-\infty$.

We now take $d\ge3$ and consider the relaxed critical energy
\begin{equation}\label{eq:critical-energy}
 \mathcal U_*(\rho)=-d\int r^{m_*}\dx,\qquad \rho=r\dx+\rho^s.
\end{equation}
Since $m_*>d/(d+2)$, \eqref{eq:internal-moment-bound} gives finite
energy on $\Ptwo$, and \cref{thm:pressure-global} gives its resolvent
$J_\tau^*$. To pass the equality conditions to arbitrary input measures,
we first establish local bounds on the critical minimizers.

\begin{lemma}
\label[lemma]{lem:critical-local-bounds}
Let $d\ge3$ and $\tau>0$. Every $J_\tau^*\eta$, $\eta\in\Ptwo(\R^d)$,
has a strictly positive locally Lipschitz density $r$. Its first variation
$q=-(d-1)r^{-1/d}$ is locally Lipschitz and satisfies
\begin{equation}\label{eq:critical-weak-Euler}
 q+\frac{|x|^2}{2\tau}\ \text{is convex},\qquad
 (\Id+\tau\nabla q)_\#(r\dx)=\eta.
\end{equation}
The map in \eqref{eq:critical-weak-Euler} is optimal and need not be
invertible.

More precisely, let $R_j\to\infty$, $\varepsilon_j\downarrow0$, and let
$\eta_j$ have a smooth strictly positive density on $\overline B_{R_j}$,
with $\sup_j\int|x|^2\dd\eta_j\le M$.
Let $r_j\dx$ minimize
\begin{equation}\label{eq:critical-regularized-objective}
 \int_{B_{R_j}}\bigl[-dr^{m_*}+\varepsilon_jr\log r\bigr]\dx
       +\frac{W_2^2(r\dx,\eta_j)}{2\tau}.
\end{equation}
For every $L<\infty$, there are $c_L,C_L>0$, depending only on
$d,\tau,M,L$, such that, for all sufficiently large $j$,
\begin{equation}\label{eq:critical-local-density-bounds}
 c_L\le r_j\le C_L,\qquad |\nabla r_j|\le C_L
       \quad\text{on }B_L.
\end{equation}
If $\eta_j\to\eta$ in $W_2$, then $r_j\dx$ converges narrowly to
$J_\tau^*\eta$, and $r_j\to r$ locally uniformly.
\end{lemma}
\begin{proof}
Put $q_j=h_j(r_j)$, with
$h_j(s)=-(d-1)s^{-1/d}+\varepsilon_j(1+\log s)$.
By \cref{lem:compact-resolvent}, $r_j$ is smooth and positive and
$\Id+\tau\nabla q_j$ is a convex gradient. Thus
\[
 q_j(y)\ge q_j(x)+\nabla q_j(x)\cdot(y-x)-\frac{|y-x|^2}{2\tau}
 \qquad(x,y\in B_{R_j}).
\]

\emph{Uniform local bounds.}
A fixed smooth compactly supported competitor bounds
\eqref{eq:critical-regularized-objective} from above. Write
$M_j=M_2(r_j\dx)$. By \eqref{eq:internal-moment-bound} the power
term is at least $-C_d(1+M_j)^{1/2}$. For $\varepsilon_j\le1$, the
Gaussian entropy bound and the transport cost give, for any $\delta>0$,
\[
 C\ge \frac{M_j/2-M}{2\tau}-C_d(1+M_j)^{1/2}
                          -\delta M_j-C_{d,\delta}.
\]
Taking $\delta<1/(4\tau)$ gives $\sup_j M_j<\infty$.

To obtain local bounds on $r_j=h_j^{-1}(q_j)$, first observe that
$q_j$ must stay uniformly below zero in the interior. Indeed, fix
$a>0$ and let $\ell=\sqrt{2\tau a}$. If $q_j(x)\ge-a$ and
$B_\ell(x)\subseteq B_{R_j}$, the preceding supporting inequality
gives $q_j(y)\ge-2a$ on the half-ball where
$\nabla q_j(x)\cdot(y-x)\ge0$. The mass of this half-ball satisfies
\[
 1\ge \frac{|B_\ell|}{2}h_j^{-1}(-2a)
 \longrightarrow \frac{|B_\ell|}{2}\left(\frac{d-1}{2a}\right)^d .
\]
The limit is proportional to $a^{-d/2}$. Choosing $a$ small excludes
such an $x$ for all large $j$, so $q_j\le-a$ on $B_{R_j-\ell}$.

For the bound in the other direction, the uniform second moment
bound puts at least half the mass in a fixed ball. We can therefore
choose $x_j$ in that ball with $r_j(x_j)\ge c>0$ and $q_j(x_j)\ge-C$.
For $x\in B_{L+1}$, the reflected point $2x_j-x$ stays in a fixed
ball, where we have already bounded $q_j$ from above. Applying
semiconvexity at the midpoint $x_j$ gives
\[
 q_j(x)\ge 2q_j(x_j)-q_j(2x_j-x)-|x-x_j|^2/\tau\ge-C_L .
\]
Now $-C_L\le q_j\le-a$ on $B_{L+1}$, so the supporting inequality
at $x\pm e_i$ bounds $|\nabla q_j(x)|$ on $B_L$. On compact subintervals of $(-\infty,0)$,
$h_j^{-1}$ converges uniformly to $s\mapsto((d-1)/(-s))^d$ and is
uniformly Lipschitz, since $h_j'(r)\ge(d-1)r^{-1-1/d}/d$.
This proves \eqref{eq:critical-local-density-bounds}.

\emph{The limiting Euler map.}
For $\eta_j\to\eta$ in $W_2$, the moment bound gives tightness and
$\liminf_j\varepsilon_j\int r_j\log r_j\dx\ge0$. Lower semicontinuity
and recovery in \cref{lem:radial-wellposed,lem:internal-recovery} identify
every narrow limit as $J_\tau^*\eta$. The local bounds and Arzel\`a--Ascoli
then give $r_j\to r$ locally uniformly, with $r$ positive and locally
Lipschitz. Thus $q_j+|\cdot|^2/(2\tau)$ converges locally uniformly to
$q+|\cdot|^2/(2\tau)$, $q=-(d-1)r^{-1/d}$, and the gradients converge
a.e.\ by convexity. For every bounded continuous $\varphi$,
\[
 \int\varphi\dd\eta
 =\lim_j\int\varphi(x+\tau\nabla q_j(x))r_j(x)\dx
 =\int\varphi(x+\tau\nabla q(x))r(x)\dx .
\]
Here dominated convergence applies on each ball and the moment bound
controls its complement. The limiting map is a convex gradient, hence optimal.
Approximation by smooth positive densities on expanding balls gives
\eqref{eq:critical-weak-Euler} for every $\eta\in\Ptwo$.
\end{proof}

The density bounds suffice to take a limit in the Euler equation,
but they do not bound the Hessians that symmetrize $\D X$. To recover
the equality condition without taking a limit of these Hessians, write
\begin{equation}\label{eq:tracefree-form}
 \mathfrak q(A)=\tr(A^2)-\frac{(\tr A)^2}{d}.
\end{equation}
This expression only involves the entries of $A$. It may have either
sign for a general matrix, but for the approximating compatible fields
\cref{lem:symmetry} gives
$\mathfrak q(\D X)=\norm{S_X^0}_{\HS}^2\ge0$. The following lemma
explains what the vanishing of these integrals implies for the limiting
velocity.
\begin{lemma} 
\label[lemma]{lem:tracefree-limit}
Let $d\ge2$, $I=(0,1)$, and $R_j\to\infty$. Suppose $X_j$ are smooth
vector fields on $I\times B_{R_j}$, with
$\mathfrak q(\D_xX_j)\ge0$. Assume that, for every compact set
$K\subseteq\R^d$,
\begin{equation}\label{eq:tracefree-limit-assumptions}
 X_j\longrightarrow w\quad\text{in }L^2(I\times K),\qquad
 \int_I\!\int_K\mathfrak q(\D_xX_j)\dx\dt\longrightarrow0.
\end{equation}
If $\partial_iw_k=\partial_kw_i$ in distributions, then
\[
 w_t(x)=a_t+b_tx
 \quad\text{for a.e. }(t,x),
 \qquad a\in L^2(I;\R^d),\quad b\in L^2(I).
\]
\end{lemma}
\begin{proof}
For smooth $X$ and $\chi\in C_c^\infty(\R^d)$, integrate the divergence
of $X\diver X-(\D X)X$ against $\chi$. A second integration by parts gives
\begin{equation}\label{eq:localized-tracefree-identity}
 \int\chi\mathfrak q(\D X)\dx
 =\frac{d-1}{d}\int\chi(\diver X)^2\dx
  +2\int(\diver X)\nabla\chi\cdot X\dx
  +\int \D^2{\chi}[X,X]\dx.
\end{equation}
Taking $\chi=\zeta^2$ and applying Young's inequality to the middle
term yields
\[
 \int\zeta^2|\diver X_j|^2\dx
 \le C_d\int\zeta^2\mathfrak q(\D X_j)\dx
     +C_{d,\zeta}\int_{\supp\zeta}|X_j|^2\dx .
\]
After integration in time, \eqref{eq:tracefree-limit-assumptions}
therefore bounds $\diver X_j$ in $L^2(I\times K)$ on each compact ball.
Its weak limit must be $\diver w$, by the local convergence of $X_j$.
Since $w$ is curl-free, the Euclidean div--curl identity applied to
mollified spatial cutoffs of $w$ gives $\D_xw\in L^2(I\times K)$.

We can now use \eqref{eq:localized-tracefree-identity} for both $X_j$
and $w$, with $\chi\ge0$. The mixed term converges because $\diver X_j$
converges weakly and $X_j$ strongly. Lower semicontinuity of the
divergence square gives the desired inequality, since for the symmetric
matrix $\D_xw$ the form on the left is the squared norm of its
trace-free part.
\[
 \int_I\!\int\chi\left\|\D_xw-\tfrac{\diver w}{d}I_d\right\|_{\HS}^2
         \dx\dt
 \le\liminf_j\int_I\!\int\chi\,\mathfrak q(\D X_j)\dx\dt=0 .
\]
Thus $\D_xw=\beta I_d$. For $i\ne k$, commutation of distributional
derivatives gives $\partial_k\beta=\partial_i\partial_k w_i=0$.
It follows that $w_t(x)=a_t+b_tx$. Integrating its squared norm on
$B_1$ shows that both coefficients belong to $L^2(I)$.
\end{proof}
\begin{theorem}\label{thm:critical-global-equality}
Let $d\ge3$, $\tau>0$, and $\mu_0,\mu_1\in\Ptwo(\R^d)$. Write
$\rho_i=J_\tau^*\mu_i$. Then
\begin{equation}\label{eq:critical-global-distance-equality}
 W_2(\rho_0,\rho_1)=W_2(\mu_0,\mu_1)
\end{equation}
if and only if
\begin{equation}\label{eq:critical-equality-outputs}
 \rho_1=(b+s\Id)_\#\rho_0
 \quad\text{for some }b\in\R^d,\ s>0.
\end{equation}
In this case $\rho_0=r_0\dx$ has the positive locally Lipschitz density
of \cref{lem:critical-local-bounds}. With $q_0=-(d-1)r_0^{-1/d}$,
\begin{equation}\label{eq:critical-equality-inputs}
 \mu_0=(\Id+\tau\nabla q_0)_\#\rho_0,\qquad
 \mu_1=(b+s\Id+\tau\nabla q_0)_\#\rho_0.
\end{equation}
\end{theorem}
\begin{proof}
Let $D$ be the common distance in \eqref{eq:critical-global-distance-equality},
$\eta_t$ a constant-speed geodesic between $\mu_0$ and $\mu_1$, and
$\rho_t=J_\tau^*\eta_t$. Nonexpansivity gives
$W_2(\rho_s,\rho_t)\le(t-s)D$ for $s<t$. Equality at the endpoints
forces equality on every subinterval: the three upper bounds on
$[0,s]$, $[s,t]$, and $[t,1]$ add up to $D$. Thus $\rho_t$ is a
constant-speed geodesic, including the constant curve when $D=0$.

\emph{The action deficit on expanding balls.}
To apply the velocity identity along this geodesic, approximate
$\eta_t$ uniformly in $W_2$ by the smooth curves $\eta_{j,t}$ on
expanding balls given by \cref{lem:smooth-paths}. Their actions
satisfy $A_j\to D^2$. For the minimizers $r_{j,t}\dx$ of
\eqref{eq:critical-regularized-objective}, let $w_j$ be the minimal
velocities and $X_j$ the compatible fields representing
$\partial_t r_j/r_j$ in \cref{prop:primal-contraction}. We compare
$A_j$ with the output actions $E_j=\int_0^1\!\int|w_j|^2r_j\dx\dt$.

For these minimizers the pressure is $r^{m_*}+\varepsilon_jr$, with
coefficient $\varepsilon_jr/d$ on the dilation term. Hence every term
of the ball formula \eqref{eq:internal-adjoint-energy}, including its
boundary term, is nonnegative. Keeping just the two terms needed to
identify the limiting field in \eqref{eq:primal-contraction-deficit}
with $\lambda=0$ gives
\begin{equation}\label{eq:critical-error-bound}
 \int_0^1\left\{\norm{w_j-X_j}_{r_j\dx}^2+
 2\tau\int_{B_{R_j}}r_j^{m_*}\mathfrak q(\D X_j)\dx\right\}\dt
 \le A_j-E_j .
\end{equation}
By \cref{lem:critical-local-bounds} the endpoint minimizers converge
narrowly. Lower semicontinuity of $W_2$ and the action bound give
\[
 D^2\le\liminf_j E_j\le\limsup_j E_j\le\lim_j A_j=D^2 .
\]
In particular, both terms on the left of \eqref{eq:critical-error-bound}
tend to zero.

\emph{Identification of the limiting velocity.}
Write $\rho_t=r_t\dx$. By \cref{lem:critical-local-bounds},
$r_{j,t}\to r_t$ locally uniformly in space for each $t$, with
upper and positive lower bounds on compact balls uniform in $j,t$.
Extend $\sqrt{r_j}\,w_j$ by zero outside $B_{R_j}$. Its squared
space-time $L^2$ norm is $E_j$, so along a subsequence
\[
 \sqrt{r_j}\,w_j\rightharpoonup\sqrt r\,w
 \quad\text{in }L^2((0,1)\times\R^d).
\]
On every compact ball, dominated convergence gives strong
space-time $L^2$ convergence of $\sqrt{r_j}$ to $\sqrt r$. This
identifies the weak limit of the fluxes $r_jw_j$ when paired with a
smooth compactly supported field. In particular we can pass to the
limit in the continuity equation for any
$\zeta\in C_c^\infty([0,1]\times\R^d)$. Tests that are nonzero at
the endpoints are allowed because the endpoint laws converge narrowly,
so the resulting identity reads
\[
 \int\zeta(1,\cdot)\dd\rho_1-\int\zeta(0,\cdot)\dd\rho_0
 =\int_0^1\!\int(\partial_t\zeta+\nabla\zeta\cdot w)\dd\rho_t\dt .
\]
In particular the limiting action bounds the squared distance between
these prescribed endpoints,
\[
 D^2\le\int_0^1\!\int |w|^2r\dx\dt\le\liminf_j E_j=D^2 .
\]
Equality identifies $w$ as the minimal velocity. Convergence of the
norms gives $\sqrt{r_j}\,w_j\to\sqrt r\,w$ strongly in space-time
$L^2$. On a fixed compact set the functions $1/\sqrt{r_j}$ are
uniformly bounded and converge pointwise to $1/\sqrt r$. Multiplying
the strong difference by $1/\sqrt{r_j}$ and applying dominated
convergence to $(\sqrt{r/r_j}-1)w$ gives $w_j\to w$ strongly in
$L^2((0,1)\times K)$ for each compact $K$. As each $w_j$ is a
gradient, $w$ is curl-free. By \eqref{eq:critical-error-bound}, $X_j\to w$ in the same
space and $\int_0^1\!\int_K\mathfrak q(\D X_j)\dx\dt\to0$.

The hypotheses of \cref{lem:tracefree-limit} are now satisfied, so
$w_t(x)=a_t+b_tx$ with $a,b\in L^2(0,1)$. Setting
$s_t=\exp(\int_0^t b_\ell\dd\ell)$, we can integrate this velocity
explicitly to give $F_t(x)=s_t(x+\int_0^t a_\ell/s_\ell\dd\ell)$.
For a compactly supported smooth $\varphi$, the test
$\varphi\circ F_t^{-1}$ is constant along the flow. It is admissible
in the continuity equation by smoothing the absolutely continuous
coefficients in time, since the supports remain in a common compact
set where the time derivatives converge in $L^1$ uniformly in space.
Consequently $(F_1)_\#\rho_0=\rho_1$, which is
\eqref{eq:critical-equality-outputs}.

\emph{The two input measures.}
Conversely, suppose $\rho_1=(b+s\Id)_\#\rho_0$. Critical scaling
gives $q_1(b+sx)=s q_0(x)$, so the Euler maps of
\cref{lem:critical-local-bounds} yield \eqref{eq:critical-equality-inputs}.
In particular, starting both maps at the same $x$ of law $\rho_0$
makes their $\tau\nabla q_0(x)$ terms cancel in the transport cost.
This coupling and nonexpansivity give
\[
 W_2^2(\rho_0,\rho_1)\le W_2^2(\mu_0,\mu_1)
 \le\int|b+(s-1)x|^2\dd\rho_0=W_2^2(\rho_0,\rho_1).
\]
The first inequality is nonexpansivity and the last equality is
optimality of $x\mapsto b+sx$.
\end{proof}

One step can therefore preserve the distance between nonradial
profiles. To construct examples, we may prescribe a smooth output
density and its translates and dilates. The only issue is whether the
corresponding Euler maps are optimal for all the proposed outputs.

\begin{corollary}\label[corollary]{thm:critical-isometries}
Let $d\ge3$ and let $\rho_0=r_0\dx$ be a smooth strictly positive,
centered probability measure such that
\begin{equation}\label{eq:critical-profile}
 M_2(\rho_0)<\infty,\qquad
 q_0:=-(d-1)r_0^{-1/d},\quad \norm{\D^2{q_0}}_\infty\le K<\infty.
\end{equation}
Fix $\tau>0$. For $b\in\R^d$ and $s>\tau K$, define
\begin{equation}\label{eq:critical-inputs}
 \rho_{b,s}=(b+s\Id)_\#\rho_0,\qquad
 \mu_{b,s}^\tau=(b+s\Id+\tau\nabla q_0)_\#\rho_0.
\end{equation}
The unique minimizer at $\mu_{b,s}^\tau$ is $\rho_{b,s}$. Moreover,
for $s,t>\tau K$ and $b,c\in\R^d$,
\begin{equation}\label{eq:critical-isometry}
 \begin{split}
 J_\tau^*\mu_{b,s}^\tau&=\rho_{b,s},\\
 W_2^2(\mu_{b,s}^\tau,\mu_{c,t}^\tau)
 &=|b-c|^2+M_2(\rho_0)(s-t)^2
 =W_2^2(\rho_{b,s},\rho_{c,t}).
 \end{split}
\end{equation}
\end{corollary}

\begin{proof}
Let $q_{b,s}=-(d-1)r_{b,s}^{-1/d}$, where $\rho_{b,s}=r_{b,s}\dx$.
Critical scaling gives
\[
 q_{b,s}(b+sx)=s q_0(x),\qquad
 T_{b,s}(b+sx)=b+sx+\tau\nabla q_0(x),
 \quad T_{b,s}=\Id+\tau\nabla q_{b,s}.
\]
Thus $(T_{b,s})_\#\rho_{b,s}=\mu_{b,s}^\tau$. Since
$\D T_{b,s}(b+sx)=I+(\tau/s)\D^2q_0(x)\succeq(1-\tau K/s)I$,
$T_{b,s}$ is optimal. The Hessian bounds and finite second moments
make the dual pairings in \eqref{eq:resolvent-support} finite.
For a competitor $\sigma$ with singular part $\sigma^s$, convexity
of $U_*$ on densities gives
\[
 \mathcal U_*(\sigma)-\mathcal U_*(\rho_{b,s})
 -\int q_{b,s}\dd(\sigma-\rho_{b,s})
 \ge -\int q_{b,s}\dd\sigma^s\ge0 .
\]
Hence \eqref{eq:resolvent-support} identifies $\rho_{b,s}$ as the
unique minimizer. The Euler map also gives
$W_2(\mu_{b,s}^\tau,\rho_{b,s})=\tau\norm{\nabla q_0}_{L^2(\rho_0)}$.

The minimizers $\rho_{b,s}$ and $\rho_{c,t}$ differ by a positive
homothety and a translation. By \cref{thm:critical-global-equality},
\[
 W_2^2(\mu_{b,s}^\tau,\mu_{c,t}^\tau)=W_2^2(\rho_{b,s},\rho_{c,t})
 =\int|b-c+(s-t)x|^2\dd\rho_0
 =|b-c|^2+M_2(\rho_0)(s-t)^2,
\]
where the last equality uses $\int x\dd\rho_0=0$.
\end{proof}

For an explicit nonradial example, take an even, nonradial
$f\in C_c^\infty(\R^d)$ and set
\begin{equation}\label{eq:nonradial-profiles}
	r_0(x)=Z^{-1}(1+|x|^2+\epsilon f(x))^{-d},
	\qquad 0<|\epsilon|\ll1.
\end{equation}
For $d\ge3$, this density has finite second moment and its first
variation has bounded Hessian. With $b=0$, the Euler maps are odd,
so two distinct admissible scales in \cref{thm:critical-isometries}
give distinct centered input measures whose distance is preserved
by one step and whose minimizers are nonradial. 

If distance is preserved a second time, the intermediate measures
are both outputs of the first step and inputs of the second. Their
two descriptions give optimal maps with the same source and target.
For a nontrivial dilation, equality of these maps will determine the
profile itself.

\begin{proof}[Proof of \cref{thm:two-step-rigidity}]
\emph{Necessity.}
Equality after two steps gives equality at each step by nonexpansivity.
Applying \cref{thm:critical-global-equality} twice yields
\[
 \rho_1=(b+s\Id)_\#\rho_0,\qquad
 \sigma_1=(c+t\Id)_\#\sigma_0,\qquad s,t>0 .
\]
When $s=1$, \eqref{eq:critical-equality-inputs} already gives
$\mu_1=(b+\Id)_\#\mu_0$. Assume $s\ne1$, write $\sigma_0=f\dx$ and
set $q=-(d-1)f^{-1/d}$. Applied to the second step, the input formula
transports $\sigma_0$ to $\rho_1$ by $c+t\Id+\tau_2\nabla q$.
Alternatively, the Euler map from $\sigma_0$ to $\rho_0$ followed by
the first homothety gives $b+s(\Id+\tau_2\nabla q)$. This composition
is still a convex gradient because $s>0$. Uniqueness identifies the
two maps $\sigma_0$-almost everywhere and $f>0$ then gives
\begin{equation}\label{eq:two-step-compatibility}
 \tau_2(1-s)\nabla q(x)=b-c+(s-t)x\qquad\text{a.e.}
\end{equation}
Thus the locally Lipschitz $q$ is a scalar quadratic plus an affine
function. Strict negativity on $\R^d$ excludes a positive quadratic
coefficient or a nonconstant affine function, and integrability of
$f$ excludes constant $q$. Completing the square gives
\[
 q(x)=-A-B|x-x_*|^2,\qquad
 f(x)=\left(\frac{d-1}{A+B|x-x_*|^2}\right)^d,\qquad A,B>0 .
\]
Thus $\sigma_0$ and its homothetic image $\sigma_1$ are Barenblatt
measures. Their Euler maps are scalar affine, so $\rho_0,\rho_1$ belong
to the same family, with positive scales by \cref{lem:critical-local-bounds}.
A second application of the Euler maps gives $\mu_0,\mu_1$, now allowing
scale zero.

\emph{Sufficiency.}
Translation invariance gives $J_\tau^*((b+\Id)_\#\mu)
=(b+\Id)_\#J_\tau^*\mu$, so translates preserve their distance at
every step. For the centered unit-scale Barenblatt profile $B$ in
\eqref{eq:intro-barenblatt-family},
\[
 \nabla q_B(x)=-a_dx,\qquad a_d=2(d-1)c_d^{-1/d}.
\]
For $R>0$ take $s=R+a_d\tau$ in \cref{thm:critical-isometries}.
The input map $b+s\Id+\tau\nabla q_B=b+R\Id$ gives
\begin{equation}\label{eq:barenblatt-step-v6}
 J_\tau^*B_{b,R}=B_{b,R+a_d\tau}\qquad(R\ge0),
\end{equation}
where the case $R=0$ follows by $W_2$ continuity.
After each common step the two outputs have positive scales and
are homothetic, so \cref{thm:critical-global-equality} gives equality
at every step.
\end{proof}

For a centered smooth profile we can also ask when its own dilates
occur as both inputs and outputs. Under the hypotheses of
\cref{thm:critical-isometries}, suppose that
$\mu_{0,s}^\tau=(t\Id)_\#\rho_0$ for one admissible $s$ and some
$t>0$. Then $s\Id+\tau\nabla q_0$ and $t\Id$ are convex gradients
from $\rho_0$ to this measure. Uniqueness gives
\[
 s x+\tau\nabla q_0(x)=t x,\qquad
 q_0(x)=q_0(0)+\frac{t-s}{2\tau}|x|^2 .
\]
Since $q_0<0$ and $r_0$ is integrable, $t<s$ and $q_0(0)<0$. Thus
\begin{equation}\label{eq:barenblatt-rigidity}
 r_0(x)=C(a+|x|^2)^{-d},\qquad a,C>0.
\end{equation}
Conversely, $q_0=-c(a+|x|^2)$ gives
$\mu_{0,s}^\tau=((s-2c\tau)\Id)_\#\rho_0$ when $s>2c\tau$.

\section{Wasserstein projections and constrained interactions}\label{sec:projections}
\label[section]{sec:projection-limit}
A projection minimizes the squared distance to the input subject to
a constraint. To apply the second-variation estimate to the projections
in \cref{thm:intro-projections}, we replace the constraint by energies that penalize violations increasingly
strongly. For a density bound in dimension $d\ge2$ we
use weighted powers. For an entropy sublevel we use increasing convex
functions of the relative entropy. In both cases the convexity condition
of the theorem makes the second variation of the penalty nonnegative
on compatible fields. Applying \cref{prop:primal-contraction} along
regularized curves and passing to the limit gives a nonexpansive
resolvent $J_n$ for each penalty.

It then suffices to prove that $J_n\mu$ converges weakly to the unique
projection $P_K\mu$ for each input. Indeed, lower semicontinuity gives
\[
 W_2(P_K\mu,P_K\nu)
 \le\liminf_nW_2(J_n\mu,J_n\nu)\le W_2(\mu,\nu).
\]
A second moment bound gives weak compactness in both problems.
Once the limiting measure is known to satisfy the constraint, comparison
with feasible competitors identifies it as the projection.

For the converse there is a geometric test that applies to both
constraints. Suppose every member of $K$ is absolutely continuous and
$P_K$ is nonexpansive. If $\gamma_t$ is the unique geodesic between
$\rho_0,\rho_1\in K$, then $P_K$ fixes its endpoints, so
\begin{equation}\label{eq:projection-geodesic}
 W_2(\rho_0,P_K\gamma_t)\le t W_2(\rho_0,\rho_1),\qquad
 W_2(P_K\gamma_t,\rho_1)\le(1-t) W_2(\rho_0,\rho_1).
\end{equation}
These two bounds add up to $W_2(\rho_0,\rho_1)$. The triangle
inequality forces equality in both, so $P_K\gamma_t$ lies at time $t$
on the unique geodesic and hence equals $\gamma_t$. This proves
displacement convexity of $K$. In particular, projecting an infeasible midpoint must increase its
distance to at least one endpoint.

\subsection{A variable density bound}
For a measurable density bound $f$, write
\[
 \cK_f=\{r\dx\in\Ptwo(\R^d):0\le r\le f\ \hbox{a.e.}\},
 \qquad P_f\mu=\argmin_{\rho\in\cK_f}\frac12W_2^2(\rho,\mu).
\]

To see which power of $f$ should enter the criterion, consider a
small mass $m$ at a point where the bound is approximately constant.
A uniform ball of this mass and density $f(x)$ has radius
$(m/|B_1|)^{1/d}f(x)^{-1/d}$. Between two such balls, optimal transport
interpolates the centers and radii linearly. Whether the midpoint still
fits under the bound is thus a convexity question for $f^{-1/d}$.

\begin{theorem}\label{thm:main-cap}
Let $d\ge1$ and let $f\ge0$ be locally integrable on $\R^d$ with
$\int f\dx>1$. Every $\mu\in\Ptwo(\R^d)$ has a unique projection
$P_f\mu$, and the following are equivalent.
\begin{enumerate}[label=\textup{(\roman*)},leftmargin=2em]
\item $P_f$ is nonexpansive on $(\Ptwo(\R^d),W_2)$.
\item $\cK_f$ is displacement convex.
\item $f^{-1/d}$ agrees a.e.\ with a lower semicontinuous
convex function $\R^d\to[0,+\infty]$.
\end{enumerate}
\end{theorem}

In particular, for any measurable $E\subseteq\R^d$ of Lebesgue measure
greater than one, $P_{\one_E}$ is nonexpansive on $\Ptwo(\R^d)$ if
and only if $E$ agrees up to a null set with a convex set. Indeed,
the effective domain of the representative in \textup{(iii)} then
agrees almost everywhere with $E$. Conversely, a convex set of
positive measure agrees almost everywhere with its closure, on which
we may take the value one and assign $+\infty$ elsewhere.

\phantomsection\label{sec:weighted-powers}
For a nonconstant example, take $\Omega=\R^d$ in
\cref{thm:inverse-power-stability} and set
$f=\lambda Z^{-1}W^{-\beta}$ with $\lambda>1$. Since $\beta\ge d$,
$f^{-1/d}$ is a positive multiple of the convex function $W^{\beta/d}$,
so projection onto this density bound is nonexpansive.

De Philippis, M\'esz\'aros, Santambrogio and Velichkov~\cite[Section~5.2]{DPMSV}
used powers of $r/f$ to approximate the constraint $r\le f$ in their
proof of BV estimates for the projection. For a smooth positive cap $f$,
consider the energy
$\mathcal E_{a,m}(r\dx)=(m-1)^{-1}\int a r^m\dx$ for $m>1$, with value
$+\infty$ on measures with a nonzero singular part. Choosing
$a_m=f^{1-m-1/d}$ gives
\[
 \mathcal E_{a_m,m}(r\dx)
   =\int \frac{f(x)^{-1/d}r(x)}{m-1}
       \left(\frac{r(x)}{f(x)}\right)^{m-1}\dx .
\]
With $f$ fixed, the integrand tends to zero where $r\le f$ and to
infinity where $r>f$ as $m\to\infty$. For this choice of $a_m$, the
following criterion for nonexpansiveness reduces to convexity of
$f^{-1/d}$ for every $m>1$.

\begin{proposition}\label[proposition]{cor:weighted-power-criterion}
Let $d\ge2$, $m>1$, and let $a\in C^\infty(\R^d)$ satisfy $a\ge a_*>0$.
Put $k=1+d(m-1)$ and $g=a^{1/k}$. For every $\tau>0$ and
$\mu\in\Ptwo(\R^d)$, the proximal problem for $\mathcal E_{a,m}$ has a
unique minimizer $J_\tau\mu$. The resolvents are nonexpansive for all
$\tau>0$ if and only if $g$ is convex.

If $\D^2g(x_0)$ has a negative direction, then for every neighborhood
$O$ of $x_0$ and every sufficiently small $\tau>0$ there are two
$C_c^\infty$ probability densities supported in $O$ whose distance
strictly increases under $J_\tau$. In particular, it suffices to test
nonexpansivity on smooth compactly supported densities along a sequence
$\tau_j\downarrow0$.
\end{proposition}

\begin{proof}
For $U(x,r)=a(x)r^m/(m-1)$, the volume transform is
$\widehat U(x,s)=g(x)^ks^{1-k}/(m-1)$. The function
$(z,s)\mapsto z^ks^{1-k}$ is jointly convex and increasing in $z>0$,
so convexity of $g$ gives joint convexity of $\widehat U$.
The proof of \cref{thm:local-volume-criterion} applies here even when
$m$ is not an integer. Indeed, the lower bound $U\ge a_*r^m/(m-1)$
gives the compactness needed for direct minimization for every such
weight, and $U+\varepsilon r\log r$ satisfies
\cref{lem:compact-resolvent}\,\textup{(i)} on a ball. Thus the speed
estimate only requires differentiation at positive densities and
yields nonexpansivity when $g$ is convex.

Conversely, choose $v$ with $\D^2g(x_0)[v,v]<0$. To isolate
this second derivative, dilate at the same relative rate at which $g$
changes along $v$. More precisely, set $b=\nabla g(x_0)\cdot v/g(x_0)$,
$X(x)=v+b(x-x_0)$ and $F_t=\Id+tX$. Then
\[
 \left.\frac{\diff^2}{\diff t^2}\right|_{t=0}
       \frac{g(F_t(x_0))^k}{(1+tb)^{k-1}}
       =k g(x_0)^{k-1}\D^2g(x_0)[v,v]<0 .
\]
The derivative stays negative in a small ball about $x_0$ contained in
$O$. Take a probability density $r_0$ proportional to
$\exp[-1/(\varepsilon^2-|x-x_0|^2)]$ on that ball and zero outside.
For $\rho_t=(F_t)_\#(r_0\dx)=r_t\dx$, a change of variables gives
\[
 A(t)=\mathcal E_{a,m}(\rho_t)
 =\frac1{m-1}\int r_0(x)^m
       \frac{g(F_t(x))^k}{(1+tb)^{k-1}}\dx .
\]
Thus $A''(0)<0$. Every positive power of $r_t$ is smooth up to and
across its support boundary, so $q_t=ma r_t^{m-1}/(m-1)$ is smooth,
compactly supported and supports the energy also at vacuum. The
endpoint coupling in the proof of \cref{thm:local-volume-criterion},
with $\lambda=0$, now gives the claimed expanding inputs in $O$.
\end{proof}

For example, take $m=2$ and $a(x)=1+|x|^2$. Then $k=d+1$ and,
for a radial unit vector $e$,
\[
 \D^2g(Re)[e,e]=\frac2k(1+R^2)^{1/k-2}
                  \left(1-\frac{k-2}{k}R^2\right).
\]
Thus the weight $a$ is convex, yet the resolvent expands local pairs
when $R^2>(d+1)/(d-1)$.

\begin{proof}[Proof of \cref{thm:main-cap}]
Choose $S$ with $\int_{B_S} f\dx>1$. The normalized restriction of
$f\dx$ to $B_S$ is a compactly supported feasible probability.
Existence and uniqueness follow from \cite[Proposition~5.2]{DPMSV},
and \textup{(i)} implies \textup{(ii)} by \eqref{eq:projection-geodesic}.
The strict mass assumption has a role in the converse: if
$\int f\dx=1$ and $f\dx\in\Ptwo$, the only feasible measure is
$f\dx$ and the projection is constant.

\emph{Sufficiency of \textup{(iii)}.}
We first show that a compactly supported input has a compactly
supported projection. Suppose $\supp\mu\subseteq\overline B_R$, write
$P_f\mu=s\dx$, and choose $S$ with $\int_{B_S}f\dx>1$. If $(x,y)$
belongs to the support of an optimal plan from $s\dx$ to $\mu$,
then $|x|>2R+S$ would imply
$B_S\subseteq B(y,|x-y|)$ since $|y|\le R$. By \cite[Lemma~5.1]{DPMSV}, the density must attain its upper bound
on this ball, giving $s=f$ a.e.\ on $B_S$,
contrary to $\int s\dx=1$. Thus $\supp(s\dx)\subseteq\overline B_{2R+S}$.

First suppose that $g=f^{-1/d}$ is smooth, convex and bounded below
by $c>0$. For $d\ge2$ and $m\ge2$, put
\[
 E_m(r)=\frac1{m-1}\int g r(r/f)^{m-1}\dx .
\]
By \cref{cor:weighted-power-criterion} with $a=g^{1+d(m-1)}$, its
unit-step resolvent $J_m$ is nonexpansive. For compactly supported
$\mu$, the support bound above and $s\le f$ give
$E_m(s)\le(m-1)^{-1}\int gs\dx\to0$.
Minimality of $r_m\dx=J_m\mu$ gives
\[
 E_m(r_m)+\tfrac12W_2^2(r_m\dx,\mu)
 \le \tfrac12W_2^2(P_f\mu,\mu)+o(1).
\]
This bounds the second moments. Moreover, for every $\varepsilon>0$,
\[
 \int_{\{r_m>(1+\varepsilon)f\}} r_m\dx
 \le \frac{(m-1)E_m(r_m)}{c(1+\varepsilon)^{m-1}}\longrightarrow0 .
\]
Every weak limit is therefore dominated by $(1+\varepsilon)f\dx$ for
all $\varepsilon>0$, hence belongs to $\cK_f$. Lower semicontinuity of
the cost and uniqueness identify the limit as $P_f\mu$. Passing to the
limit in the estimate for $J_m$ proves \textup{(i)} for compactly
supported inputs. Since $W_2$ convergence of the inputs implies weak
convergence of their projections \cite[Corollary~5.3]{DPMSV}, the
estimate extends by lower semicontinuity to all inputs.

In dimension one, quantiles identify $W_2$ with the $L^2(0,1)$
distance. For feasible densities $r_i$, let $x_i$ be their quantiles at
the same rank and set $x_t=(1-t)x_0+tx_1$. At almost every rank the
interpolated density satisfies \cite{McCann}
\[
 r_t(x_t)^{-1}=(1-t)r_0(x_0)^{-1}+t r_1(x_1)^{-1}
 \ge (1-t)g(x_0)+t g(x_1)\ge g(x_t).
\]
Thus the feasible quantiles form a convex set. It is closed since
domination by $f\dx$ passes to weak limits, so its Hilbert space
projection is nonexpansive.

For a general cap satisfying \textup{(iii)}, let $g$ be the lower
semicontinuous convex representative of $f^{-1/d}$. Set
$G_k=\max\{g,1/k\}$ and
\[
 g_{k,j}(x)=\inf_y\{G_k(y)+j|x-y|\}.
\]
For fixed $k$, $1/k\le g_{k,j}\uparrow G_k$. Convolution with a smooth
probability kernel supported in $B_{j^{-2}}$ preserves convexity and the
lower bound, and changes the $j$-Lipschitz function by at most $1/j$. By dominated
convergence the inverse $d$-th powers converge locally in $L^1$ to
$G_k^{-d}=\min\{f,k^d\}$. Since $f$ is locally integrable, a diagonal
sequence gives smooth caps $f_n\to f$ locally in $L^1$ and a.e., each
with $f_n^{-1/d}$ convex and bounded below by a positive constant.

For large $n$ the caps have $\int_{B_S}f_n\dx>1$. Their local
$L^1$ and a.e.\ convergence thus meet the hypotheses of
\cite[Lemma~5.4]{DPMSV}. Consequently $P_{f_n}\mu\rightharpoonup P_f\mu$ for every
$\mu$. Lower semicontinuity of $W_2$ gives \textup{(i)}.

\emph{The midpoint inequality under \textup{(ii)}.}
For necessity we recover the midpoint inequality for $g=f^{-1/d}$
from transport of small masses near two Lebesgue points. Set
$g=+\infty$ where $f=0$ and fix $z\in\R^d$. Since $f$ is locally
integrable with $\int f\dx>1$, we can choose $R>0$ and a bounded set
$E$ outside $B_{4R}(z)$ with $\int_E f\dx>1$. Its normalized
restriction $\lambda=f\one_E\dx/\int_E f\dx$ will hold the remaining
mass in place.

Take Lebesgue points $x_0,x_1\in B_R(z)$ with $0<f(x_i)<\infty$ and
write $g_i=g(x_i)$. To obtain the same mass at both endpoints while
respecting the local density bounds, put
\[
 h(u)=\min_{i=0,1} g_i^d f(x_i+g_i u),\qquad
 m_s=\int_{B_s(0)}h(u)\dd u .
\]
Both functions in the minimum have Lebesgue value one at the origin,
so $m_s/|B_s|\to1$. The images of $h(u)\one_{B_s}(u)\dd u$ under
$u\mapsto x_i+g_i u$ each have density at most $f$. For small $s$
they lie in $B_R(z)$ and $m_s<1$. Add $(1-m_s)\lambda$ to each
to obtain two feasible probabilities.

To identify the geodesic between these measures, we show that
optimal transport leaves their common mass on $E$ fixed. The same
argument will be needed for two moving portions of mass in the entropy
problem. Suppose each portion is supported at both endpoints in the
same ball of radius $R$, has equal mass at the two endpoints, and is
separated by at least $3R$ from the other portions and the common
background. For a given coupling, let $c_i$ be the mass leaving portion
$i$, which equals the mass entering it. Write $c$ for the total mass
transported between different portions or between a portion and the
background. These pairs cost at least $9R^2c$. We can instead keep all
pairs within each portion, couple its remaining marginals inside its
ball at cost at most $4R^2c_i$, and couple the whole background to
itself at zero cost. Since $\sum_i c_i\le c$, the replacement is
strictly cheaper whenever $c>0$.

In the present construction only the portion in $B_R(z)$ moves. Its
two endpoints are images of the same measure in $u$ under translations
and positive dilations, so matching equal $u$ is optimal there.

Put $x_*=(x_0+x_1)/2$ and $\bar g=(g_0+g_1)/2$. The midpoint puts
the mass $m_s$ in $B_{s\bar g}(x_*)$. Since it is feasible by
\textup{(ii)}, we have
\[
 \frac{m_s}{|B_s|}\le
 \bar g^d\,\frac{1}{|B_{s\bar g}|}\int_{B_{s\bar g}(x_*)}f\dx .
\]
At a Lebesgue point $x_*$, letting $s\downarrow0$ gives
$1\le\bar g^d f(x_*)$. By Fubini, almost every pair in $B_R(z)^2$
has such a midpoint. If either endpoint value of $g$ is infinite,
the extended-real midpoint inequality holds as well. We have proved
\[
 g\bigl((x_0+x_1)/2\bigr)\le\tfrac12\bigl(g(x_0)+g(x_1)\bigr)
 \qquad\text{for a.e. }(x_0,x_1)\in B_R(z)^2 .
\]

\emph{The convex representative.}
The preceding inequality holds for almost every pair, whereas a
convex representative must satisfy it everywhere on its domain.
We pass to the open set $U$ consisting of the union of balls on which
$A=\{(x,t):t>g(x)\}$ has full measure. Such balls are covered by balls
from a countable Euclidean basis, each still contained in $A$ up to a
null set. Thus $U\setminus A$ is null. To prove the reverse inclusion
up to a null set, fix one of the balls $B_R(z)$ above and two density
points $p,q$ of $A$ over this ball. For small $\delta>0$, $A$ occupies
more than three
quarters of each of $B_\delta(p)$ and $B_\delta(q)$. The convolution
of the indicators of these two portions is continuous and positive
at $p+q$. By the midpoint inequality it vanishes almost everywhere
on $\{2w:w\notin A\}$. Thus $(p+q)/2$ belongs to $U$.
Taking $p=q$ and using the density theorem gives $A=U$ a.e.
Every point of $U$ is a density point of $A$, so
$U\cap(B_R(z)\times\R)$ is open and midpoint convex, hence convex.

The set $U$ is upward closed. Its projection $D$ on $\R^d$ is open
and $D\cap B_R(z)$ is convex. By Fubini, $D=\{f>0\}$ up to a null set.
Consequently $S:=\supp(f\dx)=\overline D$. Within $B_{R/2}(z)$ this
closure agrees with $\overline{D\cap B_R(z)}$, which is convex.
It follows on the smaller ball that $S$ is locally convex and that
$D=\operatorname{int}S$ there, since an open convex set is the interior
of its closure. These smaller balls cover $\R^d$.

To pass from local to global convexity, it remains to see that
$S$ is connected. If it were the disjoint union of two nonempty
relatively open sets $S_0,S_1$, then along the geodesic between any
two feasible measures, \textup{(ii)} would put almost every transport
segment in $S$ at every rational time. Closedness extends this to all
times, so each segment is contained in one of the two pieces. All
feasible measures would then assign the same mass to $S_0$. Choose
a bounded ball $B$ meeting both
pieces in positive $f\dx$-mass, with $\int_B f\dx>1$. Multiplying
$f\one_B/\int_B f\dx$ by constants close to one on the two pieces
transfers mass between them while preserving feasibility, a contradiction.
Thus $S$ is connected.

The closed set $S$ is connected and locally convex. The
Tietze--Nakajima theorem~\cite[Theorem~1]{BjorndahlKarshon} therefore gives its convexity, and $D=\operatorname{int}S$
is convex as well.
On $D$ the lower endpoint of the vertical section of $U$ is finite and
locally convex, so it is convex. The closure of its epigraph defines a
lower semicontinuous convex function on $\R^d$. Since $U=A$ a.e.\
and $\partial S$ is null, this function agrees with $g$ almost everywhere.
\end{proof}

\subsection{Relative-entropy constraints}\label{sec:entropy-projections}
For a sublevel of relative entropy the relevant convexity condition
falls on the reference measure. Since the density $e^{-V}$ considered
here is positive almost everywhere on $\R^d$, log-concavity means that
$V$ agrees Lebesgue almost everywhere with a finite convex function.

\begin{theorem}\label{thm:entropy-projection}
Let $m=e^{-V}\dx\in\Ptwo(\R^d)$ have density positive Lebesgue
almost everywhere on $\R^d$. Fix $c>0$ and set
\[
 \cK_c=\{\rho\in\Ptwo(\R^d):\Ent_m(\rho)\le c\},\qquad
 P_c\mu=\argmin_{\rho\in\cK_c}\frac12W_2^2(\rho,\mu).
\]
The projection exists and is unique for every $\mu\in\Ptwo(\R^d)$.
The following conditions are equivalent:
\begin{enumerate}[label=\textup{(\roman*)},leftmargin=2em]
\item $P_c$ is nonexpansive in $W_2$;
\item $\cK_c$ is displacement convex;
\item $m$ is log-concave.
\end{enumerate}
When these conditions hold, the best global Lipschitz constant of $P_c$
is one. Otherwise there are compactly supported absolutely
continuous $\mu,\nu$ such that
$W_2(P_c\mu,P_c\nu)>W_2(\mu,\nu)$.
At $c=0$ the projection is the constant map with value $m$.
\end{theorem}

The necessity of log-concavity is a statement about each individual
level $c>0$. A local failure of convexity of $V$ must therefore produce
a geodesic that crosses that prescribed level, with both endpoints
below it.

\begin{lemma}\label[lemma]{lem:entropy-negative-curvature}
Let $m=e^{-V}\dx\in\Ptwo(\R^d)$ with density positive almost everywhere.
If $m$ is not log-concave, then for every $c>0$ there are compactly supported
absolutely continuous $\rho_-,\rho_+$ with
$\Ent_m(\rho_-)<c$ and $\Ent_m(\rho_+)<c$ whose unique Wasserstein midpoint
$\rho_0$ satisfies $\Ent_m(\rho_0)>c$.
\end{lemma}

\begin{proof}
Suppose $2V(x)\le V(x-h)+V(x+h)$ for a.e.\ $(x,h)$ with $|h|<r$
for some $r>0$. To obtain a convex representative, fix
$B=B_s(z)$ with $2s<r$ and choose $N$ so that
$E=B\cap\{V\le N\}$ has measure greater than $3|B|/4$.
The convolution $\mathbf 1_E*\mathbf 1_E$ is positive near $2z$,
so the midpoint inequality gives $V\le N$ almost everywhere near $z$.
Since $V^-\le e^{-V}$, it follows that $V\in L^1_{\rm loc}$.
Its mollifications satisfy the same midpoint inequality everywhere
for $|h|<r$ and hence are convex. Their local $L^1$ limit therefore
has a finite convex representative.

Consequently, if $m$ is not log-concave there are pairs violating
the midpoint inequality at arbitrarily small separations. Choose $R>0$
so small that any two balls of radius $4R$ have combined $m$-mass less
than $1-e^{-c}$. By Fubini's theorem, we can find such a pair
$x_-,x_+$ in some $B_{R/2}(z)$ for which both endpoints and their
midpoint $x_0$ are points of approximate continuity of $V$. Thus
\[
 v_0>\tfrac12(v_-+v_+),\qquad v_j=V(x_j)\quad(j=-,0,+).
\]
The two endpoint values of $V$ need not agree. To make the
midpoint entropy exceed both endpoint entropies, we balance this
difference by dilating a second portion of mass. Fix a point $y$ of
approximate continuity with $|y-z|>10R$. Place mass $p$ near $x_j$
and mass $1-p$ near $y$, with $0<p<1$. Dilating the latter by $b$
changes its Lebesgue entropy by $-(1-p)d\log b$, so to compensate
for the difference $p(v_+-v_-)$ we choose
\[
 b_-=e^{-t},\qquad b_+=e^t,\qquad
 t=\frac{p(v_+-v_-)}{2d(1-p)}.
\]
The midpoint dilation is $b_0=(b_-+b_+)/2=\cosh t$.
Fix $\varepsilon>0$. For all small $s>0$ there is a set
$F\subseteq B_s(0)$ of positive volume such that
\[
 |V(x_j+u)-v_j|<\varepsilon,\qquad
 |V(y+b_ju)-V(y)|<\varepsilon
 \qquad(u\in F,\ j=-,0,+).
\]
Indeed, each of the six conditions holds on a set of density one at
the origin. For $U$ uniform on $F$, let $\pi_j$ assign masses $p$
and $1-p$ to the laws of $x_j+U$ and $y+b_jU$, respectively. Their
supports lie in the disjoint balls $B_R(z)$ and $B_R(y)$. By change of
variables,
\[
 \Ent_m(\pi_j)=C+p v_j-(1-p)d\log b_j+e_j,\qquad |e_j|\le\varepsilon ,
\]
where $C$ is independent of $j$. The choice of $b_\pm$ thus gives
\[
 \Ent_m(\pi_0)-\max_\pm\Ent_m(\pi_\pm)
 \ge p\bigl(v_0-\tfrac12(v_-+v_+)\bigr)
      -(1-p)d\log\cosh t-2\varepsilon .
\]
The midpoint defect in parentheses is strictly positive and
$\log\cosh t=O(p^2)$ as $p\downarrow0$. Choose $p$ first and then
$\varepsilon$ so that the right-hand side is positive. Finally take
$s$ sufficiently small for $\Ent_m(\pi_0)>c$. Indeed, Jensen's
inequality bounds this entropy below by $-\log m(\supp\pi_0)$, and
the mass of the two balls carrying $\pi_0$ tends to zero with $s$.

We can now move the midpoint entropy through the prescribed level
$c$ while keeping it strictly above both endpoint entropies.
Choose a compact $A$ outside $B_{4R}(z)\cup B_{4R}(y)$ with $m(A)>e^{-c}$,
put $m_A=m|_A/m(A)$, and set $\rho_j^a=(1-a)m_A+a\pi_j$ for $a\in[0,1]$.
The two moving balls are more than $8R$ apart and each lies at least
$3R$ from $A$. By the cost comparison in the proof of \cref{thm:main-cap},
optimal transport leaves $m_A$ fixed and matches the mass in each ball
separately. The translation and dilation are optimal on their respective
balls, so $\rho_0^a$ is the midpoint of $\rho_-^a$ and $\rho_+^a$.

Finally, disjointness yields 
\[
 \Ent_m(\rho_j^a)
 =(1-a)\log\frac{1-a}{m(A)}+a\log a+a\Ent_m(\pi_j). 
\]
These entropies are continuous in $a$ and, for $a>0$,
the midpoint entropy strictly exceeds both endpoint entropies. It
starts at $-\log m(A)<c$ and ends above $c$. At
$a_*=\inf\{a:\Ent_m(\rho_0^a)>c\}$ it equals $c$, while both
endpoint entropies are below $c$. They remain below $c$ for some
$a>a_*$ with $\Ent_m(\rho_0^a)>c$, as required.
\end{proof}

\begin{proof}[Proof of \cref{thm:entropy-projection}]
\emph{Existence and uniqueness.}
The reference measure $m$ is feasible. Comparison with $m$ therefore
bounds the second moments of a minimizing sequence. Lower semicontinuity
of relative entropy and of the transport cost gives a minimizer.
If $c=0$, then $\cK_0=\{m\}$.

Suppose two distinct projections lie at distance $D$ from $\mu$.
Then $D>0$, and their average $\bar\rho$ satisfies
\[
 \Ent_m(\bar\rho)<c,\qquad W_2^2(\bar\rho,\mu)\le D^2
\]
by strict convexity of entropy and convexity of the cost under mixtures.
The strict entropy inequality leaves room to perturb $\bar\rho$
towards $\mu$. To do so with finite entropy, approximate $\mu$ by a
finitely supported measure and replace each atom by the normalized
restriction of $m$ to a small ball about it. Each ball has positive
$m$-mass, so the resulting mixture $\eta$ has finite entropy by
convexity and can be chosen with $W_2(\eta,\mu)<D$. For small
$\delta>0$, the mixture $(1-\delta)\bar\rho+\delta\eta$ remains
feasible, but its squared distance to $\mu$ is at most
$(1-\delta)D^2+\delta W_2^2(\eta,\mu)<D^2$. This contradicts
minimality and proves uniqueness.

\emph{Nonexpansivity for a smooth log-concave reference.}
Suppose $\lambda I\preceq \D^2V\preceq\Lambda I$ with
$0\le\lambda\le\Lambda<\infty$ and $V$ smooth.
Let $\Theta$ be smooth, convex and increasing on a neighborhood
of $[0,\infty)$, with $\Theta'\ge a_0>0$.
Write $J_{\tau,\Theta}$ for the resolvent of $\Theta\circ\Ent_m$.
At a smooth positive density $\rho=r\dx$ on a ball, write
$q_{\Ent_m}=1+\log r+V$ and evaluate the derivatives of $\Theta$ at
$\Ent_m(\rho)$. The first variation is $q_\rho=\Theta' q_{\Ent_m}$, and
\begin{equation}\label{eq:entropy-rank-one}
 \cK_\rho h=\Theta' h+\Theta'' q_{\Ent_m}\inner{q_{\Ent_m}}{h}_\rho.
\end{equation}
To estimate the output speed we need the complete second variation
after composition with $\Theta$. By the chain rule in
\eqref{eq:complete-second-variation}, the complete form of a composite
energy is the form of $\cE$ multiplied by $\Theta'(\cE(\rho))$, plus
$\Theta''(\cE(\rho))(\D\cE_\rho[X])^2$. Here
$H=I+\tau\D^2q_\rho\succ0$, and for $X=H^{-1}\nabla u$ with zero
normal flux we have $-\diver_\rho X=Bu$. The derivative of entropy
in the extra square is therefore $\inner{q_{\Ent_m}}{Bu}_\rho$.
Since \cref{prop:volume-variation} gives
$\mathscr Q_{\Ent_m,\rho}(X)\ge\lambda\norm X_\rho^2$, including its
nonnegative boundary term, we obtain
\begin{equation}\label{eq:entropy-composition-energy}
 \mathscr Q_\rho(X)
 =\Theta'\mathscr Q_{\Ent_m,\rho}(X)+\Theta''\inner{q_{\Ent_m}}{Bu}_\rho^2
 \ge a_0\lambda\norm X_\rho^2.
\end{equation}
Take $Bu=\partial_t r/r$ along the ball minimizers supplied by
\cref{lem:compact-resolvent}. By \cref{prop:primal-contraction,lem:smooth-paths},
their resolvents have contraction coefficient $(1+\tau a_0\lambda)^{-1}$.

For $\Theta\circ\Ent_m$, the compact recovery formula
\eqref{eq:entropy-restriction-recovery} allows us to apply
\cref{lem:resolvent-contraction-limit} and conclude that
\begin{equation}\label{eq:entropy-composition-contraction}
 W_2(J_{\tau,\Theta}\mu,J_{\tau,\Theta}\nu)
 \le(1+\tau a_0\lambda)^{-1}W_2(\mu,\nu)
 \qquad(\mu,\nu\in\Ptwo).
\end{equation}

To impose the entropy constraint, choose
\[
 \Theta_n(s)=\frac{s+e^{n(s-c)}}{n},
 \qquad \rho_n^\mu=J_{1,\Theta_n}\mu .
\]
This convex penalty has slope at least $1/n$ and is bounded by
$(c+1)/n$ on $[0,c]$. Comparing with $P_c\mu$ gives
\[
 \Theta_n(\Ent_m(\rho_n^\mu))+\tfrac12 W_2^2(\rho_n^\mu,\mu)
 \le (c+1)/n+\tfrac12 W_2^2(P_c\mu,\mu).
\]
The cost bound gives tightness. For every $\varepsilon>0$ the penalty
at $c+\varepsilon$ tends to infinity, so
$\limsup_n\Ent_m(\rho_n^\mu)\le c$. Lower semicontinuity of entropy
and cost now shows that every narrow limit is feasible and minimizes
the cost. Thus $\rho_n^\mu\rightharpoonup P_c\mu$ by uniqueness.
Passing to the lower limit in \eqref{eq:entropy-composition-contraction}
proves nonexpansivity.

\emph{Approximation and the converse.}
For a general log-concave $m$, Gaussian convolution gives reference
measures to which the preceding estimate applies. Indeed, put
$\gamma_t=\N(0,tI)$ and $m_t=m*\gamma_t=e^{-V_t}\dx$. Convolution
preserves log-concavity, and differentiating the kernel gives
$\D^2V_t\preceq t^{-1}I$ as in the proof of \cref{mm:lem:smoothing}.
We may thus take $\lambda=0$ in the smooth case. For a fixed input
$\mu$, let $\rho_t$ be its projection at level $c$ for $m_t$. Convolving
$\sigma=P_c\mu$ with the same Gaussian gives a feasible competitor for
$\rho_t$, by monotonicity of relative entropy under convolution. Hence
\[
 \Ent_{m_t}(\sigma*\gamma_t)\le\Ent_m(\sigma)\le c,\qquad
 W_2(\mu,\rho_t)\le W_2(\mu,\sigma*\gamma_t)\longrightarrow W_2(\mu,\sigma)
 \quad(t\downarrow0).
\]
The second moments of $\rho_t$ are bounded. Since $m_t\to m$ in $W_2$,
joint lower semicontinuity of entropy makes every narrow limit feasible.
The cost bound then identifies it with $\sigma$ by uniqueness. Passing
to the lower limit in the smooth-reference estimate gives
\textup{(iii)}$\Rightarrow$\textup{(i)}.

Finite relative entropy makes every feasible measure absolutely
continuous, so \eqref{eq:projection-geodesic} proves
\textup{(i)}$\Rightarrow$\textup{(ii)}.
\Cref{lem:entropy-negative-curvature} proves
\textup{(ii)}$\Rightarrow$\textup{(iii)}.
Its infeasible midpoint and the same triangle argument give strict
expansion of a pair of compactly supported absolutely continuous
measures whenever $m$ is not log-concave.

Finally, take a bounded nonzero $b$ with $\int b\dd m=0$.
For small $\epsilon$, $(1+\epsilon b)m$ lies in $\Ptwo$ and has
relative entropy tending to zero. Thus $\cK_c$ contains two
distinct measures for every $c>0$. The projection fixes them, so
the Lipschitz constant under the equivalent conditions is exactly one.
\end{proof}

\subsection{A Gaussian interaction under a density bound}\label{sec:capped-gaussian}

Craig~\cite{CraigNonconvex} studied gradient flows of nonlocal interaction
energies under a density constraint, including kernels with singular
attraction. Her estimates control convergence of the time-discrete
scheme to the flow. Craig and Topaloglu~\cite[Theorems~1.3 and~1.10]{CraigTopaloglu}
obtained the unit density cap from a power diffusion penalty as its
exponent tends to infinity, proving subsequential convergence of energy
minimizers up to translations and of gradient flows for well-prepared
initial data under their respective hypotheses. Here we ask how far
two outputs of a single step can move apart when the density bound
varies in space. For a Gaussian interaction kernel, its positive
Fourier transform controls the part of the second variation involving
the change of density. The other part involves the spatial Hessian of
the convolution, whose least eigenvalue can be bounded using the cap.
Perturbations supported near a point that realizes this lower bound
will determine when the resulting estimate is sharp.

Let $\Omega=B_R(0)\subset\R^d$, $d\ge1$, and suppose that $g$ is smooth,
positive and convex on a neighborhood of $\overline\Omega$. Set $f=g^{-d}$
on $\Omega$, with $\int_\Omega f\dx>1$, and write $\cK_f$ for the
probabilities on $\overline\Omega$ dominated by $f\dx$. For $\ell,\gamma>0$
let $W_\ell(z)=\exp(-|z|^2/(2\ell^2))$. If $V$ is smooth near
$\overline\Omega$, denote by $J_\tau^{f,V}$ the proximal map on
$\mathcal P(\overline\Omega)$ of
\[
 \cE_{f,V}(\rho)=\int V\dd\rho+
 \frac\gamma2\iint W_\ell(x-y)\dd\rho(x)\dd\rho(y)
 \qquad(\rho\in\cK_f),
\]
with value $+\infty$ elsewhere. Write $J_\tau^f=J_\tau^{f,0}$. The
size and curvature of the convolution on the feasible set enter through
\begin{equation}\label{eq:gaussian-capacities}
 \alpha_{f,\ell}=\max_{x,\rho}(W_\ell*\rho)(x),\qquad
 \Lambda_f=\min_{x,e,\rho}\D^2(W_\ell*\rho)(x)[e,e],
\end{equation}
where $x\in\overline\Omega$, $|e|=1$ and $\rho\in\cK_f$.

\begin{theorem}\label[theorem]{thm:capped-gaussian-bound}
The constants in \eqref{eq:gaussian-capacities} satisfy
$0<\alpha_{f,\ell}<1$ and $-\alpha_{f,\ell}/\ell^2<\Lambda_f<0$.
For every $\tau>0$ and every input the minimizer is unique. If
$\D^2 V\succeq\kappa I$ and $1+\tau(\kappa+\gamma\Lambda_f)>0$, then
\begin{equation}\label{eq:capped-gaussian-bound}
 \Lip_{W_2}(J_\tau^{f,V})\le
 \frac{1}{1+\tau(\kappa+\gamma\Lambda_f)}.
\end{equation}
\end{theorem}

In particular, $-\gamma\Lambda_f<\kappa<\gamma/\ell^2$ gives strict
contraction for every step size. The pointwise kernel bound
$\D^2W_\ell\succeq-\ell^{-2}I$ alone gives no contraction in this interval.

\begin{proof}
Weak compactness of $\cK_f$ and continuity of the kernels give attainment
in \eqref{eq:gaussian-capacities}. A feasible density gives zero mass to
each singleton, so $W_\ell\le1$ with equality only at the origin gives
$\alpha_{f,\ell}<1$. Moreover
\[
 \D^2 W_\ell(z)[e,e]=
 W_\ell(z)\big((z\cdot e)^2/\ell^4-1/\ell^2\big).
\]
The integral of $W_\ell(x-y)((x-y)\cdot e)^2$ against any feasible
measure is positive since a hyperplane has measure zero.
At a minimizing triple this gives $\Lambda_f>-\alpha_{f,\ell}/\ell^2$.
Take $\rho_*=f\dx/\int_\Omega f\dx$. On the sphere the outer normal
derivative of $W_\ell*\rho_*$ is strictly negative, since
$(x-y)\cdot n(x)>0$ for $y\in\Omega$. The divergence theorem gives
$\int_\Omega\Delta(W_\ell*\rho_*)\dx<0$, hence $\Lambda_f<0$.

The Fourier transform of $W_\ell$ is strictly positive. Thus
$\iint W_\ell(x-y)\dd\zeta(x)\dd\zeta(y)>0$ for every nonzero finite
signed measure $\zeta$ on the closed ball, by Fourier inversion.
The energy is strictly convex under mixtures on $\cK_f$ and the squared
cost to a fixed input is convex under mixtures of transport plans.
Compactness gives existence and strict convexity gives uniqueness.

To differentiate the minimizer along an input curve, we again
replace the density constraint by a smooth penalty. For integers
$m\ge2$ use $\int_\Omega U_m(x,r)\dx$, where
\[
 U_m(x,r)=\frac{g(x)r}{m-1}\left(\frac r{f(x)}\right)^{m-1}
              +\frac1m r\log r .
\]
Set the regularized energy equal to $+\infty$ on measures with a singular
part. Its integral term is superlinear and strictly convex in $r$.
Compactness, lower semicontinuity and comparison with the minimizer at
the limiting input give a unique continuous resolvent $J_{\tau,m}$.
We first obtain regularity for an individual minimizer $r\dx$.
Along the segment from $r\dx$ to any competitor, the right derivative
of the interaction agrees with that of the linear potential
$\gamma W_\ell*(r\dx)$. The remaining terms are convex under mixtures,
so their increments dominate their right derivatives. It follows that
$r\dx$ also minimizes the local problem with integrand
$U_m(x,s)+(V+\gamma W_\ell*(r\dx))s$. For a smooth positive input,
\cref{lem:compact-resolvent}\,\textup{(i)} now gives a smooth positive
output and a smooth transport diffeomorphism with $H\succ0$.

Along an input curve the potential $W_\ell*\rho$ changes with the
minimizer. Differentiating it as well as the local term in the
direction $\dot r=rh$ gives
\[
 \cK_\rho h=a h+\gamma W_\ell*(h\rho),\qquad
 a=r(U_m)_{rr}(x,r)>0 .
\]
In \eqref{eq:coupled-linearized-system} the conormal operator
$I+\tau\cK_\rho B$ is a compact perturbation of $I+\tau aB$ on the
H\"older spaces of \cref{lem:compact-resolvent}, since convolution by
$W_\ell$ is smoothing. The latter operator is an isomorphism. Testing
the homogeneous equation against $Bz$ gives
$\int H^{-1}[\nabla z,\nabla z]\dd\rho+
\tau\inner{Bz}{\cK_\rho Bz}_\rho=0$.
Both terms are nonnegative by Fourier positivity, so $z$ is constant
and the equation forces $z=0$. Fredholm invertibility, the same lifting
of the conormal datum, and the implicit function argument of
\cref{lem:compact-resolvent} now give a smooth branch along the given
input curve. Transport duality and mixture convexity identify its
solutions as the unique minimizers.

A minimizer of the regularized problem may exceed the cap, so the
definition of $\Lambda_f$ does not yet give a lower bound for its
convolution Hessian. To control this error uniformly in the input,
compare with $\rho_*$, whose regularized energy stays bounded. The
cost, interaction and $V$ are bounded on the ball, and
$\int r\log r\dx\ge-\log|\Omega|$. Thus each output $r\dx$ satisfies,
with $C$ independent of $m$ and of the input,
\[
 \int_{\{r>(1+\delta)f\}}r\dx
 \le \frac{C(m-1)}{(\min g)(1+\delta)^{m-1}}
 \qquad(\delta>0).
\]
For every $\delta>0$, any weak limit as $m\to\infty$ is therefore
dominated by $(1+\delta)f\dx$, even if the inputs vary with $m$. Hence
every such limit is feasible. If the least eigenvalue stayed a fixed amount below $\Lambda_f$
along a subsequence, compactness would give limiting points, directions
and a feasible limiting output. Continuity of $\D^2W_\ell$ would contradict the definition of
$\Lambda_f$ at the feasible limit. Consequently there are
$\varepsilon_m\downarrow0$ such that at every output
$\D^2(W_\ell*(r\dx))\succeq(\Lambda_f-\varepsilon_m)I$.

The volume transform of $U_m$ is
$g(x)^ks^{1-k}/(m-1)-(d/m)\log s$, where $k=1+d(m-1)$.
It is jointly convex and nonincreasing in $s$ by the perspective
convexity used in \cref{cor:weighted-power-criterion}. Thus
\cref{prop:volume-variation} gives a nonnegative local form on compatible
fields. In \eqref{eq:complete-second-variation} the density term of
the interaction is nonnegative by Fourier positivity applied to
$(\diver_\rho X)\rho$, while its spatial Hessian contributes at least
$\gamma(\Lambda_f-\varepsilon_m)\norm X_\rho^2$.
Including $V$, we obtain the lower bound
$\lambda_m\norm X_\rho^2$ with
$\lambda_m=\kappa+\gamma(\Lambda_f-\varepsilon_m)$.
For all sufficiently large $m$, \cref{prop:primal-contraction} and
\cref{lem:smooth-paths} give
$\Lip_{W_2}(J_{\tau,m})\le(1+\tau\lambda_m)^{-1}$.

Finally, for each $r\dx\in\cK_f$ the power penalty is at most
$\max g/(m-1)$ and its entropy divided by $m$ tends to zero. Comparison
with every such competitor identifies each subsequential limit of
$J_{\tau,m}\mu$ as the unique constrained minimizer, using the lower
entropy bound and continuity of the cost, potential and interaction.
Passing to the limit at both inputs proves \eqref{eq:capped-gaussian-bound}.
\end{proof}

\begin{theorem}
\label{thm:capped-gaussian-sharp}
Under the assumptions of \cref{thm:capped-gaussian-bound}, let $V=0$ and
$0<\tau\gamma\alpha_{f,\ell}<\ell^2$. On $\Ptwo(\overline\Omega)$,
\[
   \Lip_{W_2}(J_\tau^f)=\frac{1}{1+\tau\gamma\Lambda_f}>1.
\]
\end{theorem}

\begin{proof}
The upper estimate follows from \cref{thm:capped-gaussian-bound}.
For the reverse estimate we prescribe the output and use transport
duality to certify it as a minimizer. Localized oscillations of this
output then detect the least eigenvalue of the convolution Hessian,
while the compact integral term in the input velocity tends to zero.

We first realize any smooth probability
$\rho=r\dx$ with $0<r<f$ on $\overline\Omega$ as an output. Write
$q_\rho=\gamma W_\ell*\rho$ and $S_\rho=\Id+\tau\nabla q_\rho$.
The potential $\Phi(x)=|x|^2/2+\tau q_\rho(x)$ is uniformly convex
since $D S_\rho\succeq(1+\tau\gamma\Lambda_f)I$. Moreover,
\[
 S_\rho(x)=(1-a(x))x+
 \frac{\tau\gamma}{\ell^2}\int W_\ell(x-y)y\dd\rho(y),
 \qquad a(x)=\frac{\tau\gamma}{\ell^2}(W_\ell*\rho)(x).
\]
The step restriction $\tau\gamma\alpha_{f,\ell}<\ell^2$ ensures that
$a(x)<1$, so this is a convex combination of $x$ and a weighted barycenter in
$\overline\Omega$. Hence $\mu_\rho=(S_\rho)_\#\rho$ is an admissible
input. For transport from $\rho$ to $\mu_\rho$ with cost $|x-y|^2/2$,
a source potential is $|x|^2/2-\Phi(x)=-\tau q_\rho(x)$ by convex
duality on $\overline\Omega$. For any $\sigma\in\cK_f$, it follows that
\[
 \frac{W_2^2(\sigma,\mu_\rho)-W_2^2(\rho,\mu_\rho)}{2\tau}
       \ge -\int q_\rho\dd(\sigma-\rho).
\]
In the difference of the interaction energies, the linear term is
$\int q_\rho\dd(\sigma-\rho)$. It cancels the right-hand side above,
leaving $\gamma/2$ times the nonnegative Gaussian quadratic form of
$\sigma-\rho$. Thus $\rho$ minimizes the objective at $\mu_\rho$ and
$J_\tau^f\mu_\rho=\rho$.

We can therefore vary the output directly. Fix $w=\nabla\psi$
with $\psi\in C_c^\infty(\Omega)$ and put $F_t=\Id+tw$. For small
$|t|$, this is a convex-gradient diffeomorphism of $\overline\Omega$
onto itself. The density of $\rho_t=(F_t)_\#\rho$ converges uniformly
to $r$ as $t\to0$, so it remains strictly below $f$. The preceding
construction gives an input $\mu_t=(S_{\rho_t})_\#\rho_t$ with
$J_\tau^f\mu_t=\rho_t$. We know the exact output distance,
$W_2(\rho_t,\rho)=|t|\norm w_{L^2(\rho)}$, by optimality of $F_t$.
To bound the input distance, couple $\mu_\rho$ and $\mu_t$ by
$(S_\rho,S_{\rho_t}\circ F_t)_\#\rho$ and differentiate
$\int\nabla W_\ell(F_t(x)-F_t(y))\dd\rho(y)$ at $t=0$. This yields
\begin{align*}
 \limsup_{t\to0}\frac{W_2(\mu_t,\mu_\rho)}{|t|}
 &\le \norm{A_\rho w-\tau\gamma \mathsf C_\rho w}_{L^2(\rho)},\\
 A_\rho(x)&=I+\tau\gamma D^2(W_\ell*\rho)(x),\qquad
 \mathsf C_\rho w(x)=\int D^2W_\ell(x-y)w(y)\dd\rho(y).
\end{align*}
Here $\gamma(D^2(W_\ell*\rho)w-\mathsf C_\rho w)$ is the interaction
Hessian computed in \cite[Appendix~C.2, Example~2, Eq.~(36)]{BonetEtAl}.
The continuous kernel makes $\mathsf C_\rho$ compact on $L^2(\rho;\R^d)$.

To eliminate the compact term while retaining the least eigenvalue
of $A_\rho$, fix $x_0\in\Omega$ and let $e$ be a unit eigenvector
for that eigenvalue at $x_0$. Choose a nonzero cutoff
$\chi\in C_c^\infty(B_\delta(x_0))$ with $B_\delta(x_0)\subset\Omega$.
The gradients $w_n=\nabla(\chi\sin(n e\cdot x)/n)$ oscillate in the
direction $e$ and converge weakly to zero in $L^2(\rho)$ by the
Riemann--Lebesgue lemma. Compactness then gives
$\mathsf C_\rho w_n\to0$ strongly. Applying the speed bound for each
fixed $n$ before sending $n\to\infty$ gives
\[
  \lim_{n\to\infty}
  \frac{\norm{A_\rho w_n-\tau\gamma \mathsf C_\rho w_n}_{L^2(\rho)}^2}
       {\norm{w_n}_{L^2(\rho)}^2}
   =\frac{\int \chi^2|A_\rho e|^2\dd\rho}{\int\chi^2\dd\rho}.
\]
As $\delta\downarrow0$, continuity makes the right-hand side tend to
$\lambda_{\min}(A_\rho(x_0))^2$. Interior points approximate boundary
points, giving $\Lip_{W_2}(J_\tau^f)\ge
(\min_{\overline\Omega}\lambda_{\min}(A_\rho))^{-1}$.

Finally, smooth densities strictly between zero and $f$ approximate
every $r_0\dx\in\cK_f$ weakly. Indeed, mix $r_0$ with $f/\int f$ in a
positive proportion $\varepsilon$. The ratio of the mixture to $f$ then
lies in a closed subinterval of $(0,1)$. Extend this ratio outside the
ball by a value in that subinterval, mollify, multiply by $f$, and
normalize the mass. For fixed $\varepsilon$ the normalization factor
tends to one, preserving the strict margins. Letting $\varepsilon$
decrease to zero gives the required approximation. The Hessians of the
convolutions converge uniformly, by continuity of $D^2W_\ell$ on the
compact difference set. Minimizing over the limiting densities gives the
claimed lower bound. The strict inequality follows from $\Lambda_f<0$.
\end{proof}

\begin{example}
For $\Omega=(-1,1)$, $f=1$, $\ell=3$, and $\gamma=1$, the function
$W_3''(z)=(z^2-9)e^{-z^2/18}/81$ increases with $|z|$ for $|z|\le2$.
At fixed $x$, the minimizing density fills the unit interval of points
in $(-1,1)$ closest to $x$. The integral of $W_3''(x-y)$ over this
interval is smallest when it is centered at $x$, as occurs for $|x|\le1/2$. Hence
$\Lambda_f=\int_{-1/2}^{1/2}W_3''(z)\dd z=-e^{-1/72}/9$.
At $\tau=1$, \cref{thm:capped-gaussian-sharp} gives the expansion
factor $(1-e^{-1/72}/9)^{-1}$.
\end{example}

The same argument also compares two steps with a common density
bound and different interaction strengths. Write $J_{\tau,\gamma}^{f,V}$
for the map at strength $\gamma$. If $\gamma_0,\gamma_1>0$ and
$\D^2V\succeq\kappa I$, then $\lambda=\kappa+\max(\gamma_0,\gamma_1)\Lambda_f$
gives a common lower bound throughout the interval of strengths.
For $1+\tau\lambda>0$ and arbitrary inputs on the closed ball,
\[
 W_2(J_{\tau,\gamma_0}^{f,V}\mu_0,J_{\tau,\gamma_1}^{f,V}\mu_1)
 \le \frac{ W_2(\mu_0,\mu_1)+
      \tau|\gamma_1-\gamma_0|\,\|\nabla W_\ell\|_\infty}
      {1+\tau\lambda}.
\]
Indeed, interpolate the strengths linearly along the smooth input paths
of \cref{lem:smooth-paths}. For each power penalty the implicit function
argument above permits this smooth coefficient, and the excess mass
estimate is uniform on the compact interval of strengths. The derivative
at fixed density in \cref{rem:varying-energy} is
$p_t=(\gamma_1-\gamma_0)W_\ell*\rho_t$. Its gradient is bounded by
$|\gamma_1-\gamma_0|\,\|\nabla W_\ell\|_\infty$. Integrating the speed
bound, then passing to the path endpoints and the hard constraint as in
\cref{thm:capped-gaussian-bound}, proves the estimate.

\section{The entropy step on Ricci-flat manifolds}
\label{sec:ricci-resolvent}
On $\R^d$ the pressure of entropy is $p(r)=r$, so
\cref{thm:pressure-global} gives nonexpansiveness for every time step.
On a closed Ricci-flat manifold the heat flow is still nonexpansive
\cite{vonRenesseSturm}. We ask whether its variational time steps
retain this property.

The distinction between the flow and a single step appears when we
compare measures concentrated near two nearby points. Gigli and
Mantegazza used heat kernels centered at the points to obtain a metric
variation tangent to Ricci flow~\cite{GigliMantegazza}. Here we choose
a quadratic potential with Hessian $C$ near one point and move its
center in a direction $v$. In normal coordinates the minimizer has leading
covariance $\tau C^{-1}$, while its image under the transport map to
the input has covariance $\tau(I-C)C^{-1}(I-C)$. Subtracting the output covariance
from the input covariance thus gives $\tau(C-2I)$ at leading order. The scalar part has
zero pairing with $R(\,\cdot\,,v)v$ by Ricci-flatness. If the curvature
tensor is nonzero, we can choose $C$ so that the remaining pairing is
strictly positive. This will give a pair of inputs whose distance
increases after one step.

Let $(M,g)$ be a closed connected Riemannian manifold, and let $m$
be its normalized Riemannian volume. We use the geodesic distance on $M$ to define $W_2$. For $\tau>0$ put
\[
 J_\tau\mu=\argmin_{\rho\in\mathcal P(M)}
 \left\{\Ent_m(\rho)+\frac{1}{2\tau}W_2^2(\rho,\mu)\right\}.
\]
Compactness and lower semicontinuity give a minimizer. It is unique
because entropy is strictly convex on its finite domain under
ordinary mixtures, while mixing two couplings gives convexity of the
transport cost under the same mixtures. We call $J_\tau$ nonexpansive if
$W_2(J_\tau\mu,J_\tau\nu)\le W_2(\mu,\nu)$ for all $\mu,\nu\in\mathcal P(M)$.

\begin{theorem}\label{thm:ricci-resolvent}
Suppose $\operatorname{Ric}_g\equiv0$. The following are equivalent.
\begin{enumerate}[label=\textup{(\roman*)},leftmargin=*]
\item The Riemann curvature tensor vanishes identically.
\item $J_\tau$ is nonexpansive for every $\tau>0$.
\item There exists a sequence $\tau_j\downarrow0$ for which every
      $J_{\tau_j}$ is nonexpansive.
\end{enumerate}
If $M$ is nonflat, then for each sufficiently small $\tau>0$ there
are probability measures $\mu_\tau,\nu_\tau$ with smooth strictly
positive densities such that
\[
 W_2(J_\tau\mu_\tau,J_\tau\nu_\tau)>W_2(\mu_\tau,\nu_\tau).
\]
\end{theorem}

There is also a pointwise obstruction to carrying the Euclidean trace
estimate \eqref{eq:trace-positive} over to a general positive tensor.
Indeed, let $G$ be a smooth positive self-adjoint endomorphism on $M$,
$GX=\nabla u$, and $C=G^{1/2}(\nabla X)G^{-1/2}$. At a point, in an
orthonormal frame put $N_{ij}=\langle(\nabla_{e_j}G)X,e_i\rangle$.
Differentiating $GX=\nabla u$ gives $G\nabla X=\nabla^2u-N$. Since
$\tr(C^2)=\|\operatorname{sym}C\|_{\HS}^2-\|\operatorname{skew}C\|_{\HS}^2$,
\[
 \tr((\nabla X)^2)
 =\|\operatorname{sym}C\|_{\HS}^2
  -\tfrac14\|G^{-1/2}(N-N^{\mathsf T})G^{-1/2}\|_{\HS}^2 .
\]
The last square vanishes for a Codazzi tensor and measures the loss in
the derivative term of the Bochner identity for a general $G$. The
geometric obstruction in \cref{thm:ricci-resolvent} is detected below
by the anisotropic part of the Riemann curvature tensor.

\begin{lemma}\label[lemma]{lem:flat-entropy-resolvent}
Let $M$ be a closed connected flat Riemannian manifold with normalized
volume $m$. For every $\tau>0$ and all $\mu,\nu\in\mathcal P(M)$,
the entropy resolvent satisfies
\[
 W_2(J_\tau\mu,J_\tau\nu)\le W_2(\mu,\nu).
\]
\end{lemma}

\begin{proof}
\emph{The estimate on a flat torus.}
We first show how to apply the velocity estimate to arbitrary measures
on $M=\R^d/\Lambda$. Choose an optimal plan $\pi$ between $\mu$ and
$\nu$. For each $(x,y)$ choose a shortest displacement $\xi=\xi(x,y)$
on the cover, ordering the lattice to break ties measurably in a fixed
Borel fundamental cell. Smoothing the resulting geodesic will give a
smooth input curve with the same bound on its length.
Let $K_\epsilon$ be the strictly positive torus heat kernel of mass one
with respect to normalized volume. For $0\le s\le1$ put
\[
 a_s(z)=\int K_\epsilon(z-x-s\xi)\dd\pi,\qquad
 j_s(z)=\int \xi K_\epsilon(z-x-s\xi)\dd\pi .
\]
The displacements are bounded, so derivatives in $z,s$ pass under the
integrals. Thus $a_s>0$, $\partial_s a_s+\diver j_s=0$, and
Cauchy--Schwarz gives
\[
 \int_M |j_s|^2/a_s\dd m\le\int |\xi|^2\dd\pi=W_2^2(\mu,\nu).
\]
Thus $\eta_s=a_s m$ joins the heat convolutions of $\mu$ and $\nu$
with action at most $W_2^2(\mu,\nu)$. To estimate the length of
$s\mapsto J_\tau\eta_s$, we need regularity of the minimizer and its
dependence on $s$.

Fix $s$ and write $J_\tau\eta_s=r\,m$. The duality argument for
\eqref{eq:compact-Euler} gives a positive density and a Lipschitz cost
potential $-\tau q$, where $q=\log r$.
Lifting its cost-transform representation expresses
$\Phi(x)=|x|^2/2+\tau q(x)$ as a supremum of affine functions with
slopes the lifts of target points. Thus $\Phi$ is convex, its sections
are bounded, and $\nabla\Phi(x+k)=\nabla\Phi(x)+k$ for $k\in\Lambda$
at differentiability points. On the torus the optimal plan is supported
on this gradient modulo the lattice. Summing a compactly supported test
function over lattice translates shows that $\nabla\Phi$ transports
the periodic lifts of $r\,m$ and $a_s m$. The inverse convex potential
gives the reverse transport almost everywhere. The positive upper and
lower density bounds therefore bound the Alexandrov measure of $\Phi$
above and below by positive multiples of volume. The interior argument
in \cite[Section~4, p.~202]{CorderoPeriodic} applies on bounded
sections: local $C^{1,\alpha}$ regularity of $\Phi$ gives a
H\"older right-hand side $r/(a_s\circ\nabla\Phi)$, and interior
Monge--Ampère regularity gives $\Phi\in C^{2,\beta}_{\rm loc}$ for
some $\beta>0$. In dimension one the increasing transport relation
gives the same regularity directly. The equation
\[
 e^{q(x)}=a_s(x+\tau\nabla q(x))\det(I+\tau\D^2 q(x))
\]
now gives $H=I+\tau\D^2q\succ0$, uniformly on a fundamental cell.
Elliptic bootstrapping gives smoothness.

For uniqueness, let $q_1,q_2$ be periodic solutions with
$I+\tau \D^2q_i\succ0$. At a maximum $x_*$ of $q_1-q_2$, the gradients
agree and $\D^2q_1\preceq \D^2q_2$. Evaluating the equation there gives
\[
 e^{q_1(x_*)-q_2(x_*)}
 =\frac{\det(I+\tau \D^2q_1(x_*))}{\det(I+\tau \D^2q_2(x_*))}\le1 .
\]
Hence $q_1\le q_2$, and interchanging the two solutions gives equality.
As in \eqref{eq:coupled-linearized-system}, the linearization of the
logarithmic equation is $I-\tau\diver_{\rho_s}(H^{-1}\nabla)$.
The zeroth-order term makes this operator invertible on periodic
functions. The implicit function theorem gives a smooth local branch,
which uniqueness identifies with the log densities of
$\rho_s=J_\tau\eta_s=r_s m$. Since $\partial_s r_s/r_s$ is centered,
the periodic Poisson equation $Bu=\partial_s r_s/r_s$ has a smooth
centered solution by ellipticity and the Poincar\'e inequality on
the torus. For $X=H^{-1}\nabla u$, periodic integration by parts gives
$\mathscr Q_{\rho_s}(X)=\int\tr((\D X)^2)\dd\rho_s\ge0$ by
\cref{lem:symmetry}. With input velocity $j_s/a_s$, the periodic
form of \cref{prop:primal-contraction}
bounds the output length by
$\int_0^1(\int |j_s|^2/a_s\dd m)^{1/2}\dd s\le W_2(\mu,\nu)$.
As $\epsilon\downarrow0$, the endpoints converge to $\mu,\nu$ in $W_2$.
On the compact torus the squared costs converge uniformly over
competitors, so every limit point of the minimizers minimizes the
limiting objective. Uniqueness now gives the estimate for
the original endpoints.

\emph{Passage to a flat quotient.}
By Bieberbach's theorem there is a finite normal Riemannian cover
$\pi:\mathbb T\to M$ by a flat torus. For normalized volumes the
uniform lift $L$ of $r\,m$ has density $r\circ\pi$. The lift of a
singular measure is singular, so $L$ preserves relative entropy. To see
that it also preserves $W_2$, first project a coupling from the torus
to $M$. Since projection is $1$-Lipschitz, the distance downstairs is
at most the distance upstairs. For the reverse bound, start with an
optimal plan on $M$ and measurably select a minimizing geodesic for each
of its pairs.
Such a selection exists since the minimizing geodesics parametrized
with constant speed form a compact family with closed endpoint relation.
The lifts of each selected geodesic give a bijection of the two
fibres, with inverse obtained by reversing the geodesic. A uniform
choice of starting lift therefore has a uniform terminal lift.
The resulting coupling of the lifted measures has the original cost.

Every deck transformation preserves the objective for $L\mu$, so
uniqueness makes its torus minimizer deck-invariant. It is therefore the
uniform lift of its projection. Since lifting preserves both terms of
the objective, comparison with the lift of each competitor on $M$ gives
$J_\tau^{\mathbb T}L\mu=LJ_\tau^M\mu$. Applying the torus estimate
to the lifts of $\mu,\nu$ proves the claim.
\end{proof}

\begin{proof}[Proof of \cref{thm:ricci-resolvent}]
By \cref{lem:flat-entropy-resolvent}, \textup{(i)} implies
\textup{(ii)}, which implies \textup{(iii)}. For the converse, suppose
the curvature tensor is nonzero at $p\in M$. Choose a unit vector
$v\in T_pM$ for which the operator $R_vz=R(z,v)v$ is nonzero. We use
the convention that the sphere has positive sectional curvature. The operator $R_v$
is self-adjoint and $\tr R_v=\operatorname{Ric}(v,v)=0$. For small
$\theta>0$ we have $I+\theta R_v\succ0$. Put
\[
 C=\delta(I+\theta R_v),\qquad B_C=(I-C)C^{-1}(I-C).
\]
For $\delta>0$ sufficiently small, $0\prec C\prec I$ and
\[
 B_C=C^{-1}-2I+C,\qquad \tr(R_vC)=\delta\theta\tr(R_v^2)>0 .
\]

\emph{Construction of the minimizing curve.}
We prescribe the minimizers near a moving center and use the entropy
Euler equation to find the corresponding inputs. Write $p(h)=\exp_p(hv)$, let $P_h$ be parallel transport along $p(h)$
and put $E_h(z)=\exp_{p(h)}(P_hz)$. In normal coordinates at $p$,
cut off $\inner{(I+\theta R_v)z}{z}/2$ to a positive constant outside
a larger normal ball, keeping $p$ as its unique minimum. Multiplying
by $\delta$ gives a smooth function $\psi_0$ with
\[
 \psi_0(E_0(z))=\tfrac12\inner{Cz}{z}\quad\text{near }z=0,
 \qquad \|\psi_0\|_{C^2}=O(\delta).
\]
Let $F_h$ be the flow, starting at $\Id$, of a smooth cutoff of
$\partial_h E_h\circ E_h^{-1}$. For $|h|$ small,
$F_h(E_0(z))=E_h(z)$ in a smaller normal ball. Define
\[
 \rho_{\tau,0}=Z_\tau^{-1}e^{-\psi_0/\tau}m,\qquad
 \rho_{\tau,h}=(F_h)_\#\rho_{\tau,0}=r_{\tau,h}m .
\]
Set $\psi_{\tau,h}=-\tau\log r_{\tau,h}$. The change of variables
formula gives, up to a spatial constant,
\[
 \psi_{\tau,h}=(\psi_0+\tau\log\operatorname{Jac}F_h)\circ F_h^{-1},
\]
where the Jacobian is taken with respect to $m$. Choose $\delta$ small,
then restrict $|h|$ and $\tau$ so that $\psi_{\tau,h}$ is uniformly
small in $C^2$ modulo constants. For $c(x,y)=d_g(x,y)^2/2$, it is
then $c$-concave. To see the supporting inequality, let
$\psi=\psi_{\tau,h}$ and $y=\exp_x(-\nabla\psi(x))$. A minimizer of
$c(z,y)-\psi(z)$ satisfies $d_g(z,y)^2/2\le\operatorname{osc}\psi$
by comparison with $y$. Choose a fixed strongly convex normal radius
on the compact manifold on which $\nabla_z^2 c(z,y)\succeq g/2$.
For sufficiently small $C^2$ norm modulo constants this ball contains
$x$ and every minimizer, and $c(\cdot,y)-\psi$ is strictly convex
there. Its gradient at $x$ vanishes, so $x$ is its unique minimum,
as in \cite[Corollary~3.4]{GlaudoConcavity}.
The resulting map $S_{\tau,h}(x)=y$ is $C^1$-close to the identity,
hence a local diffeomorphism and, by compactness, a covering. Scaling
$\nabla\psi$ gives a homotopy to the identity, so the induced map on
the fundamental group is onto and the covering has one sheet. In
particular the proposed input $\mu_{\tau,h}=(S_{\tau,h})_\#\rho_{\tau,h}$
has a smooth positive density.

For every $\sigma$ of finite entropy, Kantorovich duality and convexity
of entropy under mixtures give
\begin{align*}
 \tfrac12 W_2^2(\sigma,\mu_{\tau,h})-\tfrac12 W_2^2(\rho_{\tau,h},\mu_{\tau,h})
 &\ge \int\psi_{\tau,h}\dd(\sigma-\rho_{\tau,h}),\\
 \Ent_m(\sigma)-\Ent_m(\rho_{\tau,h})
 &\ge \int(1+\log r_{\tau,h})\dd(\sigma-\rho_{\tau,h}).
\end{align*}
Dividing the first inequality by $\tau$ and adding the second gives
$J_\tau\mu_{\tau,h}=\rho_{\tau,h}$, since
$\psi_{\tau,h}=-\tau\log r_{\tau,h}$ and the masses agree.

\emph{Comparison of distances.}
Since $B_C-C^{-1}=C-2I$ and $\tr R_v=0$, \cref{lem:ricci-distance-quotients}
gives
\begin{equation}\label{eq:ricci-distance-gap}
 \liminf_{h\to0}
 \frac{W_2^2(\rho_{\tau,h},\rho_{\tau,0})
       -W_2^2(\mu_{\tau,h},\mu_{\tau,0})}{h^2}
 \ge \tau\delta\theta\tr(R_v^2)+o(\tau).
\end{equation}
Fix $\tau>0$ small enough that the right-hand side is positive.
Then for all sufficiently small nonzero $h$,
\[
 W_2(J_\tau\mu_{\tau,h},J_\tau\mu_{\tau,0})
 =W_2(\rho_{\tau,h},\rho_{\tau,0})
 >W_2(\mu_{\tau,h},\mu_{\tau,0}).
\]
This proves strict expansion and hence \textup{(iii)}$\Rightarrow$\textup{(i)}.
\end{proof}

\section{Propagation of functional inequalities and Gaussian rigidity}\label{sec:talagrand}
Fix a time step $\tau>0$ and a reference probability measure
$m=Z^{-1}e^{-V}\dx$ with $V$ finite and $\kappa$-strongly convex,
$\kappa>0$. If the input satisfies a Poincar\'e, Talagrand or
log-Sobolev inequality with constant $C$, then in each case the
sharp universal bound on the output constant is the positive solution
$L$ of \eqref{eq:propagation-root}. For Poincar\'e, a test
function at the output gives, through the adjoint equation of
\cref{prop:abstract-contraction}, a test function for the input
inequality. For Talagrand and log-Sobolev the relevant tests are measures. We tilt the minimizer and push the tilted measure forward by the
optimal map from the minimizer to the original input. Taking this
image as the new input, we compare the value of the objective at
the tilted measure with its minimum.

The equality question leads us first to the Talagrand bound along a
displacement geodesic. The endpoint coupling of Gozlan,
Malamut and Waldspurger~\cite{GMW} allows us to keep the remainder from
the change of variables at each endpoint. Together with the convexity
remainder of the reference potentials, this identifies the translations
that attain equality. Both remainders will also enter the equality
analysis for the entropy step.

\subsection{Talagrand interpolation}
\label{sec:talagrand-interpolation-bound}\label{sec:talagrand-rigidity}
Gozlan, Malamut and Waldspurger~\cite[Theorem~2.1]{GMW} proved the
following bound with its optimal coefficient. We use their coupling
with the endpoints interchanged.

\begin{theorem}\label{thm:talagrand}
Let $\mu_i=e^{-V_i}\dx$ be probability measures on $\R^d$, where
$V_i$ is finite and $\kappa_i$-strongly convex, $\kappa_i>0$.
Let $T$ be the Brenier map from $\mu_0$ to $\mu_1$,
$F_t=(1-t)\Id+tT$ and $\mu_t=(F_t)_\#\mu_0$. Then
\begin{equation}\label{eq:T2-interpolation}
 C_T(\mu_t)\le a_t^2,\qquad
 a_t=\frac{1-t}{\sqrt{\kappa_0}}+\frac{t}{\sqrt{\kappa_1}},
 \qquad 0\le t\le1.
\end{equation}
The bound is optimal for every $t,\kappa_0,\kappa_1$.
\end{theorem}

In the notation of \cref{sec:finite-comparisons}, the finite remainder is
$\cJ(A)=\tr A-d-\log\det A$. When $A$ is similar to a symmetric
positive definite matrix, it vanishes precisely at $A=I$.
For a Brenier map from a log-concave reference measure,
Cordero-Erausquin's entropy identity~\cite[Eq.~(2.2)]{CorderoQuantitative}
retains this trace--log-determinant remainder together with the convexity
remainder of the reference potential. Here the maps at the two endpoints
are compositions with a common coupling and may have singular
distributional derivatives. The next estimate retains their contribution.

\begin{proposition}\label[proposition]{prop:endpoint-coupling}
Under the assumptions of \cref{thm:talagrand}, write
$W_i=V_i-\kappa_i|\cdot|^2/2$. Let $\nu_0\in\Ptwo(\R^d)$ have finite entropy relative to $\mu_0$ and
suppose $\nu_1=T_\#\nu_0\in\Ptwo(\R^d)$. Let $G$ be the Brenier map
from $\mu_1$ to $\nu_0$, and put $R_0=G\circ T$, $R_1=T\circ G$.
Then $(R_i)_\#\mu_i=\nu_i$ and $R_1\circ T=T\circ R_0$. For $i=0,1$,
with $A_i=\D^{\mathrm{ap}}R_i$ the approximate derivative,
\begin{equation}\label{eq:two-deficits}
 \begin{aligned}
 \Ent_{\mu_i}(\nu_i)-\frac{\kappa_i}{2}\int|R_i-\Id|^2\dd\mu_i
 \ge \int\!\bigl\{\Dcal_{W_i}(R_i x,x)+\cJ(A_i)\bigr\}\dd\mu_i +\int e^{-V_i}\dd(\diver R_i)^s\ge0 .
 \end{aligned}
\end{equation}
Here $(\diver R_i)^s$ is the singular part of the distributional
divergence, and $\Dcal_{W_i}$ is taken at differentiability points of $W_i$.
If the left-hand side vanishes, $R_i$ is a translation.
\end{proposition}

\begin{proof}
Because the convex potentials are finite on $\R^d$, their densities
are positive and locally Lipschitz, hence locally H\"older. On each
fixed ball, upper and positive lower bounds on the densities give
doubling on ellipsoids contained in the ball.
Absolute continuity also gives zero mass to Lipschitz hypersurfaces.
As both domains are $\R^d$, the surjectivity and local regularity in
\cite[Theorems~1.6--1.7]{FigalliJhaveri} give a global $C^1$
diffeomorphism $T$ with $\D T\succ0$.
In particular $\nu_1\ll\mu_1$ and
$\Ent_{\mu_1}(\nu_1)=\Ent_{\mu_0}(\nu_0)<\infty$.

The map $G$ is the gradient of a finite convex function, since $\mu_1$
is equivalent to Lebesgue measure on $\R^d$. It is locally bounded and locally
$BV$. Precomposition by the $C^1$ diffeomorphism $T$ preserves local
$BV$ for $R_0$, and the $BV$ chain rule applies to $R_1$ since $G$
has locally bounded range~\cite[Theorems~3.16 and~3.96]{AFP}.
Their approximate derivatives are
\[
 A_0=(\D^{\mathrm{ap}}G\circ T)\,\D T,\qquad
 A_1=(\D T\circ G)\,\D^{\mathrm{ap}}G .
\]
The absolutely continuous part of the Hessian measure gives
$\D^{\mathrm{ap}}G\succeq0$. Up to a Lebesgue null set, the approximately
differentiable map $G$ agrees on countably many measurable pieces with
Lipschitz maps. On each piece, the image of the set
$\{\det \D^{\mathrm{ap}}G=0\}$ is null by the area formula. Since
$G_\#\mu_1=\nu_0\ll\dx$, that set is also $\mu_1$-null.
Thus $\D^{\mathrm{ap}}G\succ0$ almost everywhere. The two products satisfy
\[
 \D T\,A_0=\D T\,(\D^{\mathrm{ap}}G\circ T)\,\D T\succ0,\qquad
 (\D T\circ G)^{-1}A_1=\D^{\mathrm{ap}}G\succ0 .
\]
Hence $\cJ(A_i)\ge0$ at both endpoints.

The approximate derivatives do not account for the singular part of
$\D R_i$. To bound this part by its trace, mollify the convex potential
of $G$ and use its gradient $G_\varepsilon$ in both
compositions. On each compact set the $G_\varepsilon$ are uniformly
bounded, and $\D T$ and $(\D T)^{-1}$ are bounded on the compact sets involved.
For $cI\preceq H\preceq CI$ and $P\succeq0$,
\[
 \norm{PH}_{\HS}=\norm{HP}_{\HS}
 \le C\tr P\le(C/c)\tr(PH).
\]
The regularized maps converge in $L^1_{\rm loc}$ by local boundedness
and almost everywhere convergence, so their divergences converge in
distributions. The divergences are nonnegative. Testing on a compactly
supported nonnegative function and using lower semicontinuity of total
variation gives $|\D R_i|\le C_K\diver R_i$ locally on each compact set $K$.
In particular $(\diver R_i)^s\ge0$, and its vanishing forces $\D^sR_i=0$.

Write $h_i=\diff\nu_i/\diff\mu_i$. In the area formula for the pushforward
on the disjoint Lipschitz pieces, each preimage contributes its source
density divided by its positive Jacobian. Retaining just the
contribution at $x$ gives
\[
 e^{-V_i(x)}
 \le h_i(R_i x)e^{-V_i(R_i x)}\det A_i(x)
 \qquad \mu_i\text{-a.e.}
\]
Taking logarithms and using $V_i=W_i+\kappa_i|\cdot|^2/2$ gives
\[
 \begin{aligned}
 \log h_i(R_i x)-\tfrac{\kappa_i}{2}|R_i x-x|^2
 \ge \Dcal_{W_i}(R_i x,x)-\log\det A_i(x)+\nabla V_i(x)\cdot(R_i x-x).
 \end{aligned}
\]
Take a smooth cutoff $0\le\chi_N\uparrow1$ with $|\nabla\chi_N|\le C/N$.
Integration by parts in the last term gives
\[
 \begin{aligned}
 \int\chi_N\nabla V_i\cdot(R_i-\Id)\dd\mu_i
 &=\int\chi_N(\tr A_i-d)\dd\mu_i
    +\int\chi_N e^{-V_i}\dd(\diver R_i)^s\\
 &\quad+\int\nabla\chi_N\cdot(R_i-\Id)\dd\mu_i .
 \end{aligned}
\]
The cutoff error tends to zero since $R_i-\Id\in L^1(\mu_i)$.
Combining the trace and determinant terms, we pass to the limit by
$\log h_i\in L^1(\nu_i)$ and monotone convergence of the three
nonnegative remainders. This gives \eqref{eq:two-deficits}.

If the left-hand side of \eqref{eq:two-deficits} vanishes, so does the
Jacobian remainder. Since $A_i$ is
similar to a positive definite matrix, this gives $A_i=I$ almost
everywhere. The singular divergence also vanishes, so the measure
bound above gives $\D^sR_i=0$. Together these statements give
$\D(R_i-\Id)=0$ in distributions, which makes $R_i$ a translation.
\end{proof}

\begin{proof}[Proof of \cref{thm:talagrand}]
Fix $0<t<1$ and $\sigma\in\Ptwo(\R^d)$ of finite entropy relative to
$\mu_t$. Put $\nu_0=(F_t^{-1})_\#\sigma$ and $\nu_1=T_\#\nu_0$.
The bound $\D F_t\succeq(1-t)I$ makes $F_t$ a global diffeomorphism
with $(1-t)^{-1}$-Lipschitz inverse. Since $tT=F_t-(1-t)\Id$,
both measures have finite second moment.
Invariance of relative entropy gives
$\Ent_{\mu_i}(\nu_i)=\Ent_{\mu_t}(\sigma)$ for both endpoints.
To estimate the distance at time $t$, couple $F_t(x)$ with
$F_t(R_0(x))$ for $x$ distributed according to $\mu_0$. Their
displacement is the weighted sum of the two endpoint displacements
because $F_t\circ R_0=(1-t)R_0+tR_1\circ T$. Minkowski's inequality
and \eqref{eq:two-deficits} therefore give
\[
 \begin{aligned}
 W_2(\mu_t,\sigma)
 &\le \norm{F_t\circ R_0-F_t}_{L^2(\mu_0)}\\
 &\le(1-t)\norm{R_0-\Id}_{L^2(\mu_0)}
       +t\norm{R_1-\Id}_{L^2(\mu_1)}
 \le a_t\sqrt{2\Ent_{\mu_t}(\sigma)} .
 \end{aligned}
\]
At either endpoint apply \eqref{eq:two-deficits} with both reference
measures equal to that endpoint. For $\mu_i=\N(0,\kappa_i^{-1}I)$,
the interpolant is $\N(0,a_t^2I)$ and every nonzero translation
attains the bound.
\end{proof}

To determine which translations can give equality, let $E$ be the
subspace of vectors $b$ for which $\D^2W_i\,b=0$ in
distributions for both $i$. For each $b\in E$, the function $b\cdot\nabla W_i$ is constant almost
everywhere. Integrating along $E$ gives the common Gaussian factors
of \cite[Proposition~4.2]{HanZhu},
\begin{equation}\label{eq:talagrand-gaussian-factors}
 \mu_i=\N_E(P_E b_{\mu_i},\kappa_i^{-1}I_E)\otimes\bar\mu_i ,
\end{equation}
where $b_\xi$ denotes the barycenter of a measure $\xi$, $P_E$ is
orthogonal projection, and $\bar\mu_i$ is the marginal on $E^\perp$.

\begin{corollary} 
\label[corollary]{cor:talagrand-rigidity}
Under the assumptions of \cref{thm:talagrand}, fix $0<t<1$ and
$\sigma\in\Ptwo(\R^d)$ with $\sigma\ne\mu_t$ and finite relative
entropy. Then
\begin{equation}\label{eq:talagrand-attained-equality}
 W_2^2(\sigma,\mu_t)=2a_t^2\Ent_{\mu_t}(\sigma)
 \quad\Longleftrightarrow\quad
 \sigma=(\Id+c)_\#\mu_t\quad\text{for some }c\in E\setminus\{0\}.
\end{equation}
If $E\ne\{0\}$, translations along $E$ attain the bound at every
$s\in[0,1]$, and hence $C_T(\mu_s)=a_s^2$.
\end{corollary}

\begin{proof}
In the coupling used for \cref{thm:talagrand}, both weights are
positive at an interior time. Equality therefore forces the two endpoint
deficits in \eqref{eq:two-deficits} to vanish, so $R_i=\Id+b_i$.
Equality in Minkowski's inequality makes the two translation vectors
point in the same direction. Their lengths are fixed by the two endpoint
bounds, which have the same relative entropy. Thus
\[
 \sqrt{\kappa_0}\,b_0=\sqrt{\kappa_1}\,b_1\ne0,\qquad
 \sigma=(\Id+(1-t)b_0+tb_1)_\#\mu_t .
\]
For each $i$, the vanishing potential remainder in \eqref{eq:two-deficits}
gives $\Dcal_{W_i}(x+b_i,x)=0$ almost everywhere. On almost every
line parallel to $b_i$, $W_i$ is therefore affine on $[x,x+b_i]$
for almost every $x$, hence on the whole line by overlap.
Thus $b_i^{\mathsf T}\D^2W_i\,b_i=0$ in distributions, and positivity
of the Hessian measure gives $\D^2W_i\,b_i=0$.
The $b_i$ are parallel, so the translation lies in $E$.

Conversely, \eqref{eq:talagrand-gaussian-factors} splits the optimal
transport. The Gaussian factor of $\mu_s$ on $E$ has covariance
$a_s^2I_E$, so for $c\in E\setminus\{0\}$ the translate has squared
distance $|c|^2$ from $\mu_s$ and relative entropy $|c|^2/(2a_s^2)$.
\end{proof}

\subsection{Poincar\'e and Talagrand constants after one entropy step}
\label{sec:one-step-propagation}
Fix $\tau>0$ and $m=Z^{-1}e^{-V}\dx$, where $V$ is finite and
$\kappa$-strongly convex, $\kappa>0$. Write $J_\tau^m$ for the
resolvent of $\Ent_m$. For each of the three inequalities, an input bound $C$ gives a sharp
upper bound $L$ on the corresponding output constant. It is the positive
solution of
\begin{equation}\label{eq:propagation-root}
 a=1+\kappa\tau,\qquad
 aL=\tau+\sqrt{CL},\qquad
 \sqrt L=\frac{\sqrt C+\sqrt{C+4a\tau}}{2a}.
\end{equation}

\begin{theorem}\label{thm:poincare-jko}
If $\mu\in\Ptwo$ satisfies $\Var_\mu(f)\le C\int|\nabla f|^2\dd\mu$
for $f\in C_c^\infty(\R^d)$, $C\ge0$, then $\rho=J_\tau^m\mu$
satisfies the same inequality with constant $L$ given by
\eqref{eq:propagation-root}. The bound is optimal in each dimension.
\end{theorem}

\begin{proof}
First work on a ball with smooth positive input and smooth potential.
For a smooth $f$ of $\rho$-mean zero solve $(I+\tau B)u=f$ with
$X=H^{-1}\nabla u$ tangent to the boundary, using
\cref{lem:compact-resolvent}\,\textup{(i)}. Integration gives
$\int u\dd\rho=0$. The entropy form in the adjoint estimate of
\cref{prop:abstract-contraction} is bounded below by
$\kappa\norm X_\rho^2$, by the compatible derivative square and
convexity of the boundary. Hence
$a\norm X_\rho^2\le\inner{X}{\nabla f}_\rho$ and
$\norm X_\rho\le a^{-1}\norm{\nabla f}_\rho$.
The point of solving for $u$ is that $u\circ T^{-1}$ is a centered
test function for the input measure $\mu$, where $T$ is the Euler map
from $\rho$ to $\mu$. Since $\nabla(u\circ T^{-1})\circ T=X$, the input
inequality gives $\norm u_\rho\le\sqrt C\norm X_\rho$. Pairing the
equation for $u$ with $f$ now estimates the output variance,
\[
 a\norm f_\rho^2
 =a\inner f u_\rho+a\tau\inner{\nabla f}X_\rho
 \le\sqrt C\norm f_\rho\norm{\nabla f}_\rho
             +\tau\norm{\nabla f}_\rho^2 .
\]
For nonconstant $f$, division by $\norm{\nabla f}_\rho^2$ gives a
quadratic inequality for $\norm f_\rho/\norm{\nabla f}_\rho$. Its positive
root is $\sqrt L$ by \eqref{eq:propagation-root}, as required.

The contracting maps and Neumann heat kernels in
\cref{lem:lsi-approximation}\,\textup{(i)} approximate the input with constants
tending to $C$. Indeed the contraction preserves the Poincar\'e bound.
For its image $\lambda$ and the heat semigroup $P_t$, variance
decomposition and the gradient estimate on the convex ball give
\[
 \Var_{\lambda P_t}(f)
 =\Var_\lambda(P_tf)
   +2\int_0^t\!\int P_s|\nabla P_{t-s}f|^2\dd\lambda\dd s
 \le(C+2t)\int |\nabla f|^2\dd(\lambda P_t).
\]
The minimizers converge in $W_2$ by
\cref{lem:lsi-approximation}\,\textup{(ii)}, including the approximation of $V$. Testing against $C_c^\infty$ functions passes the variance bound
to the limit. For quadratic $V$ and centered isotropic Gaussian input
of variance $C$, the output has variance $L$: its Euler map has slope
$a-\tau/L=\sqrt{C/L}$. For $C=0$, the input is $\delta_0$ in this
example and the output has density proportional to
$\exp[-(\kappa+1/\tau)|x|^2/2]$. A coordinate function, approximated by
cutoffs, gives sharpness.
\end{proof}

For the transport and log-Sobolev inequalities the test at the output
is a probability measure $\nu$ obtained by tilting the minimizer $\rho$.
Pushing $\nu$ forward by the optimal map from $\rho$ to the original
input defines a new input $\eta$. The graph of this map still couples
$\nu$ optimally to $\eta$, but the minimizer for
$\eta$ can differ from $\nu$. We estimate the difference of their
objective values, first on a closed ball $\overline\Omega$ with $V$
smooth near the ball, $\D^2V\succeq\kappa I$ and $m_\Omega$ the
normalized restriction of $m$. Write
\[
 F_\xi(\zeta)=\Ent_{m_\Omega}(\zeta)+\frac{1}{2\tau}W_2^2(\zeta,\xi),
 \qquad \cE_\tau(\xi)=\min_\zeta F_\xi(\zeta),\qquad
 J\xi=\argmin_\zeta F_\xi(\zeta).
\]
For an arbitrary input $\mu$, put $\rho=J\mu$ and let $T_0$ be the
optimal transport map from $\rho$ to $\mu$. We normalize its convex potential
and the corresponding input cost potential by
\begin{equation}\label{eq:entropy-potential-gauge}
 \Phi_0=|\cdot|^2/2+\tau\log(\diff\rho/\diff m_\Omega),\qquad
 v_0=(|\cdot|^2/2-\Phi_0^*)/\tau ,
\end{equation}
taking the conjugate over $\overline\Omega$. Indeed
\eqref{eq:compact-Euler}, applied to the entropy integrand in
\cref{lem:compact-resolvent}\,\textup{(i)}, gives a positive
Lipschitz density and $\log(\diff\rho/\diff m_\Omega)=c-\phi/\tau$
for a source Kantorovich potential $\phi$. Then $T_0=\nabla\Phi_0$
$\rho$-a.e. and $v_0=c+\phi^c/\tau$.

\begin{lemma}\label[lemma]{lem:lsi-finite-defect}
Let $h$ be Lipschitz on $\overline\Omega$, with $\int e^h\dd\rho=1$,
and set $\nu=e^h\rho$, $\eta=(T_0)_\#\nu$, $\sigma=J\eta$. Denote 
\begin{equation}\label{eq:lsi-envelope-defect}
 E=\Ent_\rho(\nu),\quad
 R=\cE_\tau(\eta)-\cE_\tau(\mu)-\int v_0\dd(\eta-\mu),\quad
 G=F_\eta(\nu)-F_\eta(\sigma).
\end{equation}
Then $R\ge0$ and $G=E-R$. There is a common optimal lift from $\eta$
to $\nu$ and $\sigma$ for which, writing $\pi$ for its endpoint
coupling and $D^2=\int|y-x|^2\dd\pi(x,y)$,
\begin{equation}\label{eq:lsi-finite-gap}
 \frac{a}{2\tau}D^2\le G
 \le D\norm{\nabla h}_\nu-\frac{a}{2\tau}D^2
 \le\frac{\tau}{2a}\norm{\nabla h}_\nu^2 .
\end{equation}
In particular $W_2(\nu,\sigma)\le D\le (\tau/a)\norm{\nabla h}_\nu$.
\end{lemma}

\begin{proof}
Fenchel's inequality for $\Phi_0,\Phi_0^*$ gives, for every pair
of measures on the ball,
\begin{equation}\label{eq:tilted-fenchel}
 F_\xi(\zeta)\ge\Ent_\rho(\zeta)+\int v_0\dd\xi .
\end{equation}
Equality holds at $(\zeta,\xi)=(\rho,\mu)$ and at $(\nu,\eta)$
because the graph of $T_0$ still carries an optimal coupling after the tilt.
It follows that $\cE_\tau(\mu)=\int v_0\dd\mu$,
$R=\cE_\tau(\eta)-\int v_0\dd\eta\ge0$ and $G=E-R$.
Subtract $\int h\dd\zeta$ from \eqref{eq:tilted-fenchel} with
$\xi=\eta$. Its right side becomes
$\Ent_\nu(\zeta)+\int v_0\dd\eta$, which attains its minimum at
$\zeta=\nu$, where the bound is exact. Thus $\nu$ minimizes the
tilted objective.

For a smooth positive base measure take the two optimal transport maps from
the base to the fixed endpoints $\nu,\sigma$, and use their common
lift. Along $\lambda_t=((1-t)x+ty)_\#\pi$, entropy relative to Lebesgue
measure is convex \cite[Proposition~9.3.9]{AGS}. The potential contributes
$\kappa$ to the convexity modulus. Interpolating the two transports to
the same base bounds the squared cost by its endpoint chord minus
$t(1-t)D^2$. For general $\eta$, approximate the base in $W_2$ by
smooth positive densities and keep the two endpoints fixed. A weak
limit of the three-coordinate lifts is optimal from $\eta$ to each
endpoint. Compactness of the ball gives convergence of all quadratic
costs in the lifts, and lower semicontinuity of entropy gives
\[
 F_\eta(\lambda_t)\le (1-t)F_\eta(\nu)+tF_\eta(\sigma)
                    -\frac{a}{2\tau}t(1-t)D^2 .
\]
At $t\uparrow1$, minimality of $\sigma$ yields the lower bound on $G$.
At $t\downarrow0$, minimality of the tilted objective at $\nu$ gives
\[
 G+\frac{a}{2\tau}D^2
 \le-\int\nabla h(x)\cdot(y-x)\dd\pi(x,y).
\]
Indeed the difference quotients of $h$ converge $\nu$-a.e. in the
first coordinate and are bounded by $\Lip(h)|y-x|$.
Cauchy--Schwarz and completion of the square give the upper bounds.
\end{proof}

\begin{theorem}\label{thm:talagrand-jko}
If $\mu\in\Ptwo$ satisfies
$W_2^2(\eta,\mu)\le2C\Ent_\mu(\eta)$ for every $\eta\in\Ptwo$ of
finite relative entropy, $C\ge0$, then $J_\tau^m\mu$ satisfies this
inequality with $L$ from \eqref{eq:propagation-root}.
The bound is optimal in each dimension.
\end{theorem}

\begin{proof}
On the ball use the normalized tilt and the measures
$\rho,\nu,\eta,\sigma$ of \cref{lem:lsi-finite-defect}. Its distance
estimate controls $W_2(\nu,\sigma)$. To control $W_2(\sigma,\rho)$,
use the ball contraction preceding \eqref{eq:entropy-composition-contraction}
with $\Theta(s)=s$, the input Talagrand inequality and monotonicity
of entropy under $T_0$ to obtain
\[
 aW_2(\sigma,\rho)\le W_2(\eta,\mu)
 \le\sqrt{2C\Ent_\mu(\eta)}\le\sqrt{2C\Ent_\rho(\nu)} .
\]
Thus both terms in the triangle inequality
$W_2(\nu,\rho)\le W_2(\nu,\sigma)+W_2(\sigma,\rho)$ are controlled, and
\begin{equation}\label{eq:T2-tilted-step}
 aW_2(\nu,\rho)\le \sqrt{2C\Ent_\rho(\nu)}+\tau\norm{\nabla h}_\nu .
\end{equation}

To remove the gradient term, use the variational argument of
Fontbona, Gozlan and Jabir~\cite{FontbonaGozlanJabir}. If the bound
fails for $k>L$, the functional
$\gamma\mapsto\Ent_\rho(\gamma)-W_2^2(\gamma,\rho)/(2k)$ attains a
negative minimum at some $\nu$ on the ball. A source cost potential
$p$ from $\nu$ to $\rho$ gives an affine lower support for the cost.
Replacing the cost by this support leaves the same minimizer, so
$\nu=e^{p/k}\rho/\int e^{p/k}\dd\rho$. For this tilt the optimal map to $\rho$ is $\Id-\nabla p$. Hence
$h=p/k-\log\int e^{p/k}\dd\rho$ satisfies
$\norm{\nabla h}_\nu=W_2(\nu,\rho)/k$. Substituting this and
$2\Ent_\rho(\nu)<W_2^2(\nu,\rho)/k$ in \eqref{eq:T2-tilted-step}
and dividing by the positive distance gives $a\le\sqrt{C/k}+\tau/k$, contrary to $k>L$.
Letting $k\downarrow L$ proves the ball inequality.

Mollifying $V-\kappa|\cdot|^2/2$ on each fixed ball as in
\cref{lem:lsi-approximation}\,\textup{(ii)} preserves the curvature bound and gives
uniform convergence of the objectives there, up to an additive constant.
For general $\mu$, let $P_R$ be the metric projection onto
$\overline B_R$ and set $\mu_R=(P_R)_\#\mu$. The transport inequality
holds for $\mu_R$ with constant $C$ because composing a density tilt
with $P_R$ preserves its relative entropy, and $P_R$ is a contraction.
Since $W_2^2(\mu_R,\mu)\le\int_{|x|>R}|x|^2\dd\mu\to0$,
\cref{lem:lsi-approximation}\,\textup{(ii)} gives $W_2$ convergence
of the ball minimizers to $J_\tau^m\mu$.
In both limits the transport inequality is closed under narrow convergence. Indeed a bounded positive continuous
tilt, normalized separately along the sequence, has convergent relative
entropy, and the squared distance is lower semicontinuous. Truncation of
the density followed by approximation in $L^1$ of bounded densities by
continuous ones gives the inequality for every finite-entropy tilt.
The Gaussian examples in \cref{thm:poincare-jko} attain $L$ on translations.
\end{proof}

\subsection{The log-Sobolev constant after one entropy step}
\label{sec:lsi-jko}\label{sec:lsi-rigidity}
For the same reference measure, we use the normalization
\begin{equation}\label{eq:lsi-normalization}
 \operatorname{LSI}(C):\qquad
 \int f^2\log\frac{f^2}{\int f^2\dd\mu}\dd\mu
 \le 2C\int|\nabla f|^2\dd\mu,\qquad f\in C_c^\infty(\R^d).
\end{equation}
In the discrete Bakry--\'Emery argument of Caillet and Coudreuse
\cite[Section~3]{CailletCoudreuse}, the entropy and Fisher information
of successive iterates are measured against the same reference law.
To prove an inequality for $J_\tau^m\mu$, we take the minimizer itself
as the reference law in the entropy--information bound.
For smooth positive densities the entropy--information form is
\[
 \Ent_\mu(\eta)\le\frac C2 I_\mu(\eta),\qquad
 I_\mu(\eta)=\int|\nabla\log(\diff\eta/\diff\mu)|^2\dd\eta .
\]
In the ball approximation below we test the same inequality on functions
that are smooth up to the boundary.

\begin{theorem}\label{thm:lsi-jko}
Suppose $V:\R^d\to\R$ is finite and $\kappa$-strongly convex for
some $\kappa>0$. If $\mu\in\Ptwo(\R^d)$ satisfies $\operatorname{LSI}(C)$, $0\le C<\infty$,
then $\rho=J_\tau^m\mu$ satisfies $\operatorname{LSI}(L)$ for every
$\tau>0$, where $L$ is the positive solution of
\begin{equation}\label{eq:lsi-jko-constant}
 (1+\kappa\tau)L=\tau+\sqrt{CL}.
\end{equation}
The constant is optimal in every dimension over this class of reference and input measures. In particular, the step preserves $\operatorname{LSI}(\kappa^{-1})$.
\end{theorem}

\begin{proof}[Proof of \cref{thm:lsi-jko}]
Write $a=1+\kappa\tau$ as in \eqref{eq:propagation-root} and assume
first that $C>0$. We work on the ball of \cref{lem:lsi-finite-defect}, with smooth
positive input $\mu$. Then $T_0$ is a smooth diffeomorphism of the
closed ball by \cref{lem:compact-resolvent}\,\textup{(i)}. For each
probability $\eta$ of finite $E=\Ent_\mu(\eta)$, set
$\nu=(T_0^{-1})_\#\eta$. Invertibility gives
$E=\Ent_\rho(\nu)$. Choose
\begin{equation}\label{eq:lsi-theta}
 \theta=\frac1a\sqrt{\frac CL},\qquad
 0<\theta<1,\qquad a\tau\theta^2=C(1-\theta),\qquad
 (1-\theta)L=\tau/a .
\end{equation}
By \cref{lem:lsi-finite-defect}, the objective gap $G=E-R$ for a Lipschitz tilt
has an upper bound in terms of its Fisher information at the output.
To turn this into a log-Sobolev inequality we need a lower bound on
$G$ in terms of $E$. We will prove that
\begin{equation}\label{eq:lsi-envelope-bound}
 R\le\theta E .
\end{equation}
On the compact ball the continuity of $W_2^2(\zeta,\eta)$ in $\eta$
is uniform over $\zeta$, so the minimum value $\cE_\tau$ is continuous.
Consequently $R-\theta E$ attains a maximum at some $\eta_*$ of finite
entropy, and the maximum is nonnegative since we can test it at $\mu$.
To identify $\eta_*$, use the dual maximizer from the proof of
\eqref{eq:compact-Euler}. Its target cost potential, divided by $\tau$,
gives a Lipschitz supporting potential $v_*$ of $\cE_\tau$ at $\eta_*$
because the input occurs linearly in the dual objective.
With $w_*=v_*-v_0$, the affine
support to $R$ at $\eta_*$ shows that $\eta_*$ also maximizes
$\int w_*\dd\eta-\theta\Ent_\mu(\eta)$. Thus Gibbs' formula gives
\[
 \frac{\diff\eta_*}{\diff\mu}
 =\frac{e^{w_*/\theta}}{\int e^{w_*/\theta}\dd\mu},\qquad
 \nabla w_*=\theta\nabla\log(\diff\eta_*/\diff\mu)\quad\text{a.e.}
\]

The density of $\eta_*$ relative to $\mu$ is positive and Lipschitz.
Hence $h_*=\log(\diff\eta_*/\diff\mu)\circ T_0$ is a Lipschitz
normalized tilt on the output, and \cref{lem:lsi-finite-defect} applies
to $\nu_*=(T_0^{-1})_\#\eta_*$. Since $\eta_*$ is absolutely
continuous, the two optimal maps out of $\eta_*$ determine their common
lift. The map to $\nu_*$ is $S_0=T_0^{-1}=\Id-\tau\nabla v_0$.
The map to $J\eta_*$ is $S_*=\Id-\tau\nabla v_*$ $\eta_*$-a.e.,
because $\tau v_*$ is a target Kantorovich potential up to a constant.
Their difference gives
\[
 D^2=\int|S_*-S_0|^2\dd\eta_*
     =\tau^2\int|\nabla w_*|^2\dd\eta_* .
\]
To apply the input inequality, extend the Lipschitz function
$\sqrt{\diff\eta_*/\diff\mu}$ across the ball, truncate the extension
to its positive range and mollify. The approximants converge uniformly and
in $W^{1,2}(\mu)$, so their entropies and Dirichlet integrals converge.
If $R(\eta_*)>\theta E(\eta_*)$, the lower bound on $G$ in
\eqref{eq:lsi-finite-gap}, the preceding formula for $D^2$ and the
input log-Sobolev inequality would give
\[
 \frac{a\tau\theta^2}{2}I_\mu(\eta_*)
 \le E(\eta_*)-R(\eta_*)<(1-\theta)E(\eta_*)
 \le\frac{C(1-\theta)}2 I_\mu(\eta_*).
\]
This contradicts \eqref{eq:lsi-theta} and proves \eqref{eq:lsi-envelope-bound}.

For smooth positive $\nu$, the objective gap is $G=E-R$ by
\cref{lem:lsi-finite-defect}. Apply \eqref{eq:lsi-envelope-bound} and
\eqref{eq:lsi-finite-gap} with $h=\log(\diff\nu/\diff\rho)$ and
$\eta=(T_0)_\#\nu$ to obtain
\[
 (1-\theta)\Ent_\rho(\nu)\le G\le\frac{\tau}{2a}I_\rho(\nu).
\]
By \eqref{eq:lsi-theta} this is $\operatorname{LSI}(L)$, testing with
densities proportional to $(f^2+\varepsilon)\rho$ and letting
$\varepsilon\downarrow0$. The approximation in \cref{lem:lsi-approximation}
extends the estimate to the whole space with the same $L$.
If $\mu$ satisfies $\operatorname{LSI}(0)$, it satisfies
$\operatorname{LSI}(\varepsilon)$ for every $\varepsilon>0$.
Letting $\varepsilon\downarrow0$ in the resulting output inequality
gives $L=\tau/a$ by \eqref{eq:lsi-jko-constant}.

For $V=\kappa|\cdot|^2/2$ and $\mu=\N(0,CI)$ (with $\mu=\delta_0$ at $C=0$),
$\rho=\N(0,LI)$ has Euler
map $x+\tau\nabla\log(\diff\rho/\diff m)(x)=(a-\tau/L)x=\sqrt{C/L}\,x$.
Thus \eqref{eq:resolvent-support} gives $\rho=J_\tau^m\mu$, with sharp
Gaussian constant $L$ \cite[Example~5]{WibisonoGaussian}. Substitution
of $C=\kappa^{-1}$ gives $L=\kappa^{-1}$, and $L\to C$ as
$\tau\downarrow0$.
\end{proof}

To determine when the output attains its log-Sobolev bound, we return to
the two estimates for $G$ in \eqref{eq:lsi-finite-gap}. We first compute
the difference between $G$ and each bound for smooth data on the ball.
For
$\eta_0=\mu$, $\eta_1=\eta$ with smooth positive densities, put
$\sigma_i=J\eta_i$, $\rho=\sigma_0$. Define $\Phi_i,v_i$ from
$\sigma_i$ as in \eqref{eq:entropy-potential-gauge} and write
\[
 T_i=\nabla\Phi_i,\qquad S_i=T_i^{-1}=\Id-\tau\nabla v_i,\qquad
 M_i=\D S_i,\qquad w=v_1-v_0 .
\]
Let $\nu=(S_0)_\#\eta$, $h=\log(\diff\nu/\diff\rho)$ and
$X_i=(\nabla w)\circ T_i$. The Euler maps are diffeomorphisms of the closed ball.
Here $D^2=\tau^2\int|\nabla w|^2\dd\eta$ in the common lift, and
$R=\int w\dd\eta$ by \eqref{eq:tilted-fenchel} at each minimizing pair.

\begin{lemma}\label[lemma]{lem:entropy-jacobian-deficits}
Put $V_\kappa=V-\kappa|\cdot|^2/2$ and, for $i=0,1$, define
\begin{equation}\label{eq:lsi-matrix-remainder}
 A_i=\Dcal_{V_\kappa}(S_{1-i},S_i)+\cJ(M_i^{-1}M_{1-i}).
\end{equation}
Both terms are nonnegative. With $r_\xi$ denoting the Lebesgue
density of $\xi$ and $n$ the outward unit normal, we have
\begin{align}
 \frac{\tau}{2a}I_\rho(\nu)-G
 &=\frac{a\tau}{2}\int |X_0-a^{-1}\nabla h|^2\dd\nu
   +\int A_0\dd\eta
   +\tau\int_{\partial\Omega}r_\nu X_0\cdot n\dd S,
   \label{eq:lsi-upper-defect}\\
 G-\frac{a\tau}{2}\int|\nabla w|^2\dd\eta
 &=\int A_1\dd\eta
   -\tau\int_{\partial\Omega}r_{\sigma_1}X_1\cdot n\dd S .
   \label{eq:lsi-lower-defect}
\end{align}
\end{lemma}

\begin{proof}
The Euler relation and change of variables give
\[
 \log r_{\eta_i}
 =v_i-|\Id-S_i|^2/(2\tau)-V(S_i)+\log\det M_i-\log Z_\Omega .
\]
Expand the potential term and the negative log determinant at the
first endpoint. Their remainders after removal of the linear terms
are $\Dcal_{V_\kappa}(S_1,S_0)$ and $\cJ(M_0^{-1}M_1)$, in addition to
the explicit quadratic from the lower curvature bound. Since
$S_1-S_0=-\tau\nabla w$ and $M_1-M_0=-\tau\D^2w$, the linear terms
combine using $\nabla\log r_\rho(S_0)=(\Id-S_0)/\tau-\nabla V(S_0)$.
This gives the pointwise identity
\begin{equation}\label{eq:lsi-log-jacobian}
 \log(\diff\eta/\diff\mu)-w
 =-\tau(\diver_\rho X_0)\circ S_0-\tfrac{a\tau}{2}|\nabla w|^2-A_0 .
\end{equation}
Expanding at the other endpoint instead gives the same left side
as $-\tau(\diver_{\sigma_1}X_1)\circ S_1+a\tau|\nabla w|^2/2+A_1$.
Each $M_i^{-1}M_{1-i}$ is similar to a positive definite matrix,
so the log-Jacobian remainders are nonnegative.

We integrate each identity against $\eta$. At the first endpoint the
weight after the change of variables is $\nu=e^h\rho$. Integration by
parts therefore differentiates $e^h$, giving
$\tau\int X_0\cdot\nabla h\dd\nu$ in addition to the boundary flux.
Combining this integral with the quadratic terms gives
\eqref{eq:lsi-upper-defect}. At the other endpoint the weight is
$(S_1)_\#\eta=\sigma_1$ itself, so integration of the weighted
divergence gives only the flux in \eqref{eq:lsi-lower-defect}.
On the convex boundary the vectors $\tau X_0=x-S_1(T_0(x))$ and
$-\tau X_1=x-S_0(T_1(x))$ have nonnegative outward components. Thus
both displayed boundary contributions are nonnegative.
\end{proof}

To use this calculation on $\R^d$, we first establish the Euler
equation there for an arbitrary input. In particular it will give
a positive density at every point where we take its logarithm.

\begin{lemma}\label[lemma]{lem:entropy-arbitrary-euler}
Let $m=e^{-V}\dx$ be a probability measure on $\R^d$ with $V$ finite
and strongly convex. For every $\mu\in\Ptwo(\R^d)$ and $\tau>0$,
$\rho=J_\tau^m\mu$ has a positive locally Lipschitz density $r$.
Its optimal transport map to $\mu$ satisfies
\begin{equation}\label{eq:whole-entropy-euler}
 T(x)=x+\tau\big(\nabla V(x)+\nabla\log r(x)\big)
 \qquad\text{for Lebesgue-a.e. }x.
\end{equation}
\end{lemma}
\begin{proof}
If a minimizing density $r$ vanished on a bounded Borel set $A$ of
positive volume, then $m(A)>0$. Mixing in a fraction $t$ of
$m|_A/m(A)$ changes the entropy by $t\log t+O(t)$ since the two
components are mutually singular and both have finite entropy.
Convexity of the squared transport cost under mixtures bounds its
increase by $O(t)$. This contradicts minimality for small $t>0$. Thus for
any input in $\Ptwo$ the minimizing density $r$ is positive a.e., and
its Brenier potential $\Phi$ to the input is finite everywhere.
For $Z\in C_c^\infty(\R^d;\R^d)$ compare $\rho$ with
$(\Id+tZ)_\#\rho$. The entropy difference is the integral against
$\rho$ of $V(x+tZ)-V(x)-\log\det(I+t\D Z)$. Local Lipschitz
continuity of $V$ and absolute continuity of $\rho$ allow
differentiation at zero. Moving the source end of the optimal coupling
$(\Id,T)_\#\rho$ gives an upper bound on the change of squared cost
with derivative $2\int(x-T)\cdot Z\dd\rho$. Minimality for both signs
of $t$ yields, in distributions,
\[
 \nabla r=r\nabla\bigl((\Phi-|\cdot|^2/2)/\tau-V\bigr).
\]
The function $q=(\Phi-|\cdot|^2/2)/\tau-V$ is locally Lipschitz by
convexity of $\Phi$ and $V$. The right-hand side of the weak
equation belongs to $L^1_{\rm loc}$, so $r\in W^{1,1}_{\rm loc}$
and the product rule gives $\nabla(re^{-q})=0$. Thus $r$ has a
positive locally Lipschitz representative.
The weak equation and $T=\nabla\Phi$ give \eqref{eq:whole-entropy-euler}.
\end{proof}

Ohta and Takatsu~\cite[Theorem~1.1]{OhtaTakatsu} showed that a complete
connected weighted Riemannian manifold of dimension at least two with
$\operatorname{Ric}_\infty\ge K>0$ has a one-dimensional Gaussian
factor if a nonconstant density attains the log-Sobolev bound $1/K$.
We now use the Jacobian remainders on $\R^d$ to identify the regular
densities attaining the log-Sobolev bound.

\begin{theorem}\label{thm:lsi-jko-rigidity}
Under the assumptions of \cref{thm:lsi-jko}, suppose $C>0$ and
$\mu=e^{-U}\dx$ with $U\in C^\infty(\R^d)$. Define the common Gaussian
subspace by the following identities in the sense of distributions,
\begin{equation}\label{eq:lsi-gaussian-subspace}
 \mathcal G=\{b\in\R^d: \D^2(V-\kappa|\cdot|^2/2)b=0,\quad
                       \D^2(U-|\cdot|^2/(2C))b=0\}.
\end{equation}
Let $\nu\in\Ptwo(\R^d)$ have a positive locally Lipschitz density
relative to $\rho$, with $0<\Ent_\rho(\nu)<\infty$ and
$I_\rho(\nu)=\int|\nabla\log(\diff\nu/\diff\rho)|^2\dd\nu<\infty$. Then
\begin{equation}\label{eq:lsi-attained-equality}
 \Ent_\rho(\nu)=\frac L2 I_\rho(\nu)
 \quad\Longleftrightarrow\quad
 \frac{\diff\nu}{\diff\rho}(x)=
 \exp\!\left(b\cdot(x-b_\rho)-\frac L2|b|^2\right)
 \quad\text{for some }b\in\mathcal G\setminus\{0\},
\end{equation}
where $b_\xi$ denotes the barycenter of $\xi$. The measure $\rho$ has an
independent Gaussian factor on $\mathcal G$ with covariance $LI_{\mathcal G}$
and mean $(P_{\mathcal G}b_\mu+\kappa\tau P_{\mathcal G}b_m)/a$.
\end{theorem}

\begin{proof}
\emph{Euler maps and moments.}
By \cref{lem:entropy-arbitrary-euler}, the minimizing density is positive and locally Lipschitz. When the input density
is also positive and locally Lipschitz, both measures have full
convex support, are locally doubling and have locally H\"older
densities bounded away from zero on compact sets. The same application
of \cite[Theorems~1.6--1.7]{FigalliJhaveri} as in
\cref{prop:endpoint-coupling} gives inverse $C^1$ diffeomorphisms
$T,S$ of $\R^d$.

Put $F=\log(\diff\nu/\diff\rho)$ and $\eta=(T_0)_\#\nu$, where $T_0$ is the
Euler map of $\rho$. Then $\diff\eta/\diff\mu=\exp(F\circ S_0)$ and
$\Ent_\mu(\eta)=E:=\Ent_\rho(\nu)<\infty$. Before applying the entropy step to $\eta$, we check that it has
finite second moment. The restrictions $\eta|_{B_N}/\eta(B_N)$ have
bounded support, and their entropies tend to $E$ by
\eqref{eq:entropy-restriction-recovery}. The Talagrand inequality for
$\mu$ bounds their distances from $\mu$ uniformly, hence also their
second moments. Lower semicontinuity now gives $\eta\in\Ptwo$.
Since $\diff\eta/\diff\mu=\exp(F\circ S_0)$ is positive and locally
Lipschitz, the maps and potentials of \cref{lem:entropy-jacobian-deficits} for
$\sigma_0=\rho$ and $\sigma_1=J_\tau^m\eta$ are available on $\R^d$. Put
$w=v_1-v_0$ and $p=\nabla w$. Then $\tau p=S_0-S_1\in L^2(\eta)$, since the two images of $\eta$ are
$\nu$ and $\sigma_1$. Also $v_i\in L^1(\eta)$, since the supporting planes
for $\Phi_i^*$ at $0$ and $y$ bound $\Phi_i^*(y)-\Phi_i^*(0)$
between $S_i(0)\cdot y$ and $S_i(y)\cdot y$ and both bounds are
integrable by H\"older's inequality. The same argument with
$(S_0)_\#\mu=\rho$ gives $v_0\in L^1(\mu)$.

\emph{Equality in the bounds on relative entropy.}
Normalize the potentials by \eqref{eq:entropy-potential-gauge} on
$\R^d$. The negative part of $\log(\diff\sigma_i/\diff m)$ is bounded
by a constant times $1+|x|^2$, by an affine minorant of $\Phi_i$.
Thus the entropy chain rule and the pointwise Fenchel inequality give
\eqref{eq:tilted-fenchel} for competitors in $\Ptwo$ whenever the
input potential is integrable. For the two optimal pairs and the
tilted graph $(\nu,\eta)$ the Fenchel bound is an equality.
In the latter case $\Phi_0(x)\le\Phi_0(0)+x\cdot T_0(x)$ is integrable
under $\nu$ since $(T_0)_\#\nu=\eta\in\Ptwo$. Hence
$R=\int w\dd\eta$ and $G=E-R$ are finite.
Equation \eqref{eq:lsi-theta} gives
\begin{equation}\label{eq:lsi-rigidity-deficits}
 (1-\theta)\left(\frac L2 I_\rho(\nu)-E\right)
 =\left(\frac{\tau}{2a}I_\rho(\nu)-G\right)+(\theta E-R).
\end{equation}
The second term is nonnegative by \cref{lem:lsi-envelope-whole-space}.
For the first term put $M_i=\D S_i$ and $X_0=p\circ T_0$. The
pointwise calculation \eqref{eq:lsi-log-jacobian} uses only first
derivatives of $V$ at $S_0(y)$ and of the inverse maps. It therefore
holds $\eta$-a.e. for finite strongly convex $V$, since
$(S_0)_\#\eta=\nu$ is absolutely continuous. With
$V_\kappa=V-\kappa|\cdot|^2/2$, its remainder is
$A_0=\Dcal_{V_\kappa}(S_1,S_0)+\cJ(M_0^{-1}M_1)$.
First multiply the pointwise identity by $\chi_N\circ S_0$ and
integrate against $\eta$, where $0\le\chi_N\uparrow1$ are smooth
compactly supported cutoffs with $|\nabla\chi_N|\le c/N$.
The weak divergence identity for $X_0$ against $\chi_N e^F\rho$
and completion of the square with the same cutoff give
\[
\begin{aligned}
 &\frac{\tau}{2a}\int\chi_N |\nabla F|^2\dd\nu
   -\int (\chi_N\circ S_0)\bigl[\log(\diff\eta/\diff\mu)-w\bigr]\dd\eta\\
 =&\frac{a\tau}{2}\int\chi_N|X_0-a^{-1}\nabla F|^2\dd\nu
   -\tau\int\nabla\chi_N\cdot X_0\dd\nu+\int(\chi_N\circ S_0)
        A_0\dd\eta .
\end{aligned}
\]
The cutoff error is at most $(c\tau/N)\int|X_0|\dd\nu$. The terms on
the left converge by finite Fisher information, finite relative entropy
and $w\in L^1(\eta)$. Monotone convergence on the two nonnegative
terms gives their finiteness and the global identity
\[
 \frac{\tau}{2a}I_\rho(\nu)-G
 =\frac{a\tau}{2}\int |X_0-a^{-1}\nabla F|^2\dd\nu
 +\int A_0\dd\eta .
\]

\emph{Necessity.}
If $\nu$ attains equality, both nonnegative terms on the right of
\eqref{eq:lsi-rigidity-deficits} vanish. In the first, vanishing of
$\cJ(M_0^{-1}M_1)$ gives $M_1=M_0$, so $p=(S_0-S_1)/\tau$ is a
constant vector. Since $X_0=p\circ T_0=p$, vanishing of the integral
of $|X_0-a^{-1}\nabla F|^2$ gives $\nabla F=ap$. This vector
is nonzero since a constant normalized tilt has entropy zero. With this
choice of $p$, vanishing of the potential remainder reads
\[
 \Dcal_{V_\kappa}(x-\tau p,x)=0\qquad\text{for a.e. }x.
\]
Applying the line-restriction argument from the proof of \cref{cor:talagrand-rigidity}
to the convex function $V_\kappa$ and the displacement $-\tau p$ gives
$\D^2V_\kappa\,p=0$ in distributions.

It remains to relate the direction $p$ to the input. The second
equality in \eqref{eq:lsi-rigidity-deficits} gives
$R(\eta)=\theta\Ent_\mu(\eta)$, so their first variations agree along
bounded perturbations of both signs. For
$h\in C_c^\infty$ with $\int h\dd\eta=0$, put $\eta_t=(1+th)\eta$.
The Fenchel bound at $\eta$ and \cref{lem:lsi-envelope-whole-space} give
\[
 t\int wh\dd\eta
 \le R(\eta_t)-R(\eta)
 \le\theta\bigl(\Ent_\mu(\eta_t)-\Ent_\mu(\eta)\bigr).
\]
For $|t|\norm{h}_\infty<1$, these are positive probability
measures and the entropy derivative is
$\int h\log(\diff\eta/\diff\mu)\dd\eta$, since $h$ is bounded and
centered and $\log(\diff\eta/\diff\mu)\in L^1(\eta)$.
Taking both signs of $t$ gives
$w=\theta\log(\diff\eta/\diff\mu)+\mathrm{const}$. Since
$F=\log(\diff\eta/\diff\mu)\circ T_0$ and $\D T_0$ is symmetric, differentiation
and $\nabla F=ap$ yield
\[
 \D T_0\,p=a\theta p=\sqrt{C/L}\,p.
\]
Writing $\rho=e^{-W}\dx$, the distributional Euler identity
$\D T_0=I+\tau(\D^2V-\D^2W)$ and $\D^2V_\kappa\,p=0$ give
$\D^2W\,p=L^{-1}p$, since $a-\sqrt{C/L}=\tau/L$.

The identity for $\D^2W\,p$ splits off a Gaussian of variance $L$
along $p$. Since $\D T_0$ is symmetric with $\D T_0\,p=\sqrt{C/L}\,p$,
$T_0$ splits along the same direction and has slope $\sqrt{C/L}$
there. Thus $\mu=(T_0)_\#\rho$ has a Gaussian factor of variance $C$,
so $p\in\mathcal G$. Normalizing $e^F$ gives
\eqref{eq:lsi-attained-equality} with $b=ap$.

\emph{The converse.}
By \eqref{eq:lsi-gaussian-subspace}, both $m$ and $\mu$ split on
$\mathcal G\oplus\mathcal G^\perp$. Relative entropy and squared
transport cost each dominate the sum of their marginal values. The
product of the marginal minimizers attains both bounds, so uniqueness
gives the same splitting for $\rho$. On $\mathcal G$, a Gaussian of
covariance $LI_{\mathcal G}$ and mean $z$ has Euler map of slope
$\sqrt{C/L}$. Its image has mean $az-\kappa\tau P_{\mathcal G}b_m$.
Taking
\[
 z=(P_{\mathcal G}b_\mu+\kappa\tau P_{\mathcal G}b_m)/a
\]
gives the required Gaussian factor by \eqref{eq:resolvent-support}.
A tilt by $b\in\mathcal G\setminus\{0\}$ now has entropy $L|b|^2/2$
and Fisher information $|b|^2$, hence attains equality.
\end{proof}

\begin{remark}\label[remark]{rem:talagrand-nonattained}
Let $\diff\rho=Z^{-1}e^{-\kappa x^2/2}(\cosh x)^{-1}\dx$ on $\R$,
with $\kappa>0$. Its potential has second derivative
$\kappa+(\cosh x)^{-2}>\kappa$, so $C_{\mathrm{LS}}(\rho)\le1/\kappa$ by the Bakry--\'Emery criterion
\cite{BakryEmery}. For the translate $\rho_c=(\Id+c)_\#\rho$,
symmetry and the Lipschitz bound on $\log\cosh$ give
\[
 W_2^2(\rho_c,\rho)=c^2,\qquad
 \Ent_\rho(\rho_c)=\kappa c^2/2+O(|c|).
\]
Letting $|c|\to\infty$ gives $C_T(\rho)=C_{\mathrm{LS}}(\rho)=1/\kappa$.
The stationary interpolation with both endpoints $\rho$ thus has the
sharp constant and no nontrivial equality measure by
\cref{cor:talagrand-rigidity}.

Now take $\kappa=1$, $m=\N(0,1)$ and $0<\tau<1$. The Euler map $T(x)=x-\tau\tanh x$ is
increasing and $1$-Lipschitz. By \eqref{eq:resolvent-support},
$\rho=J_\tau^m(T_\#\rho)$. Lipschitz pushforward gives
$C_{\mathrm{LS}}(T_\#\rho)\le1$, while \cref{thm:lsi-jko} rules out
strict inequality since \eqref{eq:propagation-root} gives $L<1$ for $C<1$.
Hence $C_{\mathrm{LS}}(T_\#\rho)=C_{\mathrm{LS}}(\rho)=1$.
Since $\rho$ has no Gaussian factor, \cref{thm:lsi-jko-rigidity}
excludes a nonconstant smooth positive extremizer.
\end{remark} 

\subsection{Translation equality for the flow and for one step}

For the Gaussian reference law, equality in the propagated Talagrand
bound can be tested directly on translations. Even a single nonzero
translation attaining the bound determines a Gaussian factor of the
input. We then compare this conclusion with what happens along the
Ornstein--Uhlenbeck flow.

\begin{theorem}
\label{thm:gaussian-translation-equality}
Let $\mu\in\Ptwo(\R^d)$ satisfy $C_T(\mu)\le C$ for some $C>0$.
For $\tau>0$, set $\rho=J_\tau^\gamma\mu$, where $\gamma=\N(0,I_d)$,
and let $L$ be the positive root in \eqref{eq:propagation-root} for
$\kappa=1$. For $b\ne0$, write $S_b(x)=x+b$ and $e=b/|b|$. Then
\[
 W_2^2((S_b)_\#\rho,\rho)=2L\Ent_\rho((S_b)_\#\rho)
\]
if and only if, in the decomposition $\R e\oplus e^\perp$,
$\mu=\N(m,C)\otimes\bar\mu$ for some $m\in\R$ and
$\bar\mu\in\Ptwo(e^\perp)$. In this case $\rho$ has the independent
factor $\N(m/(1+\tau),L)$ in direction $e$, and every translation
parallel to $e$ attains equality.
\end{theorem}
\begin{proof}
For any law $\alpha\in\Ptwo$ with $C_T(\alpha)\le K$, $K>0$,
comparison of barycenters gives
$\Ent_\alpha(\beta)\ge|b_\beta-b_\alpha|^2/(2K)$.
For each $v\in\R^d$, apply this to the normalized tilt
$\one_{B_n}\exp(v\cdot(x-b_\alpha))\alpha$ and complete the square
in the mean. Letting $n\to\infty$ gives
$\log\int e^{v\cdot(x-b_\alpha)}\dd\alpha\le K|v|^2/2$.
Since $W_2((S_b)_\#\alpha,\alpha)=|b|$, equality for the translate
gives equality in the entropy variational bound at $v=b/K$. Its
density relative to $\alpha$ is therefore
\begin{equation}\label{eq:translated-tilt}
 \frac{\diff (S_b)_\#\alpha}{\diff\alpha}(x)
 =\exp\big(b\cdot(x-b_\alpha)/K-|b|^2/(2K)\big).
\end{equation}

Let $p$ be the Lebesgue density of $\rho$ and $T$ its optimal transport map to $\mu$,
as in \cref{lem:entropy-arbitrary-euler}. By \cref{thm:talagrand-jko},
$C_T(\rho)\le L$. Suppose the displayed
equality of the theorem holds and set $s=1+\tau-\tau/L=\sqrt{C/L}>0$.
Applying \eqref{eq:translated-tilt} to $p(x-b)/p(x)$ gives
$\nabla\log p(x+b)=\nabla\log p(x)-b/L$ a.e., where positivity
and local Lipschitz continuity follow from \cref{lem:entropy-arbitrary-euler}.
The Euler equation now gives $T(x+b)=T(x)+sb$ a.e. Consequently $T$
sends $(S_b)_\#\rho$ to $(S_{sb})_\#\mu$. The input bound and
monotonicity of entropy under pushforward give
\[
 \frac{|sb|^2}{2C}\le \Ent_\mu((S_{sb})_\#\mu)
 \le \Ent_\rho((S_b)_\#\rho)=\frac{|b|^2}{2L}.
\]
Since $s^2/C=1/L$, the two ends coincide. In particular the input
translation attains the Talagrand bound and pushforward by $T$
preserves the relative entropy of this tilt. To use the second
equality, write $f=\diff (S_b)_\#\rho/\diff\rho$ and
$g=\diff (S_{sb})_\#\mu/\diff\mu$. Then
$g\circ T=\E_\rho[f\mid T]$. Equality in conditional Jensen for
$z\log z$ gives $f=g\circ T$ $\rho$-a.e., including when $T$ is
many-to-one. Substituting \eqref{eq:translated-tilt} on both sides
and using $s^2/C=1/L$ yields
\[
 e\cdot(T(x)-b_\mu)=s\,e\cdot(x-b_\rho).
\]
In coordinates $x=ze+y$, \eqref{eq:whole-entropy-euler}
therefore gives $\partial_e\log p(ze+y)=-z/L+k$ with $k$ independent
of $y$. Integration on lines shows that
$p(ze+y)=\exp(-z^2/(2L)+kz)h(y)$. In particular the output
coordinate has variance $L$ and mean $kL$. The Euler map acts
affinely on this coordinate, with slope $s$, and its orthogonal
component depends only on $y$. Its image thus has an independent
Gaussian coordinate of variance $s^2L=C$ and mean $(1+\tau)kL$.

Conversely, for a product input, relative entropy against $\gamma$
and the squared transport cost both dominate the sums of their
marginal values. The product of the marginal minimizers attains both
bounds and is therefore the unique minimizer. For the one-dimensional
Gaussian input $\N(m,C)$, minimization over the mean gives
$m/(1+\tau)$, and the centered Gaussian update in
\cref{thm:poincare-jko} gives variance $L$. Translation by any
multiple of $e$ now has the asserted Gaussian entropy and distance.
\end{proof}

For the Ornstein--Uhlenbeck flow a nonzero translation can attain the
Talagrand bound even when the initial law on the line is non-Gaussian.
A periodic perturbation of the Gaussian density gives an example in
which this difference persists at every finite time.

\begin{proposition}\label[proposition]{prop:periodic-talagrand}
Let $\gamma=\N(0,1)$ and $\diff\mu_\varepsilon
=Z_\varepsilon^{-1}e^{\varepsilon\cos x}\dd\gamma$ on $\R$.
For all sufficiently small $\varepsilon>0$, the Ornstein--Uhlenbeck
flow $\sigma_t$ starting at $\mu_\varepsilon$ has $C_T(\sigma_t)=1$
at every finite $t\ge0$. It is non-Gaussian, and translation by
$b_t=2\pi e^t$ attains equality. For the same initial measure and
any $\tau>0$, no nonzero translation of $J_\tau^\gamma\mu_\varepsilon$
attains equality in the propagated bound with constant one.
\end{proposition}

\begin{proof}
First we prove $C_T(\mu_\varepsilon)\le1$. Write $\E$ for
expectation under $X\sim\mu_\varepsilon$. If $T(x)=x+u(x)$ is an
increasing smooth diffeomorphism equal to $\Id$ off a compact
interval, put $F(s)=s-\log(1+s)$ and $A=\E F(u'(X))$.
The change of variables and integration by parts give
\[
 \Ent_{\mu_\varepsilon}(T_\#\mu_\varepsilon)
 -\tfrac12\E|u|^2
 =A+\varepsilon\E\big[\cos X-\cos(X+u)-u\sin X\big].
\]
The integral $A$ depends only on $u'$, so it controls the variation
of $u$ while allowing an arbitrary translation. We therefore write
$m=\E u$ and $v=u-m$ and first estimate $v$. With
$N(z)=\min(|z|^2,|z|)$ and $K$ independent of $0<\varepsilon\le1$,
we claim that
\[
 \E N(v)\le K A,\qquad \E|v|\le K(\sqrt A+A).
\]
Indeed, write $v=v_1+v_2$ with centered summands and with
$v_1'=u'\one_{\{|u'|\le1\}}$, $v_2'=u'-v_1'$.
The density ratio to $\gamma$ is bounded above and below uniformly.
Comparison with the Gaussian Poincar\'e inequality gives
$\E|v_1|^2\le K\E|v_1'|^2$. Also $\E|v_2|\le K\E|v_2'|$:
integrate $|v_2(x)-v_2(0)|\le|\int_0^x |v_2'(z)|\dd z|$ against
the even density $p_\varepsilon$, using
$\int_x^\infty p_\varepsilon(z)\dd z\le Kp_\varepsilon(x)$ for
$x\ge0$ by the same Gaussian comparison. Centering costs at most
a factor of two. Since $N(s)\le K F(s)$ for $s>-1$ and
$N(v_1+v_2)\le 4(v_1^2+|v_2|)$, these bounds prove the claim.

Put $d_m=|\sin(m/2)|$. Expanding $-\cos(X+m+v)$ at $X+m$
gives a remainder of absolute value at most $2N(v)$. Evenness of the
law gives $\E\sin X=0$, so the constant contribution in the square
brackets is $2d_m^2\E\cos X$. The linear contribution is bounded in
absolute value by $2d_m\E|v|$. The left-hand side of the entropy identity is therefore at least
\[
 (1-K\varepsilon)A+2\varepsilon d_m^2\E\cos X
       -K\varepsilon d_m\sqrt A .
\]
Since $\E\cos X\to e^{-1/2}$, Young's inequality makes this
nonnegative for all sufficiently small $\varepsilon>0$.
To extend to arbitrary finite-entropy competitors in $\Ptwo$, first
truncate their densities relative to $\mu_\varepsilon$ above and
away from zero and normalize. For each truncation smooth its
difference from one after cutting it off on a large interval.
Correct the mass by a vanishing multiple of a fixed compact smooth
bump of integral one against $\mu_\varepsilon$. The densities stay
uniformly bounded above and away from zero at this stage and
converge in $L^1((1+x^2)\,\diff\mu_\varepsilon)$. Their monotone
transports are smooth and equal to $\Id$ on both tails. These two
approximations give convergence of entropy and $W_2$.

With $a=e^{-t}$, $\sigma_t$ is the law of $aX+\sqrt{1-a^2}Z$ with
independent $Z\sim\gamma$. Tensorization and the linear contraction
$(x,z)\mapsto ax+\sqrt{1-a^2}z$ give $C_T(\sigma_t)\le1$. Here
tensorization follows by coupling the first marginals and then the
conditional second marginals, and using the chain rule for entropy.
To see why this bound is attained, use the Mehler formula to write the
density of $\sigma_t$ relative to $\gamma$ as
$\E_Z e^{\varepsilon\cos(ay+\sqrt{1-a^2}Z)}/Z_\varepsilon$.
This is even and $b_t$-periodic, so translating by $b_t$ leaves only
the Gaussian factor in the ratio of Lebesgue densities. In other words,
if $p_t$ denotes the Lebesgue density of $\sigma_t$, then
$p_t(y-b_t)/p_t(y)=\exp(b_ty-b_t^2/2)$. Since $\sigma_t$ is centered,
the translate has relative entropy $b_t^2/2$ and squared distance $b_t^2$.
As $a=e^{-t}>0$ for every finite $t$, division by the Gaussian
characteristic function of the independent noise shows that a Gaussian
output would require $X$ itself to be Gaussian.

Finally, $C_T(\mu_\varepsilon)=1$ by the case $t=0$.
In dimension one \cref{thm:gaussian-translation-equality} with
$C=L=1$ would force $\mu_\varepsilon$ to be Gaussian if a nonzero
translation attained the propagated bound. This proves the last assertion.
\end{proof}

\subsection{Endpoint curvature and log-Sobolev bounds}
With a $1$-strongly log-concave reference measure, the entropy step
preserves $\operatorname{LSI}(1)$ by \cref{thm:lsi-jko}. In the proof,
the Euler equation and transport duality express the relative entropy of
a tilt of the minimizer as the sum of the objective gap and a nonnegative
remainder from the proximal value function (\cref{lem:lsi-finite-defect}).

What happens to the log-Sobolev constant along a Wasserstein geodesic if
both endpoints are $1$-strongly log-concave? For these endpoints the
interpolants satisfy $C_P(\mu_t)\le1$ and $C_T(\mu_t)\le1$~\cite{HanZhu,GMW}.
In every dimension $d\ge2$ and at each fixed $t\in(0,1)$, however, one
can choose such endpoints with $C_{\mathrm{LS}}(\mu_t)$ arbitrarily
large~\cite{HanXuZhu}. For Wasserstein barycenters of three measures,
even $C_P$ can be arbitrarily large, with equal weights and smooth
$1$-strongly log-concave measures on $\R^2$~\cite{HanXuZhu}.

\clearpage
\appendix
\part*{Appendices}
\phantomsection
\addcontentsline{toc}{part}{Appendices}

\section{Approximation for moment measures}
\label[appendix]{app:source-approximation}

To pass the inverse stability estimates to nonsmooth measures, we need
an approximation of the source for which its moment measure also converges.
Gaussian convolution gives such an approximation. For reconstruction
with a confinement, we will instead approximate the prescribed targets
and prove convergence of the corresponding minimizers.

\begin{lemma}\label[lemma]{mm:lem:smoothing}
Let $\rho=e^{-\psi}\dx$ be a centered probability measure with an
essentially continuous $\kappa$-strongly convex potential, $\kappa>0$,
and $\int|\nabla\psi|^2\dd\rho<\infty$. For $t>0$, let
$\rho_t=\rho*\N(0,tI)=e^{-\psi_t}\dx$. Then $\rho_t$ is centered,
$\psi_t$ is smooth and finite on $\R^d$, and
\begin{equation}\label{mm:eq:smoothing-curvature}
 \frac{\kappa}{1+\kappa t}I\preceq \D^2{\psi_t}\preceq t^{-1}I.
\end{equation}
Moreover,
\begin{equation}\label{mm:eq:both-converge}
 \rho_t\longrightarrow\rho,\qquad
 \Moment\rho_t\longrightarrow\Moment\rho
 \quad\text{in }W_2\quad(t\downarrow0).
\end{equation}
\end{lemma}

\begin{proof}
Couple $X\sim\rho$ to $Y_t=X+\sqrt t\,G$ with an independent
$G\sim\N(0,I)$. The law $\rho_t$ of $Y_t$ is centered and has a smooth
positive density. The coupling gives $W_2(\rho_t,\rho)\le\sqrt{dt}$,
and the Gaussian convolution estimate before dilation
in \cite[proof of Lemma~3.2, Step~2]{HanZhu} gives
\eqref{mm:eq:smoothing-curvature}.

To compare the moment measures, we use the same coupling for the
gradients of the two potentials. This requires a formula for the
gradient after convolution, including when the original density
vanishes. By essential continuity, $p=e^{-\psi}$ is continuous along almost
every coordinate line and vanishes at finite endpoints of its
positivity interval there \cite[Lemma~3]{CorderoKlartag}.
On each such line $p$ is locally absolutely continuous in the
interior of that interval. Its endpoint values eliminate the boundary
terms, and $p|\nabla\psi|\in L^1$ allows us to integrate the
one-dimensional identity over the remaining coordinates. We obtain
\begin{equation}\label{mm:eq:weak-score}
 \nabla p=-p\nabla\psi\quad\text{in distributions}.
\end{equation}
Convolving this identity with the Gaussian density and dividing by
the positive density of $Y_t$ gives
\begin{equation}\label{mm:eq:conditional-score}
 \nabla\psi_t(Y_t)=\E[\nabla\psi(X)\mid Y_t].
\end{equation}
Thus the coupled gradients have exactly the required moment laws.
Among functions of $Y_t$, the conditional expectation is the closest
to $\nabla\psi(X)$ in $L^2$. For any $g\in C_c^\infty(\R^d;\R^d)$
we therefore have
\[
 W_2(\Moment\rho_t,\Moment\rho)
 \le\norm{\nabla\psi(X)-g(Y_t)}_2
 \le\norm{\nabla\psi(X)-g(X)}_2+\Lip(g)\sqrt{dt}.
\]
Letting $t\downarrow0$ for a fixed $g$ and then approximating
$\nabla\psi$ in $L^2(\rho;\R^d)$ proves \eqref{mm:eq:both-converge}.
\end{proof}

\begin{lemma} 
\label[lemma]{mm:lem:regularized-model}
Fix $\alpha>0$ and a finite $\alpha$-strongly convex potential $V$ on $\R^d$.
\begin{enumerate}[label=\textup{(\roman*)},leftmargin=*]
\item
Let $\Omega$ be a ball or $\R^d$, and let $\nu\in\Ptwo(\R^d)$.
The minimizer $\rho$ of $\mathcal F_V(\cdot;\nu)$ over
$\Ptwo(\overline\Omega)$ admits the representation
\[
 \dd\rho=e^{-V-\Phi}\one_\Omega\dx,\qquad
 (\nabla\Phi)_\#\rho=\nu
\]
for some proper lower semicontinuous convex function $\Phi$.

\item
Suppose $\Omega$ and $\Lambda$ are balls, $V$ is smooth near
$\overline\Omega$, and the target probabilities $\nu_t=m_t\dy$ on $\Lambda$
have smooth positive densities on $[0,1]\times\overline\Lambda$.
Then the normalized potentials $\Phi_t$ and the source densities
$r_t=e^{-V-\Phi_t}$ depend smoothly on $t$. For each $t$,
$\nabla\Phi_t$ is a smooth diffeomorphism from $\overline\Omega$
onto $\overline\Lambda$.

\item
For any $\nu_0,\nu_1\in\Ptwo(\R^d)$, there is a sequence of ball
problems as in \textup{(ii)}, with smooth $\alpha$-strongly convex
confinements, such that at each endpoint $i=0,1$ the targets and
their reconstructions converge in $W_2$ to $\nu_i$ and
$\Reconstruct_V\nu_i$, respectively. The actions of the target
curves converge to $W_2^2(\nu_0,\nu_1)$. If the two targets have
the same barycenter, the approximating target curves can be chosen
to have constant barycenter.
\end{enumerate}
\end{lemma}

\begin{proof}
\emph{(i) The Gibbs density.}
Let $\Phi$ be a Brenier potential from the minimizer $\rho$ to $\nu$.
On an optimal plan, $\Phi(x)+\Phi^*(y)=x\cdot y$. Affine minorants
and finite second moments give $\Phi\in L^1(\rho)$ and $\Phi^*\in L^1(\nu)$.
If $\supp\nu\subseteq\overline B_L$, taking the supremum of
$x\cdot y-\Phi^*(y)$ over $\supp\nu$ gives a finite $L$-Lipschitz
choice of $\Phi$.

The optimal potential suggests a comparison density for which
the entropy and the integral of $\Phi$ can be combined into a
relative entropy. To this end set $Z=\int_\Omega e^{-\Phi}\dd m_V$ and
$\bar\rho=Z^{-1}e^{-\Phi}\one_\Omega m_V$. Since $\rho\ll m_V$ and
$\Phi$ is finite $\rho$-a.e., $Z>0$. An affine minorant for $\Phi$
and the quadratic lower bound on $V$ give $Z<\infty$,
$\bar\rho\in\Ptwo$ and $\Phi^-\in L^1(\bar\rho)$. The bound
$s e^{-s}\le1/e$ for $s\ge0$ gives $\Phi^+\in L^1(\bar\rho)$.
Fenchel's inequality bounds $\mathcal T(\bar\rho,\nu)$ by
$\int\Phi\dd\bar\rho+\int\Phi^*\dd\nu$, with equality for $\rho$ in
place of $\bar\rho$. Minimality therefore gives
\[
 0\le\Ent_{\bar\rho}(\rho)
 =\Ent_{m_V}(\rho)+\int\Phi\dd\rho+\log Z
 \le\Ent_{m_V}(\bar\rho)+\int\Phi\dd\bar\rho+\log Z=0 .
\]
Thus $\rho=\bar\rho$, and an additive constant in $\Phi$ gives \textup{(i)}.

\emph{(ii) Smooth dependence on the target measure.}
Since $\Lambda$ is bounded, the Lipschitz potential in \textup{(i)}
gives a positive Lipschitz density $r_t$ on $\overline\Omega$. For each
fixed $t$, the boundary regularity and bootstrapping in the proof of
\cref{lem:compact-resolvent} apply to the two balls, with
$r_t=e^{-V-\Phi_t}$. Thus $T_t=\nabla\Phi_t$ is a smooth
diffeomorphism of $\overline\Omega$ onto $\overline\Lambda$ and
\[
 \log\det\D^2\Phi_t+\log m_t(T_t)+V+\Phi_t=0,\qquad
 T_t(\partial\Omega)\subseteq\partial\Lambda .
\]
Put $H_t=\D T_t$ and $B_t=-\diver_{\rho_t}(H_t^{-1}\nabla)$.
Writing $w_t$ for the potential variation, we obtain from the Piola
identity and \eqref{eq:linear-conormal}
\[
 (B_t-I)w_t=(\partial_t\log m_t)\circ T_t,\qquad
 (H_t^{-1}\nabla w_t)\cdot n_\Omega=0\quad\text{on }\partial\Omega .
\]
The fixed-time Bochner identity \eqref{mm:eq:regularized-bochner}
forces $H_t^{-1}\nabla z=0$ whenever $B_tz=z$ with homogeneous conormal
data, so $z=0$. Since $B_t$ has compact resolvent, $1$ lies outside its
Neumann spectrum. The boundary lifting and Schauder estimates in the
proof of \cref{lem:compact-resolvent}
then give invertibility of $B_t-I$ with prescribed conormal data.
The implicit function theorem then gives a smooth local solution,
and the degree argument in the same proof shows that its maps are
diffeomorphisms onto $\overline\Lambda$. The Jacobian equation ensures
that each map sends its source to the prescribed target. As in the
argument following \eqref{mm:eq:regularized-selection}, this gives a
nonnegative first variation towards every competitor. Displacement
convexity then identifies the source with the unique minimizer for
that target. Hence $w_t=\partial_t\Phi_t$ and
$\partial_t r_t=-w_t r_t$.

\emph{(iii) Approximation of the endpoints.}
We first keep the source ball $B_R$ fixed while approximating the
targets and the confinement. Gluing optimal plans and applying
Cauchy--Schwarz gives, uniformly for $\supp\sigma\subseteq\overline B_R$,
\[
 |\mathcal T(\sigma,\nu)-\mathcal T(\sigma,\tilde\nu)|
 \le R W_2(\nu,\tilde\nu).
\]
Use the target curves of \cref{lem:smooth-paths} on expanding balls.
Their endpoints converge in $W_2$ and their actions converge to
$W_2^2(\nu_0,\nu_1)$. For equal endpoint barycenters their total flux
vanishes at each time, so the barycenter is constant. Mollifications
of $V$ retain the curvature bound and converge uniformly on $\overline B_R$.
The estimate above gives uniform convergence of the endpoint objectives
up to additive constants. Compactness and uniqueness then give $W_2$
convergence of the minimizers on the fixed ball.

To remove the source ball, write $\rho=\Reconstruct_V\nu$, let
$\rho_R$ be the ball minimizer and set $\sigma_R=\rho|_{B_R}/\rho(B_R)$.
By \eqref{eq:entropy-restriction-recovery} and $W_2$-continuity of
$\mathcal T$, $\mathcal F_V(\sigma_R;\nu)\to\mathcal F_V(\rho;\nu)$.
Minimality on the ball and \eqref{mm:eq:regularized-growth} give
\[
 \frac{\alpha}{2}W_2^2(\rho_R,\rho)
 \le\mathcal F_V(\sigma_R;\nu)-\mathcal F_V(\rho;\nu)\longrightarrow0 .
\]
For each source radius we may now choose a single target approximation
accurate at both endpoints and in action. Letting the source radius
increase along this sequence proves \textup{(iii)}.
\end{proof}

\section{Regularity and approximation of proximal maps}
\label[appendix]{app:analysis}

The proofs for proximal maps differentiate the minimizer as the input
measure varies. We first justify this differentiation on a ball, where
the transport boundary condition gives a conormal condition for the
linearized equation. We then approximate paths of input measures and
pass the estimates to the whole space. For the log-Sobolev estimates
we also need to control the constants during approximation.

\begin{lemma}\label[lemma]{lem:compact-resolvent}
Let $\Omega=B_R$ and $\tau>0$. The resolvent $J_\tau^{\cE}$ on
$\mathcal P(\overline\Omega)$ is single-valued and $W_2$-continuous
for either of the following energies, extended by $+\infty$ off the
absolutely continuous measures.
\begin{enumerate}[label=\textup{(\roman*)},leftmargin=2em]
\item $\cE(r\dx)=\int_\Omega F(x,r)\dx$. For $x$ near $\overline\Omega$,
 $F$ is continuous up to $r=0$, smooth for $r>0$, $F(x,0)=0$, and
 $F_r(x,\cdot)$ is an increasing diffeomorphism from $(0,\infty)$ onto $\R$.
\item $\cE=\Theta\circ\Ent_m$, where $m=Z^{-1}e^{-V}\dx\in\Ptwo(\R^d)$,
 $V$ is smooth and convex with bounded Hessian, and $\Theta$ is smooth
 and convex near $[0,\infty)$ with $\Theta'(0)>0$.
\end{enumerate}
If the input $\eta$ has a smooth positive density on $\overline\Omega$,
so does $\rho=J_\tau^{\cE}\eta$, and its density varies smoothly along
smooth curves of such inputs. The optimal map from $\rho$ to $\eta$ is
$T=\Id+\tau\nabla q_\rho$, where $q_\rho$ is the first variation of
$\cE$ at $\rho$. It is a smooth diffeomorphism of $\overline\Omega$
onto itself with $H=\D T\succ0$. With the notation of
\eqref{eq:resolvent-operators}, for every $f\in C^\infty(\overline\Omega)$
the conormal problem
\[
 (I+\tau\cK_\rho B)u=f \quad\text{in }\Omega,\qquad
 (H^{-1}\nabla u)\cdot n=0 \quad\text{on }\partial\Omega
\]
has a unique solution $u\in C^\infty(\overline\Omega)$.
\end{lemma}

\begin{proof}
In \textup{(i)}, $F_r(x,r)\to\infty$ uniformly in $x$ by monotonicity
and compactness, so $F(x,r)\ge Nr-C_N$ for every $N>0$. The convex
integral is lower semicontinuous with infinite recession cost. In
\textup{(ii)}, lower semicontinuity follows from that of entropy and
monotonicity of $\Theta$. In both cases compactness gives a minimizer,
since the uniform density has finite energy. Strict convexity under
mixtures gives uniqueness, and uniform convergence of the squared cost
in the other marginal gives $W_2$-continuity in the input.

To obtain the regularity of the minimizer, we first express its
density in terms of a Kantorovich potential. In \textup{(i)}, fix
the input $\eta$ and let $G(x,\cdot)$ be the
inverse of $F_r(x,\cdot)$. For the cost $|x-y|^2/2$, the compact convex
class $\phi(x_0)=0$, $\Lip(\phi)\le\diam\Omega$ contains an optimal
source potential for each pair of marginals. Kantorovich duality and
convex minimax give a saddle point for
\[
 (h,\phi)\longmapsto\int_\Omega F(x,h)\dx+
 \frac1\tau\left(\int_\Omega\phi h\dx+\int\phi^c\dd\eta\right).
\]
Indeed the functional is convex and lower semicontinuous in the probability
density $h$ and continuous and concave in $\phi$. Minimizing in $h$
at the saddle point gives
\begin{equation}\label{eq:compact-Euler}
 r(x)=G\bigl(x,c-\phi(x)/\tau\bigr),
\end{equation}
where $\int r\dx=1$ fixes $c$ by strict monotonicity and the endpoint
limits of $G$. In \textup{(ii)}, put $k=\Theta'(\Ent_m(\rho))>0$.
Testing minimality along $(1-s)\rho+s\sigma$, using convexity of
entropy and of the squared cost under mixtures, and letting $s\downarrow0$
gives
\[
 k\Ent_m(\rho)+W_2^2(\rho,\eta)/(2\tau)
 \le k\Ent_m(\sigma)+W_2^2(\sigma,\eta)/(2\tau).
\]
Thus the same $\rho$ minimizes the proximal objective for
$k\Ent_m$ with $k$ held fixed, and \eqref{eq:compact-Euler} applies
to that energy. In both cases the bounded Lipschitz potential $\phi$ and the
smooth positive inverse $G$ give a positive Lipschitz density on the closed ball.

With this density bound we can apply boundary regularity for
optimal transport. For a smooth positive input density $a$, write $T=\Id+\tau\nabla q$
and $H=\D T$. The densities $r,a$ are positive and H\"older on the
closed balls. For the equivalent cost $-x\cdot y$, Caffarelli's
boundary regularity theorem \cite{Caffarelli}, in the form of
\cite[Corollary~2.4]{ChenFigalliBoundary}, gives a $C^{2,\beta}$ convex
potential for some $\beta>0$.
Since $\det H=r/a(T)>0$, $H$ is uniformly positive definite.
Let $b$ define the target ball, with $b<0$ inside. On $\partial\Omega$,
$H\nabla b(T)=\nabla(b\circ T)$ is a positive multiple of the outer
normal $n$. Thus $b(T)=0$ linearizes in a gradient direction to
\begin{equation}\label{eq:linear-conormal}
 (H^{-1}\nabla z)\cdot n=0 .
\end{equation}
The Euler relation expresses $r$ smoothly in terms of $x$ and the
potential. Oblique boundary regularity for $r=a(T)\det H$ then gives
a smooth diffeomorphism up to the boundary
\cite[Remark~2.5]{ChenFigalliBoundary}. The same argument applies to
distinct source and target balls, as in \cref{mm:lem:regularized-model}\,\textup{(ii)}.

It remains to show that the map and density vary smoothly with the
input. To apply the implicit function theorem we allow nearby positive
densities of arbitrary mass and take
$q_\rho(x)=F_r(x,r(x))$ in \textup{(i)}. In \textup{(ii)}, writing
$m_0=Z^{-1}e^{-V}$ for the whole-space reference density, take
\[
 E(r)=\int_\Omega r\log(r/m_0)\dx,\qquad
 q_\rho=\Theta'(E(r))\bigl(1+\log(r/m_0)\bigr).
\]
Since $\D E(r)[rh]=\int_\Omega (1+\log(r/m_0))h\dd\rho$,
\eqref{eq:entropy-rank-one} holds for arbitrary $h$ with
$q_{\Ent_m}=1+\log(r/m_0)$.
Set $B=-\diver_\rho(H^{-1}\nabla)$. For $h=\dot r/r$ and $z=\dot q$,
the Euler relation and the Jacobian equation linearize to
\begin{equation}\label{eq:coupled-linearized-system}
 z-\cK_\rho h=f_1,\qquad h+\tau Bz=f_2 .
\end{equation}
The second identity follows from $rH^{-1}=a(T)\cof H$ and the Piola
identity. Eliminating $h$ gives
$(I+\tau\cK_\rho B)z=f_1+\cK_\rho f_2$.

In \textup{(i)}, $\cK_\rho$ is multiplication by $a_\rho=rF_{rr}(x,r)>0$.
In \textup{(ii)}, put $a_\rho=k=\Theta'(E(r))>0$. By
\eqref{eq:entropy-rank-one}, $\cK_\rho-kI$ is nonnegative and has rank
at most one. In either case, testing $(I+\tau\cK_\rho B)z=0$ with $Bz$
under the homogeneous conormal condition gives
\[
 0=\int_\Omega\inner{H^{-1}\nabla z}{\nabla z}\dd\rho
   +\tau\inner{\cK_\rho Bz}{Bz}_\rho .
\]
Both terms are nonnegative, so $z$ is constant and the kernel equation
gives $z=0$.

The conormal form of $a_\rho^{-1}+\tau B$ is coercive on $H^1(\Omega)$.
By Lax--Milgram and the boundary Schauder estimates
\cite[Chapter~6, Section~7]{GT}, $I+\tau a_\rho B$ is an isomorphism
from $C^{s+2,\beta}(\overline\Omega)$ with \eqref{eq:linear-conormal}
onto $C^{s,\beta}(\overline\Omega)$ for $s\in\mathbb N_0$ and $0<\beta<1$.
This proves surjectivity in \textup{(i)}. In \textup{(ii)},
$I+\tau\cK_\rho B$ is a finite-rank perturbation of this isomorphism,
so its Fredholm index is zero and injectivity gives surjectivity.

For $g\in C^{s+1,\beta}(\partial\Omega)$, the standard
$C^{s+2,\beta}$ extension with zero trace and normal derivative
$g/(n\cdot H^{-1}n)$ lifts the conormal datum $g$.
Extend $a$ smoothly and positively across the boundary. For $(\log r,q)$ in
$C^{s,\beta}\times C^{s+2,\beta}$, the residuals $q-q_\rho$,
$\log r-\log a(T)-\log\det H$ and $b(T)$ depend smoothly on the
unknowns with values in $(C^{s,\beta})^2\times C^{s+1,\beta}(\partial\Omega)$.
Their derivative is the isomorphism just proved, so the implicit
function theorem gives a smooth local branch. Strict convexity and the
boundary diffeomorphism persist on this branch. Degree one gives
surjectivity onto the target ball, and the Jacobian equation gives the
prescribed pushforward. Since $-\tau q$ is a source cost potential up
to a constant, transport duality and the supporting inequality for
$\cE$ identify the solution as the minimizer. Uniqueness joins the
local branches along the input curve.
\end{proof}

\begin{lemma}\label[lemma]{lem:smooth-paths}
For any $\mu_0,\mu_1\in\mathcal P(\overline B_R)$, there are smooth curves
of strictly positive densities $h_t$ on the closed ball, with smooth fluxes
$j_t$ vanishing near the boundary, whose endpoints converge in $W_2$ to
$\mu_0,\mu_1$ and whose actions have upper limit at most
$W_2^2(\mu_0,\mu_1)$.

For $\mu_0,\mu_1\in\Ptwo(\R^d)$, let $\pi$ be an optimal plan,
$D=W_2(\mu_0,\mu_1)$ and $\eta_t=((1-t)x+ty)_\#\pi$.
There are balls $B_{R_j}$, $R_j\to\infty$, and curves $\eta_{j,t}$ of
smooth positive densities with smooth velocities $v_{j,t}$ vanishing
near the boundary, such that
\[
 \sup_{t\in[0,1]}W_2(\eta_{j,t},\eta_t)\longrightarrow0,
 \qquad
 \norm{v_{j,t}}_{L^2(\eta_{j,t})}^2\le D^2,
 \qquad
 \int_0^1\norm{v_{j,t}}_{L^2(\eta_{j,t})}^2\dt\longrightarrow D^2.
\]
In both assertions, if the endpoints have the same barycenter, the curves
can be chosen so that the spatial integral of the flux vanishes at every time.
\end{lemma}

\begin{proof}
We approximate the paths of particles in an optimal coupling before
smoothing them. This keeps control of the action and, when the
barycenters agree, of the integral of the flux. Let $(X,Y)$ have law
$\pi$, choosing any optimal coupling in the fixed-ball case, and set
$D^2=\E|Y-X|^2$. Choose increasing finite $\sigma$-algebras $\mathcal F_n$
generating $\sigma(X,Y)$. For $n\ge2$, set
\[
 X_n=(1-1/n)\E[X\mid\mathcal F_n],\qquad
 Y_n=(1-1/n)\E[Y\mid\mathcal F_n].
\]
Then $(X_n,Y_n)\to(X,Y)$ in $L^2$, and conditional expectation gives
\[
 \E|Y_n-X_n|^2\le D^2,\qquad
 \E(Y_n-X_n)=(1-1/n)\E(Y-X).
\]
On a fixed ball, convexity and the factor $1-1/n$ place both endpoints
strictly inside $B_{R_n}$ with $R_n=R$. On $\R^d$, choose $R_n\uparrow\infty$
so that $B_{R_n}$ contains all values of $X_n,Y_n$ in its interior.

Write $B=B_{R_n}$ and $Z_t=(1-t)X_n+tY_n$. Choose a smooth probability
kernel $\zeta$ supported in $B_{\varepsilon_n}$, with $\varepsilon_n\downarrow0$
small enough that every segment $[X_n,Y_n]$ stays more than $\varepsilon_n$
from $\partial B$. Smoothing these segments gives a density supported away from the
boundary. To make it positive on the closed ball, mix it with a
stationary uniform density of mass $0<\theta=\theta_n<1$, chosen so
that $\theta_n(1+R_n^2)\to0$. The density and flux are then
\begin{equation}\label{eq:critical-data-approximation}
 \begin{aligned}
 h_t(z)&=\theta/|B|+(1-\theta)\E[\zeta(z-Z_t)],\\
 j_t(z)&=(1-\theta)\E[(Y_n-X_n)\zeta(z-Z_t)].
 \end{aligned}
\end{equation}
Then $h_t>0$, $\partial_th_t+\diver j_t=0$, and $v_t=j_t/h_t$ is smooth
and vanishes near $\partial B$. The total flux is $(1-\theta)\E(Y_n-X_n)$,
which vanishes for endpoints of equal barycenter. Pointwise Cauchy--Schwarz gives
\begin{equation}\label{eq:critical-data-action}
 \int_B h_t|v_t|^2\dx\le(1-\theta)\E|Y_n-X_n|^2\le D^2 .
\end{equation}
The interpolation of $(X_n,Y_n)$ converges in $L^2$ uniformly in $t$
to that of $(X,Y)$. Smoothing moves each particle by at most $\varepsilon_n$,
and coupling the background mass $\theta_n$ across $B$ costs at most
$4\theta_n R_n^2$ in squared distance. Thus $h_t\dx\to\eta_t$ uniformly
in $W_2$. The action lies between $W_2^2(h_0\dx,h_1\dx)$ and $D^2$,
so it converges to $D^2$.
\end{proof}

\begin{lemma}\label[lemma]{lem:resolvent-contraction-limit}
Fix $\tau>0$ and $C\ge0$. Suppose a proper energy $\cE$ on $\Ptwo$
and its restrictions to closed balls have unique resolvents at step $\tau$.
Assume the following:
\begin{enumerate}[label=\textup{(\roman*)},leftmargin=2em]
\item $\cE$ is narrowly lower semicontinuous along sequences with bounded
second moments. Each finite-energy measure is the limit in $W_2$ and in
energy of compactly supported measures.
\item For the point mass $\delta_0$, the minimizers on the closed balls
$\overline B_R$ have uniformly bounded second moments as $R\to\infty$.
\item On each fixed ball, the restricted energy is the $\Gamma$-limit in $W_2$
of energies whose resolvents at step $\tau$ exist, are unique and are
$C$-Lipschitz in $W_2$.
\end{enumerate}
Then $J_\tau^{\cE}$ is $C$-Lipschitz in $W_2$ on $\Ptwo$.
For $\cE=\Theta\circ\Ent_m$ as in \cref{lem:compact-resolvent}\,\textup{(ii)},
conditions \textup{(i)} and \textup{(ii)} hold.
\end{lemma}

\begin{proof}
Write $J_R$ for the resolvent on $\overline B_R$. Adding the continuous
squared transport cost preserves the $\Gamma$-limit in \textup{(iii)}.
Compactness and uniqueness give convergence of the minimizers on each
fixed ball, so $J_R$ is $C$-Lipschitz.

For $\mu\in\Ptwo$, let $\mu_R$ be its image under metric projection
onto $\overline B_R$. Then $\mu_R\to\mu$ in $W_2$. Since the projection
is a contraction fixing the origin,
\[
 W_2(J_R\mu_R,\delta_0)
 \le C W_2(\mu,\delta_0)+W_2(J_R\delta_0,\delta_0).
\]
By \textup{(ii)}, the right-hand side is bounded uniformly in $R$
for each fixed $\mu$. Thus the minimizers $J_R\mu_R$ have uniformly
bounded second moments. This yields narrow compactness, which is sufficient for the
Lipschitz bound once we identify the limit. If $\rho$ is a narrow
subsequential limit, lower semicontinuity and comparison with a
compactly supported competitor $\sigma$ give
\[
 \cE(\rho)+\frac{W_2^2(\rho,\mu)}{2\tau}
 \le \cE(\sigma)+\frac{W_2^2(\sigma,\mu)}{2\tau},
\]
since $\sigma$ is admissible on every sufficiently large ball.
Recovery in \textup{(i)} extends the comparison to all finite-energy
competitors. Thus $J_R\mu_R$ converges narrowly to $J_\tau^{\cE}\mu$.
Applying this to both input measures and using joint lower semicontinuity gives
\[
 W_2(J_\tau^{\cE}\mu,J_\tau^{\cE}\nu)
 \le\liminf_{R\to\infty} W_2(J_R\mu_R,J_R\nu_R)
 \le C W_2(\mu,\nu).
\]

For $\cE=\Theta\circ\Ent_m$, the bound $\cE\ge\Theta(0)$ and
comparison with a fixed compactly supported density of finite entropy
give \textup{(ii)}, since $W_2^2(\rho,\delta_0)=\int|x|^2\dd\rho$.
For general inputs the same cost bound gives existence by lower
semicontinuity, and strict convexity under mixtures gives uniqueness.
For recovery in \textup{(i)}, set $a_R=\sigma(B_R)$ for $\sigma=hm$
of finite entropy. The normalized restrictions converge in $W_2$, and
\begin{equation}\label{eq:entropy-restriction-recovery}
 \Ent_m(\sigma|_{B_R}/a_R)
 =a_R^{-1}\int_{B_R}h\log h\dd m-\log a_R
 \longrightarrow\Ent_m(\sigma).
\end{equation}
Continuity of $\Theta$ gives recovery in energy.
\end{proof}

\begin{lemma}
\label[lemma]{lem:lsi-approximation}
The following results hold on any sequence of balls
$B_{R_j}$ with $R_j\uparrow\infty$.
\begin{enumerate}[label=\textup{(\roman*)},leftmargin=*]
\item \emph{Approximation of the inputs.}
For any finite family of measures in $\Ptwo(\R^d)$, there are Markov
kernels $K_j$ and times $t_j\downarrow0$, common to the family, such
that for each input $\mu$
the measure $\mu K_j$ has a smooth positive density on
$\overline B_{R_j}$ and converges to $\mu$ in $W_2$. If $\mu$ satisfies
$\operatorname{LSI}(C)$, $0<C<\infty$, then $\mu K_j$ satisfies
$\operatorname{LSI}(C+2t_j)$.

\item \emph{Convergence of the minimizers.}
Fix $\tau>0$ and $m=Z^{-1}e^{-V}\dx$, where $V$ is finite and
$\kappa$-strongly convex for some $\kappa>0$. For inputs
$\mu_j\in\mathcal P(\overline B_{R_j})$ converging to $\mu$ in $W_2$,
let $\rho_j$ minimize
\[
 F_{\mu_j}(\zeta)=\Ent_m(\zeta)+\frac1{2\tau}W_2^2(\zeta,\mu_j)
 \quad\text{over }\mathcal P(\overline B_{R_j}).
\]
Then $\rho_j\to\rho=J_\tau^m\mu$ in $W_2$ and
$F_{\mu_j}(\rho_j)\to F_\mu(\rho)$.
\end{enumerate}
In particular, if $\rho_j$ satisfies $\operatorname{LSI}(L_j)$ and
$L_j\to L<\infty$, then $\rho$ satisfies $\operatorname{LSI}(L)$.
\end{lemma}

\begin{proof}
\emph{(i)} The contractions $F_R(x)=Rx/\sqrt{R^2+|x|^2}$ converge to
$\Id$ in $L^2(\mu)$ and preserve the log-Sobolev bound by pullback.
Let $P_t^{B_R}$ be the Neumann heat kernel on $B_R$ with generator
$\Delta$. On each fixed ball, the heat flow smooths any probability
to a positive density and converges to it in $W_2$ as $t\downarrow0$.
Choose $t_j\downarrow0$ so that the heat step moves every
$(F_{R_j/2})_\#\mu$ in the given family by at most $1/j$ in $W_2$.
The required kernels are
\[
 K_j(x,\diff y)=P_{t_j}^{B_{R_j}}(F_{R_j/2}(x),\diff y).
\]

Write $P_t=P_t^{B_R}$. For smooth $g>0$, Wang's gradient and local
log-Sobolev estimates on the convex ball \cite[Theorem~1.1]{WangNeumann}
give
\[
 \frac{|\nabla P_tg|^2}{P_tg}\le P_t\!\left(\frac{|\nabla g|^2}{g}\right),
 \qquad
 P_t(g\log g)-(P_tg)\log(P_tg)
 \le tP_t\!\left(\frac{|\nabla g|^2}{g}\right).
\]
If $\lambda$ is a contracted input satisfying $\operatorname{LSI}(C)$,
integrate the second bound against $\lambda$ and apply the input
inequality to $\sqrt{P_tg}$. The first bound then gives
$\operatorname{LSI}(C+2t)$ for $\lambda P_t$.

\emph{(ii)} Comparison with a fixed compactly supported competitor of
finite entropy bounds $W_2(\rho_j,\mu_j)$ and hence the second moments
of $\rho_j$. Lower semicontinuity and recovery by normalized restrictions
in \eqref{eq:entropy-restriction-recovery} give convergence of the minima
and identify every narrow limit with the unique minimizer $\rho$.

To upgrade this convergence to $W_2$ we must also recover the
second moments. We use part of the quadratic confinement for this
purpose. Fix $0<\delta<\kappa/2$. Strong convexity gives
$Z_\delta=\int e^{\delta|x|^2}\dd m<\infty$. For the probability
$m_\delta=Z_\delta^{-1}e^{\delta|\cdot|^2}m$ we have
\[
 F_\eta(\zeta)-\delta\int |x|^2\dd\zeta
 =\Ent_{m_\delta}(\zeta)+\frac1{2\tau}W_2^2(\zeta,\eta)-\log Z_\delta .
\]
The right side is jointly narrowly lower semicontinuous. Since the
minimum values converge, this gives
$\limsup_j\int|x|^2\dd\rho_j\le\int|x|^2\dd\rho$. Lower
semicontinuity gives the reverse bound for the lower limit, proving
convergence of second moments and hence in $W_2$.

On each ball mollification approximates $V$ uniformly and preserves
its curvature bound. After normalizing both references on the ball,
an error $\varepsilon$ in $V$ changes all finite objectives by at most
$2\varepsilon$. On a fixed ball compactness and uniqueness then give
convergence of the minimizers. On expanding balls the normalization of
$m$ adds $\log m(B_{R_j})\to0$, so the same comparisons apply when
the approximation errors tend to zero. Testing \eqref{eq:lsi-normalization} with $f\in C_c^\infty$
passes the output bound to the limit.
\end{proof}

\begin{lemma}\label[lemma]{lem:lsi-envelope-whole-space}
Let $V:\R^d\to\R$ be finite and $\kappa$-strongly convex, $\kappa>0$,
and let $\mu\in\Ptwo$ have a smooth positive density and satisfy
$\operatorname{LSI}(C)$ for some $0<C<\infty$. Suppose the Euler map
$T_\rho=\nabla\Phi_\rho$ of $\rho=J_\tau^m\mu$ is a $C^1$
diffeomorphism of $\R^d$.
Let $\Psi_\rho=\Phi_\rho^*$ and
$v_\rho=(|y|^2/2-\Psi_\rho)/\tau$, with the normalization in
\eqref{eq:entropy-potential-gauge}. For every
$\eta\in\Ptwo(\R^d)$ with $\Ent_\mu(\eta)<\infty$,
$v_\rho\in L^1(\eta)$ and
\begin{equation}\label{eq:lsi-envelope-global}
 0\le\cE_\tau(\eta)-\cE_\tau(\mu)
          -\int v_\rho\dd(\eta-\mu)
       \le\theta\Ent_\mu(\eta),
\end{equation}
where $\theta$ is given by \eqref{eq:lsi-theta} and
$\cE_\tau$ now denotes the whole-space minimum.
\end{lemma}

\begin{proof}
To pass the ball estimate to the limit, we need both the relative
entropy of the two inputs and the integral of the transport potential.
For the entropy bound we apply the same kernel from
\cref{lem:lsi-approximation}\,\textup{(i)} to $\mu$ and $\eta$. Then $\mu_j=\mu K_j\to\mu$ and
$\eta_j=\eta K_j\to\eta$ in $W_2$, with
\[
 \mu_j\in\operatorname{LSI}(C+2t_j),\qquad
 \Ent_{\mu_j}(\eta_j)\le\Ent_\mu(\eta)
\]
by entropy monotonicity under $K_j$. Choose smooth $\kappa$-strongly
convex $V_j$ within $1/j$ of $V$ on $\overline B_{R_j}$.
By \cref{lem:lsi-approximation}\,\textup{(ii)} and its potential
approximation, the two minima converge to their whole-space values
and the minimizers at $\mu_j$ converge to $\rho$ in $W_2$.

We next control the potential on bounded sets. Let $\Psi_j$ be the Brenier potential from $\mu_j$ to $\rho_j$ and
write $S_j=\nabla\Psi_j$. Subtract the values at $0$ from both $\Psi_j$
and $\Psi_\rho$. Both marginals converge in $W_2$, and $\mu$ has a
positive density on $\R^d$. Applying \cite[Theorem~3.4(i)]{delBarrioGonzalezSanzLoubes}
to the quadratic cost potentials $|\cdot|^2/2-\Psi_j$ gives
$\Psi_j\to\Psi_\rho$ locally uniformly with this normalization.
Convexity then gives $\sup_j|S_j(0)|<\infty$.

To integrate these locally uniform limits against the varying inputs,
we still need to bound the contribution from large $|y|$. The ball estimate in the proof of \cref{thm:lsi-jko} gives a common
log-Sobolev bound for $\rho_j$. The Talagrand inequality on
the convex balls \cite{OttoVillani,GigliLedoux} and monotonicity of relative
entropy under pushforward by $S_j$ therefore give, with $A$ independent of $j$,
\[
 W_2^2((S_j)_\#\eta_j,\rho_j)
 \le A\Ent_{\rho_j}((S_j)_\#\eta_j)\le A\Ent_\mu(\eta).
\]
As $(S_j)_\#\mu_j=\rho_j\to\rho$ in $W_2$, this bounds
$\int|S_j|^2\dd\lambda_j$ for both $\lambda_j=\mu_j$ and
$\lambda_j=\eta_j$. The supporting planes at $0$ and $y$ give
\[
 \nabla\Psi_j(0)\cdot y\le\Psi_j(y)\le S_j(y)\cdot y.
\]
The bound on $S_j(0)$ and Cauchy--Schwarz now give, for either sequence,
\[
 \int_{|y|>R}|\Psi_j(y)|\dd\lambda_j
 \le C'\left(\int_{|y|>R}|y|^2\dd\lambda_j\right)^{1/2}.
\]
The right side tends to zero uniformly in $j$ by $W_2$ convergence.
Together with local uniform convergence this gives $v_\rho\in L^1(\eta)$
and
\[
 \int \frac{|y|^2/2-\Psi_j(y)}{\tau}\dd(\eta_j-\mu_j)(y)
 \longrightarrow\int v_\rho\dd(\eta-\mu).
\]
The additive constants cancel in both pairings. The constants in the
ball estimates converge to $\theta$ by \eqref{eq:lsi-theta}.
Passing to the limit in those estimates proves
\eqref{eq:lsi-envelope-global}.
\end{proof}

The proof on a Ricci-flat manifold requires a different limiting
argument. The next lemma estimates the two distance quotients by
following particles in normal coordinates. For the lower bound on
the output distance we also need a test function whose gradient
approximates the particle velocity.

\begin{lemma}\label[lemma]{lem:ricci-distance-quotients}
For the curves in the proof of \cref{thm:ricci-resolvent}, with
$\operatorname{Ric}\equiv0$ and $p,v,C$ fixed, one has
\[
 \begin{aligned}
 \liminf_{h\to0}\frac{W_2^2(\rho_{\tau,h},\rho_{\tau,0})}{h^2}
   &\ge |v|^2-\tau\tr(R_v C^{-1})+o(\tau),\\
 \limsup_{h\to0}\frac{W_2^2(\mu_{\tau,h},\mu_{\tau,0})}{h^2}
   &\le |v|^2-\tau\tr(R_v B_C)+o(\tau)
 \end{aligned}
 \qquad(\tau\downarrow0).
\]
\end{lemma}

\begin{proof}
In the moving chart $E_h(z)=\exp_{p(h)}(P_hz)$ let $j_h$ be the
volume density. Near $z=0$, the change of variables by $F_h$ gives,
up to a spatial constant,
\[
 \psi_{\tau,h}(E_h(z))
 =\tfrac12\inner{Cz}{z}+\tau\bigl(\log j_h(z)-\log j_0(z)\bigr).
\]
Writing $u=P_hz$ at $p(h)$, the normal-coordinate expansion is
\[
 \log j_h(z)=-\tfrac16\operatorname{Ric}(u,u)
              -\tfrac1{12}(\nabla_u\operatorname{Ric})(u,u)
              +O(|z|^4).
\]
Both displayed coefficients vanish for every center along $p(h)$.
Smooth Taylor expansion in $(h,z)$ thus gives $\D_z\log j_h=O(|z|^3)$,
also after one $h$ derivative. Since the inverse metric in this chart is $I+O(|z|^2)$, it follows
that $\nabla\psi_{\tau,h}=Cz+O(|z|^3)$ in these coordinates.
The geodesic defining $S_{\tau,h}$ starts with the negative of this
velocity. Its Christoffel symbols are $O(|z|)$ and its coordinate
acceleration is $O(|z|^3)$ along the unit time interval. Hence
\[
 y_h(z):=E_h^{-1}(S_{\tau,h}(E_h(z)))=(I-C)z+O(|z|^3),
 \qquad \partial_h y_h(z)=O(|z|^3).
\]
The second bound follows by smooth Taylor expansion since all
coefficients through degree two are independent of $h$. The bounds
are uniform on a fixed smaller ball for small $h$ and $\tau$.

Under $\rho_{\tau,0}$ the local density in $z$ is proportional to
$e^{-\inner{Cz}{z}/(2\tau)}j_0(z)$. Outside a fixed smaller ball the potential is bounded away from its minimum.
Dividing by the normalization of order $\tau^{d/2}$ still leaves mass
$O(e^{-c/\tau})$ for a smaller $c>0$, and the particle velocities there are uniformly bounded.
Rescaling by $\sqrt\tau$ gives
$\E[z\otimes z]=\tau C^{-1}+o(\tau)$ and $\E|z|^k=O(\tau^{k/2})$ for
$k\ge1$, taking the coordinates to be zero off the ball.
Consequently $\E[y_0\otimes y_0]=\tau B_C+o(\tau)$.

For fixed $y$, $\partial_h E_h(y)|_{h=0}$ is the endpoint of the
Jacobi field with initial data $(v,0)$. Radial parallel transport back to $p$
gives $v-\tfrac12 R(v,y)y+O(|y|^3)$. The particles of $\rho_{\tau,h}$
follow $E_h(z)$, and those of $\mu_{\tau,h}$ follow $E_h(y_h(z))$,
with the additional velocity $\D E_0(y_0)\,\partial_h y_h|_{h=0}=O(|z|^3)$.
For either curve $\sigma_h$ with particle velocity $V$ at zero we obtain
\[
 \int |V|^2\dd\sigma_0
 = |v|^2-\E\inner{R_vy}{y}+O(\tau^{3/2}),
 \qquad y=z\ \text{or}\ y=y_0(z).
\]
For fixed $\tau$, coupling the same particles at $0$ and $h$ bounds the upper limit
of $W_2^2(\sigma_h,\sigma_0)/h^2$ by this integral, proving the
assertion for $\mu_{\tau,h}$.

For $\sigma_h=\rho_{\tau,h}$ we need a lower bound on the cost of
every coupling. The particle velocity has leading term $v$ in normal
coordinates, so we choose a smooth test function $\ell$ that agrees with $\ell(E_0(z))=\inner{v}{z}$ near $p$.
Taylor's inequality along a minimizing geodesic from $x$ to $x'$
and Young's inequality give
\[
 2h\bigl(\ell(x')-\ell(x)\bigr)
 \le h^2|\nabla\ell(x)|^2+
      (1+|h|\norm{\nabla^2\ell}_\infty)d_g(x,x')^2 .
\]
Integrate against an optimal plan from $\sigma_0$ to $\sigma_h$,
divide by $h^2$, and use
$\left.\partial_h\int\ell\dd\sigma_h\right|_{0}
 =\int\inner{\nabla\ell}{V}\dd\sigma_0$. This yields
\[
 \liminf_{h\to0}\frac{W_2^2(\sigma_h,\sigma_0)}{h^2}
 \ge \int\bigl(|V|^2-|V-\nabla\ell|^2\bigr)\dd\sigma_0.
\]
Near $p$ the transported gradient is $v+O(|z|^2)$, so the last
error is $O(\E|z|^4)+O(e^{-c/\tau})=O(\tau^2)$.
Inserting $\E[z\otimes z]=\tau C^{-1}+o(\tau)$ proves the lower bound.
\end{proof}

\section{Relaxed internal energies}
\label[appendix]{app:resolvent-tools}

When the internal energy has finite recession slope, a measure with
a singular part can have finite energy. The proof of the pressure
criterion therefore needs two facts about the relaxed problem. Its
minimizer must be uniquely determined even for singular inputs,
and the energy estimates for smooth densities must pass to every
competitor of finite energy.

\begin{lemma}\label[lemma]{lem:radial-wellposed}
Let $d$ and $U$ satisfy the hypotheses of \cref{thm:pressure-global}.
The relaxed proximal problem has
a unique minimizer for every $\mu\in\Ptwo$ and $\tau>0$. The same holds
when the input measure and the competitors are supported on a fixed closed ball.
If $U_\infty<\infty$ and the minimizer is $r\dx+\rho^s$, every
optimal plan from that minimizer to $\mu$ transports $\rho^s$ along
the diagonal. In particular $\rho^s\le\mu^s$, where $\mu^s$ is the
singular part of the input.
For the uniform probability measure $\beta_R$ on $B_R$ and the
envelope $e_\tau(R)$ in \eqref{eq:radial-envelope},
\[
 e_\tau(R)\longrightarrow \mathcal U(\beta_R)
       \quad\text{as }\tau\downarrow0.
\]
\end{lemma}

\begin{proof}
Write $\rho=r\dx+\rho^s$ and $M_2(\rho)=\int|x|^2\dd\rho$.
By convexity, $U(r)\le rU(1)$ for $r\le1$ and $U(r)\ge rU(1)$
for $r\ge1$. The assumed bound $U(r)\ge-Cr^q$ on $(0,1]$ then gives
continuity at zero and $U(r)\ge-C'(r+r^q)$ on $[0,\infty)$.
Adding a linear term to $U$ shifts the energy of every probability
by the same constant, so we may assume $U_\infty\ge0$ when it is finite.

To compare the negative part of the energy with the quadratic
transport cost, split the integral of $r^q$ at radius
$R=(1+M_2(\rho))^{1/2}$. We claim that
\begin{equation}\label{eq:internal-moment-bound}
 \int r^q\dx\le C_{d,q}(1+M_2(\rho))^{d(1-q)/2}.
\end{equation}
Indeed, $\int_{B_R}r^q\dx\le |B_R|^{1-q}$ since $\int r\dx\le1$.
For every $R>0$, H\"older's inequality and $q>d/(d+2)$ give
\begin{equation}\label{eq:internal-tail-bound}
 \int_{|x|>R}r^q\dx
 \le M_2(\rho)^q
       \left(\int_{|x|>R}|x|^{-2q/(1-q)}\dx\right)^{1-q}
 \le C_{d,q}M_2(\rho)^qR^{d-(d+2)q}.
\end{equation}
Since the exponent of $R$ is negative and $\rho(B_R^c)\le M_2(\rho)/R^2$,
the negative energy on $B_R^c$ tends uniformly to zero on families with
bounded second moments. On balls, the relaxed integral is weakly lower
semicontinuous \cite[Chapter 9]{AGS}. Exhaustion by balls whose boundaries
have zero mass for the limit measure gives narrow lower semicontinuity
on $\R^d$ along sequences with uniformly bounded second moments.

By \eqref{eq:internal-moment-bound},
\[
 \mathcal U(\rho)\ge-C'(1+M_2(\rho))^{d(1-q)/2},\qquad
 W_2^2(\rho,\mu)\ge\tfrac12 M_2(\rho)-M_2(\mu).
\]
Since $d(1-q)/2<1$, comparison with a fixed uniform-ball competitor
bounds the second moments of a minimizing sequence. Tightness and
lower semicontinuity give existence. On a closed ball the same
conclusion follows by compactness.

For uniqueness we first determine what an optimal plan can do with
the singular part of a minimizer. If $U_\infty=+\infty$, finite energy gives $\rho^s=0$. Otherwise,
decompose an optimal plan from the minimizer to $\mu$ according to
its first marginal $r\dx+\rho^s$. Let $\pi^s$ be the part with first
marginal $\rho^s$ and let $\beta$ be its second marginal. Since $U(a+b)\le U(a)+U_\infty b$, replacing
$\rho^s$ by $\beta$ and $\pi^s$ by $(\Id,\Id)_\#\beta$ gives
\[
 \begin{aligned}
 \mathcal U(r\dx+\beta)
 &\le \int U(r)\dx+U_\infty\beta(\R^d)=\mathcal U(\rho),\\
 W_2^2(r\dx+\beta,\mu)
 &\le W_2^2(\rho,\mu)-\int|x-y|^2\dd\pi^s(x,y).
 \end{aligned}
\]
Thus $\int |x-y|^2\dd\pi^s=0$ by minimality. Both marginals of
$\pi^s$ are therefore $\rho^s$. Since the second marginal is dominated
by $\mu$ and $\rho^s$ is singular, $\rho^s\le\mu^s$.

Now suppose there are two minimizers and average their optimal
plans. The objective at the average of the first marginals is bounded
above by the common minimum, so the convexity inequalities for the
energy and the cost must both be equalities. Strict convexity of $U$ gives
a common density $r$, and the averaged plan is optimal. Off a common
null set carrying the singular parts of the first marginals, its
restriction is optimal with first marginal $r\dx$. When $r\dx\ne0$
this restriction is induced by a map \cite{Brenier}, so the original
plans agree there. The remaining second marginal is common to both plans
and, by the diagonal property, equals each singular first marginal.
If the common density has zero mass, the diagonal property already
identifies each minimizer with $\mu$.

It remains to identify the limit at the uniform input $\beta_R$.
Comparing two step sizes bounds how far its minimizer can move. Fix
$\tau_0>0$ and let $\rho_\tau$ minimize at $\beta_R$. Testing
$e_{\tau_0}(R)$ with $\rho_\tau$ and $e_\tau(R)$ with $\beta_R$ gives,
for $0<\tau<\tau_0$,
\[
 \left(\frac1{2\tau}-\frac1{2\tau_0}\right)W_2^2(\rho_\tau,\beta_R)
 \le e_\tau(R)-e_{\tau_0}(R)\le\mathcal U(\beta_R)-e_{\tau_0}(R).
\]
Thus $\rho_\tau\to\beta_R$ in $W_2$. Lower semicontinuity gives
$\mathcal U(\beta_R)\le\liminf_{\tau\downarrow0}\mathcal U(\rho_\tau)$.
Since $\mathcal U(\rho_\tau)\le e_\tau(R)\le\mathcal U(\beta_R)$,
the claimed limit follows.
\end{proof}

\begin{lemma}\label[lemma]{lem:internal-recovery}
Under the hypotheses of \cref{thm:pressure-global}, the following
approximations hold.
\begin{enumerate}[label=\textup{(\roman*)},leftmargin=2em]
\item On a fixed closed ball, any probability measure of finite relaxed
 energy is a $W_2$ limit of smooth strictly positive densities on the closed
 ball whose internal energies converge to its energy. For any sequence
 $\varepsilon_j\downarrow0$, these densities can be chosen so that their
 $\mathcal U+\varepsilon_j\Ent_{\mathcal L^d}$ energies have the same limit.
\item On $\R^d$, any finite-energy member of $\Ptwo$ is a $W_2$ and
 energy limit of smooth compactly supported probability densities.
\end{enumerate}
\end{lemma}

\begin{proof}
Write $\sigma=f\dx+\sigma^s$ and smooth both parts with the same
operator $P_h$. On the ball we use the Neumann heat semigroup, and
on $\R^d$ convolution with a smooth probability kernel supported in
$B_h$. Both preserve mass and Lebesgue measure. When $U_\infty<\infty$,
apply $U(a+b)\le U(a)+U_\infty b$ to the two smoothed densities.
Jensen's inequality then gives
\[
 \mathcal U(P_h\sigma)
 \le \int U(P_h f)\dx+U_\infty\sigma^s(\R^d)
 \le \mathcal U(\sigma).
\]
On the ball this includes singular mass carried by the boundary.
If $U_\infty=+\infty$, finite energy gives $\sigma^s=0$, and Jensen's
inequality alone gives the bound.

On the ball, Neumann heat smoothing~\cite{Davies} gives a smooth
density up to the boundary for $h>0$. Conservation of mass, the strong
maximum principle and the boundary point lemma give strict positivity
also on the boundary. For a smooth test function $v$ with zero normal
derivative, the semigroup identity gives
$\|P_hv-v\|_\infty\le h\|\Delta v\|_\infty$.
Such tests are uniformly dense in the continuous functions on the
closed ball (flatten a smooth test in the normal direction in a thin
boundary collar). By symmetry and compactness, $P_h\sigma\to\sigma$
narrowly and in $W_2$, including for initial mass on the boundary.
On $\R^d$ convolution gives $W_2(P_h\sigma,\sigma)\le h$.
In both cases lower semicontinuity gives convergence of the energies.
For \textup{(i)}, each smoothed ball density has finite entropy, so we
can choose $h_j\downarrow0$ slowly enough that
$\varepsilon_j|\Ent_{\mathcal L^d}(P_{h_j}\sigma)|\to0$.

For \textup{(ii)}, smooth before truncating. The smoothed density
is bounded, so renormalizing a cutoff only evaluates $U$ on a bounded
interval. Fix $h>0$ and write $P_h\sigma=g\dx$. Let
$0\le\chi_R\le1$ be smooth with compact support and equal to one on
$B_R$. Then $g_R=\chi_R g/\int\chi_R g\dx$ converges pointwise to
$g$ and $g_R\le2g$ for large $R$. Dominated convergence of second
moments gives $g_R\dx\to g\dx$ in $W_2$. Since $g$ is bounded,
convexity and $U(0)=0$ bound $U$ above by a linear function on
$[0,2\|g\|_\infty]$. The lower bound in the proof of
\cref{lem:radial-wellposed} therefore gives
\[
 |U(g_R)|\le C_h(g+g^q)\in L^1(\R^d)
\]
by \eqref{eq:internal-moment-bound}. Dominated convergence gives
$\int U(g_R)\dx\to\int U(g)\dx$. A diagonal choice as $h\downarrow0$
yields the required smooth compactly supported densities.
\end{proof}

\clearpage
\phantomsection
\addcontentsline{toc}{section}{References}

\begingroup
\raggedright

\endgroup

\end{document}